\documentclass[preprint,sort&compress,12pt]{elsarticle}
\usepackage{bbm}
\usepackage{amstext}
\usepackage{amsmath}
\usepackage{amssymb}
\usepackage{mathrsfs}
\usepackage{stmaryrd}
\usepackage{bm}
\usepackage{pstricks}
\usepackage{pst-coil}
\usepackage{pst-3d}
\usepackage{amsfonts}
\usepackage{amsthm}

\def\journal{
 }
\def\ftimes{\hfill Completed on November 20, 2018}
\newcommand{\mf}{\mathbf}
\newcommand{\mm}{\mathrm}
\newcommand{\ml}{\mathcal}

\newcommand{\be}{\begin{equation}}
\newcommand{\bea}{\begin{equation}\begin{aligned}}
\newcommand{\beas}{\begin{equation*}\begin{aligned}}
\newcommand{\eeas}{\end{aligned}\end{equation*}}
\newcommand{\eea}{\end{aligned}\end{equation}}
\newcommand{\ee}{\end{equation}}

\renewcommand{\div}{{\rm div }}

\begin{document}
\begin{frontmatter}
\title{
 Instability  of the Abstract Rayleigh--Taylor\\
 Problem and Applications}

\author[FJ]{Fei Jiang}
\ead{jiangfei0591@163.com}
\author[sJ]{Song Jiang}
\ead{jiang@iapcm.ac.cn}
\author[WCZ,WCZ1]{Weicheng Zhan}
\ead{zhanfzu@163.com}
\address[FJ]{College of Mathematics and
Computer Science, Fuzhou University, Fuzhou,   China.}
\address[sJ]{Institute of Applied Physics and Computational Mathematics,   Beijing,  China.}
\address[WCZ]{Institute of Applied Mathematics, Chinese Academy of Sciences, Beijing, China}
\address[WCZ1]{University of Chinese Academy of Sciences, Beijing, China}

\begin{abstract}
We prove the existence of a unique unstable strong solution in the sense of $L^1$-norm for an abstract Rayleigh--Taylor (RT) problem arising from stratified viscous fluids in Lagrangian coordinates based on a bootstrap instability method. In the proof, we develop a method  to modify the initial data of the linearized abstract RT problem based on an existence theory of unique solution of stratified (steady) Stokes problem and an iterative technique, so that the obtained modified initial data satisfy  necessary compatibility conditions of the (original) abstract RT problem.
Applying an inverse transformation of Lagrangian coordinates to the obtained unstable solution, and then taking proper values of parameters, we can further get unstable solutions for the RT problems in viscoelastic fluids, magnetohydrodynamics (MHD) fluids with zero resistivity  and pure viscous fluids (with or without interface intension) in Eulerian coordinates. Our results can be also extended to the corresponding inhomogeneous case (without interface).
\end{abstract}
\end{frontmatter}


\newtheorem{thm}{Theorem}[section]
\newtheorem{lem}{Lemma}[section]
\newtheorem{pro}{Proposition}[section]
\newtheorem{concl}{Conclusion}[section]
\newtheorem{cor}{Corollary}[section]
\newproof{pf}{Proof}
\newdefinition{rem}{Remark}[section]
\newtheorem{definition}{Definition}[section]

\section{Introduction}\label{introud}
\numberwithin{equation}{section}

Considering two completely plane-parallel layers of stratified (immiscible) pure fluids, the heavier one
on top of the lighter one and both subject to the earth's gravity, it is well-known that such equilibrium state is unstable to sustain small disturbances, and this unstable disturbance will grow and lead to a release of potential energy, as the heavier fluid moves down under the gravitational force, and the lighter one is displaced upwards. This phenomenon was first studied by Rayleigh \cite{RLIS} and then Taylor \cite{TGTP}, and is called therefore the Rayleigh--Taylor (RT) instability. In the last decades, this phenomenon has been extensively investigated from mathematical, physical and numerical aspects, see \cite{CSHHSCPO,WJH,GYTI1} for examples. It has been also widely investigated how the RT instability evolves under the  effects of other physical factors, such as elasticity \cite{JFJWGCOSdd,FJWGCZXOE,HGJJJWWWOJM,WWWYYZ2018}, rotation \cite{CSHHSCPO,DRJFJSRS}, internal surface tension \cite{GYTI2,WYJTIKCT,JJTIWYHCMP}, magnetic fields \cite{JFJSWWWOA,JFJSJMFM,JFJSJMFMOSERT,JFJSWWWN,WYC,WYTIVNMI} and so on.

In this article, we are interested in the mathematical proof for the existence of (nonlinear) RT instability solutions in Eulerian coordinates. In 2003, Hwang--Guo  made the first breakthrough to prove the existence of a unique RT instability solution in the sense of $L^2$-norm (of Lebesgue space) for an inhomogeneous (inviscid) pure fluid (without interface). Jiang--Jiang further proved the corresponding viscous case in a general bounded domain \cite{JFJSO2014}.  Later Pr\"uess--Simonett obtained a RT instability solution in some norm of Sobolev--Slobodeckii space for stratified pure viscous fluids (with or without interface intension) in $\mathbb{R}^3$ in 2010 \cite{PJSGOI5x}. Recently,  Wike further extended Pr\"uess--Simonett's result to the cylindrical domain with finite height \cite{WMRTITNS2017}. It should be noted that the value of norm of Sobolev--Slobodeckii space considered by Pr\"uss--Simonett is larger than the one of $L^2$-norm. Hence it is still an open problem that whether exists a RT instability solution  in the sense of $L^p$-norm for some $p\geqslant 1$. In this article, we indeed prove the existence of  a RT instability solution in the sense of $L^1$-norm for stratified pure viscous fluids defined on a horizontally periodic layer by further developing new mathematical analysis techniques. Obviously,  $L^1$-norm is the lowest norm in Lebesgue spaces, since the norm is defined on a bounded periodic cell. Moreover, we can obtain similar results in the other two types of fluids, i.e., viscoelastic fluids and magnetohydrodynamics (MHD) fluids with zero resistivity.

To uniformly investigate the RT instability solutions in the three types of fluids, we   consider the existence of RT instability solutions of a so-called abstract (stratified)  RT problem, which includes uniform equations of the three type of fluids in Lagrangian coordinates. The abstract RT problem defined on a horizontal periodic domain  reads as follows:
\begin{equation}\label{n0101nn}\left\{\begin{array}{ll}
\eta_t=u &\mbox{ in } \Omega,\\[1mm]
{\rho} u_t+\mm{div}_{\ml{A}}\mathcal{S}_{\mathcal{A}}(q,u,\eta) =\partial_{\bar{M}}^2\eta &\mbox{ in }  \Omega,\\[1mm]
\mathrm{div}_{\mathcal{A}} {u} = 0& \mbox{ in }  \Omega ,\\
 \llbracket  (\mathcal{S}_{\mathcal{A}} (q,u,\eta)- g
\rho \eta_3I)\mathcal{A}e_3-\bar{M}_3 \partial_{\bar{M}}\eta \rrbracket = \vartheta H \mathcal{A}e_3&\mbox{ on }\Sigma,\\
\llbracket \eta  \rrbracket=  \llbracket u  \rrbracket=0 &\mbox{ on }\Sigma,\\
(\eta, u)=0 &\mbox{ on }\Sigma_-^+,\\
(\eta,u)|_{t=0}=(\eta^0,u^0)  &\mbox{ in }  \Omega.
\end{array}\right.\end{equation}
Next we explain the notations in the above ART problem, and the other relevant notations.

(1) Domains:
\begin{align}
&\Omega_+:=\{(y_{\mm{h}},y_3)\in \mathbb{R}^3~|~y_{\mm{h}}:=(y_1,y_2)\in \mathbb{T},\;\; 0<y_3<\tau\},\nonumber \\  &\Omega_-:=\{(y_{\mm{h}},y_3)\in \mathbb{R}^3~|~y_{\mm{h}} \in \mathbb{T},\;\; -l<y_3<0\},\ \Omega:=\Omega_+\cup \Omega_-,\nonumber  \\
&\Sigma_+:=  \mathbb{T} \times \{\tau\},\  \Sigma_-:=  \mathbb{T} \times \{-l\},\ \Sigma:= \mathbb{T} \times \{0\},\ \Sigma_-^+:=\Sigma_+\cup \Sigma_-, \nonumber
\end{align}
where $l$, $\tau>0$, $\mathbb{T}:=\mathbb{T}_1\times \mathbb{T}_2$, $ \mathbb{T}_i=2\pi L_i(\mathbb{R}/\mathbb{Z})$, and $2\pi L_i$ ($i=1$, $2$) are the periodicity lengths.

(2) Unknowns: $\eta:=\eta(y,t):\Omega\times [0,T]\to \mathbb{R}^3$, and $u:=u(y,t):\Omega\times [0,T]\to \mathbb{R}^3$ are the unknowns. Moreover, $\zeta:=\eta+y:\Omega\times [0,T]\to \Omega$. In addition, we denote the functions $f|_{\Omega_\pm}$ by $f_\pm$ for $f=\eta$, $u$, $\zeta$ and so on.

(3) Known parameters and vectors: $\rho_\pm:=\rho|_{\Omega_\pm}>0$, $\mu_\pm:=\mu|_{\Omega_\pm}>0$,  $\kappa_\pm:=\kappa|_{\Omega_\pm}\geqslant 0$,  $\vartheta\geqslant 0$ and  $g>0$ are constants. $\bar{M}$ is a constant vector. $e_i$ denotes the unit vector, the $i$-th component in which is $1$.  In addition, $\rho_+$ and $\rho_-$ satisfy the RT condition $\rho_+>\rho_-$.

(4) Jump notation: For $f$ defined on $\Omega$, the notation $\llbracket f \rrbracket :=f_+|_{\Sigma }
-f_-|_{\Sigma }$, where $f_\pm|_{\Sigma }$ are the traces of the functions (or parameters) $f_\pm$ on $\Sigma$.

(5) Notations involving $\mathcal{A}$: The matrix $\mathcal{A}:=(\ml{A}_{ij})_{3\times 3}$ is defined via
\begin{equation*}
\ml{A}^{\mm{T}}=(\nabla \zeta)^{-1}:=(\partial_j \zeta_i)^{-1}_{3\times 3},
\end{equation*} where the subscript $\mm{T}$ denotes the transposition, and  $\partial_{j}$ denote the partial derivative with respect to the $j$-th components of variables $y$. $\tilde{\mathcal{A}}:= \mathcal{A}-I$, and $I$ is the $3\times 3$ identity matrix.
The differential operator $\nabla_{\ml{A}}$ is defined by $$\nabla_{\ml{A}}w:=(\nabla_{\ml{A}}w_1,\nabla_{\ml{A}}w_2,\nabla_{\ml{A}}w_3)^{\mm{T}}
\;\mbox{ and }\;\nabla_{\ml{A}}w_i:=(\ml{A}_{1k}\partial_kw_i,
\ml{A}_{2k}\partial_kw_i,\ml{A}_{3k}\partial_kw_i)^{\mm{T}}$$
for vector function $w:=(w_1,w_2,w_3)$, and the differential operator $\mm{div}_\ml{A}$ is defined by
\begin{equation*}
\mm{div}_{\ml{A}}(f^1,f^2,f^3)=(\mm{div}_{\ml{A}}f^1,\mm{div}_{\ml{A}}f^2,
\mm{div}_{\ml{A}}f^3)^{\mm{T}}
\mbox{ and }\mm{div}_{\ml{A}}f^i:=\ml{A}_{lk}\partial_k f_{l}^i
\end{equation*}
for vector function $f^i:=(f_{1}^i,f_{2}^i,f_{3}^i)^{\mm{T}}$. It should be noted that we have used the Einstein convention of summation over repeated indices. In addition, we define   $\Delta_{\mathcal{A}}X:=\mm{div}_{\ml{A}}\nabla_{\ml{A}}X$.

We shall remark some properties of $\mathcal{A}$. We can deduce from \eqref{n0101nn}$_1$ and \eqref{n0101nn}$_3$ that $\partial_t \det \nabla \zeta=0 $, where $\det$ denotes the determinant of matrix $\nabla \zeta$. Hence, if the initial date $\eta^0$ satisfies $\det\nabla(\eta^0+y)=1$, then $\det \nabla \zeta=1$. Thus, by the definition of $\mathcal{A}$, we see that
\begin{equation}
\label{06131422}
\mathcal{A}=(A^*_{ij})_{3\times 3},
\end{equation} where
$A^{*}_{ij}$ is the algebraic complement minor of $(i,j)$-th entry of matrix $(\partial_j \zeta_i)_{3\times 3}$.
Moreover, it is easy to check that $\mathcal{A}e_3=\partial_1\zeta\times \partial_2 \zeta$.
In addition,\begin{equation}
\label{06132125} \partial_kA^{*}_{ik}=0\mbox{ or } \partial_k\mathcal{A}_{ik}=0,  \end{equation}
   which implies some important relations
\begin{equation}\label{diverelation}
\mm{div}_{\tilde{\mathcal{A}}}u=\mm{div} (\tilde{\mathcal{A}}^{\mm{T}}u)
\end{equation}
 and
\begin{equation}
\mm{div}_{\partial_t^i\mathcal{A}}\partial_t^ju=\mm{div} (\partial_t^i\mathcal{A}^{\mm{T}}\partial_t^j u)\mbox{ for }i,\ j\geqslant 0.\label{20178011456115}\end{equation}

(6) Compound symbols, and operators:
\begin{align}   \nonumber
&\mathcal{S}_{\mathcal{A}}(q,u,\eta):=  \mathfrak{S}_{\mathcal{A}}(q,u,\eta)  - {\kappa}\rho ( \nabla \eta\nabla \eta^{\mm{T}}),
\   \mathfrak{S}_{\mathcal{A}}(q,u,\eta) :=qI-\mu \mathbb{D}_{\mathcal{A}}u-{\kappa}\rho  \mathbb{D}\eta,\nonumber
\\
&
\mathbb{D}_{\mathcal{A}}u:=\nabla_{\mathcal{A}}u +\nabla_{\mathcal{A}}u^{\mm{T}},
\
\mathbb{D}\eta: =
\nabla \eta +\nabla \eta^{\mm{T}},\  \partial_{\bar{M}}:=\bar{M}\cdot \nabla,\ \partial_{\bar{M}}^2:=(\partial_{\bar{M}})^2,\nonumber \\
&H^{n}:=|\partial_1 \zeta|^2\partial_2^2\zeta-2(\partial_1\zeta\cdot \partial_2\zeta)\partial_1\partial_2\zeta+
|\partial_2\zeta|^2\partial_1^2\zeta,\nonumber \\
&H^{\mm{d}}:= |\partial_1\zeta|^2|\partial_2\zeta|^2-|\partial_1\zeta\cdot \partial_2\zeta|^2,   \   H:= {H^{\mm{n}}}\cdot \vec{n}/{H^{\mm{d}}} ,\ \vec{n}:= {\mathcal{A}e_3}/{|\mathcal{A}e_3|},\nonumber \\
&\mathcal{S}  (q,u,\eta) :=qI- \mathbb{D} (\mu u+\kappa \rho\eta),\ \tilde{n}:=\vec{n}-e_3.\nonumber
\end{align}
Since we consider the strong solution  of the ART problem, we have $$\llbracket\partial_i\eta \rrbracket=\llbracket \partial_i^2\eta \rrbracket=0\mbox{ for }i=1\mbox{ and }2$$
 due to $\llbracket \eta\rrbracket=0$. Thus $\llbracket \varphi\rrbracket =0$ for
$\varphi=\mathcal{A}e_3$, $\vec{n}$, $H^{n}$, $H^{\mm{d}}$ or $H$. In addition, the operator $\mathcal{S}_{\mathcal{B}}$ is defined as $\mathcal{S}_{\mathcal{A}}$ with $\mathcal{B}$ in place of $\mathcal{A}$ for some matrix $\mathcal{B}$. Similar definitions also are used in the operators $\mm{div}_{\mathcal{B}}$, $\nabla_{\mathcal{B}}$, $\mathfrak{S}_{\mathcal{B}}$ and $\mathbb{D}_{\mathcal{B}}$.

For the derivation of \emph{a priori} estimate of the ART problem in Section \ref{201805081731}, next we shall further rewrite \eqref{n0101nn} as a nonhomogeneous form, in which the unknowns at the left hands of the identities are linear.

Let $f$ and $r\neq 0$ belong to $\mathbb{R}^3$, we define that $\Pi_{r}f:=f-(f\cdot r)r/|r|^2$. It should be noted that $\Pi_{r}f=0$, if only if $f  \sslash r$.
Then, applying the operator $\Pi_{\vec{n}}$ to \eqref{n0101nn}$_4$ yields that
$$ \Pi_{\vec{n}} \llbracket\mathcal{S}_{\mathcal{A}}  (q,u,\eta) \mathcal{A} e_3  -  \bar{M}_3 \partial_{\bar{M}}\eta\rrbracket =0,$$
which can be rewritten as a nonhomogeneous form:
\begin{align}
\Pi_{e_3} \llbracket\mathcal{S}  (q,u,\eta)  e_3  -  \bar{M}_3 \partial_{\bar{M}}\eta\rrbracket =\mathfrak{N},
\label{201806211105}
\end{align}
where we have defined that
$$
\begin{aligned}
&\begin{aligned}
\mathfrak{N}:=&
\llbracket (\mathcal{S}  (q,u,\eta)   e_3 -\bar{M}_3\partial_{\bar{M}}\eta )\cdot  \tilde{n}  \vec{n}   +    (\mathcal{S}  (q,u,\eta)    e_3  -\bar{M}_3\partial_{\bar{M}}\eta ) \cdot e_3\tilde{n}
\\
& - \Pi_{ \vec{n} }(  \mathcal{S}  (q,u,\eta)  \tilde{\mathcal{A}}e_3-  (\mu \mathbb{D}_{\tilde{\mathcal{A}}}u+{\kappa}\rho \nabla \eta\nabla \eta^{\mm{T}})\mathcal{A}e_3) \rrbracket.
\end{aligned}
\end{aligned}$$

Multiplying \eqref{n0101nn}$_4$ by $\vec{n}/|\mathcal{A}e_3|$ (i.e., taking scalar product), we have
 $$
 \llbracket( \mathcal{S}_{\mathcal{A}}(q,u,\eta)  \vec{n} - \bar{M}_3 \partial_{\bar{M}}\eta /|\mathcal{A}e_3|
)\cdot  \vec{n}-g  \rho  \eta_3 \rrbracket =  \vartheta H ,$$
which can be rewritten as a nonhomogeneous form:
\begin{equation}
\label{201806202203}
\llbracket q -2(\mu \partial_3u_3+{\kappa}\rho \partial_3\eta_3) - g \rho  \eta_3 -  \bar{M}_3 \partial_{\bar{M}}\eta_3\rrbracket-\vartheta\Delta_{\mm{h}}\eta_3= \mathfrak{M},\end{equation}
where we have defined that
$$\begin{aligned}
&\mathfrak{M}:= \mathfrak{M}_1+\llbracket\mathfrak{M}_2\rrbracket,\
\mathfrak{M}_1:=\vartheta(H^{\mm{n}} \cdot \vec{n}(1-H^{\mm{d}})/H^{\mm{d}}+H^{\mm{n}}\cdot \tilde{n}+H^{\mm{n}}\cdot e_3 -\Delta_{\mm{h}}\eta_3), \\
&\begin{aligned}\mathfrak{M}_2:=&
 \bar{M}_3 \partial_{\bar{M}}\eta \cdot \tilde{n} +(1-|\mathcal{A}e_3|) \bar{M}_3 \partial_{\bar{M}}\eta\cdot \vec{n}/|\mathcal{A}e_3| + (\mu\mathbb{D}_{\tilde{\mathcal{A}}}u+{\kappa}\rho \nabla \eta\nabla \eta^{\mm{T}})    \vec{n} \cdot  \vec{n} \\
&-  \mathcal{S} (q,u,\eta) \tilde{n}\cdot \vec{n}-  \mathcal{S} (q,u,\eta)  e_3\cdot \tilde{n} \mbox{ and } \Delta_{\mm{h}}:=\partial_1^2+\partial_2^2.
\end{aligned}
\end{aligned}
$$
It is easy to observe that
\eqref{n0101nn}$_4=\eqref{201806211105}+\eqref{201806202203}\mathcal{A}e_3$.

Now we define two nonlinear mappings of $(\eta,u,q)$ as follows:
\begin{align}
&\mathcal{N}^1(q,u,\eta):=\mu \mm{div}\mathbb{D}_{\tilde{\mathcal{A}}}u-\mm{div}_{\tilde{\ml{A}}}
\mathcal{S}_{\mathcal{A}}(q,u,\eta )
+\kappa\rho\mm{div}(\nabla \eta\nabla \eta^{\mm{T}}), \nonumber \\
 & \mathcal{N}^2(q,u,\eta):=(\mathfrak{N}_1, \mathfrak{N}_2,\mathfrak{M}), \nonumber
\end{align}
then, using \eqref{201806211105} and \eqref{201806202203},  we can easily get a nonhomogeneous form of the ART problem \eqref{n0101nn}:
\begin{equation}\label{s0106pnnnn}
\left\{
\begin{array}{ll}
 \eta_t=  u &\mbox{ in }\Omega,  \\[1mm]
 \rho u_{t}+\mm{div}\mathcal{S}  (q,u,\eta) =\partial_{\bar{M}}^2\eta+
 \mathcal{N}^1 &\mbox{ in }\Omega, \\[1mm]
\div u  =-  \mm{div}_{\tilde{\mathcal{A}}} u&\mbox{ in }\Omega, \\
 \llbracket  \mathcal{S}  (q,u,\eta)  e_3-g  \rho \eta_3 e_3  - \bar{M}_3 \partial_{\bar{M}}\eta
\rrbracket-\vartheta\Delta_{\mm{h}}\eta_3  e_3 =  \mathcal{N}^2  &\mbox{ on }\Sigma, \\[1mm]
\llbracket \eta  \rrbracket=\llbracket u  \rrbracket=0&\mbox{ on }\Sigma,\\
(\eta,u)=0&\mbox{ on }\Sigma_{-}^+,\\[1mm]
(\eta,u)|_{t=0}=(\eta^0,u^0)&\mbox{ on }\Omega,
\end{array}\right.
\end{equation}
It should be noted that we have used the simplified notations $\mathcal{N}^i:= \mathcal{N}^i(q,u,\eta)$ for $i=1$ and $2$ in the above ART problem for simplicity.
The problem obtained by omitting the nonlinear terms $\mathcal{N}^1$ and $\mathcal{N}^2$ in \eqref{s0106pnnnn}  is called the linearized ART problem. We will use unstable solutions of the linearized ART problem to construct unstable solutions of the (nonlinear or original) ART problem.

Before stating our main result, we shall introduce some simplified notations throughout this article.

(1) Basic notations:

 $\overline{D}$ denotes a closure of a point set $D\subset \mathbb{R}^N$. $ \overline{\Omega_-}=\mathbb{R}^2\times [-l,0]$, $\overline{\Omega_+}=\mathbb{R}^2\times [0,\tau]$, $\Omega_\pm':= \overline{\Omega_\pm}\setminus\partial\overline{\Omega_\pm} $,  $\overline{\Omega}=\mathbb{R}^2\times [-l,\tau]$, $ {\Omega}_{-l}^\tau:=\mathbb{R}^2\times (-l,\tau)$, $\Omega\!\!\!\!\!-:=\Omega\cup  \Sigma$, $\mathcal{T}_{m,n}^\varepsilon:=(- (m+\varepsilon) L_1\pi, (m+\varepsilon) L_1\pi)\times (-(n+\varepsilon)L_2\pi,(n+\varepsilon)L_2\pi )$, $\mathcal{T}:=\mathcal{T}_{2,2}^0$,  $\mathcal{T}^\varepsilon:=\mathcal{T}_{2,2}^\varepsilon$, $\Sigma_0^+:=\Sigma \cup \Sigma_+$,   $\mathcal{T}_{-l,\tau}^\varepsilon:=\mathcal{T}^\varepsilon\times ((-l,0)\cup (0,\tau))$,
$\mathcal{T}_{-l}^\tau:=\mathcal{T}_{-l,\tau}^0$, $\int:=\int_{\mathcal{T}_{-l}^\tau}$. $\int_\Sigma:=\int_{\mathcal{T}}$. $I_T:=(0,T)$ is a times interval, ${\mathbb{R}^2_{T}}:=\mathbb{R}^2\times I_T$, $\Omega_\pm^T:=\Omega_\pm'\times I_T$, $ \Omega^T:=\Omega_{-l}^{\tau}\times I_T$.

The $j$-th difference quotient of size $h$ is $D_j^h w:=(w(y+h e_j)-w(y))/h$ for $j=1$ and $2$, and $D^h_{\mm{h}}w:=(D_1^hw_1,D_2^hw_2)$, where $|h|\in (0,1)$, and $w$ is defined on $\Omega$ and a locally summable function.
Similarly, we define  the fractional differential operator $\mathfrak{D}_{\mf{h}}^{3/2}w:= (w(y+\mf{h}  )-w(y))/|\mf{h}|^{3/2}$ for $\mf{h}\in \mathbb{R}^2\times \{0\}$.
 $\nabla_{\mm{h}}:=(\partial_1,\partial_2)^{\mm{T}}$. $f_{\mm{h}}:=(f_1,f_2)$ and $\mm{div}_{\mm{h}}f_{\mm{h}}:=\partial_1 f_1+\partial_2 f_2$ for $f=(f_1,f_2,f_3)^{\mm{T}}$.  $\partial_{\mm{h}}^\alpha  $ denote $\partial_1^{\alpha_1}\partial_2^{\alpha_2}   $ for some multiindex of order $\alpha:=(\alpha_1,\alpha_2)$;
$\partial_\mm{h}^j$ denotes $\partial_{\mm{h}}^{\alpha}$ for any $\alpha$ satisfying $|\alpha|=\alpha_1+\alpha_2=j$;  $ \mathfrak{D}_{\mf{h}}^{3/2,\alpha} :=\mathfrak{D}_{\mf{h}}^{3/2} \partial_{\mm{h}}^{\alpha}$; $\nabla_{\mm{h}}^kw\in X$ denotes that $\partial_{\mm{h}}^{\alpha}w\in X$ for any  $\alpha$ satisfying $|\alpha|=k$, where $X$ denotes a function space.

Let $w(x,t)$ be a function defined on a set $D\times \overline{I_T}$ for some $T$, then $w(D,t):=\cup_{x\in D}\{w(x,t)\}$ for given $t\in \overline{I_T}$. If $w$ is independent of $t$, we denote $w(D)= \cup_{x\in D}\{w(x)\}$.
Let $f(x,t):\overline{\Omega^T} \to \mathbb{R}^3$. If $f(y_{\mm{h}},0,t):\mathbb{R}^2 \to \mathbb{R}^2$ is a diffeomorphic mapping for given $t\in \overline{I_T}$, then we denote the inverse function of $f_{\mm{h}}(y_{\mm{h}},0,t)$ with respect to $y_{\mm{h}}$ by $(f_{\mm{h}})^{-1}(x_{\mm{h}},t)$. $(f_{\mm{h}})^{-1}(\mathcal{T},t):=\cup_{ x_{\mm{h}}\in \mathcal{T}}\{ (f_{\mm{h}})^{-1}(x_{\mm{h}},t)\}$.  If $f$ is independent of $t$, we can omit $t$ in the definitions $(f_{\mm{h}})^{-1}(x_{\mm{h}},t)$ and $(f_{\mm{h}})^{-1}(\mathcal{T},t)$.

We define a $\beta$ function such that $\beta(a)=0$ if $\vartheta=0$, else $\beta(a)=a$, where $\vartheta$ is the parameter in the ART problem. $a\lesssim b$ means that $a\leqslant cb$ for some constant $c>0$, where the positive constant $c$ may depend on the domain $\Omega$, and known parameters such as $\rho_\pm$, $\mu_\pm$, $g$, $\vartheta$ and $\bar{M}$, and may vary from line to line.

(2) Simplified notations of function spaces：
\begin{align}
& L^p:=L^p (\Omega)=W^{0,p}(\Omega),\
W^{i,2}:=W^{i,2}(\Omega ),\  {H}^i:=W^{i,2} ,\nonumber \\
&   H^\infty:=\cap_{j=1}^\infty H^j,\ \underline{H}^i:=\left\{w\in {H}^i~\left|~\int_{\mathcal{T}^1}w\mm{d}y=0\right\}\right.,\nonumber \\
& H^{i+\epsilon}:=W^{i+\epsilon,p}(\mathbb{T})\mbox{ denotes Sobelve--Sobodetskii spaces}, \nonumber \\
&H^1_0:=\{w\in {H}^1(\Omega\!\!\!\!-)~|~w|_{\Sigma_-^+}=0\mbox{ in the sense of trace}\},\nonumber \\
  &H_0^i:=H_0^1 \cap H^i,\  H^1_\sigma:=\{w\in H^1_0  ~|~\div w=0\mbox{ in }\Omega\}, \nonumber\\
 &H^i_\sigma:=H^1_\sigma \cap H^i,\ H_{\sigma,\Sigma}^1:=
 \left\{w\in {H}_\sigma^1~\left|~ w_3|_{\Sigma}\in H^1(\mathbb{T})\right.\right\}, \nonumber \\
 &H^1_{\sigma,3}:= \{ w\in H^1_{\sigma,\Sigma}~|~w_3|_{\Sigma}\neq 0\},\ H_{\sigma, \vartheta}^1=\left\{
                            \begin{array}{ll}
                            H_{\sigma,\Sigma}^1   & \hbox{ if }\vartheta\neq 0, \\
                         H_{\sigma}^1  &\hbox{ if }\vartheta= 0,
                            \end{array}
                          \right.\nonumber \\
&   \mathcal{A}: = \left\{w\in {H}_{\sigma, \vartheta}^1~\left|~\|\sqrt{\rho}w\|^2_{L^2}=1\right.\right\},\nonumber
\end{align}
where $1< p\leqslant \infty$, $\epsilon\in (0,1)$, and  $k\geqslant 0$, $i\geqslant 1$ are integers.

To prove the existence of unstable classical solutions of linearized ART problem, we shall introduce a function space
$$H_{\sigma, \vartheta}^{1,k}:=\left\{
                            \begin{array}{ll}
                            \{w\in H_{\sigma,\Sigma}^1 ~\left|\nabla_{\mm{h}}^{j} w \in H^1\mbox{ and }w_3|_{\Sigma}\in H^{k+1}(\mathbb{T})\mbox{ for } j\leqslant k\right. \} & \hbox{ if }\vartheta\neq 0, \\
           \{w\in H_{\sigma }^1 ~\left|\nabla_{\mm{h}}^{j}w \in H^1 \mbox{ for } j\leqslant k \}\right. &\hbox{ if }\vartheta= 0,
                            \end{array}
                          \right.$$
where $k\geqslant 0$ is integer. It should be noted that $H_{\sigma, \vartheta}^{1,0}=H_{\sigma, \vartheta}^{1}$.

To make sure that the ART problem can be reformulated in  Eulerian coordinates, we shall consider the solution $\eta$ belonging to the function space
$ {H}_{0,*}^{3,1}\subset H^3_0$, in which  all elements denoted by $\varpi $ enjoy the following properties: for $\phi:=\varpi+y$,
\begin{align}
&\det\nabla \phi(y)=1,\label{201805012028}\\
&  \det\nabla_{\mm{h}}\phi_{\mm{h}}(y_{\mm{h}},0)
\geqslant 1/2\mbox{ for any }y_{\mm{h}}\in \mathbb{R}^2, \nonumber
 \\
& \phi_{\mm{h}}(y_{\mm{h}},0):\mathbb{R}^2 \to \mathbb{R}^2\mbox{ is a }C^1\mbox{-diffeomorphic mapping,}\nonumber  \\
 & \phi (y): \overline{\Omega } \mapsto  \overline{\Omega } \mbox{ is a homeomorphism mapping},\nonumber \\
 & \phi_\pm (y): \Omega_\pm' \mapsto  \phi_\pm(\Omega_\pm') \mbox{ are } C^1\mbox{-diffeomorphic  mappings},\nonumber \\
&\mathcal{T}^{-1/2}_{-l,\tau} \subset \phi^{-1}(\mathcal{T}_{-l}^\tau)\subset \mathcal{T}^{1/2}_{-l,\tau},\ \mathcal{T}^{-1/2}\subset  ( \phi_{\mm{h}})^{-1}(\mathcal{T})\subset \mathcal{T}^{1/2}.  \nonumber
\end{align}
We call the property  \eqref{201805012028}  the keeping volume condition.

(3) Simplified Sobolev  norms:
$$
\begin{aligned}
&\|\cdot \|_{i} :=\|\cdot \|_{W^{i,2}(\mathcal{T}_{-l}^\tau)},\ |\cdot|_{s} := \|\cdot|_{\mathcal{T}\times\{0\}} \|_{H^{s}(\mathcal{T})}, \  \|\cdot\|_{i,j}^2:=\sum_{|\alpha|=i} \|\partial^{\alpha}_{\mm{h}}\cdot\|_{j}^2, \\
& \|\cdot \|_{L^1}:=\|\cdot \|_{L^1(\mathcal{T}_{-l}^\tau)},\ |\cdot |_{L^1}:=\|\cdot |_{\mathcal{T}\times\{0\}}   \|_{L^1(\mathcal{T})},
\\
& \|\cdot\|^2_{\underline{i},j}:=\sum_{0\leqslant k\leqslant i}\|\cdot\|_{i,j}^2,\ |\cdot|_{-s}\mbox{ denotes the norm of the dual space of }H^s(\mathbb{T}),
\end{aligned}$$
where  $s$ is a real number, and $i$, $j$ are non-negative integers. It should be noted that
$$\|f|_{\mathcal{T}_{m,n}}\|_{H^{i+1/2}(\mathcal{T}_{m,n})}\leqslant c|f|_{i+1/2}$$ for any $f\in H^{i+1} $, where the constant $c$ depends on $m$ and $n$.
In addition, we define a so-called stratified Stokes norm:
$$\|(u,q)\|_{\mm{S},k}:=\sqrt{\|u\|_{k+2}^2+\| \nabla q\|_k^2
+ |\llbracket q   \rrbracket |_{k+1/2}^2}.$$

(4) Energy functionals:
\begin{equation*}  \begin{aligned}
& \mathcal{E}(t):=\|\eta(t)\|_3^2+\|(u,q)(t)\|_{\mm{S},0}^2 ,\  \mathcal{D}(t):= \|\eta(t)\|_3^2+\|(u,q)(t)\|_{\mm{S},1}^2
, \\
 & E(w,\vartheta,\kappa,\bar{M}):= \mathcal{I}(w,\vartheta,\kappa,\bar{M})-g\llbracket\rho\rrbracket |w_3|_{0}^2,\\
& \mathcal{I}(w,\vartheta,\kappa,\bar{M}):=\vartheta|\nabla_{\mm{h}}w_3|_{0}^2
+{\|\sqrt{\kappa\rho} \mathbb{D}w\|^2_0/2}+
\|\partial_{\bar{M}}w\|^2_0,\end{aligned} \end{equation*}
where $E$ is used to discriminate the stability and instability of the ART problem, and  $\mathcal{I}$ is called the  stabilizing effect term.
Sometimes we also denote $E(w,\vartheta,\kappa,\bar{M})$, reps. $\mathcal{I}(w,\vartheta,\kappa,\bar{M})$ by $E(w)$, reps. $\mathcal{I}(w)$ for the simplicity.

Now we introduce our main result.

\begin{thm}\label{thm3}
Under the instability condition
\begin{equation}
\label{201806071047}
 E(w,\vartheta,\kappa,\bar{M})<0\mbox{ for some }w\in H_{\sigma,\vartheta}^1,
\end{equation}
a zero solution to the ART problem \eqref{n0101nn} is unstable in the Hadamard sense, that is, there are positive constants $\Lambda$, $m_0$, $\epsilon$ and $\delta_0$,
and a vector function  $(\tilde{\eta}^0,\tilde{u}^0,\tilde{q}^0,\eta^\mm{r},u^\mm{r},q^\mm{r})\in H^4_\sigma\times  H^2_\sigma\times H^1\times H^4_0\times H_0^2\times H^1$,
such that, for any
\begin{equation*}
(\eta^0, u^0,q^0):=\delta(\tilde{\eta}^0,\tilde{u}^0,\tilde{q}^0)  + \delta^2(\eta^{\mm{r}},u^\mm{r},q^\mm{r})\in H^{4,1}_{0,*}\times H^2\times H^1\mbox{ with }\delta\in (0,\delta_0),\end{equation*}
there is a unique strong solution
$(\eta,u)\in C^0(\overline{I_T},H^{3,1}_{0,*}\times  H^2_0) $ to the ART problem with initial data $(\eta^0,u^0)$ and  with a unique (up to a constant) associated pressure
$q$ (with initial data $q^0$). Moreover the solution satisfies
\begin{align}
&|\chi_3(T^\delta)|_{L^1},\ \|\chi_3(T^\delta)\|_{L^1},\ \|\partial_3\chi_3(T^\delta)\|_{{L^1}},\ \|\chi_{\mm{h}}(T^\delta)\|_{L^1},\nonumber  \\
& \|\partial_3\chi_{\mm{h}}(T^\delta)\|_{{L^1}},\ \|\mm{div}_{\mm{h}}\chi_{\mm{h}}(T^\delta)\|_{L^1},\  \|\mathcal{A}_{3k}\partial_k\chi_3(T^\delta)\|_{{L^1}}, \ \nonumber \\
& \|\mathcal{A}_{3k}\partial_k\chi_{\mm{h}}(T^\delta)\|_{{L^1}},\  \|(\mathcal{A}_{1k}\partial_k \chi_{1}+ \mathcal{A}_{2k}\partial_k \chi_{2} )(T^\delta)\|_{L^1}\geqslant 4{\epsilon}\label{201806012326}
\end{align}
and
\begin{align}
\label{2018109121834}
\|\partial_{\bar{M}}\chi_3(T^\delta)\|_{L^1},\ \|\partial_{\bar{M}}\chi_{\mm{h}}(T^\delta)\|_{L^1}\geqslant 4{\epsilon}\mbox{ if }
\bar{M}_3\neq 0
\end{align}
for some escape time $T^\delta:=\frac{1}{\Lambda}\mm{ln}\frac{8\epsilon}{m_0\delta}\in I_T$, where $\chi=\eta$ or $u$, $\zeta:=\eta+y$, and $T$ denotes some time of existence of the solution $(\eta ,u)$. Moreover, the initial data $(\eta^0,u^0,q^0)$ satisfies necessary compatibility conditions:
\begin{alignat}{2}
&\mm{div}_{\mathcal{A}^0}u^0=0& &\ \mbox{ in }\Omega, \label{201809101928} \\
& \llbracket  \mathcal{S}_{\mathcal{A}^0} (q^0,u^0,\eta^0)\mathcal{A}^0e_3- g
\rho \eta_3^0  \mathcal{A}^0 e_3-\bar{M}_3 \partial_{\bar{M}}\eta^0 \rrbracket  =\vartheta H^0  \mathcal{A}^0 e_3 & &\ \mbox{ on }\Sigma, \label{xx201809101926xxsdfs}
\end{alignat}
where $H^0$, $\mathcal{A}^0$, and $\vec{n}^0$ denote the initial data of $H$, $\mathcal{A}$,  and $\vec{n}$, respectively, and are defined by $\eta^0$.
\end{thm}
\begin{rem}
\label{201809121705}
It should be noted that the solution $(u,q)$ in the above theorem enjoys the additional regularity:
\begin{align}
& (u_t,\nabla q, \llbracket q\rrbracket ) \in C^0(\overline{I_T}, L^2 \times L^2\times H^{1/2}(\mathbb{T})), \ q\in H^1\mbox{ for each }t\in \overline{I_T}, \nonumber
\\
&(u,\nabla q, u_t,\llbracket q\rrbracket) \in L^2(I_T, H^3\times H^1\times H^1\times H^{3/2}(\mathbb{T}) ). \label{201811142107}
\end{align}\end{rem}
\begin{rem}\label{thm3xx}
Noting that (see Lemma A.4 in \cite{WYTIVNMI})
\begin{equation}
\label{201809222059}
\|w\|_0^2\leqslant (l^2+\tau^2)\|\partial_{\vec{n}}w\|_0^2\mbox{ for any }w\in H^1_0,
\end{equation}
where $\vec{n}$ is a constant vector with $\vec{n}_3=1$,
thus, making use of  \eqref{201809222059}, \eqref{201808161915}, \eqref{2018060102215680}, \eqref{201806101508} and Korn's inequality (see Lemma \ref{xfsddfsf201805072234}), we can see that, for  $(\vartheta,\kappa,\bar{M}_3)\neq 0$,
\begin{align}
 \mathcal{I}(w,\vartheta,\kappa,\bar{M}) =0\mbox{ for some }w\in H_{\sigma,\vartheta}^1,\mbox{ if and only if }
\left\{
                          \begin{array}{ll}
                        w=0   & \hbox{ for }(\kappa,\bar{M}_3)\neq 0; \\
                          | w_3|_{0}^2=0 & \hbox{ for }(\kappa,\bar{M}_3)=0
                          \end{array}
                        \right.
\nonumber
\end{align}
and
$$\infty>\mm{Dis}(\vartheta,\kappa,\bar{M}):=\left\{
                          \begin{array}{ll}
\sup_{0\neq w\in H_{\sigma,\vartheta}^1} { g\llbracket\rho\rrbracket |w_3|_{0}^2}/{ \mathcal{I}(w,\vartheta,\kappa,\bar{M})}  &\mbox{ for }(\kappa,\bar{M}_3)\neq 0;\\
\sup_{ w\in H_{\sigma,3}^1} { g\llbracket\rho\rrbracket |w_3|_{0}^2}/\vartheta|\nabla_{\mm{h}}w_3|_{0}^2 & \mbox{ for }(\kappa,\bar{M}_3)=0.     \end{array}
                        \right.$$
Thus, we  see that  the instability condition \eqref{201806071047} is equivalent to
$$
1<\mm{Dis}(\vartheta,\kappa,\bar{M}) \mbox{ for } (\vartheta,\kappa,\bar{M}_3)\neq 0.
 $$
Of course, if $ (\vartheta,\kappa,\bar{M}_3)= 0$, we always have the instability condition \eqref{201806071047} (the case of $(\vartheta,\kappa,\bar{M}_3)= 0$ with $\bar{M}_{\mm{h}}\neq0$ can be observed from Lemma \ref{201806151836}).
\end{rem}
\begin{rem}
By modifying the derivation of Growall-type energy inequality \eqref{2016121521430850} and multi-energy method in \cite{WYTIVNMI}, we can verify the global (asymptotic) stability of the ART problem  under the stability condition $\mm{Dis}(\vartheta,\kappa,\bar{M})<1$ with $(\kappa,\bar{M}_3)\neq 0$. We will present the verification of the stability result in a separate article.
\end{rem}

The proof of Theorem \ref{thm3} is based on a so-called bootstrap instability method. The bootstrap instability method has its origin in \cite{GYSWIC,GYSWICNonlinea}. Later, various versions of bootstrap approaches were presented by many authors, see \cite{FSSWVMNA,GYHCSDDC} for examples. In this article, we adapt the version of bootstrap instability method in \cite[Lemma 1.1]{GYHCSDDC} to prove Theorem \ref{thm3} by further introducing some new mathematical techniques. Our proof procedure can be divided into five steps. Firstly, we shall construct unstable solutions to the linearized ART problem, this can be achieved by a modified variational method as in \cite{GYTI2,JFJSO2014} and an existence theory of the stratified (steady) Stokes problem due to the presence of viscosity, see Proposition \ref{thm:0201201622}. Secondly, using energy method as in \cite{XLZPZZFGAR} with new auxiliary estimates in Sobelve--Sobodetskii spaces (see Lemmas \ref{201807312051}  and  \ref{20180812}), we can prove that the local-in-time solution of the ART problem enjoys some Gronwall-type energy inequality, which couples with an auxiliary solution of the linear problem  \eqref{s0106pndfsgnnnxx}, see Lemma \ref{pro:0301n0845}. It should be noted that our constructed Gronwall-type energy inequality is very different with the standard form in \cite[Lemma 1.1]{GYHCSDDC}. Thirdly, we want to use initial data of solutions of the linearized ART problem to construct initial data for solutions of the (nonlinear) ART problem as in \cite[Lemma 1.1]{GYHCSDDC}. However the initial data of the linearized and corresponding nonlinear ART problems have to satisfy different compatibility conditions. Thus we can not directly use the initial data of linearized ART problem as the one of nonlinear ART problem as in \cite[Lemma 1.1]{GYHCSDDC}. To circumvent this difficulty, we use existence theory of  stratified Stokes problem and iterative technique to modify initial data of solutions of the linearized ART problem, so that the obtained modified initial data belong to $H^{4,1}_{0,*}\times H^2_0\times H^1$, satisfy the compatibility condition  \eqref{201809101928}--\eqref{xx201809101926xxsdfs}, and are close to the initial data of the linearized ART problem, see Propositions \ref{lem:modfied}--\ref{201809012320}.
Fourthly, we deduce the error estimates between the solutions of the linearized and nonlinear ART problems based on a method of largest growth rate as in \cite{HHJGY}.
Noting that the error function $u^{\mm{d}}$ (i.e. the error between the nonlinear and linear solutions concerning $u$) does not enjoys divergence-free condition, thus we shall use existence theory of stratified Stokes problem again to modify
$u^{\mm{d}}$,  so that the method of largest growth rate can be used to deduce
the desired error estimates. Finally, we prove the existence of an escape time and thus obtain Theorem \ref{thm3}. The detailed proof of Theorem \ref{thm3} will be provided in Sections  \ref{201805302028}--\ref{sec:030845}.

Thanks to Theorem \ref{thm3}, we can use an inverse transformation of Lagrangian coordinates to establish the following theorem, which will implies some instability results of RT problem in viscoelastic fluids, magnetohydrodynamics (MHD) fluids with zero resistivity  and pure viscous fluids (with or without interface intension) in Eulerian coordinates by taking proper values of parameters.
\begin{thm}\label{201806012301xx}
Let $\lambda\geqslant 0$, $(\eta,u,q)$ be constructed in Theorem \ref{thm3} with sufficiently small $\delta$  and $\sqrt{\lambda}\bar{M}$ in place of $\bar{M}$, and $\zeta:=\eta+y$. We define that, for any $t\in \overline{I_T}$,
\begin{align}
& \label{usedopso}
\left\{
    \begin{array}{l}
d:=\eta_3((\zeta_{\mm{h}})^{-1}(x_{\mm{h}},t),0,t) \in (-l,\tau),\  {v}:={u}({\zeta}^{-1}(x,t),t),\\
V:=\nabla\zeta({\zeta}^{-1}(x,t),t)-I,\ N:= \partial_{\bar{M}}\eta,\  \sigma:=q({\zeta}^{-1}(x,t),t),\\
\nu:=(-\partial_1 d,-\partial_2 d,1 )^{\mm{T}}/\sqrt{1+|\nabla_{x_{\mm{h}}} d|^2},
    \end{array}
  \right. \\
&
\Sigma(t):=\{(x_{\mm{h}},x_3)~|~x_{\mm{h}}\in \mathbb{T},\ x_3:=d(x_{\mm{h}},t)\},\label{201809221932}\\
&\Omega_+(t):=\{(x_{\mm{h}},x_3)~|~x_{\mm{h}}\in \mathbb{T}  ,\ d(x_{\mm{h}},t)< x_3<\tau\},\label{201809221932xx}\\
&\Omega_-(t):=\{(x_{\mm{h}},x_3)~|~x_{\mm{h}}\in \mathbb{T} ,\ -l< x_3<d(x_{\mm{h}},t)\}, \label{201809221933} \\
 &  \Omega_\pm'(t):= \zeta_\pm(\Omega_\pm',t),\  \Omega(t):=\Omega_+(t)\cup \Omega_-(t),\\
& \Omega^{\mm{p}}(t):=\{(x_{\mm{h}},x_3)~|~x_{\mm{h}}\in \mathcal{T},\ -l<x_3<\tau,\  x_3 \neq d(x_{\mm{h}},t) \}, \label{201808160832}\\
&  (d^0,v^0,V^0,N^0,\sigma^0,\nu^0):=(d,v,V,N,\sigma,\nu )|_{t=0}.
\end{align}

Then we have the following conclusions:
\begin{enumerate}[\quad (1)]
\item for each $t\geqslant 0$,
\begin{align}
&(v,V, N, \sigma)(t)\in H_\sigma^2(\Omega(t))\times  H^2(\Omega(t)) \times H^2(\Omega(t)) \times H^1(\Omega(t)), \label{201808160842} \\
& (v_t ,V_t,N_t)\in L^2(\Omega(t))\times H^1(\Omega(t))\times H^1(\Omega(t)),\label{201819292017} \\
& d \in H^{5/2}(\mathbb{T}),\  d_t,\ \nu(x_{\mm{h}},t)\in H^{3/2}(\mathbb{T}), \label{201809271720},
\end{align}
where we have defined that
\begin{align}
& H^k(\Omega(t)):=\{\omega~|~\omega_\pm:=\omega|_{\Omega_\pm'(t)}\in H^k_{\mm{loc}}(\Omega_\pm'(t)),\ \omega(y_1,y_2,y_3)=\nonumber\\
&\qquad \qquad\qquad  \omega(y_1+2\pi m L_1,y_2 +2\pi m L_1,y_3)\mbox{ for any integers }m\mbox{ and }n \},  \label{201809222025} \\
&  H_\sigma^k(\Omega(t)):=H_\sigma^1\cap H^k(\Omega(t)). \label{201809222025x}
\end{align}
Moreover,
\begin{align}
& \nabla_{\mm{h}}d,\ \nu\mbox{ are continuous on }\overline{\mathbb{R}^2_T}, \label{2018090414874}\\
& v\mbox{ is continuous on } \overline{\Omega^T},\ N_\pm,\  V_\pm\mbox{ are continuous on }\tilde{\zeta}_\pm(\overline{\Omega_\pm^T}) \label{2018090714874},
\end{align}
where $\tilde{\zeta}_\pm(y,t):=(\zeta_\pm(y,t),t)$. In addition,
\begin{equation}
d \in C(\overline{I_T},H^{3}(\mathbb{T}))\cap L^2(I_T,H^{7/2}(\mathbb{T})),\ d_t\in L^2(I_T,H^{5/2}(\mathbb{T}))\mbox{ if }\vartheta\neq 0. \label{201811141303}
\end{equation}
  \item $(v,V,N,\sigma, d)$ solves a so-called mixed RT problem:
\begin{equation}\label{0103nxx}\left\{{\begin{array}{ll}
\rho  v_t+\rho v\cdot\nabla v+\mm{div}\mathcal{S}^{\mm{V},\mm{M}}(\sigma, v,V +I,N+\bar{M})=0 &\mbox{ in }  \Omega(t),\\[1mm]
 V_t + v\cdot\nabla V=\nabla v (V +I)&\mbox{ in } \Omega(t),\\[1mm]
N_t+v\cdot\nabla N=\nabla v (N+\bar{M})^{\mm{T}} &\mbox{ in } \Omega(t),
 \\[1mm]
\mathrm{div} v=\mathrm{div} N= 0& \mbox{ in } \Omega(t),\\
 d_t+v_1 \partial_1d+v_2 \partial_2d=v_3 &\mbox{ on }\Sigma(t),\\
   \llbracket \mathcal{S}^{\mm{V},\mm{M}}(\sigma,v,V +I,N+\bar{M})\nu -g \rho d \nu\rrbracket =\vartheta \mathcal{C} \nu  \ &\mbox{ on }\Sigma(t) ,  \\
  \llbracket v  \rrbracket=0&\mbox{ on }\Sigma,\\
v =0 &\mbox{ on }\Sigma_-^+,\\
(v,V,N)|_{t=0}=(v^0,V^0,N^0) &\mbox{ in } \Omega(0),\\
d|_{t=0}=d^0  &\mbox{ on } \Sigma(0),
\end{array}}\right.\end{equation}
where
$ \mathcal{S}^{\mm{V},\mm{M}}(\sigma,v,V +I,N+\bar{M}):=p^g I-\mu \mathbb{D}  v   - { {\kappa}\rho}(UU^{\mm{T}}-I) -\lambda(N+\bar{M})\otimes(N+\bar{M})$.
\item $(v,V,N, d)$ satisfies the instability relations:
\begin{align}
&  |d(T^\delta)|_{L^1} ,\ \|\varpi(T^\delta)\|_{L^1_\delta},\ \|\partial_3\varpi (T^\delta)\|_{L^1_\delta},\ \|\mm{div}_{\mm{h}}v_{\mm{h}}(T^\delta)\|_{L^1_\delta}\geqslant {\epsilon},\label{201808160850} \\
&\|V_{33}(T^\delta)\|_{L^1_\delta},\  \|(V_{13},V_{23})(T^\delta)\|_{L^1_\delta} ,\  \|(V_{11}+V_{22})(T^\delta)\|_{L^1_\delta} \geqslant {\epsilon},\\
 &\|N_3(T^\delta)\|_{L^1_\delta},\  \|N_{\mm{h}}(T^\delta)\|_{L^1_\delta}\geqslant  \epsilon,\mbox{ if }
\bar{M}_3\neq0
\label{201805201416}
\end{align}
for some escape time $T^\delta\in I_T$, where $\varpi=v_3$ or $v_{\mm{h}}$, and $L^1_\delta:=L^1(\Omega^{\mm{p}}(T^\delta))$.
 \item the initial data $(v^0,V^0,N^0,\sigma^0,d^0,\nu^0) $ satisfies
\begin{equation}
 \|(v^0,V^0,N^0)\|_2+\| \sigma^0\|_1+|d^0|_{5/2}+\beta(|d^0|_{3})\leqslant C \delta \label{2018090414874xx}
\end{equation}
and the  compatibility conditions:
\begin{alignat}{2}
&\mm{div} v^0=0 & & \ \mbox{ in }\Omega(0),\label{2018091019281x} \\
&\llbracket \mathcal{S}^{\mm{V},\mm{M}}(\sigma^0,v^0,V^0+I,N+\bar{M} )\nu^0-g
\rho d^0\nu^0\rrbracket = \vartheta \mathcal{C}^0 \nu^0& &\ \mbox{ on }\Sigma(0),\label{xx2018091019261x}
\end{alignat}
where the constant $C$ depends on the norms of functions $\tilde{\eta}^0$, $\tilde{u}^0$, $\tilde{q}^0$,
$\eta^\mm{r}$, $u^\mm{r}$ and $q^\mm{r}$ in Theorem \ref{thm3}.
\end{enumerate}
\end{thm}

\begin{rem}\label{rem:201810011518}
It should be noted that, since  $\eta(t)\in H^{4,1}_{0,*}$ for each $t\in \overline{I_T}$,  the domains defined in \eqref{201809221932}--\eqref{201809221933} automatically satisfy, for each $t\geqslant 0$,
\begin{align}\Omega\!\!\!\!-=\Omega (t)\cup \Sigma(t),\ \Omega_+(t)\cap \Omega_-(t)=\emptyset, \mbox{ and }\Omega_\pm(t)\cap \Sigma(t)=\emptyset.\nonumber
\end{align}
In addition, by the regularity of $(\eta,u)$, we further  derive many  additional conclusions, for examples,
\begin{align}&
\tilde{\zeta} (y,t) :\overline{\Omega ^T}\to  \overline{\Omega ^T},\ \tilde{\zeta}_\pm(y,t):\overline{ \Omega_\pm^T} \to \tilde{\zeta}_\pm( \overline{\Omega_\pm^T})\nonumber \\
&\mbox{ and }\bar{\zeta}(y_{\mm{h}},t): \overline{\mathbb{R}^2_{T}} \to   \overline{\mathbb{R}^2_{T}} \mbox{ are homeomorphism mappings},\label{201808151856} \\
&\tilde{\zeta}_\pm(y,t): \Omega_\pm^T \to \tilde{\zeta}_\pm( {\Omega_\pm^T})\mbox{ and } \bar{\zeta}(y_{\mm{h}},t):  {\mathbb{R}^2_{T}} \to  {\mathbb{R}^2_{T}}\nonumber \\
&\mbox{are }C^1\mbox{-diffeomorphic mappings},\label{201808151856sdfa}
\end{align}
where we have defined $\tilde{\zeta} (y,t):=(\zeta (y,t),t)$ and $\bar{\zeta}(y_{\mm{h}},t):=(\zeta_{\mm{h}}(y_{\mm{h}},0,t),t)$.
\end{rem}

The rest article are organized as follows: In Section \ref{201806101821}, by taking proper values of parameters in Theorem \ref{201806012301xx}, we further establish existence results of instability solutions for the RT problems in stratified viscoelastic fluids, stratified MHD fluids, and  stratified  pure fluids (with and without surface tension) in sequence, see Corollaries \ref{201806012301xxxxxxxx}--\ref{20180610185}. In  Sections \ref{201805302028}--\ref{sec:030845}, we will provide the detailed proof of Theorem \ref{thm3}.  Finally, we will introduce how to use Theorem \ref{thm3} to establish Theorem \ref{201806012301xx} in Section \ref{201806012303}.

\section{Applications}\label{201806101821}

In this section, we will take proper values of parameters in Theorem \ref{201806012301xx}  to  obtain the existence results of RT instability solutions of RT problems in stratified viscoelastic fluids,  stratified  MHD fluids, and   stratified  pure fluids (with or without interface intension) in Eulerian coordinates. We will also briefly introduce relevant mathematical progress of RT problems in the three types of stratified viscous fluids.

\subsection{Stratified viscoelastic fluids}\label{201809231805}

The RT instability was often investigated in various models of viscoelastic fluids from the physical point of view, see \cite{SHKCSHAR,BGMAMSVLRTI} for examples. It is well-known that viscoelasticity is a material property that exhibits both viscous and elastic characteristics when undergoing deformation. In particular, an elastic fluid strains when stretched and quickly return to its original state once the stress is removed. So the elasticity will have a stabilizing effect like internal surface tension.  Recently, Jiang--Jiang--Wu have proved that sufficiently large elasticity coefficient can contribute to the elasticity to inhibit RT instability by a mathematical model of stratified viscoelastic fluids \cite{JFJWGCOSdd}. In this article, we further provide the mathematical result of RT instability in viscoelastic fluids under an instability condition \eqref{201806071047} with $\bar{M}=0$.  To begin with, let us recall the  RT problem of stratified viscoelastic fluids in a domain $\Omega$ in \cite{JFJWGCOSdd}:
\begin{equation}\label{0101f}\left\{{\begin{array}{ll}
\rho_\pm \partial_t v_\pm+\rho_\pm v_\pm\cdot\nabla v_\pm+\mm{div}\mathcal{S}_\pm^{\mm{V}} (p_\pm^g,v_\pm, U_\pm)=0&\mbox{ in } \Omega_\pm(t),\\[1mm]
 \partial_t U_\pm+v_\pm\cdot\nabla U_\pm=\nabla v_\pm U_\pm&\mbox{ in } \Omega_\pm(t),\\[1mm]
\mathrm{div} {v}_\pm= 0& \mbox{ in } \Omega_\pm(t),\\
d_t+v_1 \partial_1d+v_2 \partial_2d=v_3 &\mbox{ on }\Sigma(t),\\
   \llbracket v_\pm  \rrbracket =0,\
   \llbracket   \mathcal{S}_\pm^{\mm{V}}( {p}^g_\pm, v_\pm,U_\pm)\nu- g d \rho_\pm \nu        \rrbracket =\vartheta \mathcal{C} \nu   &\mbox{ on }\Sigma(t) ,\\
v_\pm=0 & \mbox{ on }\Sigma_\pm,\\
(v_\pm, U_\pm)|_{t=0}=(v^0_\pm,U^0_\pm)  &\mbox{ in } \Omega_\pm (0),\\
d|_{t=0}=d^0  &\mbox{ on } \Sigma(0) .   \end{array}}  \right.
\end{equation}
Next we should explain the notations in the above (stratified) VRT problem \eqref{0101f}.

For each given $t>0$, $d:=d(x_{\mm{h}},t):\mathbb{T}\mapsto (-l,\tau)$ is a height function of a point at the interface of stratified (viscoelastic) fluids. $\Sigma(t)$ is a set of interface, and defined  by \eqref{201809221932}.

The subscripts $+$ resp. $-$ in the notation $f_+$ resp. $f_-$ mean that the functions or parameters  $f_+$ resp. $f_-$  are relevant to the upper resp. lower fluids.
$\Omega_+(t)$ and $\Omega_-(t)$ are the domains of upper and lower fluids, respectively, and defined by \eqref{201809221932xx} and \eqref{201809221933}.

For given $t>0$, $v_\pm(x,t): \Omega_\pm(t)\mapsto \mathbb{R}^3$,  $U_\pm(x,t): \Omega_\pm(t)\mapsto \mathbb{R}^9$ and $p_\pm(x,t): \Omega_\pm(t)\mapsto \mathbb{R}$ are the velocities, the deformation tensors ($3\times 3$ matrix-valued function) and the pressures of fluids.
$\rho_\pm$, $\kappa_\pm$ and $\mu_\pm$ are the density constants, viscosity coefficients and the elasticity coefficients of  fluids, where $\rho_+>\rho_-$. $g$ is the gravity constant.  $\vartheta$ the surface tension coefficient. $\mathcal{S}_\pm^{\mm{V}}(p_\pm^g, v_\pm, U_\pm):=p_\pm^g I-\mu_\pm \mathbb{D} v_\pm  - { {\kappa_\pm}\rho_\pm }(U_\pm U^{\mm{T}}_\pm-I)$, where ${p}_\pm^g: =p_\pm + g\rho_\pm x_3$.

For a function $f$ defined on $\Omega(t)$, we define $\llbracket  f_\pm  \rrbracket:=f_+|_{\Sigma(t)}
-f_-|_{\Sigma(t)}$, where $f_\pm|_{\Sigma(t)}$ are the traces of the quantities $f_\pm$ on $\Sigma(t)$.
 $\nu$ is the unit outer normal vector at boundary $\Sigma(t)$ of $\Omega_-(t)$, and $\mathcal{C}$ the twice of the mean curvature of the internal  surface $\Sigma(t)$ \cite{WRGMDSFA}, i.e.,
$$
\mathcal{C}:=\frac{\Delta_{\mm{h}}d+(\partial_1 d)^2\partial_2^2
d+(\partial_2d)^2\partial_1^2d-2\partial_1d\partial_2d\partial_1\partial_2d}
{(1+(\partial_1d)^2+(\partial_2d)^2)^{3/2}}.
$$

Finally, we briefly explain the physical meaning of each identity in \eqref{0101f}.
The equations \eqref{0101f}$_1$--\eqref{0101f}$_2$ describe the motion of the upper heavier and lower lighter fluids driven by the gravitational field along the negative $x_3$-direction, which occupy the two time-dependent disjoint open subsets $\Omega_+(t)$ and $\Omega_-(t)$ at time $t$, respectively.
We call  \eqref{0101f}$_1$ the momentum equation, and \eqref{0101f}$_2$ the deformation equation.
Since the fluids are incompressible, we naturally pose the divergence-free condition \eqref{0101f}$_3$.
The two fluids interact with each other by the motion equation of a free interface \eqref{0101f}$_4$ and the interfacial jump conditions in \eqref{0101f}$_5$. The first jump condition in \eqref{0101f}$_5$ represents that
the velocity is continuous across the interface.
The second jump in \eqref{0101f}$_5$ is equivalent to
$$   \llbracket   \mathcal{S}_\pm^{\mm{V}} (p_\pm, v_\pm,U_\pm )   \rrbracket \nu=\vartheta \mathcal{C} \nu \mbox{ on }\Sigma(t),$$
which represents that the jump in the normal stress
is proportional to the mean curvature of the surface multiplied by the normal to the surface \cite{LeMHVJ,XLZPZZFGAR}. The non-slip boundary condition of the velocities on both upper and lower fixed flat boundaries are described by \eqref{0101f}$_6$, and \eqref{0101f}$_7$--\eqref{0101f}$_8$ represent the initial status of the two fluids.

The problem \eqref{0101f} enjoys an equilibrium state (or rest) solution: $(v,U,{p}^g, d)=(0,I,\bar{p}^{g},\bar{d})$, where $\bar{d}\in (-l,\tau)$. We should point out that $\bar{p}^g$ can be uniquely computed out by hydrostatics, which depends on the variable $x_3$ and $\rho_\pm$, and is continuous with respect to $x_3\in (-l,\tau)$. Without loss of generality, we assume that $\bar{d}=0$ in this article.
If $\bar{d}$ is not zero, we can adjust the $x_3$ co-ordinate to make $\bar{d}=0$.
Thus $d$ can be regarded as the displacement away from the plane $\Sigma$.

To simply the representation of the problem \eqref{0101f}, we introduce the indicator function $\chi_{\Omega_\pm(t)}$ and denote
$$
\begin{aligned}
&\rho=\rho_+\chi_{\Omega_+(t)} +\rho_-\chi_{\Omega_-(t)},\
 \mu=\mu_+\chi_{\Omega_+(t)} +\mu_-\chi_{\Omega_-(t)},\
\kappa=\kappa_+\chi_{\Omega_+(t)} +\kappa_-\chi_{\Omega_-(t)}, \\
& v=v_+\chi_{\Omega_+(t)} +v_-\chi_{\Omega_-(t)},\
U=U_+\chi_{\Omega_+(t)} +U_-\chi_{\Omega_-(t)},\
p=p_+\chi_{\Omega_+(t)} +p_-\chi_{\Omega_-(t)}, \\
& v^0=v_+^0\chi_{\Omega_+(0)} +v_-^0\chi_{\Omega_-(0)} ,\ U^0=U_+^0\chi_{\Omega_+(0)} +U_-^0\chi_{\Omega_-(0)},\\
  & \mathcal{S}^{\mm{V}}(p^g, v,U)  :=p^g I-\mu \mathbb{D}  v   - { {\kappa}\rho}(UU^{\mm{T}}-I)  .
\end{aligned}$$
Now, we denote the perturbation  quantity around the equilibrium state $(0,0,I,\bar{p}^g)$ by
$$d=d-0,\ v=v-0,\ V=U-I \mbox{ and }\sigma={p}^g-\bar{p}^g.$$
Then we have a VRT problem in perturbation form:
\begin{equation}\label{0103nxxxxx}\left\{{\begin{array}{ll}
\rho  v_t+\rho v\cdot\nabla v+\mm{div}\mathcal{S}^{\mm{V}}(\sigma, v,V +I)=0 &\mbox{ in }  \Omega(t),\\[1mm]
 V_t + v\cdot\nabla V=\nabla v (V +I)&\mbox{ in } \Omega(t),\\[1mm]
\mathrm{div} v= 0& \mbox{ in } \Omega(t),\\
d_t+v_1 \partial_1d+v_2 \partial_2d=v_3 &\mbox{ on }\Sigma(t),\\
   \llbracket \mathcal{S}^{\mm{V}}(\sigma,v,V +I)\nu -g \rho d \nu\rrbracket =\vartheta \mathcal{C} \nu  \ &\mbox{ on }\Sigma(t) ,  \\
  \llbracket v  \rrbracket=0&\mbox{ on }\Sigma,\\
v =0 &\mbox{ on }\Sigma_-^+,\\
(v,V)|_{t=0}=(v^0,V^0) &\mbox{ in } \Omega(0),\\
d|_{t=0}=d^0  &\mbox{ on } \Sigma(0),
\end{array}}\right.\end{equation}
where we have defined that $\Omega(t):=\Omega_+(t)\cup\Omega_-(t)$, and omitted the subscript $\pm$ in $\llbracket  f_\pm  \rrbracket$ for simplicity. Thus a zero solution is an equilibrium-state solution of the above perturbation VRT problem.

It is well-known that the movement of the free interface $\Sigma(t)$ and the subsequent change of the
domains $\Omega_\pm(t)$ in Eulerian coordinates will result in severe mathematical difficulties. There exist three methods to circumvent such difficulties, see the transformation method of Lagrangian coordinates in \cite{GYTI2}, the method of translation transformation with respect to $x_3$-variable in \cite{PJSGOI5x}, and ``geometric" reformulation method in \cite{GYTIAE2}. Next we shall adopt the transformation method of Lagrangian coordinates so that the interface and the domains stay fixed in time. In addition, the VRT problem  in Lagrangian coordinate has better mathematical structure.

We  assume that there exist invertible mappings
\begin{equation*}\label{0113}
\zeta_\pm^0:\Omega_\pm\rightarrow \Omega_\pm(0),
\end{equation*}
such that
\begin{equation}
\label{05261240}
\Sigma(0)=\zeta_\pm^0(\Sigma),\
\Sigma_\pm=\zeta_\pm^0(\Sigma_\pm)
\end{equation}
and
\begin{equation}
\label{zeta0inta} \det \nabla\zeta_\pm^0= 1.
\end{equation}
We further define $\zeta:=\zeta_+\chi_{\Omega_+}+\zeta_-\chi_{\Omega_-}$, and the flow maps $\zeta$ as the solutions to
\begin{equation}
\left\{
            \begin{array}{ll}
\partial_t \zeta(y,t)=v(\zeta(y,t),t)&\mbox{ in }\Omega
\\
\zeta(y,0)=\zeta^0(y)&\mbox{ in }\Omega.
                  \end{array}    \right.
\label{201811191443}
\end{equation}
We denote the Eulerian coordinates by $(x,t)$ with $x=\zeta(y,t)$,
whereas the fixed $(y,t)\in \Omega\times \mathbb{R}^+$ stand for the
Lagrangian coordinates.

In order to switch back and forth from Lagrangian to Eulerian coordinates, we shall assume that
$\zeta_\pm(\cdot ,t)$ are invertible and $\Omega_{\pm}(t)=\zeta_{\pm}(\Omega_{\pm},t)$, and since $v_\pm$ and $\zeta_\pm^0$ are all continuous across $\Sigma$, we have $\Sigma(t)=\zeta_\pm(\Sigma,t)$. In view of the non-slip boundary condition $v|_{\Sigma_-^+}=0$, we have
$$y=\zeta(y, t)\mbox{ on }\Sigma_-^+.$$
In addition, by the incompressible condition, we have
$\det \nabla \zeta =1$ in $\Omega$
as well as the initial condition \eqref{zeta0inta},  see \cite[Proposition 1.4]{MAJBAL}.

Now, setting $\eta=\zeta-y$ and the Lagrangian unknowns
$$
(u,\tilde{U},q)(y,t)=(v,U,\sigma)(\zeta(y,t),t)\;\;\mbox{ for } (y,t)\in \Omega \times\mathbb{R}^+,
$$
then in Lagrangian coordinates $(\eta,u,q)$ satisfies a so-called transformed VRT problem (please refer to Section \ref{201806012303} for omitted derivation):
\begin{equation}\label{n0101safsnn}\left\{\begin{array}{ll}
\eta_t=u &\mbox{ in } \Omega,\\[1mm]
{\rho} u_t+\mm{div}_{\ml{A}}\mathcal{S}_{\mathcal{A}}(q,u,\eta ) =0&\mbox{ in }  \Omega,\\[1mm]
\mathrm{div}_{\mathcal{A}} {u} = 0& \mbox{ in }  \Omega ,\\
  \llbracket \eta  \rrbracket=\llbracket u  \rrbracket=0,\ \llbracket  \mathcal{S}_{\mathcal{A}} (q,u,\eta)-g
\rho  \eta_3I \rrbracket\mathcal{A}e_3=\vartheta H \mathcal{A}e_3 &\mbox{ on }\Sigma,\\
(\eta, u)=0 &\mbox{ on }\Sigma_-^+,\\
(\eta,u)|_{t=0}=(\eta^0,u^0)  &\mbox{ in }  \Omega
\end{array}\right.\end{equation}
with $\tilde{U} (y,t)=\nabla \zeta_\pm(y,t)$.
We call $\eta$ the displacement function of particle (labeled by $y$). Obviously,
the above transformed VRT problem can be deduced from the ART problem \eqref{n0101nn} by taking $\bar{M}=0$.

To focus on the elasticity effect upon the RT instability, Jiang--Jiang--Wu
omitted the surface tension in the transformed VRT problem, and  proved that
the transformed VRT problem with $\vartheta=0$ is stable \cite{JFJWGCOSdd} under the
stability condition \begin{equation*}
\mm{Dis}(0,\kappa,0)<1.
\end{equation*}
We mention that, using the energy methods in \cite{JFJWGCOSdd} and Section \ref{201805081731}, we can also verify that the transformed VRT problem is stable under the stability condition $\mm{Dis}(\vartheta,\kappa,0)<1$.

Noting that the instability condition \eqref{201806071047} with $\bar{M}=0$ and with $(\vartheta,\kappa)\neq 0$ is equivalent to
\begin{equation}
\label{201806142048}
\mm{Dis}(\vartheta,\kappa,0)>1,
\end{equation}
thus Theorem \ref{thm3} presents that the transformed VRT problem is unstable, when $(\vartheta,\kappa)= 0$, or  the instability condition \eqref{201806142048} with $(\vartheta,\kappa)\neq 0$ is satisfied.
Moreover, by Theorem \ref{201806012301xx}, we also obtain the instability result of the perturbation VRT problem \eqref{0103nxxxxx}.
\begin{cor}
\label{201806012301xxxxxxxx}
Under the instability condition \eqref{201806142048} with
 $(\vartheta,\kappa)\neq 0$ or $(\vartheta,\kappa)= 0$, the perturbation VRT problem \eqref{0103nxxxxx}  is unstable.
More precisely,  $(d,v,V,   \sigma)$ constructed by Theorem \ref{201806012301xx} with $\lambda=0$ is an unstable solution to the perturbation VRT problem.
\end{cor}

\subsection{Stratified  MHD fluids}\label{201809232038}

Now we turn to introduce the instability results in the RT problem of stratified  MHD fluids (with zero resistivity).
The topic of the effects of  magnetic fields upon the RT instability
in stratified MHD fluids goes back to the theoretical work of Kruskal and Schwarzchild \cite{KMSMSP}
in 1954. They analyzed the effect of an impressed horizontal magnetic field on the growth of
the RT instability, and pointed out that the curvature of the magnetic lines can influence the development of instability in compressible MHD fluids, but can not inhibit the RT instability growth. The effect of the impressed  vertical magnetic field upon the RT instability was also investigated by Hide for incompressible inhomogeneous MHD fluids \cite{HRWP,CSHHSCPO}. Since Kruskal and Schwarzchild's pioneering work, many physicists proceeded to develop the linear theory and nonlinear numerical simulation of the magnetic RT instability, due to its applications in a range of scales from laboratory plasma physics to atmospheric physics and astrophysics. In particular,  Hester et al. \cite{HJJSJMSPAWFPC} compared the observational characteristics of the Crab Nebula with simulations of the magnetic RT instability performed by Jun et al. \cite{JBINMLSJMA}, and found that the magnetic RT
instability could explain the observed filamentary structure, which is also observed in the Sun \cite{IHMTSKYTFstru,IHMTSKYTT}.

The (nonlinear) magnetic RT instability was first mathematically proved by Hwang in inhomogeneous inviscid MHD fluids in Eulerian coordinates in 2008. Jiang et.al. further considered the corresponding viscous case \cite{JFJSWWWN}. However, to our best knowledge, up to now, there are not any available result of mathematical proof of magnetic RT instability in the stratified  MHD fluids  in  Eulerian coordinates. Therefore, in this article, we further use Theorem \ref{201806012301xx} to establish the existence result of RT instability solutions for the following RT problem in stratified incompressible viscous MHD fluids under the instability condition \eqref{201806071047} with $\kappa=0$ in  Eulerian coordinates:
\begin{equation}\label{0101fdfadsfsd}\left\{\begin{array}{ll}
\rho_\pm \partial_t v_\pm+\rho_\pm v_\pm\cdot\nabla v_\pm+\mm{div}\mathcal{S}_\pm^{\mm{M}} (p_\pm^{g,M_\pm}, v_\pm,M_\pm)=0&\mbox{ in } \Omega_\pm(t),\\[1mm]
 \partial_t M_\pm=M_\pm\cdot \nabla v_\pm-v_\pm\cdot\nabla M_\pm &\mbox{ in } \Omega_\pm(t),\\[1mm]
\mathrm{div} {v}_\pm= 0& \mbox{ in } \Omega_\pm(t),\\
\mm{div} M_\pm=0&\mbox{ in } \Omega_\pm(t),\\[1mm]
d_t+v_1 \partial_1d+v_2 \partial_2d=v_3 &\mbox{ on }\Sigma(t),\\
   \llbracket v_\pm  \rrbracket =0,\quad
   \llbracket   \mathcal{S}_\pm^{\mm{M}} ( p^{g,M_\pm}_\pm,   v_\pm,M_\pm) \nu -g d \rho   \nu  \rrbracket =  \vartheta \mathcal{C} \nu   &\mbox{ on }\Sigma(t) ,\\
v_\pm=0 & \mbox{ on }\Sigma_\pm,\\
(v_\pm,M_\pm)|_{t=0}=(v^0_\pm,M^0_\pm) &\mbox{ in } \Omega_\pm (0),\\
d|_{t=0}=d^0  &\mbox{ on } \Sigma(0) .   \end{array}  \right.
\end{equation}
Next we shall explain new equations and  notations in the above (stratified) MRT problem.

The equations \eqref{0101fdfadsfsd}$_1$ and \eqref{0101fdfadsfsd}$_2$ describe the motion of the upper heavier and lower lighter MHD fluids with zero resistivity driven by the gravitational field along the negative $x_3$-direction.  \eqref{0101fdfadsfsd}$_2$ is called the induction equation. We remark that the resistivity is neglected in \eqref{0101fdfadsfsd}$_2$, and this arises in the physics regime with negligible electrical resistance. \eqref{0101fdfadsfsd}$_4$ represent that the magnetic fields $\bar{M}_\pm$ of the two fluids are source free, and can be seen just as restrictions on the initial value of $M_\pm$, since $(\mm{div}M_\pm)_t=0$ by \eqref{0101fdfadsfsd}$_2$. $\mathcal{S}_\pm^{\mm{M}} (p_\pm^{g,M_\pm}, v_\pm,M_\pm):= p_\pm^{g,M_\pm}I-\mu_\pm \mathbb{D}v_\pm-\lambda M_\pm \otimes M_\pm$, and $p_\pm^{g,M_\pm}:= p_\pm^{g}+\lambda  |M_\pm|^2/2$, where  $\lambda$ stands for the permeability of vacuum divided by $4\pi$.

The problem \eqref{0101fdfadsfsd} enjoys an equilibrium state solution: $(v,M,{p}^g,d)=(0,\bar{M},\bar{p}^{g},0)$, where $\bar{M}$ is a constant vector, and called an impressed field (or equilibrium magnetic field). Following the argument of \eqref{0103nxxxxx}, and further  defining that
$$
\begin{aligned}
& M:=M_+\chi_{\Omega_+(t)} +M_-\chi_{\Omega_-(t)},\ N:=M-\bar{M},\\
 & p^{g,M}:= p^{g}+\lambda  |M|^2/2,\ \sigma:= p^{g,M} - \bar{p}^g-\lambda|\bar{M}|^2/2,\\
& \mathcal{S}^{\mm{M}}(p^{g,M},v,M):=\mathcal{S}_+^{\mm{M}} (p_+^{g,M_+}, v_+,M_+)\chi_{\Omega_+(t)}+\mathcal{S}_+^{\mm{M}} (p_-^{g,M_-}, v_-,M_-)\chi_{\Omega_-(t)},
\end{aligned}$$
we can get a MRT problem in perturbation form:
\begin{equation}\label{0103nxasdfx}\left\{{\begin{array}{ll}
\rho  v_t+\rho v\cdot\nabla v+\mm{div}\mathcal{S}^{\mm{M}}(\sigma,v,N+\bar{M})=0   &\mbox{ in }  \Omega(t),\\[1mm]
N_t=(N+\bar{M})\cdot \nabla v-v\cdot\nabla N &\mbox{ in } \Omega(t),
\\
\mathrm{div} v=\mathrm{div} N= 0& \mbox{ in } \Omega(t),\\[1mm]
d_t+v_1 \partial_1d+v_2 \partial_2d=v_3 &\mbox{ on }\Sigma(t),\\
   \llbracket v  \rrbracket=0,\quad
   \llbracket \mathcal{S}^{\mm{M}} (\sigma,v ,N+\bar{M})    \nu-
 g  \rho d  \nu   \rrbracket = \vartheta \mathcal{C} \nu   &\mbox{ on }\Sigma(t) ,
 \\
v =0 &\mbox{ on }\Sigma_-^+,\\
(v,N)|_{t=0}=(v^0,N^0)  &\mbox{ in } \Omega(0),\\
d|_{t=0}=d^0  &\mbox{ on } \Sigma(0).
\end{array}}\right.\end{equation}

Let $\zeta$ be a solution of \eqref{201811191443}, $\eta:=\zeta-y$ and  the Lagrangian unknowns $(u,B,q)(y,t):=(v,M,\sigma)(\zeta(y,t),t)$. If the initial data $\mathcal{A}^0$ and $B^0$ satisfy  $
\mathcal{A}^{\mm{T}}_0B^0=\bar{M}$ and $B^0\cdot \mathcal{A}^0e_3=\bar{M}_3$,
then $(\eta,u,q)$ satisfies  a so-called transformed MRT problem \cite{WYTIVNMI}:
\begin{equation*} \left\{\begin{array}{ll}
\eta_t=u &\mbox{ in } \Omega,\\[1mm]
{\rho} u_t+\mm{div}_{\ml{A}}\mathcal{S}_{\mathcal{A}}(q,u) =\lambda\partial_{\bar{M}}^2\eta &\mbox{ in }  \Omega,\\[1mm]
\mathrm{div}_{\mathcal{A}} {u} = 0& \mbox{ in }  \Omega ,\\
  \llbracket \eta  \rrbracket=\llbracket u  \rrbracket=0,\ \llbracket  \mathcal{S}_{\mathcal{A}} (q,u) \mathcal{A}e_3- g
\rho  \eta_3  \mathcal{A}e_3- \lambda\bar{M}_3 \partial_{\bar{M}}\eta\rrbracket= \vartheta H  \mathcal{A}e_3  &\mbox{ on }\Sigma,\\
(\eta, u)=0 &\mbox{ on }\Sigma_-^+,\\
(\eta,u)|_{t=0}=(\eta^0,u^0)  &\mbox{ in }  \Omega
\end{array}\right.\end{equation*}
with $B=\partial_{\bar{M}}\zeta$.
 where $\mathcal{S}_{\mathcal{A}}(q,u ) :=q I-\mu\mathbb{D}_{\mathcal{A}}u$.
Obviously, the above transformed MRT problem can be deduced from the ART problem \eqref{n0101nn} with $\kappa=0$ and with $\sqrt{\lambda}\bar{M}$ in place of $\bar{M}$. We call the transformed MRT problem without surface tension the transformed pure-MRT problem.

Next we briefly introduce some mathematical results in the transformed pure-MRT problem.
For a special vertical (magnetic) field $\bar{M}=(0,0,\bar{M}_3)$, Wang found there exists a strength-threshold, given by a  variational formula, of $\bar{M}$ for stability/instability based on a linearized pure-MRT problem, in which he took $\tau=l=1$ for simplicity \cite{WYC}. Then Jiang--Jiang--Wang further computed out the strength-threshold obtained by Wang is equal to $\sqrt{ g  \llbracket \rho \rrbracket /2\lambda}$ in \cite{JFJSWWWOA}.  Recently,  Wang have rigorously proved that the transformed pure-MRT problem is stable under the stability condition
\begin{equation}
\label{201806082127}
\bar{M}_3>m_{\mm{S}}:= \sqrt{\frac{ g  \llbracket \rho \rrbracket}{\lambda(\tau^{-1}+l^{-1})}}.
\end{equation}
which is equivalent to
$$  \widetilde{\mm{Dis}}(0,0,\bar{M}) <1.$$
Here we have use the following definitions:
\begin{align}
& \widetilde{\mm{Dis}}(\vartheta,0,\bar{M}):=\sup_{0\neq w\in H_{\sigma,\vartheta}^1(\Omega_{-l}^\tau)} { g\llbracket\rho\rrbracket |w_3|_{0}^2}/{ \widetilde{\mathcal{I}}(w,\vartheta, \bar{M})},\nonumber \\
&   H_{\sigma, \vartheta}^1( \Omega_{-l}^\tau):=\left\{
                            \begin{array}{ll}
         H_{\sigma}^1(\Omega_{-l}^\tau):=H_0^1(\Omega_{-l}^\tau)\cap \{\mm{div} w=0\}, &\hbox{ if }\vartheta= 0;\\
           \{w\in  H_{\sigma}^1(\Omega_{-l}^\tau)~|~w_3(y_{\mm{h}},0)\in H^1(\mathbb{R}^2)\},   & \hbox{ if }\vartheta\neq 0,
                            \end{array}
                          \right.
\nonumber  \\
& \widetilde{\mathcal{I}}(w,\vartheta, \bar{M}) :=\vartheta|\nabla_{\mm{h}}w_3|_{L^2(\mathbb{R}^2)}^2+
\|\partial_{\bar{M}}w\|^2_{L^2(\Omega_{-l}^\tau )}.\nonumber
\end{align}

Since,  for any constant vector $\vec{n}$ with $\vec{n}_3=1$,
\begin{align}
\label{201808161915}
\sup_{0 \neq w\in H_0^1}
\frac{  |w_3|_{0}^2}{
\|\partial_{\vec{n}}w\|^2_0 }= \sup_{0 \neq w\in H_{\sigma}^1}
\frac{  |w_3|_{0}^2}{
\|\partial_{\vec{n}}w\|^2_0 }=  {\frac{ 1}{\tau^{-1}+l^{-1}}},
\end{align}
see Lemma 4.6 and Remark 4.3 in \cite{JFJSOMITN} for a direct proof, thus
the stability condition \eqref{201806082127} is equivalent to
$$ \mm{Dis}(0,0,\bar{M})<1.$$
We mention that, putting some techniques of energy estimates in Section \ref{201805081731} into the multi-layer energy method in \cite{WYTIVNMI}, we can also  verify that the transformed MRT problems defined on $\Omega_{-l}^\tau$, resp. $\Omega$, are stable under the stability conditions $ \widetilde{\mm{Dis}}(\vartheta,0,\bar{M})<1$, resp. $\mm{Dis}(\vartheta,0,\bar{M})<1$ with $\bar{M}_3\neq 0$.

Noting that \eqref{201806071047} with $\kappa=0$ and with $(\vartheta,\bar{M}_3)\neq 0$ is equivalent to
\begin{equation}
\label{201806142108}
\mm{Dis}(\vartheta,0,\bar{M})>1,
\end{equation}
then Theorem \ref{thm3} tells us that the transformed MRT problem  is unstable when
 $(\vartheta,\bar{M}_3)=0$ or the instability condition \eqref{201806142108} with $(\vartheta,\bar{M}_3)\neq 0$ is satisfied. Moreover, by Theorem \ref{201806012301xx}, we also obtain the instability result of the perturbation MRT problem \eqref{0103nxasdfx}.
\begin{cor}\label{201806012301} Under the instability condition \eqref{201806142108}
 with $(\vartheta,\bar{M}_3)\neq 0$ or $(\vartheta,\bar{M}_3)=0$, the perturbation MRT problem \eqref{0103nxasdfx} is unstable.
More precisely,  $(v,N,   \sigma,d)$ constructed by Theorem \ref{201806012301xx} with $\kappa=0$ is an unstable solution to the perturbation MRT problem.
\end{cor}
\begin{rem}
 Recently Jiang et.al. have proved the existence results of RT instability solutions in
 the sense of   $H^2$-norm for the non-homogenous viscoelastic fluid and the non-homogenous MHD fluid \cite{JFJSJMFMOSERT,JFJWGCOSdd}.
We mention that, by
following the derivation of Corollaries \ref{201806012301xxxxxxxx}--\ref{201806012301},
we also obtain RT instability solutions in the sense of $L^1$-norm for the corresponding non-homogenous  viscoelastic fluid and the non-homogenous MHD fluid (without interface).
\end{rem}

We further mention the transformed pure-MRT problem.
By Lemma \ref{201806151836}, we can see that, for any non-zero horizontal field $\bar{M}=(\bar{M}_{\mm{h}},0)$,  the transformed pure-MRT problem is always unstable for any horizontal field by Corollary \ref{201806012301}. This means the horizontal field can not inhibit the RT instability.
However, the transformed pure-MRT problem is always stable for any vertical field satisfying the stability condition \eqref{201806082127}. This means that a sufficiently strong vertical field can inhibit the RT instability. The reason of the inhibition of RT instability by the vertical field lies in that the non-slip velocity boundary condition is imposed in the direction of the vertical field. In \cite{JFJSJMFM}, the authors have found out that the non-slip velocity boundary condition, imposed in the direction of magnetic fields, can enhance the stabilizing effect of the magnetic fields in the linearized IMRT problem (i.e., the linearized RT problem of inhomogeneous MHD fluids. Moreover, the authors have rigorously  proved the mathematical result of  inhibition of RT instability by the horizontal field in the IMRT problem \cite{JFJSJMFMOSERT}. Such magnetic inhibition phenomenon by horizontal field was also found in the Parker instability problem \cite{JFJSSETEFP}, where the horizontal field is vertically decreasing. Recently, the authors have established an magnetic inhibition theory, which reveals why the non-slip velocity boundary condition can enhance the stabilizing effect of the magnetic fields \cite{JFJSOMITN}.
The magnetic inhibition theory can be also used to explain why the VRT problem \eqref{0103nxx} is stable when $\kappa$ is sufficiently large.

\subsection{Stratified pure viscous fluids}

Finally, we discuss the RT problem of stratified pure fluids.
Deleting the terms involving $V$ in \eqref{0103nxx} (or the terms involving $M$ in \eqref{0103nxasdfx}), we get the following  RT problem of stratified pure viscous  fluids in perturbation form:
\begin{equation}\label{0103nxxds}\left\{{\begin{array}{ll}
\rho  v_t+\rho v\cdot\nabla v+\mm{div}\mathcal{S}(\sigma, v,0 )=0 &\mbox{ in }  \Omega(t),\\[1mm]
\mathrm{div} v= 0& \mbox{ in } \Omega(t),\\
d_t+v_1 \partial_1d+v_2 \partial_2d=v_3 &\mbox{ on }\Sigma(t),\\
   \llbracket v  \rrbracket=0,\quad
   \llbracket \mathcal{S}(\sigma, v,0 ) -g \rho d I\rrbracket \nu =\vartheta \mathcal{C} \nu  \ &\mbox{ on }\Sigma(t) ,  \\
v =0 &\mbox{ on }\Sigma_-^+,\\
v|_{t=0}=v^0 &\mbox{ in } \Omega(0),\\
d|_{t=0}=d^0  &\mbox{ on } \Sigma(0).
\end{array}}\right.\end{equation}

We mention some progress in the RT problem \eqref{0103nxxds}.
If $l$ and $\tau$ are infinite, and $\mathbb{T}$ is replaced by $\mathbb{R}^2$,
Pr$\mathrm{\ddot{u}}$ess--Simonett first have proved that there exists a unique unstable solution to the RT problem in a perturbation form in 2010 \cite{PJSGOI5x}, where the solution $(d,v)$ is unstable in the norm of  Sobolev--Slobodeckii space $W^{3-2/p}_p(\mathbb{R}^2)\times W_p^{2-2/p}(\Omega(t))$ for $p>5$ based on a Henry instability method. Later, Wang--Tice--Kim proved that there exists a strength threshold $\vartheta_{\mathbb{T}}:={ g\llbracket\rho\rrbracket\max\{L_1^2,L_2^2\} }$ of the surface tension coefficient
for the instability and stability of  the RT problem \eqref{0103nxxds}  in ``geometric" coordinate. More precisely, they proved that the RT problem \eqref{0103nxx} is stable for $\vartheta >\vartheta_{\mathbb{T}} $ \cite{WYJTIKCT},
and formally verified that a RT problem in ``geometric" coordinate  is unstable in the sense of $L^2$-norm for \cite{wang2011viscous}
\begin{equation}
\label{201810011329}
 \vartheta \in [0, \vartheta_{\mathbb{T}}).
\end{equation}
In Appendix, we will provide direct proof for the relation
\begin{equation}
\label{201811151206}
\vartheta_{\mathbb{T}}=\sup_{  w\in H_{\sigma,3}^1} { g\llbracket\rho\rrbracket |w_3|_{0}^2}/{ |\nabla_{\mm{h}}w_3|_{0}^2} ,
\end{equation}
see Lemma \ref{201806091256}. Thus we see immediately that the two conditions \eqref{201810011329} and ``\eqref{201806071047} with $(\kappa,\bar{M})=0$" are equivalent. Recently,  Wilke have proved there exists a threshold $\vartheta_{\mm{c}} $ for the stability and instability of RT problem of stratified pure fluids defined on a cylindrical domain based on a rather complete analysis technique of the corresponding free boundary problem which involves a contact angle \cite{WMRTITNS2017}, where the unstable solutions is also in the sense of  norm of Sobolev--Slobodeckii space as  Pr$\mathrm{\ddot{u}}$ess--Simonett's result.

However, by \eqref{201811151206}, we can deduce from Theorem  \ref{201806012301xx} the following conclusion which presents that there exists a RT instability solution in the sense of $L^1$-norm to \eqref{0103nxxds} under the instability condition \eqref{201810011329}.
\begin{cor}\label{20180610185}
Under the instability condition \eqref{201810011329}, the perturbation RT problem \eqref{0103nxxds} is unstable.
More precisely,  $(v, \sigma,d)$ constructed by Theorem \ref{201806012301xx} with $\kappa=\lambda=0$ is an unstable solution to the perturbation RT problem.
\end{cor}

It is well-known that, in the development of the RT instability,
gravity first drives the third component of velocity $v_3$ unstable, then, the instability of the third component of the velocity further results in instability of horizontal velocity $v_{\mm{h}}$, and the height function $d_3$. Obviously, our instability result \eqref{201808160850}  accords with this instability phenomenon, and are finer than the known mathematical results of RT instability solutions in stratified pure viscous fluids \cite{PJSGOI5x,WMRTITNS2017}.

\section{Linear instability}\label{201805302028}

In this section, we will use modified variational method to construct unstable solutions for a linearized ART problem. The modified variational method  was firstly used by Guo and Tice  to  construct unstable solutions
to a class of ordinary differential equations arising from a linearized RT instability problem \cite{GYTI2}.
Later, Jiang--Jiang \cite{JFJSO2014,JFJSJMFM} further extended the modified {variational method} to construct unstable solutions to the partial differential equations (PDEs) arising from a linearized RT instability problem. Exploiting the  modified variational method of PDEs in \cite{JFJSJMFM} and an existence theory  of stratified Stokes problem (see Lemma \ref{xfsddfsf201805072212}), we can obtain a linear instability result of the ART problem.

\begin{pro}\label{thm:0201201622}
Under the assumption \eqref{201806071047}, there is an unstable solution
$$( \eta, u,{q}):=e^{\Lambda t}( w/\Lambda,w,\beta)$$
 to  the linearized ART problem
\begin{equation*}
\left\{\begin{array}{ll}
 \eta_t=  u &\mbox{ in }\Omega, \\[1mm]
\rho u_{t}+\mm{div}\mathcal{S} (q,u,\eta) =\partial_{\bar{M}}^2 \eta & \mbox{ in }\Omega,\\[1mm]
\div u  =0 &\mbox{ in }\Omega,\\
\llbracket \eta  \rrbracket= \llbracket u  \rrbracket=0,\
  \llbracket \mathcal{S}  (q,u,\eta) e_3 - g  \rho  \eta_3 e_3 - \bar{M}_3  \partial_{\bar{M}}\eta\rrbracket =\vartheta\Delta_{\mm{h}}\eta_3e_3 &\mbox{ on }\Sigma,\\
 (\eta,u)= 0 &\mbox{ on }\Sigma_-^+,\\
(\eta,u)|_{t=0}= (\eta^0,u^0)  &\mbox{ in }  \Omega,\end{array}\right.\end{equation*} where $(w,\beta)\in (\mathcal{A}\cap H^\infty) \times  ( {H}^1\cap H^\infty)$
 solves  the boundary value problem:
\begin{equation}\label{201604061413}
\left\{\begin{array}{ll}
\Lambda^2\rho w=  \partial_{\bar{M}}^2w -\mm{div}\mathcal{S}(\Lambda\beta,\Lambda w,w) &\mbox{ in }\Omega, \\ \mm{div}w=0 &  \mbox{ in }\Omega, \\
   \llbracket w \rrbracket= 0,\
   \llbracket \mathcal{S}(\Lambda \beta, \Lambda w, w)e_3
        - g \rho w_3 e_3 -  \bar{M}_3 \partial_{\bar{M}}w\rrbracket =  \vartheta\Delta_{\mm{h}}w_3  e_3  &\mbox{ on }\Sigma, \\
 w = {0}& \mbox{ on }\Sigma_-^+  \end{array} \right.
\end{equation}
with a largest growth rate  $\Lambda>0$ satisfying
\begin{equation}\label{Lambdard}
\Lambda^2=- \inf_{\varpi\in \mathcal{A} }\left({E}(\varpi)+\Lambda  \|\sqrt{\mu }\mathbb{D}\varpi\|_0^2/2\right) . \end{equation}
Moreover, for $\chi=\eta$ or $w$,
\begin{align}
& \label{201602081445MH}
\chi_3\neq 0,\ \partial_3 \chi_3\neq 0,\ \chi_{\mm{h}}\neq 0,\ \partial_3\chi_{\mm{h}}\neq 0,\ \mm{div}_{\mm{h}}\chi_{\mm{h}}\neq 0\mbox{ in }\Omega, \ |\chi_3|\neq 0\mbox{ on }\Sigma,\\
& \partial_{\bar{M}} \chi_3 \neq 0,\   \partial_{\bar{M}} {\chi}_{\mm{h}} \neq 0 \mbox{ if }
\bar{M}_3\neq0.\label{201809121901}
\end{align}
\end{pro}

\begin{pf} We divide the proof of Proposition \ref{thm:0201201622} by four steps.

(1) \emph{Existence of weak solutions to the modified problem}
 \begin{equation}
\label{2016040614130843xd}
\left\{    \begin{array}{ll}
\mm{div}\mathcal{S}(s\beta,sw, w)-\partial_{\bar{M}}^2w =   \alpha(s) \rho w &\mbox{ in }\Omega, \\ \mm{div}w=0 &  \mbox{ in }\Omega, \\
\llbracket w  \rrbracket=0,\
\llbracket   \mathcal{S}(s\beta,sw, w)e_3 - g \rho w_3  e_3- \bar{M}_3 \partial_{\bar{M}}w\rrbracket=  \vartheta\Delta_{\mm{h}}w_3  e_3&\mbox{ on }\Sigma, \\
  w = {0}& \mbox{ on }\Sigma_-^+,  \end{array} \right.
                       \end{equation}
\emph{where  $s>0$ is any given.}

To prove the existence of weak solutions of the above problem, we consider the variational problem of the functional $F(\varpi,s )$:
$$\alpha(s): =\inf_{\varpi\in\mathcal{A}} F(\varpi,s ) $$
for given $s>0$, we have defined that $F(\varpi,s):=E(\varpi )+
 s \|\sqrt{\mu}\mathbb{D}\varpi\|_0^2/2$. Sometimes, we denote $\alpha(s)$ and $F(\varpi,s)$ by $\alpha$ and $F(\varpi)$ for simplicity, resp..

Noting that
\begin{equation}
\label{201805072235}
|v|_0^2 \lesssim \|v\|_0\|\partial_3 v\|_0\mbox{ for any }v\in H_0^1,
\end{equation}
thus, by Young's inequality and Korn's inequality, we see that $F(\varpi)$ has an lower bound for any $\varpi\in \mathcal{A}$.
Hence  $F(\varpi)$ has a minimizing sequence $\{w^n\}_{n=1}^\infty\subset \mathcal{A}$, which satisfies $\alpha=\lim_{n\to\infty} F(w_n)$. Moreover, making use of \eqref{201805072235}, the fact $\|\sqrt{\rho}w^n\|_0=1$, trace estimate (see Lemma \ref{xfsddfsf201805072254}) and Young's and Korn's inequalities, we have $\|w^n\|_1+\vartheta|\nabla_{\mm{h}}w^n_3|_0+|w^n|_0\leqslant c_1$ for some constant $c_1$, which is independent of $n$. Thus, by the well-known Rellich--Kondrachov  compactness theorem,  there exist a subsequence, still labeled by  $w^n$, and a function $w\in \mathcal{A}$, such that
$$\begin{aligned}
& w^n\rightharpoonup w\mbox{ in }H^1_\sigma,\ w^n\to w\mbox{ in }L^2,\ w^n|_{y_3=0}\to w|_{y_3=0}\mbox{ in }L^2(\mathbb{T}),\\
&  w^n_3|_{y_3=0}\rightharpoonup w_3|_{y_3=0}\mbox{ in }H^1(\mathbb{T})\mbox{ if }\vartheta\neq 0.
\end{aligned}$$
Exploiting the above convergence results, and the lower semicontinuity of weak convergence, we have
\begin{equation*}
\alpha =\liminf_{n\to \infty}F(w^n)
\geqslant F(w) \geqslant \alpha .\end{equation*}
Hence $w$ is the minimum point of the functional $F(\varpi)$ with respect to $\varpi\in \mathcal{A}$.

Obviously, $w$ constructed above is also  the minimum point of the functional $F(\varpi )/  \|\sqrt{\rho }\varpi\|^2_0$ with respect to  $\varpi\in H_{\sigma, \vartheta}^1$. Moreover $\alpha= F(w)/\|\sqrt{\rho}w\|^2_0$. Thus, for any given $\varphi \in H^1_{\sigma,\vartheta} $, the point $t=0$ is the minimum point of the function
$$I(t):=F(w+t\varphi )-\int \alpha  \rho |w+t\varphi|^2 \mm{d}y\in C^1(\mathbb{R}).$$
Then, by computing out $I'(0)=0$, we have the weak form:
\begin{align}
\label{201805161108}
& \frac{1}{2}\int (s\mu +\kappa\rho)\mathbb{D} w:\mathbb{D}\varphi\mm{d}y +\vartheta  \int_\Sigma \nabla_{\mm{h}}w_3\cdot \nabla_{\mm{h}}\varphi_3 \mm{d}y_{\mm{h}} \nonumber \\
& +
\int \partial_{\bar{M}}w\cdot \partial_{\bar{M}}\varphi \mm{d}y  = g\llbracket\rho\rrbracket \int_\Sigma w_3\varphi_3 \mm{d}y_{\mm{h}}+ \alpha \int \rho w\cdot \varphi \mm{d}y.
\end{align}
Noting the fact,   for any $f^1$, $f^2\in H^1$, and any matrix $A$, $B\in \mathbb{R}^{3\times 3}$,
\begin{equation}
\label{201800515}
\frac{1}{2}\int \mathbb{D}_{A}f^1:\mathbb{D}_Bf^2\mm{d}y=\int \mathbb{D}_A f^1:\nabla_B f^2\mm{d}y,
\end{equation}
thus \eqref{201805161108} is equivalent to
\begin{align}
& \int (s\mu +\kappa\rho)\mathbb{D} w:\nabla \varphi\mm{d}y+ \vartheta  \int_\Sigma\nabla_{\mm{h}}w_3\cdot \nabla_{\mm{h}}\varphi_3\mm{d}y_{\mm{h}}\nonumber \\
&+
\int \partial_{\bar{M}}w \cdot \partial_{\bar{M}}\varphi\mm{d}y = g\llbracket\rho\rrbracket \int_\Sigma w_3\varphi_3\mm{d}y_{\mm{h}}+\alpha  \int \rho w\cdot \varphi\mm{d}y. \nonumber
\end{align}
The means that $w$ is the weak solution of the modified problem \eqref{2016040614130843xd}.

(2) \emph{Improving the regularity of the weak solution $w$.}

To begin with, we shall establish the following preliminary conclusion:

\emph{For any $i\geqslant 0$, we have}
\begin{equation}
w \in H^{1,i}_{\sigma,\vartheta}
\label{201806181000}
\end{equation}
\emph{and}
\begin{align}
&\frac{1}{2}\int (s\mu+\kappa\rho )\mathbb{D} \partial_{\mm{h}}^i w :\mathbb{D}\varphi \mm{d}y + \vartheta\int_\Sigma \nabla_{\mm{h}} \partial_{\mm{h}}^i w_3\cdot\nabla_{\mm{h}}  \varphi_3\mm{d}y_{\mm{h}}\nonumber \\
& +
\int \partial_{\bar{M}} \partial_{\mm{h}}^i w\cdot \partial_{\bar{M}}\varphi\mm{d}y  = g\llbracket\rho\rrbracket \int_\Sigma  \partial_{\mm{h}}^i w_3\varphi_3\mm{d}y_{\mm{h}}+\alpha \int \rho   \partial_{\mm{h}}^i w\cdot   \varphi \mm{d}y.
\label{201805161108sdfs}
\end{align}

Obviously, by induction, the above assertion reduces to verify the following recurrence relation:

\emph{For given $i\geqslant 0$, if $w \in H^{1,i}_{\sigma,\vartheta}$ satisfies \eqref{201805161108sdfs}
for any $\varphi \in H^1_{\sigma,\vartheta} $,  then}
\begin{equation}
\label{201806181412}
w \in H^{1,i+1}_{\sigma,\vartheta}
\end{equation} and
 $w$ satisfies
\begin{align}
&\frac{1}{2}\int (s\mu+\kappa\rho )\mathbb{D} \partial_{\mm{h}}^{i+1} w :\mathbb{D} \varphi \mm{d}y+ \vartheta \int_\Sigma \nabla_{\mm{h}}  \partial_{\mm{h}}^{i+1} w_3\cdot\nabla_{\mm{h}}  \varphi_3\mm{d}y_{\mm{h}} \nonumber \\
& +
\int \partial_{\bar{M}} \partial_{\mm{h}}^{i+1} w\cdot \partial_{\bar{M}}\varphi \mm{d}y  =
g\llbracket\rho\rrbracket \int_\Sigma  \partial_{\mm{h}}^{i+1} w_3\varphi_3 \mm{d}y_{\mm{h}}+\alpha  \int \rho \partial_{\mm{h}}^{i+1} w\cdot   \varphi \mm{d}y.\label{201805161108sdfssafas}
\end{align}
Next we verify the above recurrence relation by method of difference quotients.

Now we assume that $w \in H^{1,i}_{\sigma,\vartheta}$ satisfies \eqref{201805161108sdfs}
for any $\varphi \in H^1_{\sigma,\vartheta} $. Noting that $\partial_{\mm{h}}^iw \in H^1_{\sigma,\vartheta}$, we can deduce from \eqref{201805161108sdfs} that, for $j=1$ and $2$,
\begin{align}
& \frac{1}{2}\int (s\mu+\kappa\rho) \mathbb{D} \partial_{\mm{h}}^i w : \mathbb{D}  D_j^h \varphi
 \mm{d}y+ \vartheta\int_\Sigma \nabla_{\mm{h}} \partial_{\mm{h}}^i w_3 \cdot \nabla_{\mm{h}}  D_j^h  \varphi_3 \mm{d}y_{\mm{h}} \nonumber \\
& +
\int \partial_{\bar{M}}\partial_{\mm{h}}^i w \cdot \partial_{\bar{M}} D_j^h \varphi  \mm{d}y = g\llbracket\rho\rrbracket \int_\Sigma \partial_{\mm{h}}^i w _3  D_j^h  \varphi_3\mm{d}y_{\mm{h}}+\alpha \int \rho
  \partial_{\mm{h}}^i w\cdot  D_j^h \varphi  \mm{d}y \nonumber
\end{align}
and
\begin{align}
& \frac{1}{2}\int (s\mu+\kappa\rho) \mathbb{D} \partial_{\mm{h}}^i w : \mathbb{D} D_j^{-h} D_j^h \partial_{\mm{h}}^i w
 \mm{d}y+ \vartheta \int_\Sigma \nabla_{\mm{h}}\partial_{\mm{h}}^i w_3\cdot \nabla_{\mm{h}} D_j^{-h} D_j^h\partial_{\mm{h}}^i w_3 \mm{d}y_{\mm{h}}  \nonumber \\
& +
\int \partial_{\bar{M}}  \partial_{\mm{h}}^i w\cdot \partial_{\bar{M}} D_j^{-h} D_j^h \partial_{\mm{h}}^i w  \mm{d}y = g\llbracket\rho\rrbracket \int_\Sigma   \partial_{\mm{h}}^i w_3 D_j^{-h} D_j^h \partial_{\mm{h}}^i w_3\mm{d}y_{\mm{h}} +\alpha  \int \rho
 \partial_{\mm{h}}^i w\cdot D_j^{-h}  D_j^h  \partial_{\mm{h}}^i w  \mm{d}y, \nonumber
\end{align}
which yield that
\begin{align}
&\frac{1}{2}\int (s\mu+\kappa\rho) \mathbb{D} D_j^{-h} \partial_{\mm{h}}^i w : \mathbb{D} \varphi
 \mm{d}y +\vartheta \int_\Sigma \nabla_{\mm{h}} D_j^{-h} \partial_{\mm{h}}^i w_3\cdot \nabla_{\mm{h}}    \varphi_3 \mm{d}y_{\mm{h}}\nonumber \\
& +\int \partial_{\bar{M}}  D_j^{-h}\partial_{\mm{h}}^i w \cdot \partial_{\bar{M}} \varphi\mm{d}y = g\llbracket\rho\rrbracket \int_\Sigma D_j^{-h} \partial_{\mm{h}}^i w_3 \varphi_3\mm{d}y_{\mm{h}}+\alpha \int \rho
  D_j^{-h}\partial_{\mm{h}}^i w\cdot \varphi\mm{d}y ,
\label{20180612713459}
\end{align}
and
\begin{align}
&\|\sqrt{s\mu+\kappa\rho}\mathbb{D} D_j^h\partial_{\mm{h}}^i w\|^2_0/2+ \vartheta|D_j^h\nabla_{\mm{h}} \partial_{\mm{h}}^i w |^2_0
\nonumber \\
&+ \|\partial_{\bar{M}}  D_j^h\partial_{\mm{h}}^i w \|^2_0\lesssim  g\llbracket\rho\rrbracket |D_j^h\partial_{\mm{h}}^i w_3 |^2_0 +|\alpha|\|\sqrt{ \rho}
D_j^h\partial_{\mm{h}}^i w \|^2_0, \label{201806161512}
\end{align}
resp..

By Korn's inequality,
$$\| D_j^h \partial_{\mm{h}}^i w\|^2_{1}\lesssim \|\sqrt{s\mu+\kappa\rho}\mathbb{D} D_j^h \partial_{\mm{h}}^i w\|^2_0,$$
thus, using \eqref{201805072235}, Young's inequality, and the first conclusion in Lemma \ref{xfsddfsf2018050813379safdadf} , we further deduce from \eqref{201806161512} that
\begin{align}
  \|D^h_{\mm{h}}  \partial_{\mm{h}}^i w\|^2_1+ \vartheta| D^h_{\mm{h}}\nabla_{\mm{h}} \partial_{\mm{h}}^i w|^2_0
\lesssim  \|D^h_{\mm{h}} \partial_{\mm{h}}^i w \|^2_0\lesssim \| \nabla_{\mm{h}} \partial_{\mm{h}}^i  w \|^2_0\lesssim 1.\nonumber
\end{align}
Thus, using \eqref{201805072235}, trace estimate and the second conclusion in Lemma \ref{xfsddfsf2018050813379safdadf},  there exists a subsequence of $\{-h\}_{h\in \mathbb{R}}$ (still denoted by $-h$) such that
\begin{equation}
\label{201806127345}
 \left\{
   \begin{array}{ll}
 D^{-h}_{\mm{h}} \partial_{\mm{h}}^i  w \rightharpoonup \nabla_{\mm{h}}\partial_{\mm{h}}^i w \mbox{ in }H^1_\sigma,\  D^{-h}_{\mm{h}} \partial_{\mm{h}}^i  w \to \nabla_{\mm{h}}\partial_{\mm{h}}^i w \mbox{ in }L^2, &\\
 D^{-h}_{\mm{h}} \partial_{\mm{h}}^i  w|_{y_3=0} \to \nabla_{\mm{h}}\partial_{\mm{h}}^i w |_{y_3=0} \mbox{ in }L^2(\mathbb{T})& \\
 D^{-h}_{\mm{h}}\partial_{\mm{h}}^i w|_{\Sigma}\rightharpoonup \nabla_{\mm{h}}\partial_{\mm{h}}^i w|_{\Sigma}\mbox{ in }H^1(\mathbb{T})\mbox{ if }\vartheta\neq 0.&
   \end{array}
 \right.
\end{equation}
Using  regularity of $w$ in \eqref{201806127345} and the fact $w\in H_{\sigma,\vartheta}^{1,i}$, we have \eqref{201806181412}.
In addition, exploiting the limit results in \eqref{201806127345},  we can deduce \eqref{201805161108sdfssafas} from \eqref{20180612713459}.
This complete the proof of the recurrence relation, and thus \eqref{201806181000} holds.

With \eqref{201806181000} in hand, we can consider a stratified   Stokes problem:
 \begin{equation}\label{2016040614130843x}      \left\{  \begin{array}{ll}
s \nabla \beta^k -(s\mu +\kappa\rho+\bar{M}_3^2)\Delta\omega^k =\partial_{\mm{h}}^k\mathcal{L}^1
&\mbox{ in } \Omega,\\
\mm{div}\omega^k=0&\mbox{ in } \Omega,\\
 \llbracket \omega^k  \rrbracket=0,\
\llbracket( s  \beta^k- (s\mu+\kappa\rho +\bar{M}_3^2)\mathbb{D}\omega^k)e_3 \rrbracket= \partial_{\mm{h}}^k\mathcal{L}^2&\mbox{ on }\Sigma,
\\ \omega^k =0  &\mbox{ on }\Sigma_{-}^+,\end{array}\right.
\end{equation}
where $k\geqslant 0$ is a given integer, and we have defined that
$$
\begin{aligned}
&\mathcal{L}^1:= \partial_{\bar{M}}^2w- \bar{M}_3^2\Delta w+\alpha \rho w,\ \mathcal{L}^2 :=g \llbracket\rho \rrbracket w_3e_3
 + \vartheta\Delta_{\mm{h}}w_3 e_3.
\end{aligned}$$
Recalling the regularity \eqref{201806181000} of $w$, we see that $\partial_{\mm{h}}^k\mathcal{L}^1\in L^2$, and $ \partial_{\mm{h}}^k\mathcal{L}^2 \in H^{1}(\mathbb{T})$. Applying the existence theory of stratified   Stokes problem (see Lemma \ref{xfsddfsf201805072212}), there exists a unique strong solution $(\omega^k,\beta^k)\in H^2\times \underline{H}^1$ of the above problem \eqref{2016040614130843x}.

Multiplying \eqref{2016040614130843x}$_1$ by $\varphi\in H^1_{\sigma,\vartheta}$ in $L^2$ (i.e., taking the inner product in $L^2$), and using the integration by parts and  \eqref{2016040614130843x}$_2$--\eqref{2016040614130843x}$_4$, we have
\begin{align}
 &\frac{1}{2}\int (s\mu+\kappa\rho+\bar{M}_3^2)\mathbb{D}\omega^k:\mathbb{D}\varphi \mm{d}y \nonumber
\\
& = \int \alpha \rho \partial_{\mm{h}}^k w\varphi\mm{d}y +g\llbracket\rho\rrbracket \int_\Sigma \partial_{\mm{h}}^k w_3\varphi_3
\mm{d}y_{\mm{h}}\nonumber \\
&\quad -  \int_\Sigma \vartheta\partial_{\mm{h}}^k \nabla_{\mm{h}}w_3\cdot \nabla_{\mm{h}}\varphi_3\mm{d}y_{\mm{h}}+\int \partial_{\mm{h}}^k( \partial_{\bar{M}}^2w- \bar{M}_3^2\Delta w)\cdot \varphi\mm{d}y.
\label{201808072057}
\end{align}

Noting that
\begin{align}
\label{20180803939}
\mm{div}\partial_{\mm{h}}^i w=0\mbox{ in }\Omega,\mbox{ and } \llbracket\partial_{\mm{h}}^i w\rrbracket =0\mbox{ on }\Sigma\mbox{ for any }i\geqslant 0
\end{align}
and
\begin{align}
\bar{M}_3^2 \int\mathbb{D}\partial_{\mm{h}}^k  w:\nabla \varphi\mm{d}y -\int \partial_{\mm{h}}^k  ( \partial_{\bar{M}}^2 w- \bar{M}_3^2\Delta w)\cdot \varphi \mm{d}y=\int \partial_{\bar{M}} \partial_{\mm{h}}^{k} w\cdot \partial_{\bar{M}}\varphi \mm{d}y,\nonumber
\end{align}
we can use the integration by parts and \eqref{201800515} to rewrite \eqref{201805161108sdfs} as follows:
\begin{align}
&\frac{1}{2}\int (s\mu +\kappa\rho +\bar{M}_3^2)\mathbb{D}\partial_{\mm{h}}^k  w:\mathbb{D} \varphi\mm{d}y+  \int_\Sigma \vartheta \nabla_{\mm{h}}\partial_{\mm{h}}^k  w_3\cdot \nabla_{\mm{h}}\varphi_3\mm{d}y_{\mm{h}}\nonumber \\
&-\int \partial_{\mm{h}}^k  ( \partial_{\bar{M}}^2 w- \bar{M}_3^2\Delta w)\cdot \varphi \mm{d}y= \int \alpha \rho \partial_{\mm{h}}^k  w\cdot \varphi\mm{d}y + g\llbracket\rho\rrbracket \int_\Sigma \partial_{\mm{h}}^k  w_3\varphi_3\mm{d}y_{\mm{h}} .  \label{sdfsaf201808072057}
\end{align}
Subtracting the two identities \eqref{201808072057} and \eqref{sdfsaf201808072057} yields that
$$\int(s\mu+\kappa\rho +\bar{M}_3^2) \mathbb{D}(\partial_{\mm{h}}^k  w-\omega^k):\mathbb{D} \varphi \mm{d}y=0.
$$
Taking $\varphi:= \partial_{\mm{h}}^k w-\omega^k\in H^1_{\sigma,\vartheta}$ in the above identity, and using the Korn's inequality, we find that $\omega^k= \partial_{\mm{h}}^k w$. Thus we immediately see that
\begin{equation}
\label{201806181507}
\partial_{\mm{h}}^k  w\in H^{2}\mbox{ for any }k\geqslant 0,
\end{equation}
which implies $\partial_{\mm{h}}^k\mathcal{L}^1\in H^1$, and $\partial_{\mm{h}}^k\mathcal{L}^2\in H^{2}(\mathbb{T})$ for any $k\geqslant 0$. Thus, applying the stratified Stokes estimate \eqref{2011805302036} to \eqref{2016040614130843x}, we have
\begin{equation}
\label{201806181507fdsdsgsdfgsg}
\partial_{\mm{h}}^k  w\in H^{3} \mbox{ for any }k\geqslant 0,
\end{equation}
Obviously, by induction, we can easily follow the improving regularity method from \eqref{201806181507} to \eqref{201806181507fdsdsgsdfgsg} to deduce that $w\in H^\infty$.
In addition, we have an associated pressure $\beta:=\beta^0\in H^\infty$; moreover,  $\beta^k$ in \eqref{2016040614130843x}  is equal to $\partial^k_{\mm{h}}\beta$.

Finally, recalling \eqref{20180803939}, and the embedding inequality \eqref{esmmdforinfty}, we easily see that $(w,\beta)$ constructed above is indeed a classical solution to the modified problem \eqref{2016040614130843xd}.

(3) \emph{Some properties of the function $\alpha(s)$ on $(0,\infty)$:}
\begin{align}
\label{201702081046}
&\alpha(s_2)  >\alpha(s_1)\mbox{ \emph{ for any }}s_2  >s_1>0,\\
\label{201702081047}&
\alpha(s) \in C^{0,1}_{\mm{loc}}(0,\infty),\\
&\label{201702081122}\alpha(s)<0\mbox{ \emph{ on some interval }}(0,c_2),\\
&\label{201702081122n}
\alpha(s)>0\mbox{ \emph{ on some interval }}(c_3,\infty).
\end{align}

To being with, we  verify \eqref{201702081046}.
For given $s_2>s_1$, then there exist  $v^{s_2}\in \mathcal{A}$ such that
$\alpha(s_2)  = F(v^{s_2},s_2)$.
 Thus, by Korn's inequality and the fact $\|\sqrt{\rho}v^{s_2}\|_0=1$,
$$
 \alpha(s_1)\leqslant F(v^{s_2},s_1) = \alpha(s_2) +
 ( s_1- {s_2}) \|\sqrt{\mu }\mathbb{D}v^{s_2}\|_0^2/2< \alpha(s_2),
$$
 which yields \eqref{201702081046}.

Now we turn to prove \eqref{201702081047}. Choosing a bounded interval $[c_4,c_5]\subset (0,\infty)$, then, for any $s\in [c_4,c_5] $, there exists a function $v^s$ satisfying $\alpha(s)=E(v^{s} ) +s\|\sqrt{\mu }\mathbb{D}v^s\|_0^2 /2$. Thus, by the monotonicity \eqref{201702081046}, we have
 $$\alpha(c_5)-
c_4\|\sqrt{\mu }\mathbb{D}v^s\|_0^2/4\geqslant F(v^{s},s/2) \geqslant \alpha(s/2) \geqslant \alpha( c_4/2),$$
  which yields
$$\|\sqrt{\mu } \mathbb{D}v^s\|_0^2/2 \leqslant  2(\alpha(c_5)-\alpha(c_4/2))/c_4=:\xi\mbox{ for any }s\in [c_4,c_5].$$
Thus, for any $s_1$, $s_2\in [c_4,c_5]$,
$$
\begin{aligned}
\alpha(s_1)-\alpha(s_2)\leqslant & F(v^{s_2},s_1)-F(v^{s_2},s_2)\leqslant
\xi|  {s_2}-s_1|
\end{aligned}$$
 and $$\alpha(s_2)-\alpha(s_1)\leqslant
\xi|  {s_2}-s_1|,$$
which immediately imply $|\alpha(s_1)-\alpha(s_2)|\leqslant \xi| {s_2}-s_1|$. Hence \eqref{201702081047} holds.

Finally, \eqref{201702081122} is obvious by the definition of $\alpha$ and the assumption \eqref{201806071047}, while \eqref{201702081122n}
can be also deduced from  the definition of $\alpha$ by using  Korn's inequality and \eqref{201805072235}.

(4) \emph{Construction of an interval for fixed point}:
Let
$$\mathfrak{I}:=\sup\{\mbox{all the real constant }s, \mbox{ which satisfy that }\alpha( \tau)<0\mbox{ for any }\tau\in (0,s)\}.$$
In virtue of \eqref{201702081122} and \eqref{201702081122n}, $0<\mathfrak{I}<\infty$. Moreover, $\alpha(s)<0$ for any $s\in (0,\mathfrak{I})$, and, by the continuity of $\alpha(s)$,
\begin{equation}\label{nzerolin}
 \alpha( \mathfrak{I})=0.
\end{equation}
Using the monotonicity and the upper boundedness of $\alpha(s)$, we see that
 \begin{equation}\label{zeron}
 \lim_{s\rightarrow 0}\alpha(s)=\varsigma\mbox{ for some negative constant }\varsigma.
 \end{equation}

Now, exploiting \eqref{nzerolin},  \eqref{zeron} and the continuity of ${\alpha}(s)$ on $ (0,\mathfrak{I})$,
we find by a fixed-point argument on $(0,\mathfrak{I})$ that there is a unique $\Lambda\in(0,\mathfrak{I})$ satisfying
\begin{equation*}  \Lambda=\sqrt{-\alpha(\Lambda)}=
\sqrt{-\inf_{\varpi\in\mathcal{A}}F(\varpi, \Lambda )}\in (0,\mathfrak{I}).
\end{equation*}
Thus there is a classical solution $w \in \mathcal{A}\cap H^\infty$ to the boundary problem  \eqref{201604061413} with $\Lambda$ constructed above and with $\beta\in \underline{H}^1\cap H^\infty $.
Moreover,
\begin{equation}
\label{growthnn} \Lambda= \sqrt{-F(w, \Lambda )}>0.
\end{equation}
In addition, \eqref{201602081445MH} and \eqref{201809121901} directly follows \eqref{growthnn},  the fact
$\chi\in H_\sigma^1$ for $\eta$ or $w$, and the estimate \eqref{201809222059}.
This completes the proof of Proposition \ref{thm:0201201622}.
\hfill $\Box$
\end{pf}

\section{Growall-type energy inequality of nonlinear solutions}\label{201805081731}

In this section, we mainly derive that any solution of the ART problem enjoys a Growall-type energy inequality. We will derive such inequality by \emph{a priori} estimate method for simplicity.
Let $(\eta,u)$ be a solution
of the ART problem, such that
\begin{equation}\label{aprpiosesnew}
 \sup_{0\leqslant t < T}\sqrt{\|\eta(t)\|_3^2 +\|u(t)\|_2^2}\leqslant \delta \in (0,1)\;\;\mbox{ for some  }T,
\end{equation}
where $\delta$ is sufficiently small, and the initial data $\eta^0$ satisfies $\det(\nabla\eta^0+I)=1$. It should be noted that the smallness depends on the domain and known parameters in the ART problem, and  will be repeatedly used in what follows. Moreover, the solution enjoys fine regularity, which makes sure the procedure of formal deduction.
We will divide the derivation into three steps: first we shall introduce some preliminary estimates, then further derive energy estimates of the solution, finally establish a Growall-type energy inequality.

\subsection{Preliminaries}

This subsection is devoted to derive some preliminary estimates, which include the estimates concerning $\mathcal{A}$,   the estimate of $\partial_t^i\mm{div}\eta$, the equivalent estimate of $\mathbb{D}_{\mathcal{A}}w$, and the estimates of nonlinear terms, such as $\mathcal{N}^1$, $\mathcal{N}^2$ and so on.

\begin{lem}
\label{201805141072}
Under the assumption \eqref{aprpiosesnew}, we have
\begin{enumerate}
  \item[(1)] the estimate for $\mathcal{A}$: for $0\leqslant i\leqslant 2$,
\begin{align}
&\label{aimdse}
\|\mathcal{A}\|_{C^0(\bar{\Omega})}+ \|\mathcal{A}\|_2 \lesssim  1,\\
& \|\ml{A}_t\|_i \lesssim   \|   u\|_{i+1} , \label{06142100x}\\
&  \label{06142100}
\|\ml{A}_{tt}\|_0 \lesssim   \|  (u, u_t)\|_{1} ,\\
& |\mathcal{A}_{tt}e_3|_{-1/2}\lesssim \| u \|_3+\| u_t \|_1 , \label{2018081015456} \\
&  \label{06041533fwqg}
\|\tilde{\mathcal{A}}\|_{i}\lesssim   \| \eta\|_{i+1}.
\end{align}
  \item[(2)] the estimate of $\div { {\eta}}$: for $0\leqslant i\leqslant 2$,
\begin{align}
&\label{improtian1}
\|\div { {\eta}}\|_{i} \lesssim \|\eta\|_{3} \|\eta\|_{i+1} ,\\
& \|\div u\|_i\lesssim \|\eta\|_{3}\|u\|_{i+1}.\label{201808181500}
\end{align}
 \item[(3)] the equivalent estimate
 \begin{equation}\label{20160614fdsa1957}
\|w\|_1 \lesssim \|\mathbb{D}_{\mathcal{A}}w\|_0\lesssim \|w\|_1\quad\mbox{for any }\; w\in H^1_0.
\end{equation}
\end{enumerate}
\end{lem}
\begin{pf}
(1) We can employ \eqref{06131422}, embedding inequality \eqref{esmmdforinfty} and
 product estimate \eqref{fgestims} to have \eqref{aimdse}.
Similarly, we can further use \eqref{n0101nn}$_1$, \eqref{0604221702}, \eqref{201808121247}, trace estimate and the fact
\begin{align}
\label{2018081810135}
\mathcal{A}e_3=\partial_1(\eta+y)\times \partial_2 (\eta+y)=e_3+ e_1\times \partial_2\eta+\partial_1\eta\times e_2 +\partial_1\eta\times \partial_2\eta
\end{align}
to deduce \eqref{06142100x}--\eqref{2018081015456}.

To bound $\tilde{\mathcal{A}}$, we assume that $\delta$ is so small that the following expansion holds.
$$\mathcal{A}^{\mm{T}}=I-\nabla \eta+(\nabla\eta)^2\sum_{i=0}^\infty (-\nabla\eta)^i=I-\nabla \eta+(\nabla\eta)^2\mathcal{A}^T,$$
which implies that
$$ \tilde{\mathcal{A}}^T =(\nabla\eta)^2\tilde{\mathcal{A}}^T+(\nabla\eta)^2 -\nabla \eta. $$
Thus, using  product estimate \eqref{fgestims}, we find that, for $0\leqslant i\leqslant 2$,
\begin{equation*}
\|\tilde{\mathcal{A}}\|_{i}\lesssim  \|\nabla\eta\|^2_2\|\tilde{\mathcal{A}}\|_i+\|\nabla \eta\|_i (1+\|\nabla \eta\|_2),
\end{equation*}
which yields \eqref{06041533fwqg}.

(2)
Since $\det(\nabla\eta_0 +I)=1$, one can derive from \eqref{n0101nn}$_1$ that $\det(\nabla \eta+ I)=1$.
In addition, using determinant expansion theorem, we see that
$$ 1=\det (\nabla \eta+ I)=1+ \mm{div} \eta+r_\eta, $$
where $r_\eta:= ((\mm{div}\eta)^2-\mm{tr} (\nabla \eta)^2)/2+ \det \nabla \eta$. Consequently,
\begin{equation*} \mm{div} \eta=-r_\eta.  \end{equation*}
By \eqref{aprpiosesnew}, and  product estimate \eqref{fgestims}, we immediately get \eqref{improtian1} from the above relation.

Using \eqref{06142100x}, \eqref{06041533fwqg} and product estimate \eqref{fgestims}, we can derive \eqref{201808181500}  from \eqref{s0106pnnnn}$_3$.

(3) Using  \eqref{06041533fwqg} and  Korn's inequality, we can derive that
\begin{align}  &
 \|\mathbb{D}_{\mathcal{A}}w\|_0\lesssim \|w\|_1 , \nonumber  \\
 &  \|w\|_1 \lesssim \| \mathbb{D} w\|_0\lesssim \|\mathbb{D}_{\tilde{\mathcal{A}}}w\|_0+ \|\mathbb{D}_{{\mathcal{A}}} w\|_0\lesssim \|\eta\|_3\|w\|_1+ \|\mathbb{D}_{\mathcal{A}}w\|_0,\nonumber
\end{align}
which imply \eqref{20160614fdsa1957} for sufficiently small $\delta$. \hfill $\Box$
\end{pf}
 \begin{lem}[Estimates of nonlinear terms]
\label{201806291049}
Under the assumption \eqref{aprpiosesnew},
\begin{enumerate}
 \item[(1)] we have
\begin{align}
\label{06011733}
& 　\|\mathcal{N}^1\|_i\lesssim  \|\eta\|_3(\|(\kappa \eta,u)\|_{i+2}+\|\nabla q\|_i)\mbox{ for }0\leqslant i\leqslant  1, \\
 &  |     \mathcal{N}^2   |_{1/2}
 \lesssim \|\eta\|_3(\|((\kappa+|\bar{M}_3|) \eta,u)\|_{2}
+ \|\eta\|_{\underline{2},1}+|  \llbracket q \rrbracket  |_{1/2}),
    \label{06011734}\\
&  |    \mathcal{N}^2   |_{3/2}
 \lesssim \|\eta\|_3( \|(\eta,u)\|_3+|  \llbracket q \rrbracket  |_{3/2})\mbox{ if }\vartheta=0,
    \label{060117sdfa34}\\
&  |    \mathcal{N}^2_{\mm{h}}    |_{3/2}
 \lesssim \|\eta\|_3  (\|(\eta,u)\|_3+ |  \llbracket q \rrbracket  |_{3/2} ) ,\label{201807241645} \\
 &  |\partial_t \mathfrak{M}_1|_{1/2} \lesssim  \|\partial_t \mathfrak{M}_1\|_{1} \lesssim \|\eta\|_3\|(\eta,u)\|_3.\label{201806271120dsfs}
 \end{align}
  \item[(2)]
In addition,
$$\left\{\begin{aligned}
&\begin{aligned} \mathcal{N}^{ 3}:=& \mu \mm{div}_{\mathcal{A}}\mathbb{D}_{\mathcal{A}_t}u+\kappa\rho \mm{div}_{\mathcal{A}}\partial_t (\nabla \eta\nabla \eta^{\mm{T}})-\mm{div}_{\ml{A}_t}\mathcal{S}_{\mathcal{A}}(q,u,\eta ) ,
\end{aligned}\\
&\begin{aligned}\mathcal{N}^4:=&
\llbracket  \mu \mathbb{D}_{ {\mathcal{A}}_t}u\mathcal{A}e_3-\mathfrak{S}_{ {\mathcal{A}} } (q,u,\eta )\mathcal{A}_te_3+\partial_t( g\rho\eta_3 \tilde{\mathcal{A}}e_3  +\kappa\rho\nabla\eta\nabla\eta^{\mm{T}}\mathcal{A}e_3 )\rrbracket
 \\
&+\vartheta\partial_t( H \tilde{\mathcal{A}}e_3+ (H-\Delta_{\mm{h}}\eta_3)e_3)
    \end{aligned}\\
& \mathcal{N}^5:= (\kappa \rho  +\bar{M}_3^2) \mm{div}\eta/\mu - \mathrm{div}_{\tilde{\mathcal{A}} } {u},
\end{aligned}\right.
$$
  we have
\begin{align}
& \|\mathcal{N}^{ 3}\|_0 +|   \mathcal{N}^4  |_{1/2}
\lesssim \|\eta\|_3\|(\eta,u)\|_{ 3}+\|u\|_2(\|u\|_3+\|\nabla q\|_1) +
|  \llbracket q \rrbracket  |_{ 1/2} \|u\|_3, \label{201806291} \\
& \|\mathcal{N}^{5}\|_{j+1}\lesssim \|\eta\|_3\|(\eta,u)\|_{j+2}\mbox{ for }0\leqslant j\leqslant 1.\label{201807241959}
\end{align}
 \end{enumerate}
\end{lem}
\begin{pf}
Next we only derive  \eqref{201806271120dsfs} for example, since the derivations of the other estimates are similar, and simpler.
We mention that we shall further use product estimates \eqref{06011948} and \eqref{06041605xdsf} for the derivations of   $|\llbracket q \rrbracket |_{3/2}$ and $|\llbracket q \rrbracket |_{1/2}$ in \eqref{06011734}--\eqref{201807241645} and \eqref{201806291}.

Using  embedding inequality \eqref{esmmdforinfty}  and product estimate \eqref{fgestims}, we derive from \eqref{2018081810135} that, for sufficiently small $\delta$,
\begin{align}    & \nonumber
 \left\||\mathcal{A}e_3|^{-1}\right\|_2\lesssim 1 \mbox{ and }\|1-|\mathcal{A}e_3|\|_j\lesssim \|\eta\|_{j+1}\mbox{ for }0\leqslant j \leqslant 2.  \nonumber  \end{align}
Making use of  \eqref{06142100x}, \eqref{06041533fwqg}, product estimate \eqref{fgestims} and  the two estimates above, we get
\begin{align}\label{06012130}
 \left\|\vec{n}\right\|_2\lesssim 1,\ \left\|\tilde{n}\right\|_j\lesssim \|\eta\|_{j+1}\mbox{ and }  \|\tilde{n}_t\|_j\lesssim \|u\|_{j+1}.
\end{align}
Recalling the definitions of $H^{\mm{n}}$ and $H^{\mm{d}}$,  and then making use of \eqref{06012130}, product estimate \eqref{fgestims} and \eqref{201808121247}, we
can easily estimate, for $0\leqslant k\leqslant 1$,
\begin{equation}\left\{\begin{aligned}
& \|H^{\mm{n}}\|_{k}\lesssim  \|\eta\|_{2,k}, \ \|H^{\mm{n}}\cdot e_3-\Delta_{\mm{h}}\eta_3\|_{k} \lesssim  \|\eta\|_3\|\eta\|_{2,k},\ \|H^{\mm{n}}_t\|_k  \lesssim \|(\eta,u)\|_{2,k},
\\& \|H^{\mm{n}}_t\cdot e_3-\Delta_{\mm{h}}u_3\|_{k} \lesssim  \|\eta\|_3\|u\|_{k+2}+\|u\|_3\|\eta\|_{k+2}, \\
&\|1/H^{\mm{d}}\|_2\lesssim 1,\ \|1-H^{\mm{d}}\|_{j}\lesssim \|\eta\|_{j+1},\ \|H^{\mm{d}}_t\|_j  \lesssim \|u\|_{j+1} .\end{aligned}\right.\label{0dsafdf6012130n}
\end{equation}

Consequently, making use of \eqref{06012130}--\eqref{0dsafdf6012130n},  product estimate  \eqref{fgestims}  and trace estimate,   we can easily get \eqref{201806271120dsfs}.
\end{pf}

\subsection{Basic energy estimates}

This subsection is devoted to the derivation of energy estimates of $(\eta,u)$, which include  the $y_{\mm{h}}$-derivative  estimates of $(\eta,u)$ (see Lemmas \ref{201612132242nn}--\ref{201612132242nxsfs}), the estimate of temporal derivative of $u$ (see Lemma \ref{201612132242nxsfssdfs}), the stratified Stokes estimates of $u$ (see Lemma \ref{201612132242nx}), the hybrid derivative  estimate of $(\eta,u)$ (see Lemma \ref{lem:dfifessim2057}), and the equivalent estimate of $\mathcal{E}$ (see Lemma \ref{lem:dfifessim2057sadfafdas}). Next we will establish those estimates  in sequence.
\begin{lem}\label{201612132242nn}
Under the assumption \eqref{aprpiosesnew}, the following estimates holds:
\begin{align}
&
\frac{\mm{d}}{\mm{d}t}\int  \left( \rho\partial_\mm{h}^i \eta \cdot  \partial_\mm{h}^i u  + \frac{ \mu  }{4}| \mathbb{D}\partial_\mm{h}^i \eta |^2\right) \mm{d}y+
  \mathcal{I}(\partial_{\mm{h}}^i \eta)\nonumber \\
    & \lesssim
| \eta_3|_i^2+ \|   u\|^2_{i,0}
 + \sqrt{\mathcal{E}} \mathcal{D}\mbox{ for }\left\{
                                   \begin{array}{ll}
                         0\leqslant i\leqslant 1;\\
                              ( i,\vartheta)=(2, 0),
                                   \end{array}
                                 \right.
\label{ssebdaiseqinM0846} \\
&
\frac{\mm{d}}{\mm{d}t}\int_{\mathcal{T}_{4,4}}
\int  \left( \rho  \mathfrak{D}_{\mf{h}}^{3/2} \partial_{\mm{h}}\eta \cdot   \mathfrak{D}_{\mf{h}}^{3/2} \partial_{\mm{h}} u  + \frac{ \mu}{4}  |\mathbb{D} \mathfrak{D}_{\mf{h}}^{3/2} \partial_{\mm{h}}\eta|^2\right) \mm{d}y\mm{d}\mf{h}+
\vartheta|\nabla_{\mm{h}}\partial_{\mm{h}} \eta_3|_{1/2}^2
 \nonumber \\
    & \lesssim |\eta_3|_{2}^2+\|u\|_{1,1}^2 + \sqrt{\mathcal{E}}\mathcal{D}\label{ssebdaisedsaqinM0dsfadsff846asdfadfad},\\
 &
 \frac{\mm{d}}{\mm{d}t}\int  \left( \partial_\mm{h}^2 \eta \cdot  \partial_\mm{h}^2 u  + \frac{ \mu}{4}    | \mathbb{D}\partial_\mm{h}^2 \eta |^2\right) \mm{d}y+
 \kappa \|\partial_{\mm{h}}^2 \eta\|_1^2+\|\partial_{\bar{M}}\partial_{\mm{h}}^2  \eta\|_{0}^2\nonumber \\
  & \lesssim
  \|   u\|^2_{2,0}+ |\partial_{\mm{h}}^2 \eta_3|_{1/2}(|\nabla_{\mm{h}}\partial_{\mm{h}}(\eta_\mm{h} ,u_\mm{h})|_{1/2}+| \llbracket \nabla_{\mm{h}}q \rrbracket |_{1/2})
 + \sqrt{\mathcal{E}} \mathcal{D}.
\label{ssebdaiseqinM0asdfa846} \end{align}
\end{lem}
\begin{pf}
(1) Applying $\partial_{\mm{h}}^i$ to  \eqref{s0106pnnnn},  we have
 \begin{equation}\label{20171102237}  \left\{\begin{array}{ll}
{\rho}\partial_{\mm{h}}^i u_t =\partial_{\mm{h}}^i(\partial_{\bar{M}}^2\eta- \mm{div}\mathcal{S}  (q,u,\eta) +\mathcal{N}^1)&\mbox{ in }  \Omega,\\[1mm]
\mathrm{div}\partial_{\mm{h}}^i {u} =-\partial_{\mm{h}}^i\mathrm{div}_{\tilde{\mathcal{A}}} {u} & \mbox{ in }  \Omega,\\
 \llbracket \partial_{\mm{h}}^i \eta  \rrbracket=\llbracket \partial_{\mm{h}}^i u  \rrbracket=0&\mbox{ on }\Sigma, \\
  \partial_{\mm{h}}^i  \llbracket  \mathcal{S}  (q,u,\eta)    e_3
- g \rho   \eta_3 e_3- \bar{M}_3 \partial_{\bar{M}}\eta   \rrbracket=
 \partial_{\mm{h}}^i (  \vartheta\Delta_{\mm{h}}\eta_3 e_3+  \mathcal{N}^2    )&\mbox{ on }\Sigma ,  \\
 \partial_{\mm{h}}^i(\eta,u)=0 &\mbox{ on }\Sigma_-^+.
\end{array}\right.\end{equation}
Multiplying \eqref{20171102237}$_1$  by $\partial^i_\mm{h}\eta$ in $L^2$ yields that
\begin{align}
\frac{\mm{d}}{\mm{d}t} \int {\rho} \partial_{\mm{h}}^i\eta \cdot \partial_{\mm{h}}^i u \mm{d}y =&\int  \partial_{\mm{h}}^i(  \partial_{\bar{M}}^2\eta-\mm{div}\mathcal{S}  (q,u,\eta )) \cdot  \partial_{\mm{h}}^i\eta\mm{d}y \nonumber  \\
&+\int { {\rho} } |\partial_{\mm{h}}^i u|^2\mm{d}y +\int  \partial_{\mm{h}}^i{\mathcal{N}}^1 \cdot  \partial_{\mm{h}}^i\eta\mm{d}y =:\sum_{ k=1}^3 I_{k}. \label{201808100154628}
\end{align}

Exploiting the integration by parts, and the boundary conditions in \eqref{20171102237}, we get
 \begin{align}
I_{1}
= &\int_\Sigma
   \partial_{\mm{h}}^i\llbracket\mathcal{S}  (q,u,\eta )   e_3 -\bar{M}_3\partial_{\bar{M}}\eta\rrbracket
  \cdot \partial_{\mm{h}}^i\eta\mm{d}y_{\mm{h}} + \int( \partial_{\mm{h}}^i \mathcal{S}  (q,u,\eta )  : \nabla  \partial_{\mm{h}}^i\eta-  |\partial_{\mm{h}}^i \partial_{\bar{M}}\eta|^2)\mm{d}y
\nonumber \\
 =&   g  \llbracket \rho \rrbracket |\partial_{\mm{h}}^i \eta_3 |_0^2-  \mathcal{I}(\partial_{\mm{h}}^i \eta)-\frac{1}{4}\frac{\mm{d}}{\mm{d}t}\int \mu | \mathbb{D}\partial_{\mm{h}}^i\eta|^2\mm{d}y+ \int \partial_{\mm{h}}^i q\mm{div}\partial_{\mm{h}}^i\eta\mm{d}y +I_4,\nonumber
\end{align}
where we have defined that
\begin{equation*}
I_{4}:=\int_\Sigma \partial_{\mm{h}}^i    \mathcal{N}^2   \cdot \partial_{\mm{h}}^i\eta\mm{d}y_{\mm{h}}.
\end{equation*}
Thus we have
\begin{align}
& \frac{\mm{d}}{\mm{d}t}\int\left( {\rho} \partial_{\mm{h}}^i \eta \cdot \partial_{\mm{h}}^i u +\frac{ \mu }{4}| \mathbb{D} \partial_{\mm{h}}^i \eta|^2\right) \mm{d}y+
  \mathcal{I}(\partial_{\mm{h}}^i \eta)\nonumber  \\
  &\leqslant c\left(|\partial_{\mm{h}}^i \eta_3|_0^2+\|  \partial_{\mm{h}}^i u\|^2_{0} \right) +\left\{
              \begin{array}{ll}
  \int q\mm{div} \eta\mm{d}y  +|I_3|+|I_4|      &\mbox{ for }i= 0;\\
\| \partial_{\mm{h}}\mm{div}\eta\|_0 \| \partial_{\mm{h}} q \|_0+ |I_3|+|I_4|&\mbox{ for } i=1;\\
\|\partial_{\mm{h}}^2  \mm{div}\eta\|_0\|\partial_{\mm{h}}^2 \nabla q \|_0+ |I_3|+|I_4|&\mbox{ for } (i,\vartheta)=(2,0).\end{array}
            \right.
   \label{estimforhoedsds1stnsafdaf}
\end{align}

Exploiting the partial integrations and \eqref{201808121247}, we get
\begin{align}
|I_3|=
 \left|\int \partial_{\mm{h}}^i{\mathcal{N}^1}  \cdot  \partial_{\mm{h}}^i\eta\mm{d}y\right|\lesssim  \|\partial_{\mm{h}}^i \eta\|_0 \| \mathcal{N}^1\|_i
   \mbox{ for }0\leqslant i\leqslant 1        \nonumber \end{align}
 and
\begin{align}
|I_4|    \lesssim\left\{
                 \begin{array}{ll}
     | \eta |_0  |   \mathcal{N}^2  |_0
 &\mbox{ for } i=0; \\
 |\partial_{\mm{h}}
\eta |_{1/2}  |    \mathcal{N}^2   |_{1/2}&\mbox{ for }i=1.
                 \end{array}
               \right.
 \nonumber
\end{align}
If $(i,\vartheta)=(2,0)$, by \eqref{201808121247} and the partial integrations, we can further estimate that
\begin{align}
&|I_3|
    \lesssim  \|\nabla_{\mm{h}}\partial_{\mm{h}}^2 \eta\|_0\|\mathcal{N}^1  \|_1 \mbox{ and }|I_4|  =
    \left|\int_\Sigma\partial_{\mm{h}}^2{\mathcal{N}^2}  \cdot  \partial_{\mm{h}}^2\eta\mm{d}y_{\mm{h}}\right|
        \lesssim |\partial_{\mm{h}}^2\eta|_{1/2}  |     \mathcal{N}^2   |_{3/2}.
  \nonumber
\end{align}
In addition, we have (see (3.32) in \cite{JFJWGCOSdd} for the derivation)
  \begin{equation}
\begin{aligned}
 \int q\mm{div} \eta\mm{d}y = -\int  \psi\cdot\nabla   q\mm{d}y-\int_\Sigma  \llbracket
 q  \rrbracket e_3\cdot  \psi\mm{d}y_{\mm{h}}
\lesssim    \|\eta\|_2^2(  \|\nabla q\|_0 + | \llbracket
q  \rrbracket|_{1/2}) ,
\end{aligned} \nonumber
\end{equation}
where
 $$  \psi:=\left(\begin{array}{c}
         - \eta_1(\partial_2\eta_2+\partial_3\eta_3 )+
           \eta_1(\partial_2\eta_3 \partial_3\eta_2
- \partial_2\eta_2\partial_3\eta_3) \\
      \eta_1\partial_1\eta_2  - \eta_2\partial_3\eta_3+
        \eta_1( \partial_1\eta_2 \partial_3\eta_3
-\partial_1\eta_3\partial_3\eta_2)\\
          \eta_1\partial_1\eta_3
+\eta_2\partial_2\eta_3+\eta_1(\partial_1\eta_3 \partial_2\eta_2
-\partial_1\eta_2\partial_2\eta_3)
        \end{array}\right)\mbox{ and }\mm{div}\eta=\mm{div}\psi. $$
Thus, inserting all the estimates above into \eqref{estimforhoedsds1stnsafdaf} yields that
\begin{align}
&\frac{\mm{d}}{\mm{d}t}\int\left( {\rho} \partial_{\mm{h}}^i \eta \cdot \partial_{\mm{h}}^i u +\frac{ \mu }{4}| \mathbb{D} \partial_{\mm{h}}^i \eta|^2\right) \mm{d}y+
  \mathcal{I}(\partial_{\mm{h}}^i \eta)\nonumber \\
 & \leqslant c\left(|\eta_3|_i^2+\|  u\|^2_{i,0}\right)\nonumber \\
&\quad  + \left\{
                \begin{array}{ll}
\sqrt{\mathcal{E}}\mathcal{D} +\|\eta\|_0 \|\mathcal{N}^1  \|_0+ |
\eta|_{0 } |    \mathcal{N}^2    |_{0}& \mbox{ for }i=0 ; \\
\|\mm{div}\eta\|_1\| \partial_{\mm{h}} q\|_0+\|\partial_{\mm{h}}\eta\|_0 \|\mathcal{N}^1  \|_1+ |\partial_{\mm{h}}
\eta|_{1/2} |      \mathcal{N}^2  |_{1/2}& \mbox{ for }i=1 ; \\
\|\mm{div}\eta\|_2\|\partial_{\mm{h}}^2  q\|_0+ \|\nabla \partial_{\mm{h}}^2 \eta\|_0 \|\mathcal{N}^1  \|_1+ |\partial_{\mm{h}}^2
\eta|_{1/2} |    \mathcal{N}^2   |_{3/2}
& \mbox{ for } (i,\vartheta)=(2,0).
                \end{array}
              \right.
 \label{estimforhoedsds1stn}
\end{align}

Finally, making using of \eqref{improtian1}, \eqref{06011733}--\eqref{060117sdfa34}, trace estimate and Young's inequality, we immediately get \eqref{ssebdaiseqinM0846} from \eqref{estimforhoedsds1stn}.

(2) Applying the operator $\mathfrak{D}_{\mf{h}}^{3/2}$ to \eqref{20171102237} with $i=1$, and then augmenting as \eqref{estimforhoedsds1stnsafdaf}, we have
\begin{align}
&
\frac{\mm{d}}{\mm{d}t} \int \left(\rho\mathfrak{D}_{\mf{h}}^{3/2} \partial_{\mm{h}} \eta \cdot\mathfrak{D}_{\mf{h}}^{3/2}  \partial_{\mm{h}}  u + \frac{ \mu }{4} | \mathbb{D} \mathfrak{D}_{\mf{h}}^{3/2}    \partial_{\mm{h}}  \eta  |^2\right) \mm{d}y\nonumber \\
&+\vartheta \left( |\nabla_{\mm{h}} \partial_{\mm{h}} \eta_3|_{(y_{\mm{h}},y_3)=(\mf{h},0)} |^2 +
|\nabla_{\mm{h}}\mathfrak{D}_{\mf{h}}^{3/2}\partial_{\mm{h}}  \eta_3 |_{0 }^2\right)\nonumber \\
    & \lesssim   \vartheta  |\nabla_{\mm{h}} \partial_{\mm{h}} \eta_3 |_{(y_{\mm{h}},y_3)=(\mf{h},0)}|^2  + |\mathfrak{D}_{\mf{h}}^{3/2}  \partial_{\mm{h}} \eta_3|_{0 }^2+\|\mathfrak{D}_{\mf{h}}^{3/2} \partial_{\mm{h}}  u\|^2_{0 } +\|\mathfrak{D}_{\mf{h}}^{3/2}  \partial_{\mm{h}}  \mm{div} \eta\|_{0 } \|\mathfrak{D}_{\mf{h}}^{3/2}  \partial_{\mm{h}} q\|_{0 }\nonumber
\\
 &\quad+ \left|\int \mathfrak{D}_{\mf{h}}^{3/2} \partial_{\mm{h}}  \eta\cdot \mathfrak{D}_{\mf{h}}^{3/2} \partial_{\mm{h}}  \mathcal{N}^1  \mm{d}y  +\int_{\Sigma} \mathfrak{D}_{\mf{h}}^{3/2}  \partial_{\mm{h}} \eta \cdot\mathfrak{D}_{\mf{h}}^{3/2}    \partial_{\mm{h}}  \mathcal{N}^2 \mm{d}y_{\mm{h}}\right|
:=I_5(\mf{h}). \label{ssebdaiseqinM0846asdfadfad}
 \end{align}

Making use of \eqref{201807312054}, \eqref{201806291928}, \eqref{201808121247}  and the integration by parts,  we can estimate that
\begin{align}
\int_{\mathcal{T}_{4,4}}I_5(\mf{h})\mm{d}\mf{h}\lesssim   & |\nabla_{\mm{h}}\partial_{\mm{h}} \eta_3|_0^2+\|  u\|_{1,1}^2 +|\partial_{\mm{h}}\eta_3 |_{1/2}^2\nonumber \\
&+\|\mm{div}\eta\|_2\|\nabla q\|_1+\| \eta\|_3\| \mathcal{N}^1  \|_1+|\nabla_{\mm{h}}\partial_{\mm{h}}\eta|_{1/2}| \mathcal{N}^2 |_{1/2}.\nonumber
\end{align}
Thus, integrating \eqref{ssebdaiseqinM0846asdfadfad} with respect to $\mf{h}$ over $\mathcal{T}_{4,4}$, and making use of \eqref{improtian1}, \eqref{06011733}, \eqref{06011734}, the first inequality in  \eqref{201807312054}, trace estimate   and Young's inequality, we can deduce  \eqref{ssebdaisedsaqinM0dsfadsff846asdfadfad} from \eqref{ssebdaiseqinM0846asdfadfad}.

(3) Finally, we turn to derive \eqref{ssebdaiseqinM0asdfa846}. For $i=2$, using the integration by parts, we have
\begin{equation}
\label{201807022000}
I_{1}
=   -  {\|\sqrt{\kappa\rho} \mathbb{D} \partial_{\mm{h}}^2 \eta  \|^2_0/2}-
\|\partial_{\bar{M}}\partial_{\mm{h}}^2 \eta \|^2_0-\frac{1}{4}\frac{\mm{d}}{\mm{d}t}\int \mu | \mathbb{D}\partial_{\mm{h}}^2\eta|^2\mm{d}y+I_6 ,
\end{equation}
where we have defined that
\begin{equation*}
\begin{aligned}
I_{6}:=&-\int_\Sigma  \partial_{\mm{h}} \llbracket q-2(\mu\partial_3u_3+\kappa\rho\partial_3\eta_3)-
\bar{M}_3\partial_{\bar{M}}\eta_3 \rrbracket  \partial_{\mm{h}}^3\eta_3\mm{d}y_{\mm{h}}   \\
&+\int_\Sigma \partial_{\mm{h}}^2  \mathcal{N}^2_{\mm{h}}     \cdot \partial_{\mm{h}}^2\eta_{\mm{h}}\mm{d}y_{\mm{h}}+ \int \partial_{\mm{h}}^2 q\partial_{\mm{h}}^2\mm{div}\eta\mm{d}y .
\end{aligned}
\end{equation*}
Putting \eqref{201807022000} into \eqref{201808100154628} for $i=2$, and then using the integration by parts, we have
\begin{align}
& \frac{\mm{d}}{\mm{d}t}\int  \left( \rho \partial_\mm{h}^2 \eta \cdot  \partial_\mm{h}^2 u  + \frac{\mu }{4}  |\mathbb{D}\partial_\mm{h}^2 \eta |^2\right) \mm{d}y+
 {\|\sqrt{\kappa\rho} \mathbb{D}\partial_{\mm{h}}^2 \eta\|^2_0/2}+
\|\partial_{\bar{M}}\partial_{\mm{h}}^2 \eta  \|^2_0\nonumber \\
 & \lesssim   \|   u\|^2_{2,0}+\|\eta\|_3\|\mathcal{N}^1\|_1+ |I_6|.\label{201807022042}
\end{align}

Noting that
$$
 \llbracket \partial_{\mm{h}}^\alpha \eta  \rrbracket =\llbracket \partial_{\mm{h}}^\alpha u  \rrbracket =0 \mbox{ on }\Sigma,
$$
thus, making use of \eqref{improtian1}, \eqref{201808181500},  \eqref{201807241645},  \eqref{201808121247}, trace estimate  and Young's inequality, we can estimate that
\begin{align}
|I_6|\lesssim & |\partial_{\mm{h}}^2 \eta_3|_{1/2}(|\partial_{\mm{h}} ( \llbracket \mu(\mm{div}u-\mm{div}_{\mm{h}}u_{\mm{h}})\rrbracket ,\llbracket \kappa\rho( \mm{div}\eta-\mm{div}_{\mm{h}}\eta_{\mm{h}})\rrbracket )|_{1/2}
\nonumber \\
&+| \partial_{\mm{h}} \mm{div}\eta|_{1/2}
+|  \llbracket \partial_{\mm{h}}  q \rrbracket  |_{1/2})+ |\eta|_{5/2} |  \mathcal{N}^2_{\mm{h}}  |_{3/2}+\|\mm{div} \eta\|_2\|\nabla q\|_1
\nonumber \\
\lesssim& |\partial_{\mm{h}}^2 \eta_3|_{1/2}(|\nabla_{\mm{h}}\partial_{\mm{h}} (\eta_{\mm{h}},u_{\mm{h}})|_{1/2} +| \llbracket \nabla_{\mm{h}}q\rrbracket |_{1/2})+\sqrt{\mathcal{E}}\mathcal{D}. \label{2018051214552}
\end{align}
Finally, putting the above estimate and \eqref{06011733} into \eqref{201807022042} yields \eqref{ssebdaiseqinM0asdfa846}.
\hfill$\Box$
\end{pf}


\begin{lem}\label{201612132242nxsfs}Under the assumption \eqref{aprpiosesnew},
the following estimates hold:
\begin{align}
&
 \label{201702061418}
 \frac{\mm{d}}{\mm{d}t}(\|\sqrt{\rho} \partial_\mm{h}^i u\|^2_0+
  \mathcal{I}(\partial_{\mm{h}}^i \eta ))
+ c\|\partial_\mm{h}^i   u \|_{1}^2 \lesssim
\| \eta_3\|_{\sigma(i),1}^2+ \sqrt{\mathcal{E} }\mathcal{D}\mbox{ for }\left\{
                                   \begin{array}{ll}
                         0\leqslant i\leqslant 1 ;\\
                              ( i,\vartheta)=(2, 0),
                                   \end{array}
                                 \right.  \\
&\frac{\mm{d}}{\mm{d}t}\int_{\mathcal{T}_{4,4}}\left( \|\sqrt{\rho}   \mathfrak{D}_{\mf{h}}^{3/2} \partial_{\mm{h}}  u\|^2_0+
  \mathcal{I}( \mathfrak{D}_{\mf{h}}^{3/2} \partial_{\mm{h}}  \eta  )\right)\mm{d}\mathbf{h} \lesssim  \|\eta_3\|_{1,1} \|u_3\|_{1,1} + \sqrt{\mathcal{E}}\mathcal{D},
 \label{ssebdaisedsaqinM0dsfadsff846sdfaaasdfadfad}\\&\frac{\mm{d}}{\mm{d}t}   \|\partial_{\mm{h}}^2 (\sqrt{\rho}u, \sqrt{\kappa\rho} \mathbb{D} \eta_3,
\partial_{\bar{M}} \eta_3 )\|^2_0 + c\|\partial_\mm{h}^2 u\|_1^2
\nonumber \\
  & \lesssim
  |\partial_{\mm{h}}^2 u_3|_{1/2}
 \left( |\nabla_{\mm{h}}\partial_{\mm{h}} (\eta_{\mm{h}},u_{\mm{h}}) |_{1/2}  +|\llbracket \nabla_{\mm{h}} q \rrbracket |_{1/2}\right)+ \sqrt{\mathcal{E}} \mathcal{D} ,
 \label{2018072033}
 \end{align}
where we have defined that $\sigma(0)=0$ and $\sigma(i)=i-1$ for $i\neq0$.
\end{lem}
\begin{pf}
(1) Multiplying \eqref{20171102237}$_1$  by $\partial^i_{\mm{h}}u$,  and then
integrating the resulting  identity over $\Omega$, we have
\begin{equation}
  \frac{1}{2}\frac{\mm{d}}{\mm{d}t} \int {\rho}   | \partial_{\mm{h}}^i u|^2 \mm{d}y = \int\partial_{\mm{h}}^{i} (\partial_{\bar{M}}^2\eta- \mm{div} \mathcal{S}  (  q, u,\eta )) \cdot  \partial_{\mm{h}}^iu\mm{d}y   +\int  \partial_{\mm{h}}^i{\mathcal{N}}^1 \cdot  \partial_{\mm{h}}^iu\mm{d}y.
  \label{201807022044}
  \end{equation}
Exploiting \eqref{diverelation}, \eqref{s0106pnnnn}$_3$, \eqref{06041533fwqg}, product estimate \eqref{fgestims}, trace estimate and the integration by parts,  we have
$$
 \int  q\mm{div} u\mm{d}y=
 \int  \nabla q  \cdot (\tilde{\mathcal{A}}^{\mm{T}} {u} )\mm{d}y +\int_\Sigma   \llbracket q\rrbracket \tilde{\mathcal{A}}e_3\cdot {u}\mm{d}y
\lesssim \|\eta\|_3\|u\|_1 (\| \nabla q\|_0+|\llbracket q\rrbracket|_{1/2}) . $$Thus, following the argument of \eqref{estimforhoedsds1stn}, and using   Korn's inequality, we have
\begin{align}
& \frac{\mm{d}}{\mm{d}t} \left(\int  {\rho} |\partial_{\mm{h}}^i u|^2 \mm{d}y  +
  \mathcal{I}(\partial_{\mm{h}}^i \eta)\right)+c\| \partial_{\mm{h}}^i u\|^2_1\lesssim \left|
\int_\Sigma \partial_{\mm{h}}^i \eta_3 \partial_{\mm{h}}^iu_3\mm{d}y_{\mm{h}}
\right| \nonumber \\
&+   \left\{
                                   \begin{array}{ll}
\sqrt{\mathcal{E}}\mathcal{D}+ \|u\|_0 \|\mathcal{N}^1  \|_0+ |
u|_0|  \mathcal{N}^2  |_0  &\mbox{ for }   i=0 ;\\
 \|\mathrm{div} {u} \|_{1} \|\partial_{\mm{h}}  q\|_{0}+\|\partial_{\mm{h}}  u\|_0 \|\mathcal{N}^1  \|_1+ |\partial_{\mm{h}}
u|_{1/2} |    \mathcal{N}^2    |_{1/2} &\mbox{ for }            i=1 ;\\
  \|\mathrm{div}  {u} \|_2 \| \partial_{\mm{h}}^2 q\|_{0}+   \|\nabla_{\mm{h}}\partial^2_{\mm{h}}u\|_0 \|\mathcal{N}^1  \|_1+ |
 \partial_{\mm{h}}^2 u|_{1/2 } |      \mathcal{N}^2  |_{3/2} &                           \mbox{ for }(i,\vartheta)=(2, 0).
                                   \end{array}
                                 \right.
 \label{estimforhoedsds1stnn1524}\end{align}

Using  trace estimate and \eqref{201808121247},
we have
\begin{align}
\left|
\int_\Sigma \partial_{\mm{h}}^i \eta_3 \partial_{\mm{h}}^iu_3\mm{d}y_{\mm{h}}
\right|\lesssim &\left\{
                   \begin{array}{ll}
     | \partial_{\mm{h}}^{i-1} \eta_3|_{1/2} |\partial_{\mm{h}}^iu_3|_{1/2}\lesssim
\| \eta_3\|_{i-1,1} \|\partial_{\mm{h}}^i u_3\|_{1}   & \hbox{ for }1\leqslant i\leqslant 2; \\
  |\eta_3|_0|u_3|_0\lesssim \|\eta_3\|_1\|u_3\|_1 & \hbox{ for }i=0.
                   \end{array}
                 \right. \nonumber
\end{align}
Consequently, putting the above estimate into \eqref{estimforhoedsds1stnn1524}, and then using \eqref{201808181500}, \eqref{06011733}--\eqref{060117sdfa34}, trace estimate  and  Young's inequality, we get the desired conclusion \eqref{201702061418} from \eqref{estimforhoedsds1stnn1524}.

(2) Similar to \eqref{ssebdaiseqinM0846asdfadfad}, we can derive from \eqref{20171102237} with $i=2$  that
\begin{align}
&\frac{1}{2}\frac{\mm{d}}{\mm{d}t}\left( \int \rho
| \mathfrak{D}_{\mf{h}}^{3/2} \partial_{\mm{h}}  u|^2\mm{d}y+
  \mathcal{I}(\mathfrak{D}_{\mf{h}}^{3/2} \partial_{\mm{h}}  \eta )\right)  \nonumber
  \\  & \lesssim |\mathfrak{D}_{\mf{h}}^{3/2} \partial_{\mm{h}} \eta_3|_0 |\mathfrak{D}_{\mf{h}}^{3/2} \partial_{\mm{h}}  u_3|_0  + \|\mathfrak{D}_{\mf{h}}^{3/2} \partial_{\mm{h}} \mathrm{div} u \|_0\|\mathfrak{D}_{\mf{h}}^{3/2} \partial_{\mm{h}}   q\|_0\nonumber \\
&\qquad
+\left|\int \mathfrak{D}_{\mf{h}}^{3/2} \partial_{\mm{h}} u\cdot \mathfrak{D}_{\mf{h}}^{3/2}\partial_{\mm{h}}  \mathcal{N}^1 \mm{d}y+
\int_{\Sigma}\mathfrak{D}_{\mf{h}}^{3/2} \partial_{\mm{h}}   u\cdot  \mathfrak{D}_{\mf{h}}^{3/2} \partial_{\mm{h}}   \mathcal{N}^2 \mm{d}y_{\mm{h}}\right|.\nonumber
 \end{align}
Thus, integrating the above inequality  with respect to $\mf{h}$ over $\mathcal{T}_{4,4}$, and then following the argument of \eqref{ssebdaisedsaqinM0dsfadsff846asdfadfad} from \eqref{ssebdaiseqinM0846asdfadfad}, we easily get \eqref{ssebdaisedsaqinM0dsfadsff846sdfaaasdfadfad} from the above estimate.

(3) Finally, we turn to derive \eqref{2018072033}.
Similar to \eqref{201807022042}, using Korn's  inequality, we can derive from \eqref{201807022044} for $i=3$ that
 \begin{align}
  \frac{\mm{d}}{\mm{d}t}  \|\partial_{\mm{h}}^2 (\sqrt{\rho}u, \sqrt{\kappa\rho} \mathbb{D} \eta /\sqrt{2},
\partial_{\bar{M}} \eta )\|^2_0 + c\|\partial_\mm{h}^2 u\|_1^2   \lesssim
 \|u\|_3\|\mathcal{N}^1\|_1+|I_7|. \label{2018070220257}
 \end{align}
 where we have defined that
\begin{equation*}
\begin{aligned}
I_7:=&-\int_\Sigma \partial_{\mm{h}}  \llbracket q-2(\mu\partial_3u_3+ \kappa\rho\partial_3\eta_3)-\bar{M}_3\partial_{\bar{M}}\eta_3  \rrbracket \partial_{\mm{h}}^3 u_3\mm{d}y_{\mm{h}}\\
   &+\int_\Sigma \partial_{\mm{h}}^2    \mathcal{N}^2_{\mm{h}}   \cdot \partial_{\mm{h}}^2 u_{\mm{h}}\mm{d}y_{\mm{h}}+ \int \partial_{\mm{h}}^2 q\partial_{\mm{h}}^2\mm{div} u\mm{d}y .
   \end{aligned}
\end{equation*}
Following the argument of \eqref{2018051214552}, we get
\begin{align}
|I_7|\lesssim &   |\partial_{\mm{h}}^2 u_3|_{1/2}(|\nabla_{\mm{h}}\partial_{\mm{h}} (\eta_{\mm{h}},u_{\mm{h}})|_{1/2}+|\llbracket \nabla_{\mm{h}} q \rrbracket|_{1/2}) +\sqrt{\mathcal{E}}\mathcal{D}.
\nonumber
\end{align}
 Thus, putting the above estimate and \eqref{06011733} into \eqref{2018070220257} yields \eqref{2018072033}.
\hfill$\Box$
 \end{pf}
\begin{lem}\label{201612132242nxsfssdfs}Under the assumption \eqref{aprpiosesnew},
the following estimates hold:
\begin{align}
&  \frac{\mm{d}}{\mm{d}t}\left(\|\sqrt{\rho}  u_{t}\|_{0}^2+ \mathcal{I}(  u)
-\int \nabla q\cdot (\mathcal{A}_t^{\mm{T}}u)\mm{d}y-\int_\Sigma \llbracket
q\rrbracket \mathcal{A}_te_3\cdot u\mm{d}y_{\mm{h }}
\right)
+ c\| u_t\|_{1}^2\nonumber \\
& \lesssim  |  u_3|_{ 0}^2+ \sqrt{\mathcal{E} } \mathcal{D} . \label{Lem:030dsfafds1m0dfsf832}
 \end{align}
\end{lem}
\begin{pf}
 Applying $ \partial_t$ to \eqref{n0101nn}$_2$--\eqref{n0101nn}$_6$, and  using \eqref{s0106pnnnn}$_1$, we have that
\begin{equation}\label{06051331}
\left\{\begin{array}{ll}
{\rho}  u_{tt}+\mm{div}_{\ml{A}}\mathfrak{S}_{\mathcal{A}} ( q_t,  u_t,  u)
=\partial_{\bar{M}}^2    u+\mathcal{N}^{ 3} &\mbox{ in }  \Omega,\\[1mm]
\mathrm{div}_{\mathcal{A}}  {u}_{t} =-\mm{div}_{\mathcal{A}_t}u   & \mbox{ in }  \Omega ,\\
\llbracket  u_t  \rrbracket= 0 &\mbox{ on }\Sigma,\\
   \llbracket \mathfrak{S}_{\mathcal{A}} (  q_t, u_t, u )
    \mathcal{A}e_3-g \rho  u_3 e_3
-\bar{M}_3\partial_{\bar{M}} u\rrbracket  = \vartheta   \Delta_{\mm{h}}u_3e_3+ \mathcal{N}^4 &\mbox{ on }\Sigma,\\
  u_t=0 &\mbox{ on }\Sigma_-^+.
\end{array}\right.\end{equation}

By \eqref{20178011456115} with $(i,j)=(1,0)$,
\begin{align}
\int q_t  \mm{div}_{\mathcal{A}_t}u\mm{d}y=&-
\frac{\mm{d}}{\mm{d}t}\left( \int \nabla q\cdot (\mathcal{A}_t^{\mm{T}}u)\mm{d}y+\int_\Sigma \llbracket
q\rrbracket \mathcal{A}_te_3\cdot u\mm{d}y_{\mm{h }}
\right)\nonumber \\
& +\int  \nabla q  \cdot \partial_t(\mathcal{A}_t^{\mm{T}}u)  \mm{d}y
  +\int_\Sigma   \llbracket q  \rrbracket \partial_t\left(\mathcal{A}_t e_3\cdot u \right) \mm{d}y_{\mm{h}}.
\label{201811031645}
\end{align}
Thus, multiplying \eqref{06051331}$_1$ by $  u_{t}$ in $L^2$, and using the integration by parts, \eqref{06051331}$_2$--\eqref{06051331}$_5$ and \eqref{201811031645}, we can compute out that
\begin{align}
&\frac{1}{2}\frac{\mm{d}}{\mm{d}t}\left(\|\sqrt{\rho} u_{t}\|^2_0+\mathcal{I}( u)- \int \nabla q\cdot (\mathcal{A}_t^{\mm{T}}u)\mm{d}y+\int_\Sigma \llbracket
q\rrbracket \mathcal{A}_te_3\cdot u\mm{d}y_{\mm{h }}\right)+\frac{1}{2}\|\sqrt{\mu}\mathbb{D}_{\mathcal{A}} u_{t}\|_{0}^2\nonumber \\
&=\int \left(\mathcal{N}^{ 3}\cdot   u_t -\kappa\rho \mathbb{D}  u:\nabla_{\tilde{\mathcal{A}}}   u_{t}-  \nabla q  \cdot \partial_t(\mathcal{A}_t^{\mm{T}}u)  \right)\mm{d}y\nonumber \\
 &\quad +\int_{\Sigma}\left(  \mathcal{N}^4  \cdot  u_t+g\llbracket \rho\rrbracket u_3  \partial_t u_3 -  \llbracket q  \rrbracket \partial_t( \mathcal{A}_t e_3\cdot u)\right)\mathrm{d} y_{\mm{h}}
   = :K_1. \label{0425sdfasfa}
\end{align}

Making use of \eqref{06142100x}, \eqref{06142100}, \eqref{06041533fwqg},   product estimates \eqref{fgestims} and \eqref{0604221702}, and trace estimate, we can estimate that \begin{align}
K_1  \lesssim| u_3|_0\|  u_t\|_1+\|u\|_2| \llbracket q \rrbracket|_{1/2}|\mathcal{A}_{tt}e_3|_{-1/2}+ \|  u_t\|_0\|\mathcal{N}^{ 3}\|_0+ \|  u_t\|_1 |   \mathcal{N}^4 |_0+\sqrt{\mathcal{E}}\mathcal{D} .
\nonumber
\end{align}
Putting the above estimate into \eqref{0425sdfasfa}, and then using \eqref{2018081015456}, \eqref{201806291} and Korn's and Young's inequalities, we get  \eqref{Lem:030dsfafds1m0dfsf832}.
  \hfill$\Box$
 \end{pf}
\begin{lem}\label{201612132242nx}
Under the assumption \eqref{aprpiosesnew},
the following estimate hold:
\begin{align}
  &\|(u,q)\|_{\mm{S},0} \lesssim |(\kappa,\bar{M}_3)|\|\eta\|_{2} + \|\eta\|_{\underline{2},1}+\|u_t\|_0.\label{2017020614181721}
   \end{align}
\end{lem}
\begin{pf}
We rewrite \eqref{s0106pnnnn}$_2$--\eqref{s0106pnnnn}$_6$  as a stratified   Stokes problem:
\begin{equation}\label{n0101nn928}\left\{\begin{array}{ll}
\nabla  q-\mu\Delta  u= \mathcal{M}^1 &\mbox{ in }  \Omega, \\[1mm]
\mathrm{div}  u =\mathrm{div}  u & \mbox{ in }  \Omega,\\
  \llbracket  u \rrbracket=0,\ \llbracket  (qI-\mu \mathbb{D}u)e_3\rrbracket =  \mathcal{M}^2  &\mbox{ on }\Sigma, \\
 u=0 &\mbox{ on }\Sigma_-^+ ,
\end{array}\right.\end{equation}
where we have defined that
\begin{align}
&\mathcal{M}^1:=\kappa\rho\Delta \eta+ \partial_{\bar{M}}^2\eta-\rho u_t+\nabla  (\mu \mathrm{div} {u}+\kappa\rho  \mm{div}\eta) +\mathcal{N}^1,\nonumber\\
& \mathcal{M}^2:=\llbracket\kappa\rho\mathbb{D}\eta e_3+    g \rho \eta_3e_3+\bar{M}_3\partial_{\bar{M}}\eta\rrbracket+\vartheta \Delta_{\mm{h}}\eta_3  e_3 +  \mathcal{N}^2.
\end{align}

Applying  the  stratified  Stokes estimate   to \eqref{n0101nn928}, we have
\begin{align}
& \|  (u,q)\|_{\mm{S},0}  \lesssim
\|   \mathcal{M}^1 \|_0+\| \mathrm{div} {u}\|_{ 1}+|   \mathcal{M}^2   |_{ 1/2}.\label{20171111841}
\end{align}
Using   \eqref{improtian1}--\eqref{201808181500} and trace estimate,
we can estimate that
$$
\begin{aligned}
&\|\mathcal{M}^1 \|_0+\|  \mathrm{div} {u}\|_{1}+|  \mathcal{M}^2   |_{1/2}\\
&\lesssim |(\kappa,\bar{M}_3)|\|\eta\|_{2}  +\|\eta \|_{\underline{2},1}+\|u_t\|_0+\|\eta\|_3\|u\|_{2} + \| \mathcal{N}^1 \|_0+|  \mathcal{N}^2  |_{1/2}.
\end{aligned}$$
Putting the above  estimate  into \eqref{20171111841}, and then using \eqref{06011733}  and \eqref{06011734}, we get \eqref{2017020614181721}.
\hfill $\Box$
\end{pf}

\begin{lem}
\label{lem:dfifessim2057}
Let $(\eta,u)$ satisfy  the assumption \eqref{aprpiosesnew}.
\begin{enumerate}
 \item[(1)] It holds that
 \begin{align}
& \frac{\mm{d}}{\mm{d}t}
{\left\|\sqrt{(\kappa\rho+\bar{M}_3^2)/\mu}\eta\right\|}_{2}^2+
c\left\|  \left((\kappa\rho+\bar{M}_3^2)\eta,u\right)\right\|_{2}^2+
\|\nabla q \|_{0}^2  \nonumber \\
& \lesssim  \|(\eta,u)\|_{\underline{1},1}^2+\| u_{t}\|_{0}^2+\mathcal{E}\mathcal{D}.  \label{20808011436} \end{align}
\item[(2)] If $(\kappa,\bar{M}_3)\neq 0$, we have
 \begin{align}
& \frac{\mm{d}}{\mm{d}t}
\overline{\|\eta\|}_3^2+
\left\|  \left( \eta,u\right)\right\|_3^2+
\|\nabla q \|_1^2+| \llbracket \nabla q \rrbracket |_{1/2}^2 \lesssim  \|(\eta,u)\|_{\underline{2},1}^2+\| u_t\|_1^2,
 \label{201611111952} \end{align}
   where the norm $\overline{\|\eta\|}_{3}$ is equivalent to $\|\eta\|_{3}$.
  \item[(3)]
  If $(\kappa,\bar{M}_3)=0$, we have
 \begin{align}
&\left\|u\right\|_3+
\|\nabla q \|_1+| \llbracket \nabla q \rrbracket |_{1/2}  \lesssim \|\eta\|_{2,1}+ \|  u \|_{\underline{2},1}+\| u_{t}\|_1+\|\eta\|_3^2,\label{20161111195sdfafsd2}\\
&\left\|u\right\|_{\underline{1},2}+
\|\nabla q \|_{\underline{1},0}  +
\| \llbracket\nabla_{\mm{h}} q \rrbracket\|_{1/2}  \lesssim \|(\eta,u)\|_{\underline{2},1}+\| u_{t}\|_{1}+\sqrt{\mathcal{E}\mathcal{D}}.\label{201808292010}  \end{align}
 \end{enumerate}
 \end{lem}
\begin{pf}
We shall rewrite \eqref{n0101nn928} with $i=0$ as the two one-layer (steady) Stokes problems:
\begin{equation}\label{n0101nn928n}\left\{\begin{array}{ll}
\nabla q_\pm-\mu_\pm\Delta \omega_\pm= \mathcal{M}^3_\pm&\mbox{ in }  \Omega_\pm, \\[1mm]
\mathrm{div} \omega_\pm = \mathcal{N}_\pm^5 & \mbox{ in }  \Omega_\pm,\\
\omega_\pm=\omega_\pm|_{\Sigma}& \mbox{ on }  \Sigma, \\
 \omega_\pm=0 &\mbox{ on }\Sigma_\pm,
\end{array}\right.\end{equation}
where we have defined that
$$
\begin{aligned}
&\mathcal{M}^3_\pm=\mathcal{M}^3\chi_{\Omega_\pm},\
\mathcal{M}^3 := ( \partial_{M}^2\eta- \bar{M}_3^2\Delta\eta)-\rho u_t+\nabla(\mu\mm{div} u+\kappa\rho \mm{div}\eta) +\mathcal{N}^1,\\
&  \omega_\pm:=\omega\chi_{\Omega_\pm},\ \omega := u +(\kappa \rho  + \bar{M}_3^2)\eta /\mu,\ \mathcal{N}^5_\pm:=\mathcal{N}^5\chi_{\Omega_\pm},
\end{aligned}
$$
and  $(\mathcal{N}^5_\pm, \omega_\pm|_{\Sigma})$ satisfies
that
$$\int_{\Omega_\pm} \mathcal{N}_\pm^5\mm{d}y=\int_\Sigma \omega_\pm|_{\Sigma}\cdot (\mp e_3 )\mm{d}y_{\mm{h}}.$$

Let the integers $i$ and $j$  satisfy $0\leqslant i\leqslant 1$ and $0\leqslant j\leqslant i$.
Applying $\partial_{\mm{h}}^j $ to \eqref{n0101nn928n} yields that
\begin{equation*}\label{n0101nn928nadfa}\left\{\begin{array}{ll}
\nabla \partial_{\mm{h}}^j   q_\pm-\mu_\pm\Delta \partial_{\mm{h}}^j \omega_\pm= \partial_{\mm{h}}^j  \mathcal{M}^3_\pm&\mbox{ in }  \Omega_\pm, \\[1mm]
\mathrm{div}\partial_{\mm{h}}^j   \omega_\pm =\partial_{\mm{h}}^j  \mathcal{N}_\pm^5 & \mbox{ in }  \Omega_\pm,\\
\partial_{\mm{h}}^j \omega_\pm=\partial_{\mm{h}}^j \omega_\pm|_{\Sigma}& \mbox{ on }  \Sigma, \\
 \partial_{\mm{h}}^j  \omega_\pm=0 &\mbox{ on }\Sigma_\pm,
\end{array}\right.\end{equation*}
Applying the one-layer Stokes estimate \eqref{xfsddfsf201705141252} to the above Stokes problems,  we have
\begin{align}
 \| \omega \|_{j,i+2-j }^2+
\|\nabla    q \|_{j,i-j }^2  \lesssim |   \omega  |_{i+3/2  }^2 +
\|    \mathcal{M}^3 \|_{j,i-j }^2+\|   \mathcal{N}^5 \|_{j,i+1-j }^2.
 \label{20171111029}
 \end{align}

(1)
By \eqref{n0101nn}$_1$, we obviously see that
\begin{align}
\|  \omega \|_{j,i+2-j}^2=&  \frac{\mm{d}}{\mm{d}t}
\left\| \sqrt{(\kappa\rho  + \bar{M}_3^2)/\mu}   \eta\right\|_{j,i+2-j}^2 \nonumber \\
&+\|(\kappa\rho +\bar{M}_3^2)    \eta/\mu \|_{j,i+2-j}^2+
\|    u\|_{j,i+2-j}^2.\label{201808101923}\end{align}
Exploiting   \eqref{improtian1} and \eqref{201808181500}, we have
\begin{align}
& \|\mathcal{M}^3\|_{j,i-j}\lesssim \bar{M}_3\|\eta\|_{j+1,i+1-j}+\|\eta\|_{j+2,i-j}+\|u_t\|_{j,i-j}
+\|\eta\|_3\|(\eta,u)\|_{i+2}+\|\mathcal{N}^1\|_i. \label{201808101951}
\end{align}
Making use of \eqref{06011733}, \eqref{201807241959}, \eqref{201808101923} and \eqref{201808051835}, we can  derive \eqref{20808011436} from \eqref{20171111029} with $(i,j)=(0,0)$.

(2)  By \eqref{201808101923} and \eqref{201808101951}, we can derive from \eqref{20171111029} that
  \begin{align}
& \frac{\mm{d}}{\mm{d}t}
\left\| \sqrt{(\kappa\rho  + \bar{M}_3^2)/\mu}  \eta\right\|_{j,i+2-j}^2+\|(\kappa\rho +\bar{M}_3^2 )  \eta/\mu \|_{j,i+2-j}^2+
\|    u\|_{j,i+2-j}^2+
\|  \nabla   q\|_{j,i-j}^2\nonumber\\
&\lesssim
\|  \eta \|_{j+1,i+1-j}^2+
|  \omega  |_{i+3/2}^2+\|u_{t}\|_{j,i-j}^2 + \|\eta\|_3^2\|(\eta,u)\|_3^2 +\|\mathcal{N}^1\|_1^2+\|\mathcal{N}^5\|_{2}^2.\label{20180627sdafa1120safaf}
 \end{align}

Since $(\kappa,\bar{M}_3)\neq 0$,
we can derive from \eqref{20180627sdafa1120safaf} that
\begin{align}
& \frac{\mm{d}}{\mm{d}t}
\overline{\left\|\eta\right\|}_{j,i+2-j}^2+\|(\eta,u)\|_{j,i+2-j}^2
+\|\nabla q \|_{j,i-j}^2 \nonumber
\\
&\lesssim
\|  \eta \|_{i+1,1}^2+
| \omega  |_{i+3/2}^2+\|u_{t}\|_{j,i-j}^2 + \|\eta\|_3^2\|(\eta,u)\|_3^2 +\|\mathcal{N}^1\|_{1}^2 +\|\mathcal{N}^5\|_{2}^2,\label{201807231906sdfsaf}
 \end{align}
where $\overline{\left\|\eta\right\|}_{j,i+2-j}^2:=\sum_{0\leqslant k\leqslant i-j}a_k \left\| \sqrt{(\kappa\rho  + \bar{M}_3^2)/\mu}  \eta\right\|_{j+k,i+2-j-k}^2$, and $a_k$ are positive constants.
Thus, putting \eqref{06011733} and \eqref{201807241959} into \eqref{201807231906sdfsaf} with $(i,j)=(1,0)$, and then using trace estimate, we  get \eqref{201611111952}.

 (3) If $(\kappa,\bar{M}_3)=0$, then we can deduce from \eqref{20171111029} and \eqref{201808101951} that
 \begin{align}
 \| u \|_{j,i+2-j }+
\|\nabla   q \|_{j,i-j } \lesssim &
| u  |_{i+3/2}+\|\eta\|_{j+2,i-j}+\|u_t\|_{j,i-j}\nonumber \\
&+ \|\eta\|_3\|(\eta,u) \|_3+\| \mathcal{N}^1 \|_1+\|  \mathcal{N}^5 \|_2.
 \nonumber
 \end{align}
Making use of \eqref{06011733},  \eqref{201807241959} and trace estimate, we can easily deduce \eqref{20161111195sdfafsd2}--\eqref{201808292010}  from the above estimate. \hfill$\Box$
\end{pf}

\begin{lem}
\label{lem:dfifessim2057sadfafdas}
Under the assumption \eqref{aprpiosesnew}, we have
\begin{equation} \mathcal{E}\mbox{ is equivalent to }\|\eta\|_3^2+\|u\|_2^2 \label{201702071610nb}
\end{equation}
\end{lem}
\begin{pf}
By \eqref{2017020614181721} and Friedrichs's inequality, we see that
$$\|\eta\|_3^2+\|u\|_2^2\lesssim  \mathcal{E}\lesssim \|\eta\|_3^2+\|u\|_2^2+\| u_t\|_0^2.$$
To get \eqref{201702071610nb}, it suffices to verify
\begin{equation}
\label{201807311651}
\|  u_t\|_0^2\lesssim \|\eta\|_3^2+\|u\|_2^2.
\end{equation}

 We get from \eqref{n0101nn}$_2$ that
\begin{equation}\label{06051331asdfasfd}
{\rho}  u_{t }+\nabla q
=\partial_{\bar{M}}^2 \eta- \mm{div}_{\mathcal{A}}\mathcal{S}_{\mathcal{A}} (0,u,\eta)-\nabla_{\tilde{\mathcal{A}}} q.\end{equation}
Multiplying \eqref{06051331asdfasfd} by $\partial_t(\mathcal{A}^{\mm{T}}u)$ in $L^2$ and using the integral by parts and \eqref{20178011456115} with $(i,j)=(0,0)$, we get
\begin{align}
 \int \rho| u_{t}|^2\mm{d}y=K_6+K_7,
\label{201806271923}
\end{align}
where we have defined that
\begin{align}
K_6:=&
  \int_{\Sigma} \llbracket   q \rrbracket \partial_t(\mathcal{A}^{\mm{T}}u)\cdot e_3 \mm{d}y_{\mm{h}}
,\nonumber \\
K_7:=& \int  \left(\left(\partial_{\bar{M}}^2  \eta-
 \mm{div}_{\ml{A}}\mathcal{S}_{\mathcal{A}} (0,u,\eta)-\nabla_{\tilde{\mathcal{A}}} q
\right) \cdot  \partial_t(\mathcal{A}^{\mm{T}}u)-\rho u_t \partial_t(\tilde{\mathcal{A}}^{\mm{T}}u) \right)\mm{d}y.\nonumber
\end{align}

Exploiting \eqref{201806202203}, \eqref{06201929201811041930} and trace estimate, we can estimate that
\begin{align}
K_6
\lesssim &  |\llbracket   q \rrbracket |_{1/2}|\partial_t(\mathcal{A}^{\mm{T}}u)\cdot e_3|_{-1/2}\lesssim
  (\|u\|_1+\|u_t\|_0)  (\|\eta\|_3+\|u\|_2+|   \mathcal{N}^2_3   |_{1/2}) .
 \nonumber
\end{align}
In addition, making use of  \eqref{06041533fwqg} and product estimate \eqref{fgestims}, we can estimate that
\begin{align}
K_7 \lesssim (\|u\|_1+\|  u_t\|_0)(\|(\eta,u)\|_2+\|\eta\|_3\|\nabla q\|_0)+\|\eta\|_3\|u_t\|_0^2.\nonumber
\end{align}
Putting the above two estimates into \eqref{201806271923} yields that
\begin{align}
\| u_{t}\|^2_0\lesssim(\|u\|_1+\|  u_t\|_0)(\|\eta\|_3+\|u\|_2+\|\eta\|_3\|\nabla q\|_0+|   \mathcal{N}^2_3   |_{1/2})+ \|\eta\|_3\|  u_t\|_0^2 .\nonumber
\end{align}
Inserting  \eqref{06011734} and  \eqref{2017020614181721} into the above estimate, and then using Young's inequality, we get \eqref{201807311651}.
This completes the proof of Lemma \ref{lem:dfifessim2057sadfafdas}.
  \hfill$\Box$
\end{pf}

\subsection{Highest-order  estimates of $u_3$ at interface for $\vartheta\neq 0$}

In this subsection, we further establish the highest-order boundary estimate of $u_3$.
It is difficult to directly derive the desired estimate based on the ART problem.
Motivated by \cite{XLZPZZFGAR}, we shall divide the nonhomogeneous form of the ART problem into two subproblems.

Let $({\eta}^1,{u}^1,{q}^1):=({\eta},{u},{q})-({\eta}^2,{u}^2,{q}^2)$,
 where $({\eta}^2,{u}^2,{q}^2)$ is a solution to the linear problem
 \begin{equation}\label{s0106pndfsgnnnxx}
\left\{
\begin{array}{ll}
 \eta_t^2=  u^2 &\mbox{ in }\Omega,  \\[1mm]
 \rho u^2_{t}+\mm{div}\mathcal{S}  (q^2,u^2,\eta^2)-\partial_{\bar{M}}^2\eta^2 =\mathcal{N}^1 &\mbox{ in }\Omega, \\[1mm]
\div u^2  =   \mm{div}  u&\mbox{ in }\Omega, \\
  \llbracket \eta^2  \rrbracket=\llbracket u^2  \rrbracket=0&\mbox{ on }\Sigma,\\
 \llbracket  \mathcal{S}  (q^2,u^2,\eta^2)   e_3-  \bar{M}_3 \partial_{\bar{M}}\eta^2
\rrbracket= g \llbracket \rho\rrbracket \eta_3  e_3+
(\mathcal{N}_1^2,\mathcal{N}_2^2,\mathfrak{M}_2)^{\mm{T}}  &\mbox{ on }\Sigma, \\[1mm]
(\eta^2,u^2)=0&\mbox{ on }\Sigma_{-}^+,\\[1mm]
(\eta^2,u^2)|_{t=0}=(\eta^0,u^0)&\mbox{ on }\Omega.
\end{array}\right.
\end{equation}We will mention the existence of unique solution of the above problem \eqref{s0106pndfsgnnnxx} in next subsection.
Then $({\eta}^1,{u}^1,{q}^1)$ satisfies
\begin{equation}\label{s0106pndfsgnnn}
\left\{
\begin{array}{ll}
 \eta_t^1=  u^1 &\mbox{ in }\Omega,  \\[1mm]
 \rho u_{t}^1+\mm{div}\mathcal{S}  (q^1,u^1,\eta^1)-\partial_{\bar{M}}^2\eta^1 =0 &\mbox{ in }\Omega, \\[1mm]
\div u^1  =0&\mbox{ in }\Omega, \\
  \llbracket \eta^1  \rrbracket=\llbracket u^1 \rrbracket=0&\mbox{ on }\Sigma,\\
 \llbracket  \mathcal{S}  (q^1,u^1,\eta^1)  e_3-  \bar{M}_3 \partial_{\bar{M}}\eta^1
\rrbracket-\vartheta\Delta_{\mm{h}}\eta_3^1 e_3
=  (\mathfrak{M}_1+\vartheta\Delta_{\mm{h}}\eta_3^2)e_3=:\widetilde{\mathfrak{M}}_1e_3 &\mbox{ on }\Sigma, \\[1mm]
(\eta^1,u^1)=0&\mbox{ on }\Sigma_{-}^+,\\[1mm]
(\eta^1,u^1)|_{t=0}=0&\mbox{ on }\Omega.
\end{array}\right.
\end{equation}
Then we can use the above two auxiliary problems to derive the highest-order boundary estimate of $u_3$.
\begin{lem}\label{201807251431}
 Let $\vartheta\neq0$ and $\eta^0$ further satisfy $H^0\in H^1(\mathbb{T})$. Under the assumption \eqref{aprpiosesnew}, for any given $\Lambda>0$,
  there exists a non-negative functional $ \mathfrak{E}_L(t) $ of solutions $(\eta^2,u^2)$ and $(\eta^2,u^2)$ (see \eqref{201808301358xx} for the detailed definition), such that
\begin{align}
& \mathfrak{E}_L(t)+ \int_0^t  | \nabla_{\mm{h}}^2 u_3|^2_{1/2} \mm{d}\tau \nonumber \\
& \leqslant   c\left(\|\eta^0\|_3^2+\|u^0\|_2^2+|H^0|_1^2\right)+  {\Lambda}  \int_0^t  \mathfrak{E}_L(t) \mm{d}\tau \nonumber \\
&  \quad + c\int_0^t\left( \|\eta_3\|_{1,1}^2+| \eta_3|_2^2+ \Upsilon \left(\Upsilon +\|\partial_\mm{h}^2 q\|_{0}\right)+ \sqrt{\mathcal{E} }\mathcal{D}\right)\mm{d}\tau
\label{201807012235}
\end{align}
where we have defined that  $ | \nabla_{\mm{h}}^2 u_3|^2_{1/2}:= \sum_{|\alpha|=2}| \partial_{\mm{h}}^\alpha u_3|^2_{1/2}$.
\end{lem}
\begin{pf}
(1) Let $0\leqslant i\leqslant 1$.
Applying $\partial_{\mm{h}}^i\partial_t$ to \eqref{s0106pndfsgnnn} yields
\begin{equation}\label{sadfdxsafdsa}
\left\{\begin{array}{ll}
 \rho   \partial_{\mm{h}}^i u_{tt}^1+  \partial_{\mm{h}}^i(\mm{div}\mathcal{S}  (q^1_t,u^1_t,u^1 ) -\partial_{\bar{M}}^2u^1)=0 &\mbox{ in }\Omega, \\[1mm]
\div   \partial_{\mm{h}} u^1  =\div u_t^1  =0&\mbox{ in }\Omega, \\
  \llbracket   \partial_{\mm{h}}  u^1  \rrbracket=\llbracket    u^1_t \rrbracket=0&\mbox{ on }\Sigma,\\
   \partial_{\mm{h}}^i\llbracket  \mathcal{S}  (q^1_t,u^1_t,u^1 )    e_3 -\bar{M}_3 \partial_{\bar{M}}u^1\rrbracket
 =\partial_{\mm{h}}^i   \partial_t  \widetilde{\mathfrak{M}}_1  e_3 &\mbox{ on }\Sigma,\\
  ( \partial_{\mm{h}}u^1, u^1_t)=0&\mbox{ on }\Sigma_{-}^+.\\[1mm]
\end{array}\right.
\end{equation}

Multiplying \eqref{sadfdxsafdsa}$_1$ with $i=0$ by $  u^1_t$ in $L^2$, and then using the integration by parts, we can estimate that
\begin{align}
 \frac{1}{2}\frac{\mm{d}}{\mm{d}t}\left(\int  \rho |   u^1_t |^2\mm{d}y +\mathcal{I}(  u^1)\right)+\frac{1}{2} \|\sqrt{\mu}\mathbb{D}  \partial_{\mm{h}} u^1_t\|^2_0 \lesssim |\partial_t   u^1_3|_{1/2} | \partial_t   \widetilde{\mathfrak{M}}_1|_{-1/2} . \label{201811041954}
\end{align}
Similarly to \eqref{201807311651}, we also derive from
 \eqref{s0106pndfsgnnn} that
$$\| u_t^1\|_0^2\lesssim \|\eta^1\|_3^2+\|u^1\|_2^2  +| \widetilde{\mathfrak{M}}_1 |_{1/2}^2,$$
which, together with \eqref{06011734}, \eqref{201702071610nb}, the zero initial data \eqref{s0106pndfsgnnn}$_7$ and the fact $H^0:=H|_{t=0}= \widetilde{\mathfrak{M}}_1 |_{t=0}/\vartheta $,  yields that
\begin{equation}
\label{201811041957}
\| u_t^1|_{t=0}\|_0^2 \lesssim  |H^0|_{1/2}^2.
\end{equation}
Integrating \eqref{201811041954} with respect to $t$, and using \eqref{201811041957}, and Korn's and Young's inequalities, we further have
\begin{align}
\label{201807011935}
 \|\sqrt{\rho}   u^1_t \|^2_0 +\mathcal{I}(   u^1) +c \int_0^t \| u^1_\tau\|^2_{1}\mm{d}\tau\lesssim  |H^0|_{1/2}^2+  \int_0^t  |  \partial_\tau  \widetilde{ \mathfrak{M}}_1 |_{-1/2}^2 \mm{d}\tau.
\end{align}

Multiplying \eqref{sadfdxsafdsa}$_1$ with $i=1$ by $  \partial_{\mm{h}}  u^1$ in $L^2$, and using the integration by parts, we have
\begin{align}
&\frac{\mm{d}}{\mm{d}t}\int\left(\frac{\mu}{4}|\mathbb{D}   \partial_{\mm{h}} u^1 |^2- \rho    u_{t}^1  \partial_{\mm{h}}^2 u^1 \right)\mm{d}y
+\mathcal{I}(  \partial_{\mm{h}} u^1)\nonumber
\\
&\lesssim \int|   \partial_{\mm{h}}  u_{t}^1|^2\mm{d}y+|  \nabla_{\mm{h}} \partial_{\mm{h}}  u^1_3|_0 |  \partial_t  \widetilde{\mathfrak{M}}_1 |_0 . \nonumber
\end{align}
Integrating the above identity  with respect to $t$, and then using
 \eqref{s0106pndfsgnnn}$_7$ and \eqref{201811041957}, and
Korn's  and Young's inequalities, we further have
\begin{align}
c  \| \partial_{\mm{h}} u^1\|^2_1
+\int_0^t \mathcal{I}(  \partial_{\mm{h}} u^1(\tau))\mm{d}\tau\lesssim \int_0^t\left(\|   \partial_{\mm{h}}   u_{\tau}^1\|^2_0+
|   \partial_{\tau}  \widetilde{  \mathfrak{M}}_1 |_0^2 \right)\mm{d}\tau +\| u_t^1\|_0^2. \label{201807011913}
\end{align}
Thus we can deduce from \eqref{201807011935} and  \eqref{201807011913} that
\begin{align}
 \int_0^t |\nabla_{\mm{h}} \partial_{\mm{h}} u^1_3|_0^2\mm{d}\tau
 \lesssim  |H^0|_{1/2}^2+ \int_0^t   |   \partial_\tau   \widetilde{\mathfrak{M}}_1 |_0^2   \mm{d}\tau. \label{2001807012105}
\end{align}

Applying the operator $\mathfrak{D}_{\mf{h}}^{3/2}$ to \eqref{sadfdxsafdsa} with $i=1$, thus, following the argument of \eqref{2001807012105},  we have
\begin{align}
  \int_0^t  | \mathfrak{D}_{\mf{h}}^{3/2}\nabla_{\mm{h}}\partial_{\mm{h}} u^1_3|_0^2\mm{d}\tau
 \lesssim  |\mathfrak{D}_{\mf{h}}^{3/2}H^0|_{1/2}^2+  \int_0^t  | \mathfrak{D}_{\mf{h}}^{3/2}\partial_\tau  \widetilde{ \mathfrak{M}}_1 |_0^2 \mm{d}\tau
  .\nonumber
\end{align}
Integrating the above integral over $\mathcal{T}_{4,4}$, adding the resulting estimate and \eqref{2001807012105} together, and then using \eqref{201807312054} and \eqref{201811132119}, we get
\begin{align}
 \int_0^t | \nabla_{\mm{h}}\partial_{\mm{h}}  u^1_3|^2_{1/2}\mm{d}\tau
 \lesssim  |H^0|_{1}^2 + \int_0^t   ( | \partial_\tau   \mathfrak{M}_1 |_{1/2}^2 +
 | \Delta_{\mm{h}} u_3^2 |_{1/2}^2 )  \mm{d}\tau. \label{20018070saf121sfafsad05}
\end{align}

(2) Noting that $\mm{div}\eta^1=\mm{div}u^1=0$,
 thus, following the argument of \eqref{ssebdaiseqinM0846}--\eqref{ssebdaiseqinM0asdfa846}, we can derive from \eqref{s0106pndfsgnnn} that, for $0\leqslant  i\leqslant 1$,
\begin{align}
& \frac{\mm{d}}{\mm{d}t}\int\left( {\rho} \partial_{\mm{h}}^i \eta^1 \cdot \partial_{\mm{h}}^i u^1 +\frac{ \mu }{4}| \mathbb{D} \partial_{\mm{h}}^i \eta^1|^2\right) \mm{d}y+
  \mathcal{I}(\partial_{\mm{h}}^i \eta^1)  \lesssim \|  u^1\|^2_{i,0} +   | \eta^1_3|_{3/2} |   \widetilde{\mathfrak{M}}_1  |_{1/2}  ,\nonumber  \\
&\frac{\mm{d}}{\mm{d}t}\int_{\mathcal{T}_{4,4}}\int  \left( \rho  \mathfrak{D}_{\mf{h}}^{3/2} \partial_{\mm{h}} \eta^1 \cdot   \mathfrak{D}_{\mf{h}}^{3/2} \partial_{\mm{h}}  u ^1 + \frac{ \mu}{4}  |\mathbb{D} \mathfrak{D}_{\mf{h}}^{3/2} \partial_{\mm{h}} \eta^1|^2\right) \mm{d}y\mm{d}\mf{h}+
\vartheta|\nabla_{\mm{h}}\partial_{\mm{h}} \eta_3^1|_{1/2}^2
 \nonumber \\
    & \lesssim |  \eta_3^1|_2^2+\| u^1\|_{1,1}^2 + |\eta^1_3|_{5/2}|  \widetilde{\mathfrak{M}}_1  |_{1/2} \nonumber ,\\
 &
 \frac{\mm{d}}{\mm{d}t}\int  \left( \partial_\mm{h}^2 \eta ^1\cdot  \partial_\mm{h}^2 u^1  + \frac{ \mu}{4}    | \mathbb{D}\partial_\mm{h}^2 \eta^1 |^2\right) \mm{d}y+
 \kappa \|\partial_{\mm{h}}^2 \eta^1\|_1^2+\|\partial_{\bar{M}}\partial_{\mm{h}}^2  \eta^1\|_{0}^2\nonumber \\
  & \lesssim
  \|   u^1\|^2_{2,0}+ |\partial_{\mm{h}}^2 \eta_3^1|_{1/2}(|\nabla_{\mm{h}}\partial_{\mm{h}} (\eta^1_{\mm{h}} ,u_\mm{h}^1)|_{1/2}+| \llbracket \nabla_{\mm{h}}q ^1\rrbracket \|_{1/2}) .\nonumber
\end{align}
Similarly, following the argument of \eqref{201702061418}--\eqref{2018072033} and \eqref{Lem:030dsfafds1m0dfsf832},    we can derive from \eqref{s0106pndfsgnnn} that
\begin{align}
& \frac{\mm{d}}{\mm{d}t} \left(\|\sqrt{\rho}  \partial_{\mm{h}}^i u^1\|^2_0  +
  \mathcal{I}(\partial_{\mm{h}}^i \eta^1)\right)+c\| \partial_{\mm{h}}^i u^1\|^2_1 \lesssim | u^1_3 |_{3/2} | \widetilde{\mathfrak{M}}_1  |_{1/2}  \mbox{ for }0\leqslant  i\leqslant 1, \nonumber \\
&\frac{\mm{d}}{\mm{d}t}\int_{\mathcal{T}_{4,4}}\left( \|\sqrt{\rho}   \mathfrak{D}_{\mf{h}}^{3/2} \partial_{\mm{h}}  u^1\|^2_0+
  \mathcal{I}( \mathfrak{D}_{\mf{h}}^{3/2} \partial_{\mm{h}} \eta^1  )\right)\mm{d}\mathbf{h} \lesssim  \|   \eta_3^1\|_{1,1}\| u_3^1\|_{1,1} + |\nabla_{\mm{h}}\partial_{\mm{h}}u^1_3|_{1/2} |  \widetilde{\mathfrak{M}}_1  |_{1/2},\nonumber
\\&\frac{\mm{d}}{\mm{d}t} \|\partial_{\mm{h}}^2 (\sqrt{\rho}u^1, \sqrt{\kappa\rho} \mathbb{D} \eta_3^1,
\partial_{\bar{M}} \eta_3^1 )\|^2_0 + c\|\partial_\mm{h}^2 u^1\|_1^2
\nonumber \\
  & \lesssim
  | \partial_{\mm{h}}^2u_3^1|_{1/2}
 \left( |\nabla_{\mm{h}}\partial_{\mm{h}} (\eta^1_{\mm{h}},u_{\mm{h}}^1) |_{1/2}  +|\llbracket \nabla_{\mm{h}} q ^1 \rrbracket |_{1/2}\right)  ,\nonumber  \\
&\frac{\mm{d}}{\mm{d}t}(\|\sqrt{\rho} u_{t}^1\|_{0}^2+ \mathcal{I}(  u^1))
+ c\|  u_t^1\|_{1}^2 \lesssim  |\partial_t\widetilde{\mathfrak{M}}_1 |_{1/2}^2 .\nonumber
 \end{align}
Using trace estimate and Young's and Korn's inequalities, we can derive from the above seven estimates that, for any sufficiently large ${c}_1$,
 \begin{align}
& \frac{\mm{d}}{\mm{d}t} \mathfrak{E}_1(t)+  c{c}_1\left( \|   u_t^1\|_{ 1}^2+  \| u^1\|_{\underline{2},1}^2\right)\nonumber
\\
&\leqslant  ( | \llbracket \nabla_{\mm{h}}q ^1\rrbracket |_{1/2}^2/ {{c}_1}+\|  \eta^1\|_{2,1}^2
+{c}_1^3 \|  \eta^1\|_{\underline{1},1}^2
+{c}^5_1  |(\nabla_{\mm{h}}\partial_{\mm{h}} u_3^1,   \widetilde{\mathfrak{M}}_1  ,   \partial_t\widetilde{\mathfrak{M}}_1  )|_{1/2}^2 )/c,\label{201808301318} \\
&
c\left({c}_1  \| \eta^1 \|_{2,1}^2+ {c}_1^4  \| \eta^1 \|_{\underline{1},1}^2\right) \leqslant \mathfrak{E}_1(t), \label{201810011343}
 \end{align}
where we have defined that
\begin{align}
\mathfrak{E}_1(t):=&{c}_1^4\Bigg(\sum_{ |\alpha|  \leqslant 1}
\left(\int   {\rho}\partial_{\mm{h}}^{\alpha}\eta^1 \cdot \partial_{\mm{h}}^{\alpha} u^1  \mm{d}y +{c}_1 \mathcal{I}(\partial_{\mm{h}}^{\alpha} \eta^1)
\right)+ \frac{ 1 }{4}\|\sqrt{\mu} \mathbb{D}  \eta^1\|^2_{\underline{1},0}+{c}_1  \|\sqrt{\rho}  u^1\|^2_{\underline{1},0}\Bigg)\nonumber \\
&
 +{c}_1^3\int_{\mathcal{T}_{4,4}}  \sum_{ |\alpha|= 1}\left(\int \rho  \mathfrak{D}_{\mf{h}}^{3/2,\alpha}\eta^1 \cdot   \mathfrak{D}_{\mf{h}}^{3/2,\alpha}u ^1 \mm{d}y +{c}_1
  \mathcal{I}( \mathfrak{D}_{\mf{h}}^{3/2,\alpha}\eta^1  )\right) \mm{d}\mf{h}\nonumber
\\
&+{c}_1^3
\int_{\mathcal{T}_{4,4}}\left(\frac{ 1}{4}  \|\sqrt{\mu}\mathbb{D} \mathfrak{D}_{\mf{h}}^{3/2} \eta^1\|^2_{1,0} + {c}_1 \|\sqrt{\rho}   \mathfrak{D}_{\mf{h}}^{3/2} u^1\|^2_{1,0} \right) \mm{d}\mf{h}
\nonumber \\
&+c_1\left(  \sum_{ |\alpha| =2}\int  \partial_{\mm{h}}^{\alpha} \eta ^1\cdot  \partial_{\mm{h}}^{\alpha} u^1 \mm{d}y + \frac{1}{4}  \| \sqrt{\mu}\mathbb{D} \eta^1 \|^2_{2,0}\right) \nonumber\\
&+{c}_1 \Bigg(\|\sqrt{\rho}u_{t}^1\|_{0}^2+  \mathcal{I}(   u^1)+ c_1 \| (\sqrt{\rho}u^1, \sqrt{\kappa\rho} \mathbb{D} \eta_3^1,
\partial_{\bar{M}} \eta_3^1 )\|^2_{2,0} \Bigg).\nonumber
\end{align}

Following the argument of \eqref{201808292010} and \eqref{201807231906sdfsaf}, we  can derive from \eqref{s0106pndfsgnnn} that
\begin{align}
& \frac{\mm{d}}{\mm{d}t}
 \widetilde{\| \eta^1\|}_{\underline{1},2}^2+
\left\|  \left( (\kappa\rho+\bar{M}_3^2)\eta^1,u^1\right)\right\|_{\underline{1},2}^2+
\|\nabla q^1 \|_{\underline{1},0}^2+| \llbracket \nabla_{\mm{h}} q^1\rrbracket |_{1/2}^2 \nonumber \\
&\lesssim  \|(\eta^1,u^1)\|_{\underline{2},1}^2+\| u_t^1\|_{1}^2, \label{201808301318sss}
\end{align}
where
$$\widetilde{\|\eta^1\|}_{\underline{1},2} \left\{
                                             \begin{array}{ll}
\mbox{is equivalent to }\|\eta^1\|_{\underline{1},2} & \hbox{ for }(\kappa,\bar{M}_3)\neq 0;  \\
=0, &  \hbox{ for }(\kappa,\bar{M}_3)=0.
                                             \end{array}
                                           \right.
 $$
Then, by choosing a sufficienlty large ${c}_1$, we derive from \eqref{201808301318} and \eqref{201808301318sss}  that
 \begin{align}
&\frac{\mm{d}}{\mm{d}t} \left( \mathfrak{E}_1(t) +\widetilde{\| \eta^1\|}_{\underline{1},2}^2\right)+c  \|\nabla q^1 \|_{\underline{1},0}^2\nonumber  \\
 & \leqslant
\left( \|  \eta^1\|_{2,1}^2
+{c}_1^3 \|  \eta^1\|_{\underline{1},1}^2+ c_1^5|( \nabla_{\mm{h}}\partial_{\mm{h}}  u^1_3,\widetilde{\mathfrak{M}} ,\partial_t \widetilde{\mathfrak{M}}) |_{1/2}^2 \right)/c.  \label{201808302128} \end{align}

By \eqref{201811041957} and the zero initial data condition of $(\eta^1,u^1)$, we see that
$$\left.\left(\mathfrak{E}_1 +(\kappa\rho+\bar{M}_3^2)
\overline{\| \eta^1\|}_{\underline{1},2}^2\right)\right|_{t=0}\leqslant cc_1
 \left.\|u^1_t\|_0^2 \right|_{t=0}\leqslant cc_1( \|\eta^0\|_3^2+\|u^0\|_2^2 +|H^0|_1^2).$$
Thus, integrating \eqref{201808302128} over $(0,t)$, and then using  \eqref{20018070saf121sfafsad05}, \eqref{201810011343} and trace estimate, we get
 \begin{align}
 &  {c}_1  \| \eta^1 \|_{2,1}^2+ {c}_1^4  \| \eta^1 \|_{\underline{1},1}^2    + \int_0^t \left(c_1^5 | \nabla_{\mm{h}}^2 u^1_3|^2_{1/2}+\|\nabla q^1 \|_{\underline{1},0}^2\right)\mm{d}\tau \nonumber
\\
&  \leqslant c c_1^5\left(\|\eta^0\|_3^2+\|u^0\|_2^2+|H^0|_1^2\right)   \nonumber \\
&\quad +c\int_0^t  \left( \|  \eta^1\|_{2,1}^2
+{c}_1^3 \|  \eta^1\|_{\underline{1},1}^2 + c_1^5( |(\mathfrak{M}_1,\partial_{\tau} \mathfrak{M}_1)|_{1/2}^2+|(\eta^2, u^2 )  |_{2,1}^2 \right)\mm{d}\tau. \label{201808310859}
\end{align}

(3)
Following the argument of \eqref{estimforhoedsds1stn} with $(i,\vartheta)=(2,0)$, we
can derive from \eqref{s0106pndfsgnnnxx} that
\begin{align}
&\frac{\mm{d}}{\mm{d}t}\int  \left( \rho\partial_\mm{h}^2 \eta^2 \cdot  \partial_\mm{h}^2 u^2  + \frac{ \mu  }{4}| \mathbb{D}\partial_\mm{h}^2 \eta^2 |^2\right) \mm{d}y \nonumber  \\
& \lesssim
| \eta_3|_2| \eta_3^2|_2+ \|   u^2\|^2_{2,0} +\|\mm{div}\eta^2\|_{2}
\|\partial_\mm{h}^2   q^2\|_{0}+\|\nabla_{\mm{h}}\partial_{\mm{h}}^2 \eta^2\|_{0}\|\mathcal{N}^1\|_{1}\nonumber \\
&\quad +|\partial_{\mm{h}}^2 \eta^2|_{1/2}|(\mathcal{N}_1^2,\mathcal{N}_2^2,\mathfrak{M}_2)|_{3/2}.\nonumber
\end{align}
Using \eqref{06011733}, \eqref{060117sdfa34} and trace estimate, we further have
 \begin{align}
&\int  \left( \rho\partial_\mm{h}^2 \eta^2 \cdot  \partial_\mm{h}^2 u^2  + \frac{ \mu  }{4}| \mathbb{D}\partial_\mm{h}^2 \eta^2 |^2\right) \mm{d}y \nonumber \\
& \lesssim \int_0^t\left(| \eta_3|_2\| \eta_3^2\|_{2,1}+\|u^2\|_{2,0}^2 + \|\mm{div}\eta^2\|_2\|\partial_\mm{h}^2 q^2\|_{0}\right.\nonumber \\
&\qquad \left.+\|\eta^2\|_{2,1} \|\eta\|_3\sqrt{\mathcal{D}}\right)\mm{d}\tau+
\|\eta^0\|_3^2+\|u^0\|_2^2. \label{201811050925}
\end{align}
In addition, following the argument of   \eqref{estimforhoedsds1stnn1524} with $(i,\vartheta)=(2,0)$, we can derive from \eqref{s0106pndfsgnnnxx} that
 \begin{align}
 &\frac{\mm{d}}{\mm{d}t}(\|\sqrt{\rho} \partial_\mm{h}^2 u^2\|^2_0+
   \kappa \|\partial_{\mm{h}}^2 \eta\|_1^2+\|\partial_{\bar{M}}\partial_{\mm{h}}^2  \eta^2\|_{0}^2 )
+ c\|\partial_\mm{h}^2   u^2 \|_{1}^2\nonumber  \\
&\lesssim
\| \eta_3\|_{1,1}^2+\|\mm{div}u \|_{2}
\|\partial_\mm{h}^2   q^2\|_{0}+\|\nabla_{\mm{h}}\partial_{\mm{h}}^2 u^2\|_{0}\|\mathcal{N}^1\|_{1}+|\partial_{\mm{h}}^2 u^2|_{1/2}|(\mathcal{N}_1^2,\mathcal{N}_2^2,\mathfrak{M}_2 )|_{3/2}.\nonumber
\end{align}
Thus, similarly to \eqref{201811050925}, making use of  \eqref{201808181500}, \eqref{06011733}, \eqref{060117sdfa34}, trace estimate and Young's inequality, we further have
 \begin{align}
&\| \partial_\mm{h}^2 u^2\|^2_0
+  \int_0^t\|\partial_\mm{h}^2  u^2 \|_{1}^2\mm{d}\tau\nonumber \\
& \lesssim
\|\eta^0\|_3^2+\|u^0\|_2^2+\int_0^t\left(\| \eta_3\|_{1,1}^2+ \|\eta\|_3\|u\|_3\|\partial_\mm{h}^2 q^2\|_{0}+ \sqrt{\mathcal{E} }\mathcal{D}\right)\mm{d}\tau. \label{201811050925xyz}
\end{align}
Combining \eqref{201811050925} with \eqref{201811050925xyz} yields that, for sufficient large $c_1$,
 \begin{align}
 &c\left( c_1\left( \|  \partial_\mm{h}^2 u^2\|^2_0 +  \int_0^t\|\partial_\mm{h}^2  u^2 \|_{1}^2\mm{d}\tau\right)+
\|\eta^2\|^2_{2,1} \right)\nonumber  \\
&\leqslant   \int_0^t\left(c_1\|\eta_3\|_{1,1}^2+| \eta_3|_2\| \eta_3^2\|_{2,1}+ \left(c_1\|\eta\|_3\|u\|_3+\|\mm{div}\eta^2\|_2\right)\|\partial_\mm{h}^2 q^2\|_{0}+ c_1\sqrt{\mathcal{E} }\mathcal{D}\right. \nonumber\\
&\qquad + \left. \|\eta^2\|_{2,1} \|\eta\|_3\sqrt{\mathcal{D}}\right)\mm{d}\tau + c_1\left( \|\eta^0\|_3^2+\|u^0\|_2^2\right).  \label{xx201808310859}
\end{align}

(4) Now we let
\begin{equation}
  \mathfrak{E}_L(t):=  \| (\eta-\eta^2,\eta^2) \|_{2,1}^2+ {c}_1^3  \| \eta- \eta^2  \|_{\underline{1},1}^2 .    \label{201808301358xx}
\end{equation}
Making use of \eqref{06011734}, \eqref{201806271120dsfs}, Young inequality and trace estimate, we can derive from  \eqref{201808310859}  and  \eqref{xx201808310859}
that, for sufficiently large $c_1$,
 \begin{align}
 & c c_1 \left(\mathfrak{E}_L(t) + \int_0^t    | \nabla_{\mm{h}}^2 u_3|^2_{1/2} \mm{d}\tau \right) \nonumber
\\
&  \leqslant c_1^{10}\int_0^t\left( \|\eta_3\|_{1,1}^2+| \eta_3|_2^2+ \|\mm{div}\eta^2\|_2 \left(\|\mm{div}\eta^2\|_2 +\|\partial_\mm{h}^2 q\|_{0}\right)+ \sqrt{\mathcal{E} }\mathcal{D}\right)\mm{d}\tau \nonumber \\
&  \quad +  \int_0^t  \mathfrak{E}_L(t) \mm{d}\tau +c_1^5(\|\eta^0\|_3^2+\|u^0\|_2^2+|H^0|_1^2). \label{20180831085xdsfd9}
\end{align}

Noting that $$\mm{div}\eta^2=\mm{div}\eta^0+\int_0^t\mm{div}u\mm{d}\tau,$$
then applying the norm $\|\cdot\|_2$ to the above relation, we easily derive
$$
\begin{aligned}
\|\mm{div}\eta^2\|_2 \leqslant & \|\mm{div}\eta^0\|_2 +\int_0^t\|\mm{div}u\|_2\mm{d}\tau \\
\leqslant   &\Upsilon:= \|\eta^0\|_3^2 +\left(\int_0^t\|\eta\|_3^2 \mm{d}\tau \int_0^t\|u\|_3^2\mm{d}\tau\right)^{1/2}.
\end{aligned}$$
Putting the above estimate into \eqref{20180831085xdsfd9} yields \eqref{201807012235}.
\hfill $\Box$
\end{pf}

\subsection{Growall-type energy inequality}

With the estimates of $(\eta,u)$ in Lemmas \ref{201612132242nn}--\ref{201807251431}, we are in the position to derive  \emph{a prior} Growall-type energy inequality, which couples with the solution $(\eta^2,u^2)$ of the linear problem \eqref{s0106pndfsgnnnxx} for the case $\vartheta\neq 0$.
\begin{lem}  \label{pro:0301n0845}
For any given constant $\Lambda>0$, there exists a energy functional $\overline{\mathcal{E}}(t)$,  and constants $\delta_1\in (0,1)$, $c>0$ and $C_1>0$ such that, for any $\delta\leqslant \delta_1$, if the solution $(\eta,u)$ of the ART problem satisfies \eqref{aprpiosesnew} and $(\eta^2,u^2)$ is a solution of linear problem \eqref{s0106pndfsgnnnxx}, then $(\eta^2,u^2)$ satisfies the  Growall-type energy inequality
\begin{align}
&  \mathfrak{E}(t)+ \beta(c \mathfrak{E}_L(t) +|d(t)|_3^2) +\int_0^t\mathcal{D}(\tau)\mm{d}\tau\nonumber \\
& \leqslant \Lambda\int_0^t
(\mathfrak{E}(\tau)+\beta(c \mathfrak{E}_L(\tau))\mm{d}\tau\nonumber \\
&\qquad +  C_1\left(\|\eta^0\|_3^2+ \|u^0\|_2^2+\beta(|d^0|_{3}^2) +\int_0^t(\beta( \Upsilon^2)+
\|(\eta,u)\|_{0}^2)\mm{d}\tau\right),
\label{2016121521430850}
\end{align}
and the equivalent estimate
\begin{align}
\label{2018008121027}
 \mathcal{E}(t)\leqslant C_1\mathfrak{E}(t)\lesssim \mathcal{E}(t)
\end{align}
for any $t\in I_T$,
where
\begin{align}\label{201811191638}
& d(x_{\mm{h}},t):=\zeta_3((\zeta_{\mm{h}})^{-1}(x_{\mm{h}},t),0,t),\\
& d^0=d|_{t=0}\in H^3(\mathbb{T})\mbox{ and } \mathfrak{E}_L(t)\mbox{ is provided by Lemma \ref{201807251431}}.\nonumber
 \end{align}
It should be noted that the above three constants $\delta_1$, $C_1$ and $c$ depend on $\Lambda$, the domain $\Omega$, and  parameters in the ART problem.
\end{lem}
\begin{rem}
By Proposition \ref{201809012320}, for sufficiently small $\delta_1$, the definition \eqref{201811191638} makes sense.
\end{rem}
\begin{pf}
Making use of \eqref{ssebdaiseqinM0846} for $0\leqslant i\leqslant 1$, \eqref{201702061418} for $0\leqslant i\leqslant 1$, \eqref{Lem:030dsfafds1m0dfsf832} and \eqref{20808011436}, we can derive that, there exist a constant $c$, a (sufficiently small) constant $\tilde{\delta}_1$ and a proper large constant  $c_3$ such that,
for any sufficiently large constant  $c_2>1$ and any $\delta\leqslant \tilde{\delta}_1$,
\begin{align}
 \label{201702061553}
\frac{\mm{d}}{\mm{d}t} \mathfrak{E}_2+c \mathfrak{O}_2\leqslant c_2 \left(| \eta_3|_1^2 + \|\eta\|_{\underline{1},1}^2+|u_3|_0^2 +
 \sqrt{{\mathcal{E}}} {\mathcal{D}} \right) /c,
 \end{align}
  where we have defined that
 $$\begin{aligned}\mathfrak{E}_2:=&\sum_{|\alpha|\leqslant 1}
 \left(\int \rho \partial_{\mm{h}}^{\alpha}\eta \cdot  \partial_{\mm{h}}^{\alpha}u\mm{d}y +
  c_2 \mathcal{I}(  \partial_{\mm{h}}^{\alpha}\eta )\right) +
\frac{1}{4}\|\sqrt{\mu}  \mathbb{D}\eta\|_{\underline{1},0}
+\left\|\sqrt{(\kappa\rho+\bar{M}_3^2)/\mu}\eta\right\|_{2}^2\\
 &+c_2\|\sqrt{\rho} u\|^2_{\underline{1},0}  + c_3\left( \|\sqrt{\rho}  u_{t}\|_{0}^2  \mathcal{I}(u)  -\int \nabla q\cdot (\mathcal{A}_t^{\mm{T}}u)\mm{d}y-\int_\Sigma \llbracket
q\rrbracket \mathcal{A}_te_3\cdot u\mm{d}y_{\mm{h }}\right)
\end{aligned}$$
and
$$
\mathfrak{O}_2:= c_3\|u_t\|_{1}^2+   c_2 \| u \|_{\underline{1},1}^2+\left\|  \left((\kappa\rho+\bar{M}_3^2)\eta,u\right)\right\|_{2}^2+
\|\nabla q \|_{0}^2.
 $$

Making using of \eqref{06142100x},  \eqref{2017020614181721},
the product estimates \eqref{fgestims} and \eqref{0604221702}, Korn's inequality and trace estimate, we can choose  constants $c$ and $\tilde{\delta}_2\leqslant \tilde{\delta}_1$ such that, for any $\delta \leqslant \tilde{\delta}_2$ and any sufficiently large $c_2>1$,
\begin{equation}
\label{201809081915}
 c(\| \eta\|_{\underline{1},1}^2+\|u_t\|_0^2- \|u\|_2\|\eta\|_{ {2},1}^2) \leqslant \mathfrak{E}_2.
\end{equation}

In addition, we can further deduce that
\begin{align}
&|H^0|_1\lesssim |d^0|_3,\label{20181114101015} \\
& |d(t)|_{3}^2\lesssim |d^0|_3^2+\int_0^t \mathcal{D}(\tau)\mm{d}\tau,\label{20181114101015xx}
\end{align}
please refer to the derivation \eqref{201811141623} and \eqref{201811141623xx} in Section \ref{201811142125}. Here we does not directly provide the derivation of the above two estimates, since we shall need auxiliary properties of inverse transformation of Lagrangian coordinates, which will be introduced in Section \ref{201811191634}.

Next we derive  \eqref{2016121521430850} by four cases.

(1) Case of ``$(\kappa,\bar{M}_3)\neq 0$ and $\vartheta=0$".

We define that
$$\begin{aligned}{\mathfrak{E}}_3:=&\mathfrak{E}_2+ \sum_{|\alpha|=2} \left( \int \rho \partial_{\mm{h}}^{\alpha}\eta \cdot  \partial_{\mm{h}}^{\alpha} u\mm{d}y
+ c_2 \mathcal{I}(  \partial_{\mm{h}}^{\alpha}\eta ) \right)
+   \frac{1}{4}\|\sqrt{\mu}  \mathbb{D}\eta\|_{2,0}  +c_2  \|\sqrt{\rho} u\|^2_{2,0} ,\\
\mathfrak{O}_{3}:=& \mathfrak{O}_2 +\sum_{|\alpha|= 2}\mathcal{I}(\partial_{\mm{h}}^{\alpha}\eta) + c_2  \| u \|_{2,1}^2  . \end{aligned}$$

Since $\vartheta=0$, we can further derive from \eqref{ssebdaiseqinM0846} for $i=2$, \eqref{201702061418} for $i=2$, \eqref{201611111952},  and \eqref{201702061553} that, for any  sufficiently large $c_2>1$,
\begin{align}
 \frac{\mm{d}}{\mm{d}t}{\mathfrak{E}}_4 +  c  \mathfrak{O}_4
 \leqslant \left( c_2^2\left(\| \eta \|_{\underline{1},1}^2 +| \eta_3|_{2}^2  +|u_3|_0^2+ \sqrt{{\mathcal{E}}} {\mathcal{D}} \right) +\|\eta\|_{2,1}^2 \right)/c, \label{201702061553n}
 \end{align} where we have defined that
$$\begin{aligned}{\mathfrak{E}}_4:=&c_2{\mathfrak{E}}_3 + \overline{\|\eta\|}_{3}^2
\end{aligned}$$
and
$$\mathfrak{O}_4:= c_2 \mathfrak{O}_3+  \left\|  \left( \eta,u\right)\right\|_{3}^2+
\|\nabla q \|_1^2+| \llbracket \nabla q \rrbracket |_{1/2}^2  .$$
In addition, making use of \eqref{2017020614181721}, \eqref{201702071610nb}, \eqref{201809081915}  and  Young's and Korn's inequalities,
we can choose constant $\tilde{\delta}_3\leqslant \tilde{\delta}_2$ such that, for any $\delta \leqslant \tilde{\delta}_3$ and any sufficiently large $c_2$,
\begin{align}
\label{201808021219}
& {\mathfrak{E}}_4,\ \mathcal{E}\mbox{ and } \|\eta\|_3^2+\|u\|_2^2 \mbox{ are equivalent},\\
\label{201808safdas02asdfsadf1219}
&c\left(c_2  \|\eta\|_{\underline{2},1}^2+\|\eta\|_3^2+\|u\|_2^2\right) \leqslant  \mathfrak{E}_4,\\
\label{201808061629}&
\mathfrak{O}_4\mbox{ is equivalent to }\mathcal{D}.
\end{align}

By  \eqref{201808021219}--\eqref{201808061629}, we can derive from \eqref{201702061553n} that there exist  constants $\tilde{\delta}_4\leqslant \tilde{\delta}_3 $ and  $c_2$ such that,  for any $\delta \leqslant \tilde{\delta}_4$,
\begin{align}
 \frac{\mm{d}}{\mm{d}t}{\mathfrak{E}}_4 +  c\mathcal{D}_4 \leqslant  \left( \|\eta\|_2^2+| \eta_3|_2^2 +|u_3|_0^2 \right)/c+\frac{\Lambda}{2} \mathfrak{E}_4(\tau ) .
 \label{2017020615dfgss53n}
 \end{align}
Using \eqref{201805072235}, \eqref{201808safdas02asdfsadf1219}, interpolation inequality \eqref{201807291850} and Young's inequality, we further get from \eqref{2017020615dfgss53n} that
\begin{align}
 \frac{\mm{d}}{\mm{d}t}\mathfrak{E}_4 +  c \mathcal{D}_4 \leqslant\Lambda
\mathfrak{E}_4(\tau ) +
 \|(\eta,u_3)\|_{0}^2/c .\nonumber
 \end{align}
Integrating the above inequality with respect to $t$, and
using \eqref{201702071610nb} with $t=0$, we get \eqref{2016121521430850} by taking $\delta_1:=\tilde{\delta}_4$ and $\mathfrak{E}=\mathfrak{E}_4/c $ for some constant $c$.

(2) Case of $(\kappa,\bar{M}_3,\vartheta) =0$.

Let $0\leqslant j\leqslant 1$ and $0\leqslant k\leqslant 1-j$.
Applying $\partial_{\mm{h}}^{k}\partial_3^{j}$ to \eqref{s0106pnnnn}$_2$, and then multiplying the resulting identity by $ \partial_{\mm{h}}^{k}\partial_3^{2+j}\eta$ in $L^2$,
we have
$$
 \frac{\mm{d}}{\mm{d}t}\|\sqrt{\mu } \partial_{\mm{h}}^k\partial_3^{2+j} \eta \|_{0}^2 =-  2 \int \partial_{\mm{h}}^{k}\partial_3^j ( \partial_{\bar{M}}^2\eta +\mu\Delta_{\mm{h}}u +\mu \nabla\mm{div} u -{\rho} u_t  -\nabla q+ \mathcal{N}^1)\cdot \partial_{\mm{h}}^{k}\partial_3^{2+j}\eta \mm{d}y.
$$
Making use of \eqref{20161111195sdfafsd2}--\eqref{201808292010}, we further have
\begin{align}
 \frac{\mm{d}}{\mm{d}t}\|\sqrt{\mu } \partial_3^{2+j} \eta \|_{\underline{1-j},0}^2  \lesssim & \| \sqrt{\mu}\partial_3^{2+j} \eta \|_{\underline{1-j},0} \| \partial_3^{j}(\partial_{\bar{M}}^2\eta,\Delta_{\mm{h}}u ,u_{t},\nabla q, \nabla \mm{div} u, \mathcal{N}^1)\|_{\underline{1-j},0}\nonumber \\
\lesssim & \| \sqrt{\mu}\partial_3^{2+j} \eta \|_{\underline{1-j},0}  \left( \|\partial_3^ju_t\|_{\underline{1-j},0}
+\sqrt{\mathcal{E}\mathcal{D}}\right)\nonumber \\
& + \| \sqrt{\mu}\partial_3^{2+j} \eta \|_{\underline{1-j},0}
\left\{  \begin{array}{ll}
\|(\eta,u)\|_{\underline{2},1}+\|u_t\|_{1} &\mbox{ for }j=0;\\
  \|\eta\|_{2,1}+\|u\|_{\underline{2},1}+\|u_t\|_1  &\hbox{ for }j=1.
  \end{array}
 \right. \nonumber
\end{align}
Exploiting Young inequality  and   \eqref{201702071610nb}, we derive from the above estimate that
\begin{align}
& \frac{\mm{d}}{\mm{d}t}  {\|\sqrt{\mu}\partial_3^2\eta\|}_{\underline{1},0}^2
 \leqslant  \frac{\Lambda}{2}
  \| \sqrt{\mu}\partial_3^{2} \eta \|_{\underline{1},0}^2  + c ( \|(\eta,u)\|_{\underline{2},1}^2 + \|u_t\|_1^2+\sqrt{\mathcal{E}} \mathcal{D}),\label{201808131706} \\
 &
\frac{\mm{d}}{\mm{d}t} {\|\sqrt{\mu}\partial_3^3\eta\|}_{0}^2
\leqslant  \frac{\Lambda}{2}
 \|\sqrt{\mu} \partial_3^3 \eta \|_{0}^2+c
\big( \|\sqrt{\mu}\eta\|_{ {2},1}^2+ \|u\|_{\underline{2},1}^2+\|u_t\|_1^2  +\sqrt{\mathcal{E}} \mathcal{D}\big).\label{201808131708}
 \end{align}

Similarly to \eqref{201702061553n}, we can derive from  \eqref{ssebdaiseqinM0846} for $i=3$, \eqref{201702061418} for $i=3$,  \eqref{201702061553} and the above two estimates that, for any sufficiently small $\delta$ and any sufficiently large $c_2$,
\begin{align}
\frac{\mm{d}}{\mm{d}t}\mathfrak{E}_5+ c\mathfrak{O}_5\leqslant &
 \frac{\Lambda c_2}{2}
  \| \sqrt{\mu}\partial_3^{2} \eta \|_{\underline{1},0}^2
 +\frac{\Lambda }{2}
 \|\sqrt{\mu} \partial_3^3 \eta \|_{0}^2  \nonumber \\
& + \left( c_2\|\eta\|_{\underline{2},1}^2+ c_2^3 \left(|\eta_3|_2^2+\|\eta\|_{\underline{1},1}^2+|u_3|_0^2+\sqrt{{\mathcal{E}}} {\mathcal{D}}\right)\right)/c , \label{2018081315583}
\end{align}
where we have defined that $\mathfrak{O}_5:=c_2^2\mathfrak{O}_3$ and
\begin{align}
\nonumber
& \mathfrak{E}_5:=  c_2^2\mathfrak{E}_3 + {\|\sqrt{\mu}\partial_3^3\eta\|}_{0}^2+  c_2 {\|\sqrt{\mu}\partial_3^2\eta\|}_{\underline{1},0}^2.
 \end{align}

Similarly to \eqref{201808021219}--\eqref{201808061629},  exploiting \eqref{2017020614181721},  \eqref{20161111195sdfafsd2}, \eqref{201702071610nb}, \eqref{201809081915} and Young's and Korn's inequalities,  we have, for any sufficiently small $\delta$ and any sufficiently large $c_2$,
\begin{align}
&{\mathfrak{E}}_5 ,\ \mathcal{E}\mbox{ and } \|\eta\|_3^2+\|u\|_2^2 \mbox{ are equivalent},
\label{201808131542}
\\
&{\|\sqrt{\mu}\partial_3^3\eta\|}_{0}^2 + c_2 {\|\sqrt{\mu}\partial_3^2\eta\|}_{\underline{1},0}^2 +c \left(\|\eta\|_3^2+\|u\|_2^2 +c_2^2  \|\eta\|_{\underline{2},1}^2\right)  \leqslant {\mathfrak{E}}_5, \label{201809091222} \\
\label{2018080616dsfa29}
&   \mathfrak{O}_5 +\|\eta\|_3^2\mbox{ is equivalent to }\mathcal{D}.
\end{align}
Then, making use of \eqref{201805072235}, \eqref{201808131542}--\eqref{2018080616dsfa29}, interpolation inequality \eqref{201807291850} and Young's inequality,
we derive from \eqref{2018081315583} that, for any sufficiently small $\delta$ and some sufficiently large $c_2$,
\begin{align}
\frac{\mm{d}}{\mm{d}t}\mathfrak{E}_5+ c\mathfrak{O}_5\leqslant
 \Lambda \mathfrak{E}_5+\|(\eta_3,u_3)\|_0^2 /c    \nonumber
\end{align}
for any sufficiently mall $\delta$.
Consequently we can easily derive \eqref{2016121521430850} from the above inequality.

(3) Case of ``$(\kappa,\bar{M}_3)\neq0$ and $\vartheta\neq 0$".

We define that
$$\begin{aligned}
\mathfrak{E}_6:=&c_2 \mathfrak{E}_2 + {c_2^3}\int_{\mathcal{T}_{4,4}}\sum_{|\alpha|= 2}\left(
\int   \rho \mathfrak{D}_{\mf{h}}^{3/2,\alpha} \eta \cdot  \mathfrak{D}_{\mf{h}}^{3/2,\alpha} u \mm{d}y
+c_2\mathcal{I}( \mathfrak{O}_{\mf{h}}^{3/2,\alpha}\eta  )
\right)\mm{d}\mathbf{h}
\\
&
+ c_2^2  \|(\sqrt{\rho} u,\sqrt{\kappa\rho} \mathbb{D} \eta_3 ,\partial_{\bar{M}} \eta_3 )\|^2_{2,0} +c_2\sum_{|\alpha|=2}\int   \partial_\mm{h}^\alpha \eta \cdot  \partial_\mm{h}^\alpha u \mm{d}y\\
& + c_2^3 \int_{\mathcal{T}_{4,4}}\left(
 \frac{1}{4}  \|\sqrt{ \mu}\mathbb{D} \mathfrak{D}_{\mf{h}}^{3/2}\eta\|^2_{1,0} +
c_2 \|\sqrt{\rho}  \mathfrak{D}_{\mf{h}}^{3/2}   u\|^2_{1,0} \right)\mm{d}\mathbf{h}+ \frac{ c_2}{4}    \| \sqrt{\mu}\mathbb{D}  \eta \|^2_{2,0}
\end{aligned}$$
and
 $$\begin{aligned}
\mathfrak{O}_6:=&
c_2^2\| u\|_{2,1}^2+c_2 \mathfrak{O}_2+
{c_2^3}\vartheta|\nabla_{\mm{h}}\partial_{\mm{h}}  \eta_3|_{1/2}^2 .
\end{aligned}$$
Then, using Young's inequality and trace estimate, we can derive from \eqref{ssebdaisedsaqinM0dsfadsff846asdfadfad}, \eqref{ssebdaiseqinM0asdfa846}, \eqref{ssebdaisedsaqinM0dsfadsff846sdfaaasdfadfad}, \eqref{2018072033}, \eqref{201611111952} and \eqref{201702061553} that, for any sufficiently small $\delta$ and any sufficiently large $c_2$,
\begin{align}
\frac{\mm{d}}{\mm{d}t}\mathfrak{E}_7+c\mathfrak{O}_7\leqslant &
\left( \| \eta\|_{2,1}^2 +c_2^4\left(| \eta_3|_2^2 +\|(\eta,u)\|_2^2+  |\nabla_{\mm{h}}^2u_3|_{1/2}^2
+ \sqrt{{\mathcal{E}}} {\mathcal{D}} \right) \right)/c ,
 \label{20170206155dsfadsfas3}
 \end{align}
where we have defined that
$$\begin{aligned}
&\mathfrak{E}_7:= \mathfrak{E}_6+ \overline{\|\eta\|}_3^2 ,\\
&\mathfrak{O}_7:=\mathfrak{O}_6+
 \left\|  \left( \eta,u\right)\right\|_3^2+
\|\nabla q \|_1^2+| \llbracket\nabla q \rrbracket |_{1/2}^2
.
\end{aligned}$$

In addition, similarly to \eqref{201808021219}--\eqref{201808061629},
we also derive that,  for any sufficiently small $\delta$ and any sufficiently large $c_2$,
\begin{align}
& {\mathfrak{E}}_7,\ \mathcal{E}\mbox{ and } \|\eta\|_3^2+\|u\|_2^2 \mbox{ are equivalent},  \label{201809091655} \\
&c (\|\eta\|_3^2+\|u\|_2^2+c_2  \|\eta\|_{2,1}^2) \leqslant  \mathfrak{E}_7,  \\
&
\mathfrak{O}_7\mbox{ is equivalent to }\mathcal{D}. \label{201809091655xx}
\end{align}
Integrating \eqref{20170206155dsfadsfas3} over $(0,t)$, and then making use of  \eqref{201805072235},
interpolation inequality, Young's inequality
and the three results above, we have,  for any sufficiently   sufficiently small $\delta$ and any sufficiently  large constant $c_2$,
\begin{align}
 \mathfrak{E}_7+c\int_0^t\mathfrak{O}_7\mm{d}\tau\leqslant & (\|\eta^0\|_3^2+\|u^0\|_2^2)/c+\frac{\Lambda}{2}\int_0^t \mathfrak{E}_7\mm{d}\tau\nonumber
\\
&+\int_0^t\left( \| (\eta ,u )\|_0^2+  |\nabla_{\mm{h}}^2u_3|_{1/2}^2
  \right) \mm{d}\tau /c.  \label{201809091850}
 \end{align}
Consequently, using \eqref{201809091655}--\eqref{201809091655xx} and interpolation inequality again, we can easily derive \eqref{2016121521430850} from \eqref{201807012235}, \eqref{20181114101015}, \eqref{20181114101015xx} and \eqref{201809091850}.

(4) Case of ``$(\kappa,\bar{M}_3)=0$ and $\vartheta\neq0$".

We can derive from \eqref{ssebdaisedsaqinM0dsfadsff846asdfadfad}, \eqref{ssebdaiseqinM0asdfa846}, \eqref{ssebdaisedsaqinM0dsfadsff846sdfaaasdfadfad},   \eqref{2018072033}, \eqref{201702061553}, \eqref{201808131706} and \eqref{201808131708} that, for some sufficiently large $c_2$ and for any sufficiently small $\delta$,
\begin{align}
\frac{\mm{d}}{\mm{d}t}\mathfrak{E}_8+c\mathfrak{O}_8\leqslant &  \frac{\Lambda c_2}{2}
  \| \sqrt{\mu}\partial_3^{2} \eta \|_{\underline{1},0}^2
 +\frac{\Lambda }{2}
 \|\sqrt{\mu} \partial_3^3 \eta \|_{0}^2+ \left(| \llbracket\nabla_{\mm{h}} q \rrbracket |_{1/2}^2+ c_2\|\eta\|_{\underline{2},1}^2\right)/c    \nonumber \\
&  + c_2^6 \left(|\eta_3|_2^2+\|(\eta,u)\|_2^2+|u_3|_0^2+  |\nabla_{\mm{h}}^2u_3|_{1/2}^2+\sqrt{{\mathcal{E}}} {\mathcal{D}} \right)/c ,\nonumber
 \end{align}
  where we have defined that $\mathfrak{O}_8:=c_2\mathfrak{O}_6$ and
 $$\begin{aligned}\mathfrak{E}_8:=&c_2\mathfrak{E}_6+ {\|\sqrt{\mu}\partial_3^3\eta\|}_{0}^2+  c_2 {\|\sqrt{\mu}\partial_3^2\eta\|}_{\underline{1},0}^2 .
\end{aligned}$$
By \eqref{20161111195sdfafsd2}, we further have, for any sufficiently   sufficiently small $\delta$ and any sufficiently  large constant $c_2$,
\begin{align}
 \label{2017020615safd5dsfadsfas3}
\frac{\mm{d}}{\mm{d}t}\mathfrak{E}_8+c\mathfrak{O}_8\leqslant &  \frac{\Lambda c_2}{2}
  \| \sqrt{\mu}\partial_3^{2} \eta \|_{\underline{1},0}^2
 +\frac{\Lambda }{2}
 \|\sqrt{\mu} \partial_3^3 \eta \|_{0}^2+ c_2\|\eta\|_{\underline{2},1}^2/c \nonumber \\
& + \left(c_2^6 \left(|\eta_3|_2^2+\|(\eta,u)\|_2^2+|u_3|_0^2+  |\nabla_{\mm{h}}^2u_3|_{1/2}^2+\sqrt{{\mathcal{E}}} {\mathcal{D}}\right)\right)/c.
 \end{align}

 Similarly to \eqref{201808131542}--\eqref{2018080616dsfa29},  we also derive that,
 for any sufficiently large $c_2$, $\mathfrak{E}_8$ and $\mathfrak{O}_8$ also
enjoys the properties \eqref{201808131542}--\eqref{2018080616dsfa29} with  $(\mathfrak{E}_8,\mathfrak{O}_8)$ in place of $(\mathfrak{E}_5,\mathfrak{O}_5)$
for any sufficiently   sufficiently small $\delta$.
Consequently, similarly to \eqref{201809091850}, we can also derive from \eqref{2017020615safd5dsfadsfas3} that, for any sufficiently   sufficiently small $\delta$ and any sufficiently  large constant $c_2$,
\begin{align}
 \mathfrak{E}_8+c\int_0^t\mathfrak{O}_8\mm{d}\tau\leqslant & (\|\eta^0\|_3^2+\|u^0\|_2^2)/c+\frac{3\Lambda}{4}\int_0^t \mathfrak{E}_8\mm{d}\tau\nonumber
\\
&+\int_0^t\left( \| (\eta,u)\|_0^2+  |\nabla_{\mm{h}}^2 u_3|_{1/2}^2
  \right) \mm{d}\tau /c. \nonumber
\end{align}
which, together with \eqref{201807012235}, \eqref{20181114101015} and \eqref{20181114101015xx}, further implies   \eqref{2016121521430850}. This completes the proof of Lemma  \ref{pro:0301n0845}.
\hfill$\Box$ \end{pf}

The local existence of strong solutions to the 3D PDEs model of
stratified incompressible viscous fluids based on the Navier--Stokes equations have been
established, see \cite{WYJTIKCT,PJSGOI5x,WMRTITNS2017} for examples. Similarly to the existence results of solution of stratified viscous fluids in \cite{WYJTIKCT} or the ones of viscous fluids with upper free boundary in \cite{XLZPZZFGAR,GYTILW1,GYTIDAP}, by using a Faedo--Galerkin approximation scheme for the linearized problem and an iterative method, we can also get a unique local-in-time existence result of a unique strong   solution $(\eta,u)$ with an associated pressure $q$ to the ART problem for some $T$,  unique global-in-time strong  solution $(\eta^2,u^2)$ with an associated pressure $q^2$ to the linear problem \eqref{s0106pndfsgnnnxx} for any given $\mathcal{N}_1^2$, $\mathcal{N}_2^2$ and $\mathfrak{M}_2$ defined by $(\eta,u)$; moreover the strong  solution also satisfies the \emph{a priori }estimates in Lemma \ref{pro:0301n0845}.  Since the proof is standard in the local-wellposedness theory of PDEs, and hence we omit  the details, and only directly state the local well-posedness result, in which the solution enjoys the Growall-type energy inequality.
\begin{pro} \label{pro:0401nxd}
\begin{enumerate}
  \item[(1)] Let  $(\eta^0,u^0)\in H^3_0\times H_0^2$,  $\zeta^0:=\eta^0+y$,
  $\zeta^0_{\mm{h}}(y_{\mm{h}},0):\mathbb{R}^2 \to \mathbb{R}^2$ is a homeomorphism mapping, and $d^0:=\zeta_3((\zeta_{\mm{h}})^{-1}(x_{\mm{h}}),0) $. There exists a sufficiently small $\delta_2\in (0,1)$, such that, if $(\eta^0,u^0,d^0)$ satisfying
\begin{align}
& \det \zeta^0=1,  d^0\in (-l,\tau)\\
& \sqrt{\|\eta^0\|_3^2+\|{u}^0\|_2^2+\beta( |d^0|_3^2)}<  \delta_2
\end{align}
and necessary compatibility conditions
\begin{alignat}{2}
&\mm{div}_{\mathcal{A}^0}u^0=0& &\ \mbox{ in }\Omega , \label{201808040977}    \\
& \Pi_{\mathcal{A}^0 e_3} \llbracket   \mathcal{S}_{\mathcal{A}^0}(0, u^0,   \eta^0 )\mathcal{A}^0 e_3-   \bar{M}_3 \partial_{\bar{M}}\eta^0 \rrbracket = 0  & &\
 \mbox{ on }\Sigma,\label{201811031045}
\end{alignat}
then there is a local existence time $T^{\max}>0$, depending on $\delta_2$, the domain and the known parameters in the ART problem, and a unique local-in-time strong  solution
$(\eta, u)\in C^0([0,T^{\max})$, $H^3\times H^2 )$ with a unique (up to a constant) associated pressure $q$ to the ART problem, where $(\eta,u,q)$ enjoys the following regularity
$$
\begin{aligned}
&(u_t,\nabla q, \llbracket  q \rrbracket  )\in C^0([0,T^{\max}),L^2\times L^2\times H^{1/2}(\mathbb{T})),\\
&(u,u_t, \nabla q )\in L^2((0,T^{\max}),H^3\times H^1\times H^1),\\
& q(x,t)\in H^1\mbox{ for }t\in [0,T^{\max}).
\end{aligned}$$
Moreover, $d:=\zeta_3((\zeta_{\mm{h}})^{-1}(x_{\mm{h}},t),0,t)$ makes sense, and belongs to $C^0(\overline{I_T},H^3(\mathbb{T}))$ for $\vartheta\neq 0$.
  \item[(2)]  Using the above strong  solution $(\eta,u)$ to define $\mathcal{N}_1^2$, $\mathcal{N}_2^2$ and $\mathfrak{M}_2$,
then the linear problem \eqref{s0106pndfsgnnnxx} possesses unique local-in-time strong   solutions
$(\eta^2, u^2) $ with an associated pressure functions $q^2$. Moreover,  $(\eta^2, u^2,q^2)$ enjoys the regularity as $(\eta, u,q)$.
  \item[(3)] In addition, if $(\eta,u)$ satisfies
$$\sup_{t\in [0,T)}\sqrt{\|\eta(t)\|_3^2 +\|u(t)\|_2^2}\leqslant \delta_1\mbox{ for some }T<T^{\max},$$
 then the solution $(\eta, u)$ enjoys the Growall-type energy inequality \eqref{2016121521430850} and the equivalent estimate \eqref{2018008121027} on $(0,T)$.
\end{enumerate}
\end{pro}
\begin{rem}
In the second assertion in Proposition \ref{pro:0401nxd}, we do not mention the necessary compatibility conditions as in the first assertion, since the initial data in linear problem \eqref{s0106pndfsgnnnxx}  automatically satisfies necessary compatibility conditions, similarly to \eqref{201808040977} and \eqref{201811031045}.
\end{rem}

\begin{rem}
The initial date of the associated pressure $q$ constructed in Proposition \ref{pro:0401nxd} satisfies
\begin{equation}
\llbracket  \mathcal{S}_{\mathcal{A}^0} (q^0,u^0,\eta^0)\mathcal{A}^0e_3- g
\rho \eta_3^0  \mathcal{A}e_3-\bar{M}_3 \partial_{\bar{M}}\eta^0 \rrbracket  =\vartheta H^0  \mathcal{A}^0 e_3  \mbox{ on }\Sigma , \label{201808040977xx}
\end{equation}
and is the weak solution of the mixed boundary value problem (referring to the derivation (4.12) in \cite{WYJTIKCT} for the definition of weak solution)
\begin{equation}
\left\{
\begin{array}{ll}
\Delta_{\mathcal{A}^0} q^0 /\rho=\mm{div}_{\mathcal{A}^0} \left((\partial_{\bar{M}}^2\eta^0 -\mm{div}_{\mathcal{A}^0}\mathcal{S}_{\mathcal{A}^0}(0,u^0,\eta^0))/\rho  - u^0 \cdot \nabla_{\mathcal{A}^0} u^0  \right)  &\mbox{ in }\Omega,\\
 \llbracket \nabla_{\mathcal{A}^0}q^0 /\rho \rrbracket\cdot  \mathcal{A}^0 e_3\\
= \llbracket (\partial_{\bar{M}}^2\eta^0-\mm{div}_{\mathcal{A}^0}\mathcal{S}_{\mathcal{A}^0}(0,u^0,\eta^0 ))/\rho- u^0 \cdot \nabla_{\mathcal{A}^0} u^0  ) \rrbracket \cdot  \mathcal{A}^0 e_3 &\mbox{ on }\Sigma,
\\
 \nabla_{\mathcal{A}^0}q^0  \cdot  \mathcal{A}^0 e_3 =\left(  \partial_{\bar{M}}^2\eta^0-\mm{div}_{\mathcal{A}^0}\mathcal{S}_{\mathcal{A}^0}(0,u^0,\eta^0   )  \right) \cdot  \mathcal{A}^0 e_3 &\mbox{ on }\Sigma_-^+ ,\\
\eqref{201808040977xx}.
\end{array}
\right. \label{201808040916}
\end{equation}
It should be noted that test functions for the weak form of \eqref{201808040916}$_1$ belong  to $H^1$ (referring to (3.64) in \cite{WYJTIKCT}).
Since there exists a $\delta_3$ such that, for any $\|\eta\|_3\leqslant \delta_3$,
$$\|\nabla q\|_0 \leqslant \|\nabla_{\mathcal{A}} q\|_0\mbox{ and }\mathcal{A}e_3\neq 0,$$
thus we easily see that, for any
$\|\eta^0\|_3\leqslant \delta_3$, the weak solution solution $q^0\in H^1$ of \eqref{201808040916} is unique up to a constant for given $(\eta^0,u^0)$.
\end{rem}

\begin{rem}
\label{201811031411}
It should be noted that, for any strong solution $(\eta,u)$ with an associated pressure $q$ constructed by Proposition \ref{pro:0401nxd}, for any $t_0\in (0,T^{\max})$,
$(\eta,u,q)|_{t=t_0}$ automatically satisfies \eqref{201808040977} and the weak form of \eqref{201808040916}  with $(\eta,u,q)|_{t=t_0}$ in place of $(\eta^0,u^0,q^0)$. We take  $(\eta,u,q)|_{t=t_0}$ as a new initial data, then the new initial data can define a unique local-in-time strong solution $(\tilde{\eta},\tilde{u},\tilde{q})$ constructed by Proposition
\ref{pro:0401nxd}; moreover the initial date of $\tilde{q}$ is equal to $q|_{t=t_0}$ up to  a constant.
\end{rem}


\section{Construction of initial data for  nonlinear solutions}

For any given $\delta>0$,  let
\begin{equation}\label{0501}
\left(\eta^\mm{a}, u^\mm{a} ,q^{\mm{a}}\right)=\delta e^{{\Lambda t}} (\tilde{\eta}^0,\tilde{u}^0,\tilde{q}^0),
\end{equation}
where $(\tilde{\eta}^0,\tilde{u}^0, \tilde{q}^0):=(w/\Lambda,w,\beta)$, and $(w,\beta)\in (\mathcal{A}\cap H^\infty) \times  H^\infty$ comes from  Proposition \ref{thm:0201201622}.
Then $\left(\eta^\mm{a}, u^\mm{a}\right)$  is a solution to the linearized ART problem, and enjoys the estimates, for any $i$, $j\geqslant 0$,
 \begin{equation}
\label{appesimtsofu1857}
\|\partial_t^i(u^{\mm{a}}, q^{\mm{a}})\|_{\mm{S},j}\leqslant c(i,j) \delta e^{\Lambda t}.
 \end{equation}Moreover, by \eqref{201602081445MH} and \eqref{201809121901}, for $\tilde{\chi}^0:=\tilde{\eta}^0$ or $\tilde{u}^0$,
 \begin{align}
&
 |\tilde{\chi}_3^0|_{L^1}  \|\tilde{\chi}_{\mm{h}}^0\|_{L^1} \|\tilde{\chi}_3^0\|_{L^1}
\|\partial_3\tilde{\chi}^0_{\mm{h}}\|_{L^1}
\|\partial_3\tilde{\chi}_3^0\|_{L^1}
 \|\mm{div}_{\mm{h}}\tilde{\chi}_{\mm{h}}^0\|_{L^1} > 0, \label{n05022052} \\
&\|\partial_{\bar{M}}\tilde{\chi}^0_{\mm{h}}\|_{L^1},\  \|\partial_{\bar{M}}\tilde{\chi}^0_3\|_{L^1}>0 \mbox{ if }
\bar{M}_3\neq0. \label{n05022052xx}
\end{align}
Unfortunately, the initial data of linear  solution $\left(\eta^\mm{a}, u^\mm{a}, q^\mm{a}\right)$ does not satisfy the necessary compatibility conditions   of ART problem in general. Therefore, next we modify  the initial data of the linear solutions, so that the obtained new initial data satisfy \eqref{201808040977} and  \eqref{201811031045}.
\begin{pro}
\label{lem:modfied}
Let $(\tilde{\eta}^0,\tilde{u}^0,\tilde{q}^0)$
be the same as in \eqref{0501}, then there exists a constant $\delta_4$ (depending on $(\tilde{\eta}^0,\tilde{u}^0,\tilde{q}^0)$, the domain and the parameters in ART problem), such that for any $\delta\in(0,\delta_4)$, there is a couple $(\eta^{\mm{r}},u^{\mm{r}},q^{\mm{r}})\in H^4_0\cap H^2_0\cap H^1$ enjoying the following properties:
\begin{enumerate}[\quad (1)]
                \item The modified initial data
  \begin{equation}\label{mmmode04091215}(  {\eta}_0^\delta,{u}_0^\delta,q_0^\delta): =\delta
   (\tilde{\eta}^0,\tilde{u}^0, \tilde{q}^0) + \delta^2 (\eta^\mm{r},   u^\mm{r},q^{\mm{r}})
\end{equation}
belongs to $H^4_{0}\times H^2_0\times {H}^1$, and satisfies
\begin{align}
\det \nabla ( {\eta}_0^\delta+y)=1,\nonumber
\end{align}
the compatibility conditions  \eqref{201808040977} and  \eqref{201811031045} with $(\eta_0^\delta,u_0^\delta,q_0^\delta)$ in place of $(\eta^0,u^0$, $q^0)$.
\item Uniform estimate:
\begin{equation}
\label{201702091755}
 \sqrt{\|\eta^{\mm{r}}\|_4^2+ \|(u^{\mm{r}},q^{\mm{r}})\|_{\mm{S},0}^2 }\leqslant {C}_2,
 \end{equation}
where the constant $C_2\geqslant 1$ depends on the domain and the known parameters, but is independent of $\delta$.
 \end{enumerate}
\end{pro}
\begin{pf}
For the simplicity in the proof of Proposition \ref{lem:modfied}, we rewrite $( {\eta}_0^\delta,{u}_0^\delta,q_0^\delta,  \tilde{\eta}^0,\tilde{u}^0, \tilde{q}^0)$ by $( {\eta}^\delta,{u}^\delta,q^\delta, \tilde{\eta},\tilde{u}, \tilde{q})$.

Recalling the construction of $(\tilde{\eta} ,\tilde{u} ,\tilde{q} )$, we can see that
$(\tilde{\eta} ,\tilde{u} ,\tilde{q} )$ satisfies
\begin{equation}
\label{201705141510}
\left\{
\begin{array}{ll}
 \mm{div}\tilde{\eta} =\mm{div}\tilde{u}  =0   &\mbox{ in }\Omega,\\
 \Lambda\rho \tilde{u} +\mm{div}\mathcal{S}(\tilde{q} , \tilde{u} ,\tilde{\eta} )=\partial_{\bar{M}}^2\tilde{\eta}   & \mbox{ in }\Omega,\\
  \llbracket   \mathcal{S}(\tilde{q} , \tilde{u} ,\tilde{\eta} )e_3-
 g\rho  \eta_3   e_3- \bar{M}_3 \partial_{\bar{M}}\tilde{\eta} \rrbracket= \vartheta \Delta_{\mm{h}}\tilde{\eta}_3   e_3 & \mbox{ on }\Sigma,\\ \llbracket \tilde{\eta}   \rrbracket= \llbracket \tilde{u}   \rrbracket= 0 & \mbox{ on }\Sigma,\\
 (\tilde{\eta} ,\tilde{u} ) = {0} & \mbox{ on }\Sigma_{-}^+.
 \end{array}\right.
\end{equation}
If $(\eta^\mm{r},u^\mm{r},q^\mm{r})\in H^4\times H^2\times H^1$ satisfies
\begin{equation}
\label{201702061320}
\left\{\begin{array}{ll}
\mm{div}\eta^{\mm{r}}=O(\eta^{\mm{r}}),\ \mm{div}u^\mm{r} = - \mm{div}_{\tilde{\mathcal{A}}^\delta} ( \tilde{u}   + \delta  u^\mm{r} ) /\delta &\mbox{ in }\Omega,
\\    \llbracket
\mathcal{S}(q^{\mm{r}}, u^{\mm{r}},\eta^{\mm{r}})  e_3 - g \rho  \eta_3^{\mm{r}}e_3 - \bar{M}_3 \partial_{\bar{M}}\eta^{\mm{r}}  \rrbracket
=  \vartheta \Delta_{\mm{h}}\eta_3^{\mm{r}} e_3+ \mathcal{N}^2({\eta} ^\delta,{u} ^\delta,  {q} ^\delta )/\delta^{2}  & \mbox{ on }\Sigma,
\\
\llbracket {\eta}^{\mm{r}}  \rrbracket= \llbracket {u}^{\mm{r}}  \rrbracket=  0 & \mbox{ on }\Sigma,\\
(\eta^\mm{r}, u^\mm{r}  )=0&\mbox{ on }\Sigma_{-}^+,
\end{array}\right.\end{equation}
 where $({\eta} ^\delta,{u} ^\delta,q ^\delta)$ is given in the mode \eqref{mmmode04091215},
  $\zeta^\delta:=\eta^\delta +y$, $\mathcal{A} ^\delta:= (\nabla\zeta^{\delta})^{-\mm{T}}$, $\tilde{\mathcal{A}}^{\delta} := \mathcal{A} ^\delta-I$ and $O(\eta^{\mm{r}}):=\delta^{-2}(1+ \delta^2\mm{div}\eta^{\mm{r}}-\det\nabla  \zeta^\delta)$,
then, by \eqref{201705141510}, it is easy to check that $({\eta} ^\delta,{u} ^\delta,q ^\delta)$  belongs to $H^4_{0}\times H^2_0 \times H^1$, and satisfies
\begin{align}
\det \nabla \zeta^\delta=1 \nonumber
\end{align}
\eqref{201808040977}  and \eqref{201811031045}  with $(\eta ^\delta,u ^\delta,q ^\delta)$ in place of $(\eta^0 ,u^0 ,q^0 )$.
Next we construct such $(\eta^\mm{r},u^\mm{r},q^\mm{r} )$ enjoying \eqref{201702061320} by three steps.

(1) We begin with the construction of $\eta^\mm{r}$.
 We consider a stratified   Stokes problem:  for a given function $\xi\in H^4(\Omega)$,
 \begin{equation}    \label{201705141945}
 \left\{\begin{array}{ll}
 \nabla \varpi-\mu\Delta \eta=0,\quad  \mm{div}\eta =O(\xi) &\mbox{ in }\Omega,  \\
\llbracket \eta \rrbracket=0,\ \llbracket (\varpi I-\mathbb{D}\eta)e_3 \rrbracket=0 & \mbox{ on } \Sigma,\\
\eta =0&\mbox{ on } \Sigma_-^+,
\end{array}\right.\end{equation}
where $O(\xi)= (1+\delta^2\mm{div}\xi-\det \nabla (y+\delta\tilde{\eta} +\delta^2\xi ))/\delta^2$. By virtue of the product estimate \eqref{fgestims}, it is easy to estimate that
$$\|O(\xi)\|_3\leqslant c(1+\delta \|\xi\|_4+\delta^2\|\xi\|_4^2+\delta^3\|\xi\|_4^3)\leqslant c_1(1 +\delta^3\|\xi\|_4^3). $$
By the existence  theory on the stratified   Stokes problem, there exists a solution $(\eta,\varpi)\in H^4\times  {H}^3$ of \eqref{201705141945} with $\|\varpi\|_0^2=1$.
Moreover, by  stratified Stokes estimate, it holds that
\begin{align}
\|(\eta,\varpi) \|_{\mm{S},2} \leqslant\|O(\xi)\|_3
\leqslant c_1 \left(1+\delta^3\|\xi\|_4^3\right),
\label{201702141356nn}  \end{align}
where the letter $c_1$ denotes a fixed constant, and is dependent of    $\delta$.

Therefore, one can construct an approximate function sequence $\{(\eta^{n},\varpi^{n})\}_{n=1}^\infty$, such that for any $n\geqslant 2$,
\begin{equation}  \label{201702051905}
 \left\{\begin{array}{ll}
\nabla \varpi^n-\mu \Delta \eta^{n}=0,\quad \mm{div}\eta^{n}
=O(\eta^{n-1})&\mbox{ in }\Omega, \\
\llbracket \eta^n \rrbracket=0,\ \llbracket (\varpi^n I-\mathbb{D}\eta^n)e_3 \rrbracket=0 & \mbox{ on } \Sigma,\\
\eta^{n}=0 &\mbox{ on }\Sigma_-^+,  \end{array}\right.\end{equation}
where $\|\eta^{1}\|_4\leqslant 2c_1$. Moreover, from \eqref{201702141356nn} one gets
$$\|(\eta^n,\varpi^n) \|_{\mm{S},2} \leqslant  c_1 (1+\delta^3\|\eta^{n-1}\|_4^3)
\quad\mbox{for any }n\geqslant 2, $$
which implies
\begin{equation}   \label{201702090854}
\|(\eta^n,\varpi^n) \|_{\mm{S},2} \leqslant 2c_1
\end{equation}
for any $n\geqslant 2$, and any $\delta\leqslant 1/2c_1$. Recalling that  $\|\varpi^n\|_0^2=1$,
thus there exists a subsequence of $\{ \varpi^n\}_{n=1}^\infty$ denoted by $\{ \varpi^{n_k} \}_{k=1}^\infty$, such that
\begin{align}
& \varpi^{n_k}\to \varpi^{\mm{r}} \mbox{ weakly in }H^3, \label{201809111106} \\
&( \varpi^{n_k},\varpi^{n_k}_\pm|_{y_3=0} )\to (\varpi^{\mm{r}},\varpi^{\mm{r}}_\pm|_{y_3=0} ) \mbox{ strongly in }H^2 \times H^{3/2}(\mathbb{T}).\label{201809111106xx}
\end{align}

Next we want to show that $\{(\eta^{n},\nabla \varpi^{n},\llbracket  \varpi^{n }\rrbracket )\}_{n=1}^\infty$,   is a Cauchy sequence. Noting that
\begin{equation*}     \left\{\begin{array}{ll}
\nabla (\varpi^{n+1}-\varpi^n)-\mu \Delta( \eta^{n+1} - \eta^{n} )=0,\quad \mm{div}(\eta^{n+1}  -\eta^{n} )
=O(\eta^{n} )-O(\eta^{n-1} )&\mbox{ in }\Omega, \\
\llbracket \eta^{n+1} \rrbracket=0,\ \llbracket ((\varpi^{n+1}-\varpi^{n})  I-\mathbb{D}(\eta^{n+1}-\eta^{n} ))e_3 \rrbracket=0  & \mbox{ on } \Sigma,\\
\eta^{n+1} -\eta^{n} =0 &\mbox{ on }\Sigma_-^+,
\end{array}\right.\end{equation*}
we obtain
\begin{equation}  \label{Ellipticestimate0837}
\|(\eta^{n+1} -\eta^{n},\varpi^{n+1}-\varpi^n) \|_{\mm{S},2}    \leqslant c \|O(\eta^{n})-O(\eta^{n-1})\|_3.
\end{equation}
On the other hand, using \eqref{201702090854} and the product estimate \eqref{fgestims}, we arrive at
 $$\|O(\eta^{n} )-O(\eta^{n-1} )\|_3\leqslant c \delta\| \eta^{n}  - \eta^{n-1}  \|_4.$$
Substituting the above inequality into \eqref{Ellipticestimate0837}, one sees by taking $\delta$ appropriately small
that $\{(\eta^n,\nabla \varpi^n,\llbracket  \varpi^{n}    \rrbracket)\}_{n=1}^\infty$ is a Cauchy sequence in $H^4 \times  {H}^2 \times  H^{5/2}(\mathbb{T})$.
Thus, using \eqref{201809111106} and \eqref{201809111106xx}, we can derive that,
for some  limit function $\eta^{\mm{r}}$,
$$
\begin{aligned}
&(\eta^n,\nabla \varpi^n)\to ( \eta^{\mm{r}}, \nabla \varpi^{\mm{r}})\mbox{ in }H^3\times {H}^2,\\ &(\eta^n|_{y_3=0},\llbracket\varpi^n\rrbracket)\to ( \eta^{\mm{r}}|_{y_3=0},\llbracket\varpi^{\mm{r}}\rrbracket)\mbox{ in }H^{5/2}(\mathbb{T})\times  {H}^{3/2}(\mathbb{T}).
\end{aligned}
$$
Consequently, we can take to the limit in \eqref{201702051905} as $n\to\infty$ to see that the
limit function $(\eta^{\mm{r}},\varpi^{\mm{r}})$ solves
\begin{equation}
\label{20170514}
 \left\{\begin{array}{ll}
\nabla \varpi^{\mm{r}}-\mu \Delta \eta^{\mm{r}} =0,\quad  \mm{div}\eta^{\mm{r}} =O(\eta^{\mm{r}}) &\mbox{ in }\Omega,  \\
\llbracket \eta^{\mm{r}} \rrbracket=0,\  \llbracket ( \varpi^{\mm{r}} I-\mu\mathbb{D}\eta^{\mm{r }})e_3 \rrbracket=0& \mbox{ on } \Sigma\\
\eta^{\mm{r}} =0&\mbox{ on } \Sigma_-^+.
\end{array}\right.  \end{equation}
 Furthermore, by \eqref{201702090854},
\begin{align}
\label{201808311904}
\|\eta^\mm{r}\|_4\leqslant 2c_1.
\end{align}

(2) Now we turn to construct $(u^\mm{r},q^{\mm{r}})$, we consider a stratified Stokes problem:  for a given function $(w,p) \in H^2\times H^1 $,
\begin{equation}
\label{201702051905xfwea}
 \left\{\begin{array}{ll}
\nabla q- \mu \Delta u   =0, \
\mm{div} u
=- \mm{div}_{\tilde{\mathcal{A}}^{\delta}} ( \tilde{u}   + \delta w ) /\delta &\mbox{ in }\Omega, \\
\llbracket u\rrbracket =0, \ \llbracket(q I-\mu\mathbb{D}u)e_3\rrbracket
 =\mathcal{M}^6(w,p)&\mbox{ on } \Sigma,\\
u=0 &\mbox{ on }\Sigma_{-}^+,
\end{array}\right.\end{equation}
where $\eta^\mm{r}$ is provided by \eqref{20170514}, and we have defined that
$$
\begin{aligned}
& \mathcal{M}^6(w,p):=
\llbracket \kappa\rho\mathbb{D}\eta^{\mm{r}}+ (g \rho  \eta_3^{\mm{r}} + \vartheta \Delta_{\mm{h}}\eta_3^{\mm{r}} )e_3+ \bar{M}_3 \partial_{\bar{M}}\eta^{\mm{r}} \rrbracket\\
&\qquad \qquad\qquad +\mathcal{N}^2(\delta\tilde{\eta} +\delta^2 \eta^{\mm{r}} , \delta \tilde{u}  +\delta^2 w, \delta \tilde{q}  +\delta^2 p)/\delta^{2}.
\end{aligned}$$
Then, by the existence theory of stratified   Stokes problem, there exist a solution $(u,q)\in H^2\times  {H}^1$
 to \eqref{201702051905xfwea} with $\|q\|_0^2=1$;  moreover, similarly to \eqref{201702141356nn},
 \begin{align}
&\|(u,q)\|_{\mm{S},0}\lesssim  \| \mm{div}_{\tilde{\mathcal{A}}^{\delta}} ( \tilde{u}   + \delta w )\|_1 /\delta+| \mathcal{M}^6(w,p)|_{1/2}.
\label{Ellipticestimate0839nn}
\end{align}
Recalling \eqref{201808311904},  then, following the arguments of \eqref{06011733} and \eqref{06011734}, we can estimate that
$$  | \mathcal{M}^6(w,p)|_{1/2}\lesssim 1+\delta(\|w\|_2+\| \nabla p\|_0+
 | \llbracket p\rrbracket|_{1/2}  ) +\delta^2\|w\|_2^2.$$
Thus, by Young inequality,  we get,  for some constant $c_2$,
\begin{equation}
\label{201702141356nn123}
\|(u,q)\|_{\mm{S},0} \leqslant c_2 ( 1+\delta^2(\|w\|_2^2+\|\nabla p\|_0^2+  | \llbracket p\rrbracket|_{1/2}^2 )) .
\end{equation}

Therefore, one can construct an approximate function sequence $\{(u^{n}, q^n)\}_{n=1}^\infty$, such that, for any $n\geqslant 2$,
\begin{equation}
\label{201702051905xxxx}
 \left\{\begin{array}{ll}
\nabla q^{n}-\mu \Delta u^{n}=0,\quad
\mm{div}u^{n}=
- \mm{div}_{\tilde{\mathcal{A}}^{\delta}} ( \tilde{u}   + \delta u^{n-1} ) /\delta &\mbox{ in }\Omega, \\
\llbracket u^{n} \rrbracket=0,\
\llbracket (q^{n} I-\mu\mathbb{D} u^{n} )e_3  \rrbracket
=\mathcal{M}^6( u^{n-1}, q^{n-1})   &\mbox{ on } \Sigma,\\
u^{n}=0 &\mbox{ on }\Sigma_{-}^+,
\end{array}\right.\end{equation}
where $\| u^{1}\|_2+\|\nabla q^{1}\|_0+   | \llbracket q^1\rrbracket|_{1/2}  \leqslant c_2$.   Moreover,  by \eqref{201702141356nn123}, one has
$$ \|(u^n,q^n)\|_{\mm{S},0}   \leqslant  c_2(1+\delta^2(\|u^{n-1}\|_2^2+\| \nabla q^{n-1}\|_0^2
+ | \llbracket q^{n-1}\rrbracket|_{1/2}^2 ))$$
for any $n\geqslant 2$, which implies that  \begin{equation}
 \label{20170dfsf2090854}
\|(u^n,q^n)\|_{\mm{S},0} \leqslant 2c_2
 \end{equation} for any $n$, and any $\delta\leqslant 1/2c_2$.

 Next we further prove that  $\{(u^n,\nabla q^n,   \llbracket q^n\rrbracket)\}_{n=1}^\infty$ is a Cauchy sequence.
Noting that
$$
 \left\{\begin{array}{ll}
\nabla (q^{n+1}-q^n)-\mu \Delta (u^{n+1}-u^{n})= 0,&\mbox{ in }\Omega,
\\
\mm{div}(u^{n+1}-u^{n})
=- \mm{div}_{\tilde{\mathcal{A}}^{\delta}} ( u^{n} -   u^{n-1})  &\mbox{ in }\Omega, \\
\llbracket(u^{n+1}-u^{n})e_3\rrbracket=0 &\mbox{ on } \Sigma,\\ \llbracket ((q^{n+1}- q^n) I-\mu\mathbb{D}(u^{n+1}-u^{n}))e_3\rrbracket =
\mathcal{M}^6( u^n, q^n)-  \mathcal{M}^6( u^{n-1}, q^{n-1}) &\mbox{ on } \Sigma,\\
u^{n+1}-u^{n}=0 &\mbox{ on }\Sigma_{-}^+,
\end{array}\right.$$
thus we can use the stratified Stokes estimate  to get that
\begin{align}
&
\|(u^{n+1}-u^{n } , q^{n+1}-q^{n })\|_{\mm{S},0}  \nonumber \\
& \lesssim  \|\mm{div}_{\tilde{\mathcal{A}}^{\delta}} (   u^{n} -  u^{n-1}) \|_1 +|\mathcal{M}^6( u^n, q^n)-  \mathcal{M}^6( u^{n-1},q^{n-1})|_{1/2}=:\beth^n.
\label{Ellipticestimate0837df}
\end{align}
Following the arguments of \eqref{06011733} and \eqref{06011734}, it is easy to estimate that
$$\beth^n \leqslant  c  \delta \|(u^{n}-u^{n-1 } , q^{n}-q^{n-1 })\|_{\mm{S},0}.$$ Putting the above estimate into \eqref{Ellipticestimate0837df} yields
$$\begin{aligned}
\|(u^{n+1}-u^{n } , q^{n+1}-q^{n })\|_{\mm{S},0}  \lesssim c  \delta \|(u^{n}-u^{n-1 } , q^{n}-q^{n-1 })\|_{\mm{S},0} ,
\end{aligned}$$
which presents that $\{ (u^{ n},\nabla q^n,\llbracket q^n\rrbracket)\}_{n=1}^\infty $ is a Cauchy sequence in $H^2\times {H}^1\times H^{1/2}(\mathbb{T})$  by choosing a sufficiently small $\delta$.
Consequently, we can easily get
a limit function $(u^{\mm{r}},q^{\mm{r}})\in H^2\times {H}^1$ as in step (1), which solves
\begin{equation}
\nonumber
 \left\{\begin{array}{ll}
\nabla q^{\mm{r}}-\mu \Delta u^{\mm{r}} =0, \
\mm{div}u^{\mm{r}}=
- \mm{div}_{\tilde{\mathcal{A}}} ( \tilde{u}   + \delta u^{\mm{r}} ) /\delta &\mbox{ in }\Omega, \\
\llbracket u^{\mm{r}}  \rrbracket=0,\
\llbracket (q^{\mm{r}} I-\mu\mathbb{D} u^{\mm{r}} )e_3\rrbracket
= \mathcal{M}^6( u^{\mm{r}}, q^{\mm{r}})  &\mbox{ on } \Sigma,\\
u^{\mm{r}}=0 &\mbox{ on }\Sigma_{-}^+
\end{array}\right.\end{equation}
by \eqref{201702051905xxxx}. Moreover,  by \eqref{20170dfsf2090854},
$  \|(u^{\mm{r}},q^{\mm{r}})\|_{\mm{S},0}  \leqslant 2c_2$,
which, together with \eqref{201808311904}, yields  \eqref{201702091755}. \hfill $\Box$
\end{pf}

If $\delta$  is sufficiently small, then $\eta_0^\delta$ constructed in
 Proposition \ref{lem:modfied} belongs to $H_{0,*}^{3,1}$. This fact can be directly observed from the following conclusion:
\begin{pro}\label{201809012320}
Let $k\geqslant 3$.
There exists a constant $\delta_5$, depending on $\Omega$, such that, for any $\eta\in H_0^{k}$ satisfying
$\|\eta\|_3\leqslant \delta_5\in (0,1)$ and $\det \nabla \zeta=1$, we have
 \begin{align}
& \det\nabla_{\mm{h}}\zeta_{\mm{h}} (y_{\mm{h}},0)
\geqslant  {1}/{2}\mbox{ for any }y_{\mm{h}}\in \mathbb{T}, \label{201803121601xx21082109}
 \\
&\label{201803121601xx2108}
\zeta_{\mm{h}} (y_{\mm{h}},0):\mathbb{R}^2 \to \mathbb{R}^2\mbox{ is a }C^{k-2}\mbox{-diffeomorphic mapping},
\\
&\label{201803121601xx}
\zeta  : \overline{\Omega}\to \overline{\Omega} \mbox{ is a homeomorphism mapping},  \\
& \zeta _\pm : \Omega_\pm' \to \zeta _\pm(\Omega_\pm')\mbox{ are }C^{k-2}\mbox{-diffeomorphic mapping} , \label{2018031adsadfa21601xx}\\
&\mathcal{T}^{-{1/2}}_{-l,\tau} \subset \zeta^{-1}(\mathcal{T}_{-l}^\tau)
\subset \mathcal{T}^{{1/2}}_{-l,\tau}, \ \mathcal{T}^{-{1/2}}  \subset (\zeta_{\mm{h}} )^{-1}(\mathcal{T}) \subset \mathcal{T}^{1/2},\label{2011801941224538}
\end{align}
where $\zeta:=\eta+y$, $(\zeta_{\mm{h}} )^{-1}$ denotes the inverse mapping of $\zeta_{\mm{h}} (y_{\mm{h}},0)$, and $\zeta_{\mm{h}} $ denotes the first two components of $\zeta :={\eta} +y$.
\end{pro}
\begin{rem} \label{201809273138}
Since $\eta$ is a horizontally periodic function, then, by \eqref{201803121601xx2108} and \eqref{201803121601xx}, we can verify that $(\zeta_{\mm{h}} )^{-1}(x_{\mm{h}})-x_{\mm{h}}$  and $\zeta^{-1}(x)-x$ are horizontally periodic functions as well as $\eta$.
\end{rem}
  \begin{pf}
Let $\eta\in H_0^{k}$ satisfy $\|\eta\|_3\leqslant \delta_5$ and $\det \nabla \zeta=1$, where $\zeta:=\eta+y$. By the embedding $H^k(\Omega_\pm)\hookrightarrow C^{k-2}(\overline{\Omega_\pm})$, we have $\eta_\pm \in C^{k-2}(\overline{\Omega}_\pm)$
and $\eta\in C^0(\overline{\Omega})$.
Moreover we also have, for any sufficiently small $\delta$,
\begin{align}
& \label{201809011630}
|\eta _i(y)|\leqslant L_i\pi/2\mbox{ for any }y\in \overline{\Omega},\\
&-l<\inf_{y_{\mm{h}}\in \mathbb{R}^2 }\{ \zeta_3(y_{\mm{h}},0)\}\leqslant\sup_{y_{\mm{h}}\in \mathbb{R}^2 }\{ \zeta_3(y_{\mm{h}},0)\}<  \tau ,\label{201809111914}
\end{align}
where $i=1$  or $2$.
Next we prove the results in \eqref{201803121601xx21082109}--\eqref{2011801941224538} in sequence.

(1) The infimum \eqref{201803121601xx21082109} for sufficiently small $\delta_5$ is obvious by the embedding inequality \eqref{esmmdforinfty}.

(2) We  turn to verify \eqref{201803121601xx2108} for sufficiently small $\delta_5$ by three steps.

(a)
We choose a cut-off function $\xi\in C_0^\infty(\mathbb{R}^2)$ such that
$\xi=1$ in $\mathcal{T}^{1} $, and $0\leqslant \xi\leqslant 1$ in $\mathbb{R}^2$.  Then there exists a disk $b_{R}(0)$ of radius $R$ and center $0\in   \mathbb{R}^2$, such that $\xi=0$ in $\mathbb{R}^2 \backslash b_{R}(0)$ and $\mathcal{T}^{1} \subset b_{R}(0)$. Let $\eta^\xi_-:=\eta_-  \xi$. Then $\eta^\xi_- \in C^{k-2}(\overline{\Omega_-})$, which yields that $\eta^\xi_- (y_{\mm{h}},0)\in C^{k-2}(\mathbb{R}^2)$.
By the embedding inequality \eqref{esmmdforinfty}, we easily see that $(\zeta^\xi_-)_{\mm{h}}:=(\eta^\xi_-)_{\mm{h}}+y_{\mm{h}}:b_R(0)\to \mathbb{R}^2$ is injective by Lemma \ref{xfsddfsf201803121937} for sufficiently small $\delta_5$, where $(\eta^\xi_-)_{\mm{h}}$ denotes the first two components of $\eta^\xi_-$. Noting that  $\zeta_{\mm{h}}(y_{\mm{h}},0)=(\zeta^\xi_-)_{\mm{h}}$ in $\mathcal{T}^1$, then
 \begin{equation*}
 \zeta_{\mm{h}}(y_{\mm{h}},0): \mathcal{T}^1\to  \mathbb{R}^2\mbox{ is also injective}.
 \end{equation*}
Thus, by \eqref{201809011630}  and the periodicity of $\zeta_{\mm{h}}(y_{\mm{h}},0)-y_{\mm{h}}$, we easily see that \begin{equation}
 \label{2018050051441sdfaf}
 \zeta_{\mm{h}}(y_{\mm{h}},0): \mathbb{R}^2\to  \mathbb{R}^2\mbox{ is also injective}.
 \end{equation}

(b) We will further prove that \begin{equation}
\label{2018safasdf04121545dfadf}
 \zeta_{\mm{h}}(y_{\mm{h}},0): \mathbb{R}^2\to  \mathbb{R}^2 \mbox{ is surjective}.
\end{equation}

We define that
$$\mathbb{J}:=\{x_{\mm{h}}\in \mathbb{R}^2 ~|~x_{\mm{h}}\neq \zeta_{\mm{h}}(y_{\mm{h}},0)\mbox{ for any }y_{\mm{h}}\in \mathbb{R}^2\}.$$
Obviously, to get \eqref{2018safasdf04121545dfadf}, it suffices to prove that
\begin{equation}
\label{201804151053xx}
\mathbb{J}=\emptyset.
\end{equation}
Next we verify \eqref{201804151053xx} by contradiction.

We assume that $\mathbb{J}\neq \emptyset$. Then there exists a point
\begin{equation*}
x^0_{\mm{h}}\in \partial \mathbb{J}.
\end{equation*}
Thus there exists a disk $b_{\delta_0}(x^0_{\mm{h}})\subset \mathbb{R}^2 $, such that, for any $b_{\delta}(x^0_{\mm{h}})\subset b_{\delta_0}(x^0_{\mm{h}})$,
\begin{equation}
\label{201804151502xx}
b_{\delta}(x^0_{\mm{h}})\cap \mathbb{J}\neq\emptyset
\end{equation}
and
\begin{equation}
\label{201804152200xx}
b_{\delta}(x^0_{\mm{h}})\cap (\mathbb{R}^2\setminus \mathbb{J})\neq\emptyset.
\end{equation}
By \eqref{201804152200xx}, we can choose two sequences $\{x^n_{\mm{h}}\}_{n=m}^\infty\subset (\mathbb{R}^2 \setminus \mathbb{J})$ and $\{y^n_{\mm{h}}\}_{n=m}^\infty\subset \mathbb{R}^2$ for some $m>1/\delta_0$, such that $x^n_{\mm{h}}\in b_{1/n}(x^0_{\mm{h}})\subset b_{\delta_0}(x^0_{\mm{h}}) $, $x^n_{\mm{h}}\to x^0_{\mm{h}} $ and $x^n_{\mm{h} }=\zeta(y^n_{\mm{h}})$. Since $(\zeta-y)\in  C^0(\overline{\Omega})$, then $\{y^n_{\mm{h}}\}_{n=m}^\infty$ is a bounded sequence. Therefore, there exists a subsequence (still labeled by $y^n_{\mm{h}}$)
such that $y^n_{\mm{h}}\to y^0_{\mm{h}}$ for some $y^0_{\mm{h}}\in \mathbb{R}^2$. By the continuity of $\zeta$,  $x^n_{\mm{h}}=\zeta_{\mm{h}}(y^n_{\mm{h}},0)\to x^0_{\mm{h}}=\zeta_{\mm{h}}(y^0_{\mm{h}},0)$.
Noting that $\nabla \zeta_{\mm{h}}(y_{\mm{h}},0)$ is invertible for $y^0_{\mm{h}}$ by \eqref{201803121601xx21082109},
thus, by the well-known inverse function theorem (see Lemma \ref{xfsddfsf201805081359})
there exist two open sets $U$, $V\subset \mathbb{R}^2$ such that
$y^0_{\mm{h}}\in U$, $x^0_{\mm{h}}\in V$, and  $\zeta_{\mm{h}}(U,0)=V$, which imply that there exists a ball
$b_{\delta_1}(x^0_{\mm{h}})\subset (\mathbb{R}^2\setminus \mathbb{J})$ for $\delta_1<\delta_0$, i.e. $b_{\delta_1}(x^0_{\mm{h}})\cap \mathbb{J}=\emptyset$, which contradicts with to \eqref{201804151502xx}. Hence \eqref{201804151053xx} holds.

(c) By \eqref{2018050051441sdfaf} and \eqref{2018safasdf04121545dfadf}, $\zeta_{\mm{h}}(y_{\mm{h}},0) :\mathbb{R}^2\to \mathbb{R}^2$ is bijective.
Then we can consider an inverse mapping (denoted by $(\zeta_{\mm{h}})^{-1}$) of $\zeta_{\mm{h}}(y_{\mm{h}},0) $  defined on $\mathbb{R}^2$
  by \begin{equation*}
(\zeta_{\mm{h}})^{-1}(\zeta_{\mm{h}}(y_{\mm{h}},0))=y\mbox{ for }y_{\mm{h}}\in \mathbb{R}^2.
 \end{equation*}
 Obviously, $(\zeta_{\mm{h}})^{-1} :\mathbb{R}^2\to \mathbb{R}^2$ is bijective.
 Moreover,
\begin{equation}
\label{201809110852}
(\zeta_{\mm{h}})^{-1}\in C^{k-2}(\mathbb{R}^2).
\end{equation}
In fact,  using the second conclusion of Lemma
\ref{xfsddfsf201805081359}, we see that $(\zeta_{\mm{h}})^{-1}\in C^1(\mathbb{R}^2)$ and
\begin{equation}
\label{201811132040}
\nabla_{x_{\mm{h}}} (\zeta_{\mm{h}})^{-1}(x_{\mm{h}})=(\nabla_{y_\mm{h}} \zeta_{\mm{h}}(y_{\mm{h}},0))^{-1}|_{y_{\mm{h}}=(\zeta_{\mm{h}})^{-1}(x_{\mm{h}}) }.
\end{equation}
Here $x_i$ resp. $y_i$ represent the $i$-th component of $x_{\mm{h}}$ resp. $y_{\mm{h}}$.  By  \eqref{201811132040} and the regularity $\eta\in H^k$, we immediately get \eqref{201809110852}.

Finally, using  \eqref{2018050051441sdfaf}, \eqref{2018safasdf04121545dfadf} and \eqref{201809110852}, we   get \eqref{201803121601xx2108}.

(3) Now we verify \eqref{201803121601xx} for sufficiently small $\delta$ by three steps.

(a) To begin with, we proceed to the proof
\begin{equation}
\label{201804121545}
\zeta :\overline{\Omega}\to \overline{\Omega}\mbox{ is injective for sufficiently small }\delta.
\end{equation}

Let $d(x_{\mm{h}})=\zeta_3((\zeta_{\mm{h}})^{-1}(x_{\mm{h}}),0) $. Then, by \eqref{201809111914} and \eqref{201809110852},
$d(x_{\mm{h}})\in C^{k-2}(\mathbb{R}^2)$ and
\begin{equation}
 -l<\inf_{x_{\mm{h}}\in \mathbb{R}^2 }\{d(x_{\mm{h}}) \}\leqslant\sup_{x_{\mm{h}}\in \mathbb{R}^2 }\{ d(x_{\mm{h}} )\}<  \tau ,
\label{201809231554}
\end{equation}
We denote
$$
 {\Omega^-}:=\{x\in \mathbb{R}^3~|~ -l<x_3< d(x_{\mm{h}} )\}\mbox{ and } {\Omega^+}:=\{x\in \mathbb{R}^3~|~d(x_{\mm{h}} )< x_3< \tau\}.
$$
By the properties of $d(x_{\mm{h}})$ and $\zeta(y_{\mm{h}},0)$, we can see the following properties of $\Omega^\pm$:
\begin{align}
& \overline{\Omega^-}:=\{x\in \mathbb{R}^3~|~ -l\leqslant x_3\leqslant d(x_{\mm{h}} ) \},\  \overline{\Omega^+}:=\{x\in \mathbb{R}^3~|~d(x_{\mm{h}} )\leqslant  x_3\leqslant \tau\},\nonumber \\
& \overline{\Omega^-}\cap \Omega^+=\Omega^-\cap \overline{\Omega^+}=\emptyset,\ \overline{\Omega^-}\cup \overline{ \Omega^+}=\overline{\Omega}, \nonumber\\
& \overline{\Omega^-}\cap\overline{\Omega^+}=\{x_3=d(x_{\mm{h}})\}=\zeta(\{y_3=0\}). \label{201809231444}
\end{align}
Moreover, we also see that the proof of
\eqref{201804121545} reduces to the proof of
\begin{equation}
\label{201805050954}
\zeta_\pm  :\overline{\Omega_\pm}\to \overline{\Omega^\pm}\mbox{ are injective}
\end{equation}
for sufficiently small $\delta$.
Next we only provide the proof for $\zeta_-$, since $\zeta_+$ can be similarly proved.

Let $S:=\overline{\mathcal{T}^{1}_{-l,0}}$.
We choose a cut-off function $\chi\in C_0^\infty(\mathbb{R}^3)$ such that
$\chi=1$ in $S$, and $0\leqslant \chi\leqslant 1$ in $\mathbb{R}^3$.
Then there exists a ball $B_{R}(0)$ of radius $R$ and center $0\in \mathbb{R}^3$, such that
$\chi=0$ in $\mathbb{R}^3 \backslash B_{R}(0)$ and $S\subset B_{R}(0)$.
Let $\eta^\chi_-:=\eta_- \chi$, then $\eta^\chi_- \in H^k(\Omega_-)\cap C^{k-2}(\overline{\Omega_-})$.

{Since $\Omega_-'$ is locally Lipschitz} (see \cite[Section 4.9]{ARAJJFF} for the definition),
by virtue of the well-known Stein extension theorem (see Lemma \ref{xfsddfsf2018050813379}), there is an extension operator for $\Omega_-'$, such that
\begin{equation}
\label{20180312}
E(\eta^\chi_-)=\eta^\chi_-\mbox{ a.e. in }\Omega_-' ,
\quad \|E(\eta^\chi_-)\|_{H^3(\mathbb{R}^3)}\leqslant{c}\|\eta_-^\chi\|_{H^3(\Omega_-)},
\end{equation}
where the constant ${c}$ depends on $\Omega_-'$.

 Define $\zeta^E_-:=y+\chi E(\eta^\chi_-)$, then
\begin{equation}   \label{201803121926}
\zeta^{E}_- =\zeta_- \mbox{ in } S .
\end{equation}
From  \eqref{20180312} and the periodicity of $\eta$, it follows that
 $$\| E(\eta^\chi_-)\|_{H^3(\mathbb{R}^3)}\leqslant{c}\|\chi \eta_-\|_{H^3(\Omega_-')} \leqslant c\|\eta\|_3.$$  Thus, in terms of the embedding inequality \eqref{esmmdforinfty},
 $\nabla \zeta^{E}_-$ is invertible in $\mathbb{R}^3$
for sufficiently small $\delta_5$.

 Since $\zeta^{E}_-(y)=y$ on $\partial B_R(0)$, then $\zeta^{E}_-: B_R(0)\to \mathbb{R}^N$
 is injective  by Lemma \ref{xfsddfsf201803121937}.
 Noting that $S\subset B_R(0)$, we see by \eqref{201803121926} that
 \begin{equation*}
 \zeta_- :S\to  \zeta_- (\overline{\Omega_-})\mbox{ is also injective},
 \end{equation*}
which, together with \eqref{201809011630}, yields
 \begin{equation}
 \label{2018050051441}
\zeta_- :\overline{\Omega_-}\to  \zeta_- (\overline{\Omega_-})\mbox{ is  injective}.
\end{equation}

Now, we further show that
\begin{equation}
\label{2018050541438}
\zeta_- (\overline{\Omega_-})\subset \overline{\Omega^-}.
\end{equation}
To this end, it suffices to prove that, for any given $y_{\mm{h}} \in \mathbb{R}^2$,
$$\zeta _-: l_{y_{\mm{h}}}\to \overline{\Omega^-},$$
where $l_{y_{\mm{h}}} :=\{(y_{\mm{h}},y_3)~|~y_3\in [-l,0] \}$. We prove this by contradiction.

Assume that there exists a point $y^0\in l_{y_{\mm{h}}}$, such that $\zeta_- ( y^0)\in  \!\!\!\!\!/ \  \overline{\Omega^-}$.
Obviously, $y_3^0\in (-l, 0)$, where $y^0_3$ denotes the third component of $y^0$. Without loss of generalization, we assume that the point $\zeta _-(y^0)$ is above the closed set $\overline{\Omega_-}$. Thus, the curve function $\zeta _-$ defined on $l_0:=\{y \in l_{y_{\mm{h}}}~|~y_3\in (-l,y^0_3)\}$ must go through the upper boundary of $\zeta(\{y_3=0\})$. We denote by $z$ the intersection of the curve and the upper boundary. This means that there is a point
$y\in l_0$, such that $y_3\in (-l,y^0_3)$ and $\zeta _-(y)=z$.
However, there exists a point $\sigma\in \{y_3=0\}$  such that $\zeta _-(\sigma)=z$. Thus we have $\zeta _-(y)=\zeta _-(\sigma)$ where $y$, $\sigma\in \overline{\Omega_-}$ and $y\neq \sigma$. This contradicts with injectivity of $\zeta _-$ on $\overline{\Omega_-}$. Hence \eqref{2018050541438} holds. Consequently, by \eqref{2018050051441} and \eqref{2018050541438}, we easily get \eqref{201805050954}.

(b) Now we turn to prove that
\begin{equation}
\label{201804121545dfadfxx}
\zeta:\overline{\Omega}\to \overline{\Omega}\mbox{ is surjective},
\end{equation}
which reduces to the proof of
\begin{equation}
\label{201804121545dfadf}
\zeta_\pm:\overline{\Omega_\pm}\to \overline{\Omega^\pm}\mbox{ is surjective}.
\end{equation}
Next we only provide the proof for $\zeta_-$, since the case of $\zeta_+$ can be similarly proved.

We define that
$$\mathbb{K}:=\{x\in \overline{\Omega^-} ~|~x\neq \zeta(y)\mbox{ for any }y\in \overline{\Omega_-}\ \}.$$
By \eqref{201809231444} and the fact $\zeta(y)|_{\Sigma_-}=y$, we see that $\mathbb{K}\subset  \Omega^-$. Obviously, to get the surjective conclusion of $\zeta_-$, it suffices to prove that
\begin{equation}
\label{201804151053}
\mathbb{K}=\emptyset.
\end{equation}
We  can also verify \eqref{201804151053} by contradiction.

In fact, we assume that $\mathbb{K}\neq \emptyset$. Noting that
$\mathbb{K} \cup \zeta(\Omega_-)\subset \Omega^-$ and \eqref{201809231554}, then there exists a point
\begin{equation*}
x^0\in \partial \mathbb{K}\cap \Omega^-.
\end{equation*}
Noting that $\Omega^-$ is open, then, there exists a ball $B_{\delta_0}(x^0)\subset {\Omega^-}$, such that, for any $B_{\delta}(x^0)\subset B_{\delta_0}(x^0)$,
\begin{equation}
\label{201809111954}
B_{\delta}(x^0)\cap \mathbb{K}\neq\emptyset\mbox{ and }
B_{\delta}(x^0)\cap (\Omega^-\setminus \mathbb{K})\neq\emptyset.
\end{equation}
Thus, following the argument of \eqref{201804151053xx}, we also have a result which contracts with the first relation in
\eqref{201809111954}. Hence \eqref{201804151053} holds.

(c) By \eqref{201804121545} and \eqref{201804121545dfadfxx}, $\zeta :\overline{\Omega}\to \overline{\Omega}$ is bijective.
Then we can consider an inverse mapping (denoted by $\zeta^{-1}$) of $\zeta$  defined on $\overline{\Omega}$
  by \begin{equation*}
\zeta^{-1} (\zeta (y))=y\mbox{ for }y\in \overline{\Omega}.
 \end{equation*}
 Obviously, $\zeta^{-1}  :\overline{\Omega}\to \overline{\Omega}$ is bijective.
Next we verify that
\begin{equation}
\label{201804151826}
\zeta^{-1} \mbox{ is a continuous mapping of }\overline{\Omega}\mbox{ onto } \overline{\Omega}
\end{equation}
by contradiction.

We assume that $\zeta^{-1}$ is not continuous for some $x^0\in \overline{\Omega}$. Then, there exists a constant $\varepsilon>0$, such that, for any $\iota>0$, there exists a point $x^\iota\in \overline{\Omega}$ satisfying
$|x^\iota-x^0|<\iota$ and
\begin{equation}
\label{201804151826xx}
|\zeta^{-1} (x^\iota)- \zeta^{-1} (x^0)|\geqslant  \varepsilon.
\end{equation}
Let $\iota=1/n$, we denote $y^n=\zeta^{-1} (x^{1/n})$.  Since $(\zeta -y)\in   C^0(\overline{\Omega})$, and $\{x^{1/n}\}_{n=1}^\infty\subset B_1(x^0)\cap \overline{\Omega}$, then $\{y^n\}_{n=1}^\infty$ is a bounded sequence. Thus there exists a subsequence (still labeled by $y^n$) such that $y^n \to y^0$ for some $y^0 \in \overline{\Omega}$. By the continuity of $\zeta $ and the fact $x^{1/n}\to x^0$ as $n\to \infty$, $x^{1/n}=\zeta (y^n)\to \zeta ( y^0)=x^0$ as $n\to \infty$. Thus $ \zeta^{-1} (x^{1/n})=y^n\to y^0=\zeta^{-1}(x^0)$, which contracts with \eqref{201804151826xx}. Hence \eqref{201804151826} holds. Consequently, we obtain \eqref{201803121601xx} from \eqref{201804121545}, \eqref{201804121545dfadfxx} and \eqref{201804151826}.

(4) By \eqref{201805050954} and \eqref{201804121545dfadf}, we see that
$$\zeta_\pm :\overline{\Omega_\pm}\to \overline{\Omega^\pm}\mbox{ are bijective}.$$
By \eqref{201809231444} and $\zeta(y)|_\Sigma=y$, we further have
\begin{equation}
\label{201809231510}
\zeta_\pm : {\Omega_\pm'}\to {\Omega^\pm}\mbox{ are bijective}.
\end{equation}
In addition, following the argument of \eqref{201809110852},  we also  $\zeta ^{-1}_\pm\in C^1(\Omega^\pm)$, which, together with \eqref{201809231510}, yields  \eqref{2018031adsadfa21601xx}.

(5) By the embedding inequality \eqref{esmmdforinfty}, we have, for sufficiently small $\delta$,
\begin{align}
&\zeta (\mathcal{T}^{-{1/2}}_{-l,\tau} )\subset \mathcal{T}_{-l}^\tau
\subset \zeta ( \mathcal{T}^{{1/2}}_{-l,\tau}),\nonumber \\
&\zeta _{\mm{h}}(\mathcal{T}^{-{1/2}}\times \{0\})  \subset  \mathcal{T} \subset \zeta _{\mm{h}}( \mathcal{T}^{1/2}\times \{0\}).\nonumber
\end{align}
Consequently, by \eqref{201803121601xx2108} and \eqref{201803121601xx}, we immediately get \eqref{2011801941224538} from the above relations.
This completes the proof of Proposition \ref{201809012320}. \hfill $\Box$
\end{pf}

\section{Error estimates}\label{201807111049}

This section is devoted to the derivation of error estimates  between the solutions
 of the nonlinear ART problem and the solutions of linearized ART problem.

To begin with, let us introduce the estimate
\begin{equation}
\label{201811201409}
|\phi|_{3}\leqslant |\phi|_{7/2}\leqslant C_3(\|\chi\|_4)\|\chi\|_4\mbox{ for any }\chi\in H^{4,1}_{0,*},
\end{equation}
where we have defined that $\phi(x_\mm{h}):=\tilde{\chi}_3((\tilde{\chi}_{\mm{h}})^{-1}(x_{\mm{h}}),0)$, $\tilde{\chi}:=\chi+y$, and the positive constant $C_3(\|\chi\|_4)$ increasing with respect to $\|\chi\|_4$. Please refer to \eqref{201809292130} in Section \ref{201811142125} for derivation of \eqref{201811201409}.

We define that
$$
\begin{aligned}
&C_4:=(1+C_3(\|\tilde{\eta}^0\|_4+C_2))\left(\sqrt{\|\tilde{\eta}^0\|_4^2+ \|\tilde{u}^0\|_2^2 } +C_2\right)\geqslant 1,\\
 & \delta<\delta_0:= \min\{\delta_1,\delta_2,\delta_3,C_4\delta_4,\delta_5\}/2C_4<1,
 \end{aligned}$$
and $(\eta_0^\delta,u_0^\delta)$ be constructed by Proposition \ref{lem:modfied}.
Let $\zeta^0:=\eta^0+y $ and $d(x_\mm{h}):=\zeta^0_3((\zeta_{\mm{h}}^0)^{-1}(x_{\mm{h}}),0)$. Then, by
\eqref{mmmode04091215}, \eqref{201702091755} and \eqref{201811201409},
$$|d_0^\delta|_3\leqslant C_3(\|\tilde{\eta}^0\|_4+C_2 )(\|\tilde{\eta}^0\|_4+C_2).$$
Thus, using \eqref{mmmode04091215} and \eqref{201702091755} again, we see that
$$\sqrt{\| {\eta}_0^\delta\|_3^2+ \|{u}_0^\delta\|_2^2+\beta(|d^\delta_0|_3^2) }<\delta_2,\ \delta_3 \mbox{ or }\delta_5.$$
Thus, by Proposition \ref{pro:0401nxd}, there exists a (nonlinear) solution $(\eta,u )$ with an associated  pressure $q$ of the ART problem with initial value $(  {\eta}_0^\delta ,{u}_0^\delta )$. Moreover, by Proposition \ref{201809012320},
${\eta}_0^\delta\in H_{0,*}^{4,1}$.

Now we estimate the error between the (nonlinear) solution $(\eta,u )$ and the linear solution $(\eta^\mm{a},u^{\mm{a}})$ provided by \eqref{0501}.
To this purpose, we define a error function
$(\eta^{\mathrm{d}}, u^{\mathrm{d}}):=(\eta, u)-(\eta^\mm{a},u^{\mm{a}})$.
Then we can establish the following estimate of error function.
\begin{pro}\label{lem:0401xx}
Let  $(\eta,u)$ be the strong  solution of the ART problem with initial value $(  {\eta}_0^\delta ,{u}_0^\delta )$ with an associated pressure $q$,  and $\theta\in (0,1)$  a small constant such that \eqref{aimdse}--\eqref{06041533fwqg}, \eqref{06011734} and \eqref{201806291} hold  for $(\eta,u)$ satisfying \eqref{aprpiosesnew} with $\theta$  in place of $\delta$. For given constants  $\gamma$  and $\beta$, if
\begin{align}
& \label{201702100940n} \sqrt{\|\eta\|_3^2+
\| u \|_2^2}\leqslant \theta,\ \delta e^{\Lambda t}\leqslant \beta,\\
& \|\eta(t)\|_3^2+\| (u(t), q(t))\|_{\mm{S},0}^2  +
\| u_t\|_0^2 \nonumber\\
& +
\int_0^t\left(\| u(\tau)\|_3^2+\|\nabla q(\tau)\|_{1}^2+\|u_{\tau}\|_1^2\right)\mm{d}\tau \leqslant (\gamma \delta e^{\Lambda t})^2 \label{201702100940}
\end{align}
on some interval $(0,T)$, then there exists a constant $C_5$ such that, for  any  $\delta\in (0,1)$ and any $t\in (0,T)$,
\begin{align}
\label{ereroe}
&|\chi^{\mm{d}}|_{L^1 }+ \|\chi^{\mm{d}}\|_{\mathfrak{X}}+ \| \mm{div}_{\mm{h}}\chi^{\mm{d}}_{\mm{h}}\|_{L^1}+\|u_t^{\mm{d}}\|_{L^1}   \leqslant C_5\sqrt{\delta^3e^{3\Lambda t}}, \\
&\|\partial_{\bar{M}}\chi_3^{\mm{d}} (t)\|_{L^1} , \ \|\partial_{\bar{M}}\chi_{\mm{h}}^{\mm{d}} (t)\|_{L^1} \leqslant C_5\sqrt{\delta^3 e^{3\Lambda t}},\label{201805291548}\\
&\|(\mathcal{A}_{3k}\partial_k\chi_3
-\partial_3 \chi^{\mm{a}}_3)(t)
\|_{{L^1}}\leqslant C_5\sqrt{\delta^3 e^{3\Lambda t}},\label{2018090119021} \\
&\|(\mathcal{A}_{3k}\partial_k\chi_{\mm{h}}
- \partial_3\chi^{\mm{a}}_{\mm{h}})(t)\|_{{L^1}}\leqslant C_5\sqrt{\delta^3 e^{3\Lambda t}},\label{201809011902} \\
&\|(\mathcal{A}_{1k}\partial_k\chi_{1}+ \mathcal{A}_{2k}\partial_k\chi_{2} -\mm{div}_{\mm{h}}\chi_{\mm{h}}^{\mm{a}})(t)\|_{L^1} \leqslant C_5\sqrt{\delta^3 e^{3\Lambda t}}, \label{201809231547}
\end{align}
where $\chi=\eta$ or $u$, $ \mathfrak{X}=W^{1,1}$ or $H^1$, and $C_5$ is independent  of  $T$.
\end{pro}
\begin{pf} We divide the proof of Proposition \ref{lem:0401xx} by four steps.

(1) Derivation of energy inequality.

Subtracting the ART  and the linearized ART problems, we get
\begin{equation}\label{201702052209}\left\{\begin{array}{ll}
\eta_t^{\mathrm{d}}=u^{\mathrm{d}} &\mbox{ in } \Omega,\\[1mm]
\rho u_{t}^{\mathrm{d}}+ \mathrm{div}_{\mathcal{A}}\mathcal{S}_{\mathcal{A}} ( q ,u ,\eta)-\mathrm{div}\mathcal{S} ( q^{\mm{a}},u^{\mm{a}},\eta^{\mm{a}})
= \partial_{\bar{M}}^2  \eta^{\mathrm{d}}& \mbox{ in }\Omega, \\
\div_{\mathcal{A}}  u^{\mm{d}}=-  \div_{\tilde{\mathcal{A}}}  u^{\mathrm{a}} & \mbox{ in }\Omega,\\
\llbracket u^{\mm{d}} \llbracket =0, \ \llbracket \mathcal{S}_{\mathcal{A}}( q,u,\eta)\mathcal{A}e_3 - \mathcal{S} ( q^{\mm{a}},u^{\mm{a}},\eta^{\mm{a}})e_3
&\\
-  g\rho \eta_3^{\mathrm{d}} e_3 - \bar{M}_3^2\partial_3 \eta^{\mathrm{d}} \rrbracket  =\vartheta H\mathcal{A}e_3  -\vartheta\Delta_{\mm{h}}\eta_3^{\mathrm{a}}e_3& \mbox{ on }\Sigma,\\
 (\eta^{\mathrm{d}},u^{\mathrm{d}})=0 &\mbox{ on }\Sigma_{-}^+,\\
 (\eta^{\mathrm{d}},u^{\mathrm{d}})|_{t=0}=\delta^2 (\eta^{\mm{r}},u^{\mm{r}}) &\mbox{ on }\Omega. \end{array}\right.\end{equation}
Then we formally apply $\partial_t $ to \eqref{201702052209}$_2$--\eqref{201702052209}$_5$ to get
\begin{equation}\label{201702052221}\left\{\begin{array}{ll}
   \rho u_{tt}^{\mathrm{d}}+\mm{div}_{\mathcal{A}}
\mathfrak{S}_{\mathcal{A}}(q^{\mathrm{d}}_t, u^\mm{d}_t, u^\mm{d})
= \partial_{\bar{M}}^2 u^{\mathrm{d}}+\mu \mathrm{div} \mathbb{D}_{\tilde{\mathcal{A}}}u^{\mm{a}}_t-\mathrm{div}_{\tilde{\mathcal{A}}}\mathfrak{S}_{\mathcal{A}}  ( q^{\mm{a}}_t,u^{\mm{a}}_t,u^{\mm{a}})+\mathcal{N}^3 & \mbox{ in }\Omega, \\
\div_{\mathcal{A}}  u^{\mathrm{d}}_t =-\mm{div}_{\mathcal{A}_t}u- \mm{div}_{\tilde{\mathcal{A}}}u^{\mm{a}}_t  & \mbox{ in }\Omega,\\
 \llbracket u_t^{\mm{d}} \llbracket =0,\ \llbracket \mathfrak{S}_{\mathcal{A}}(q^{\mathrm{d}}_t,u^{\mathrm{d}}_t,u^{\mathrm{d}})\mathcal{A}e_3
 - g\rho u_3^{\mathrm{d}} e_3
- \bar{M}_3 \partial_{\bar{M}}u^{\mathrm{d}}  \rrbracket &\\
= \vartheta\Delta_{\mm{h}}u_3^{\mathrm{d}} e_3 +\llbracket\mu\mathbb{D}_{\tilde{\mathcal{A}}} u^{\mm{a}}_t e_3- \mathfrak{S}_{\mathcal{A}} ( q^{\mm{a}}_t,u^{\mm{a}}_t,u^{\mm{a}}_t)\tilde{\mathcal{A}}e_3
\rrbracket +
\mathcal{N}^4. & \mbox{ on }\Sigma,\\
 u^{\mathrm{d}} _t=0 &\mbox{ on }\Sigma_{-}^+.
\end{array}\right.\end{equation}

Noting that
$$\begin{aligned}
 \int q_t  \mm{div}_{\tilde{\mathcal{A}}}u^{\mm{a}}_t\mm{d}y=& -\frac{\mm{d}}{\mm{d}t} \left(\int \nabla q \cdot( \tilde{\mathcal{A}}^{\mm{T}}u^{\mm{a}}_t)\mm{d}y+ \int_\Sigma\llbracket q \rrbracket \tilde{\mathcal{A}}e_3 \cdot u^{\mm{a}}_t\mm{d}y\right)\\
&+ \int \nabla q \cdot\partial_t( \tilde{\mathcal{A}}^{\mm{T}}u^{\mm{a}}_t)\mm{d}y+ \int_\Sigma\llbracket q \rrbracket \partial_t( \tilde{\mathcal{A}}e_3\cdot u^{\mm{a}}_t)\mm{d}y_{\mm{h}},
\end{aligned}$$
thus, following the argument of \eqref{0425sdfasfa},
we can formally deduce from \eqref{201702052221} that
\begin{align}
& \frac{1}{2}\frac{\mm{d}}{\mm{d}t}\left(\|\sqrt{\rho}  u_{t}^{\mm{d}}\|^2_0+ {E}( u^{\mm{d}})-2\int\left( \nabla q^{\mm{d}} \cdot (\mathcal{A}_t^{\mm{T}}u)+\nabla q \cdot( \tilde{\mathcal{A}}^{\mm{T}}u^{\mm{a}}_t)\right)\mm{d}y\right.\nonumber \\
& \left. -2\int_\Sigma \left(\llbracket
q^{\mm{d}} \rrbracket \mathcal{A}_te_3\cdot u+\llbracket q \rrbracket \tilde{\mathcal{A}} e_3\cdot  u^{\mm{a}}_t \right)\mm{d}y_{\mm{h}}\right)+\frac{1}{2}\|\sqrt{\mu} \mathbb{D}_{\mathcal{A}}  u_{t}^{\mm{d}} \|_{0}^2\nonumber \\
&= \int \left( \left(\mu \mathrm{div} \mathbb{D}_{\tilde{\mathcal{A}}}u^{\mm{a}}_t -\mathrm{div}_{\tilde{\mathcal{A}}}\mathfrak{S}_{\mathcal{A}} ( q^{\mm{a}}_t,u^{\mm{a}}_t,u^{\mm{a}})+\mathcal{N}^3 \right)\cdot u_t^{\mm{d}}
- \kappa\rho \mathbb{D} u^{\mm{d}} :\nabla_{\tilde{\mathcal{A}}}u_t^{\mm{d}}
\right)
\mm{d}y\nonumber  \\
&\quad
+\int_{\Sigma}  \left(\llbracket
\mu\mathbb{D}_{\tilde{\mathcal{A}}} u^{\mm{a}}_t e_3-\mathfrak{S}_{\mathcal{A}} ( q^{\mm{a}}_t,u^{\mm{a}}_t,u^{\mm{a}} )\tilde{\mathcal{A}}e_3\rrbracket+\mathcal{N}^{4} \right)\cdot  u_t^{\mm{d}}\mathrm{d} y_{\mm{h}} \nonumber \\
& \quad
- \int (\nabla q^{\mm{d}}  \cdot \partial_t(\mathcal{A}_t^{\mm{T}}u)+ \nabla q \cdot\partial_t( \tilde{\mathcal{A}}^{\mm{T}}u^{\mm{a}}_t))\mm{d}y
\nonumber \\
&\quad -\int_\Sigma \left( \llbracket q^{\mm{d}}  \rrbracket \partial_t(\mathcal{A}_t e_3\cdot u)  + \llbracket q \rrbracket \partial_t( \tilde{\mathcal{A}} e_3 \cdot u^{\mm{a}}_t)\right)\mm{d}y_{\mm{h}}=:R_1(t). \nonumber
\end{align}
Integrating the above identity in time from $0$ to $t$ yields that
\begin{align}
& \|\sqrt{\rho} u_t^\mm{d}\|^2_{0}+E(u^{\mm{d}})+ \int_0^t\|\sqrt{\mu} \mathbb{D}_{\mathcal{A} } u^{\mathrm{d}}_t\|_0^2\mm{d}\tau\nonumber \\
& =2 \left(\int_0^t {R}_1(\tau)\mm{d}\tau  +R_2(t)-R_2(0)\right)+ {E}(u^{\mm{d}}|_{t=0}+ \| \sqrt{\rho}u_t^{\mathrm{d}}\|^2_0\big|_{t=0})
, \label{0314}
\end{align}
where we have defined that
$$R_2(t):= \int\left( \nabla q^{\mm{d}} \cdot (\mathcal{A}_t^{\mm{T}}u)+\nabla q \cdot( \tilde{\mathcal{A}}^{\mm{T}}u^{\mm{a}}_t)\right)\mm{d}y+\int_\Sigma \left(\llbracket
q^{\mm{d}} \rrbracket \mathcal{A}_te_3\cdot u+\llbracket q \rrbracket \tilde{\mathcal{A}} e_3\cdot  u^{\mm{a}}_t \right)\mm{d}y_{\mm{h}} .$$
Here we have formally derived \eqref{0314} for simplicity. Of course, \eqref{0314} can be rigorously verified by a complicated regularized method.

We call $\delta^3 e^{3\Lambda t}$ the three-orders term in general. Next we see that  the integrals in the right hand of the above identity can be controlled by the three-orders term.

(2) Estimate for three-orders small terms.

Noting that $(\eta,u)$ satisfies the estimates \eqref{aimdse}--\eqref{06041533fwqg} and \eqref{201806291},
thus, by further exploiting  \eqref{appesimtsofu1857}, \eqref{201702100940}, the product estimates \eqref{fgestims} and \eqref{06041605xdsf}, and trace estimate
  we have
\begin{align}
  {R}_1   \lesssim &
 \|  u_t^{\mm{d}}\|_1\left(
\|\tilde{\mathcal{A}}\|_2(\|(u^{\mm{a}},u_t^{\mm{a}})\|_{2}+\|u^{\mm{d}}\|_1
+\| \nabla q_t^{\mm{a}} \|_0+|  \llbracket q_t^{\mm{a}} \rrbracket |_{1/2})+\|\mathcal{N}^3\|_0+ |\mathcal{N}^4|_{1/2}
\right)\nonumber \\
&+ \|\nabla q^{\mm{d}}\|_0\left(\| u \|_2 \|\mathcal{A}_{tt}\|_0+\|\mathcal{A}_t\|_1 \|u_t\|_1 \right)
+\|\nabla q\|_0\left(\|\mathcal{A}_t\|_1\|u^{\mm{a}}_t\|_1+
\|\tilde{\mathcal{A}}\|_2\|u_{tt}^{\mm{a}}\|_0\right)\nonumber \\
&+|\llbracket q^{\mm{d}} \rrbracket |_{1/2}(\|u\|_2|\mathcal{A}_{tt}e_3|_{-1/2} +
\|\mathcal{A}_t\|_2\|u_t\|_1)+|\llbracket q  \rrbracket |_{1/2}(\|\mathcal{A}_{t}\|_{1} \|u_t^{\mm{a}}\|_2+
\|\tilde{\mathcal{A}}\|_2\|u_{tt}^{\mm{a}}\|_1)  \nonumber \\
\lesssim &\delta e^{\Lambda t}( \delta^2 e^{2\Lambda t}+
 \| u \|_3^2+\|\nabla q \|_{1}^2+\|u_t\|_1^2). \label{201811081020}
\end{align}
Thus, by \eqref{201702100940}, we get
\begin{align}
\label{201808141005}
\int_0^t R_1(\tau)\mm{d}\tau\lesssim \delta^3 e^{3\Lambda \tau}.
\end{align}

Similarly, we easily estimate
\begin{align}
  {R}_2(t) \lesssim  \delta^3 e^{3\Lambda t},
\end{align}
which also yields that
\begin{align}
  {R}_2(0) \lesssim  \delta^3 e^{3\Lambda t}.
\end{align}
Noting that $u^\mm{d}(0)=\delta^2 u^{\mm{r}}$, we have
\begin{equation}
E(u^{\mm{d}}|_{t=0})\lesssim \delta^4\|u^{\mm{r}}\|_2^2\lesssim  \delta^3 e^{3\Lambda t}.
\label{2018110702038}
\end{equation}

Next we turn to estimate for $\| \sqrt{\rho}u_t^{\mathrm{d}}\|^2_0\big|_{t=0} $.
Recalling that
$$\div  u_t^{\mm{d}}=  -\mm{div}\partial_t( \tilde{\mathcal{A}}^{\mm{T }}  u),$$
thus, taking the inner product of \eqref{201702052209}$_2$ and $u_t^{\mm{d}}$ in $L^2$, and using the integration by parts,  we have
\begin{align}
\int \rho |u_{t}^{\mathrm{d}}|^2\mm{d}y=
&  \int( \partial_{\bar{M}}^2  \eta^{\mathrm{d} }
-\mathrm{div}\mathcal{S} ( 0,u^{\mm{d}},\eta^{\mm{d}})- \mathrm{div}_{\tilde{\mathcal{A}}}\mathcal{S}_{\mathcal{A}} (q ,u ,\eta)
-\mu \mathrm{div}\mathbb{D}_{\tilde{\mathcal{A}}}  u)\cdot u_t^{\mm{d}}   \mm{d}y\nonumber
 \\
& +
\int  \nabla  q^{\mathrm{d}}\partial_t({\tilde{\mathcal{A}}}^{\mm{T}} u)  \mm{d}y+ \int_\Sigma \left(\llbracket q^{\mm{d}}\rrbracket \partial_t(\mathcal{A}^{\mm{T}}u  )\cdot e_3 -  \llbracket q^{\mm{d}}\rrbracket  \partial_tu_3^{\mm{a}}  \right)\mm{d}y_{\mm{h}} .\label{201808141332}
  \end{align}
Noting that
$$
\begin{aligned}
&\llbracket q^{\mm{d}}\rrbracket=   \llbracket 2(\mu \partial_3u_3^{\mm{d}}+{\kappa}\rho \partial_3\eta_3^{\mm{d}})+ g \rho  \eta_3^{\mm{d}} +  \bar{M}_3 \partial_{\bar{M}}\eta_3^{\mm{d}}\rrbracket +\vartheta\Delta_{\mm{h}}\eta_3^{\mm{d}}+ \mathfrak{M},\\
&  \sqrt{\|\eta\|_3^2+
\| u \|_2^2}\leqslant \theta\mbox{ and }(\eta,u)\mbox{ satisfies the estimate }\eqref{06011734} ,
\end{aligned}$$
thus, following the arguments of \eqref{201807311651} and \eqref{201811081020}, we easily derive from \eqref{201808141332} that
$$
  \|  u_t^{\mathrm{d}}\|_0^2
\lesssim\| \eta^{\mm{d}}\|_3^2+\|u^{\mm{d}}\|_2^2+ \delta e^{\Lambda t}(\delta^2 e^{2\Lambda t}+\| \eta^{\mm{d}}\|_3+\|u^{\mm{d}}\|_2),
$$
which, together with the initial data \eqref{201702052209}$_6$, implies that
\begin{equation}
\label{201807112050}
\|  u_t^{\mathrm{d}}\|_0^2 |_{t=0}\lesssim  \delta^3 e^{3\Lambda t} .
\end{equation}

Consequently, putting \eqref{201808141005}--\eqref{2018110702038} and \eqref{201807112050} into \eqref{0314}  yields
\begin{align}
 \|\sqrt{\rho} u_t^\mm{d}\|^2_{0}+ \int_0^t\|\sqrt{\mu} \mathbb{D}_{\mathcal{A} } u^{\mathrm{d}}_{\tau}\|_0^2\mm{d}\tau\leqslant - {E}( u^{\mm{d}})+  c\delta^3 e^{3\Lambda t}.\label{0314xxdfafdss}
\end{align}

(3) Application of the largest growth rate $\Lambda$.

Since $u^{\mm{d}}$ does not satisfies the divergence-free condition (i.e.,
 $\mm{div}u^{\mm{d}}=0$), thus we can not use \eqref{Lambdard} to deal with $ -E(u^{\mm{d}})$.
To overcome this trouble, we consider the following stratified   Stokes problem
\begin{equation*}
 \left\{\begin{array}{ll}
\nabla \varpi-\Delta \tilde{u}=0,\quad
 \mm{div}\tilde{u} =-  \div_{\tilde{\mathcal{A}}}
  u  &\mbox{ in }\Omega,  \\
\llbracket (\varpi I-\mu\mathbb{D}\tilde{u})e_3\rrbracket=0&\mbox{ on } \Sigma,\\
\tilde{u}   =0&\mbox{ on } \Sigma_{-}^+.
\end{array}\right.\end{equation*}
Then, by the existence  theory of  stratified   Stokes problem, we have
a solution $(\tilde{u},\varpi)\in H_0^2\times  {H}^1$. Moreover,
\begin{equation}
\label{201705012219xx}
\|\tilde{u}\|_2\lesssim \|\mm{div}_{\tilde{\mathcal{A}}}u\|_1 \lesssim \delta^2 e^{2\Lambda t}.
\end{equation}
It is easy to check that $v^{\mm{d}}:=u^{\mm{d}}-
\tilde{u}\in H_\sigma^2$.

Now, we can apply \eqref{Lambdard} to $ -E(v^{\mathrm{d}})$, and get
  \begin{equation*}\begin{aligned}
-E(v^{\mathrm{d}})\leqslant  \Lambda^2{\|v^{\mathrm{d}}\|^2_0}+ \Lambda\| \sqrt{\mu} \mathbb{D}  v^{\mathrm{d}}\|_0^2/2  ,
\end{aligned}\end{equation*}
which, together with \eqref{201702100940n} and \eqref{201705012219xx}, immediately implies
\begin{equation*}\begin{aligned}
-E(u^{\mathrm{d}})\leqslant  \Lambda^2{\|\sqrt{\rho}u^{\mathrm{d}}\|^2_0} + \Lambda\|\sqrt{\mu} \mathbb{D}  u^{\mathrm{d}}\|_0^2/2+c\delta^3 e^{3\Lambda t}.
\end{aligned}\end{equation*}
Inserting it into \eqref{0314xxdfafdss}, we arrive at
\begin{equation}
\label{201702232221}
\begin{aligned}  \| \sqrt{\rho}u_t^\mm{d}\|^2_{0}+\int_0^t\|\sqrt{\mu}\mathbb{D}_{\mathcal{A}}
u^{\mathrm{d}}_{\tau}\|_0^2\mm{d}\tau \leqslant {\Lambda^2} \|  \sqrt{\rho}u^\mm{d}\|_{0}^2 + \frac{\Lambda}{2}\|\sqrt{\mu}\mathbb{D}u^{\mathrm{d}}\|_0^2 +  c \delta^3 e^{3\Lambda t}.
\end{aligned}\end{equation}
In addition,
$$
\begin{aligned}
\int_0^t\|\sqrt{\mu}\mathbb{D}u_{\tau}^{\mathrm{d}}\|_0^2\mm{d}\tau
\leqslant &\int_0^t(\|\sqrt{\mu}\mathbb{D}_{\mathcal{A}}u_{\tau}^{\mathrm{d}}\|_0^2+
\|\tilde{\mathcal{A}}\|_2\|u_{\tau}^{\mathrm{d}}\|_1^2)\mm{d}\tau
\leqslant \int_0^t\|\sqrt{\mu}\mathbb{D}_{\mathcal{A}}u_{\tau}^{\mathrm{d}}\|_0^2\mm{d}\tau+  c \delta^3 e^{3\Lambda t}.
\end{aligned}$$
Thus we further deduce from \eqref{201702232221} that \begin{equation} \label{new0311}
\begin{aligned}  \| \sqrt{\rho} u_t^\mm{d}\|^2_{0}+ \int_0^t\|\sqrt{\mu} \mathbb{D}u^{\mathrm{d}}_{\tau}\|_0^2\mm{d}\tau \leqslant {\Lambda^2} \|\sqrt{\rho} u^\mm{d}\|_{0}^2 + \frac{\Lambda  }{2}\|\sqrt{\mu} \mathbb{D}u^{\mathrm{d}}\|_0^2+  c \delta^3 e^{3\Lambda t}.
\end{aligned}\end{equation}

(4) Derivation of \eqref{ereroe}--\eqref{201809231547}:

Recalling the initial data $u^{\mm{d}}(0)=\delta^2 u^{\mm{r}}$, we apply Newton--Leibniz's formula and Young's inequality to find that
 \begin{align}
 \Lambda  \|\sqrt{\mu} \mathbb{D}u^{\mathrm{d}}\|_0^2
 = &2\Lambda  \int_0^t \int  {\mu} \mathbb{D}^{\mm{d}}u:\mathbb{D}u^{\mathrm{d}}_\tau \mm{d}y\mathrm{d}\tau  +\delta^4\Lambda \|\sqrt{\mu} \mathbb{D}u^{\mm{r}}\|_0^2\nonumber \\
 \leqslant& \Lambda^2\int_0^t \|\sqrt{\mu} \mathbb{D}u^{\mathrm{d}}\|_0^2\mm{d}\tau+\int_0^t
\|\sqrt{\mu} \mathbb{D} u^{\mathrm{d}}_\tau\|_0^2\mathrm{d}\tau+c\delta^3 e^{3\Lambda t}. \nonumber
 \end{align}
Then we derive from \eqref{new0311} that
 \begin{equation}\label{inequalemee}\begin{aligned}
 \frac{1}{\Lambda}\|   u_t^{\mathrm{d}}\|^2_{0 }+
\frac{1}{2}\|\sqrt{\mu}\mathbb{D}u^{\mathrm{d}}\|_0^2  \leqslant   {\Lambda}\|  u^{\mathrm{d}}\|^2_{0 }+  {\Lambda}\int_0^t
\|\sqrt{\mu} \mathbb{D}u^{\mathrm{d}}\|_0^2\mm{d}\tau + c\delta^3 e^{3\Lambda t}.
\end{aligned}\end{equation}
In addition,
\begin{equation*}\begin{aligned}\label{}
\frac{\mm{d}}{\mm{d}t}\|u^\mm{d} \|^2_{0}=&2\int
u^\mm{d} \cdot  u^\mm{d}_t \mm{d}y
\leqslant \frac{1}{\Lambda}\| u_t^\mm{d} \|^2_{0}
+\Lambda\|u^\mm{d} \|^2_{0} .
\end{aligned}\end{equation*}
If we put the previous two estimates together, we get the differential inequality
\begin{equation}\label{growallsinequa}
\begin{aligned}
 \frac{\mm{d}}{\mm{d}t} \|  u^\mm{d}\|^2_{0}+\frac{1}{2}\|\sqrt{\mu} \mathbb{D}u^{\mathrm{d}} \|_0^2\leqslant 2\Lambda\left( \|  u^\mm{d}(t)\|^2_{0}
 +\frac{1 }{2}\int_0^t\|\sqrt{\mu} \mathbb{D}(u^{\mathrm{d}})\|_0^2\mathrm{d}\tau\right) +c\delta^3e^{3\Lambda t}.
\end{aligned}
\end{equation}
Recalling $u^{\mm{d}}(0)=\delta^2u^{\mm{r}}$, one can apply Gronwall's inequality to \eqref{growallsinequa} to conclude that
 \begin{equation}\label{estimerrvelcoity}
\begin{aligned}
\| u^{\mathrm{d}}\|^2_{0}+ \frac{1 }{2} \int_0^t\|\sqrt{\mu} \mathbb{D}u^{\mathrm{d}}\|_0^2\mm{d}\tau \lesssim   e^{2\Lambda t}\left(\int_0^t   \delta^3 e^{3\Lambda t}e^{-2\Lambda\tau}\mm{d}\tau
+\delta^4\|  u^{\mathrm{r}}\|^2_{0}\right) \lesssim \delta^3e^{3\Lambda t}.
 \end{aligned}  \end{equation}
Moreover,  we can further deduce from \eqref{new0311}, \eqref{inequalemee},  \eqref{estimerrvelcoity} and Korn's inequality that
\begin{eqnarray}\label{uestimate1n}
\|u^{\mathrm{d}}\|_{1 }^2+\| u_t^{\mathrm{d}}\|^2_0 +\|u^{\mathrm{d}}_\tau\|^2_{L^2((0,t),H^1)}\lesssim \delta^3e^{3\Lambda t}.
\end{eqnarray}

We turn to derive the error estimate for $\eta^{\mathrm{d}}$.
It follows from \eqref{201702052209}$_1$ that
\begin{equation*}\begin{aligned}
 \frac{\mm{d}}{\mm{d}t}\|\eta^{\mathrm{d}}\|_{1}^2
\lesssim   \| u^{\mathrm{d}}\|_1 \|\eta^{\mathrm{d}}\|_{1}.
\end{aligned}\end{equation*}
Therefore, using \eqref{uestimate1n} and  the condition $\eta^{\mm{d}}|_{t=0}=\delta^2 \eta^{\mm{r}}$,  it follows that
\begin{equation}\begin{aligned}\label{erroresimts}
 \|\eta^{\mathrm{d}}\|_{1 }\lesssim  &  \int_0^t \|u^{\mathrm{d}}\|_{1} \mm{d}\tau
+\delta^2\|\eta^{\mm{r}}\|_1\lesssim  \sqrt{ \delta^3e^{3\Lambda t}}.
\end{aligned}\end{equation}
Noting that $L^2(\mathbb{T})\hookrightarrow L^1(\mathbb{T})$ and $H^1\hookrightarrow W^{1,1}$, then we can derive
 \eqref{ereroe} and \eqref{201805291548} from \eqref{uestimate1n},  \eqref{erroresimts} and trace estimate.

Using \eqref{201702100940} and \eqref{ereroe}, we can estimate
\begin{align}
 \|\mathcal{A}_{3k}\partial_k \chi_{3} - \partial_3\chi^{\mm{a}}_3\|_{L^1}  \leqslant  \|(\mathcal{A}_{3k}-\delta_{3k}) \partial_k\chi_3 + \partial_3 \chi^{\mm{d}}_3\|_{L^1} \lesssim \sqrt{\delta^3 e^{3\Lambda t}},\nonumber
\end{align}
which yields \eqref{2018090119021}.
Similarly, we also can derive  \eqref{201809011902} and \eqref{201809231547}.
This completes the proof of Proposition \ref{lem:0401xx}. \hfill$\Box$
\end{pf}

\section{Existence of escape times}\label{sec:030845}

Let $\delta<\delta_0$, $(\eta^{\mm{d}},u^{\mm{d}})$ be defined  in Section \ref{201807111049},  and $(\eta,u)$ be the strong  solution constructed in Section \ref{201807111049} with an existence time $[0,T^{\max})$.
Let $\epsilon_0\in (0,1)$ be a constant, which will be defined in \eqref{defined}.
 We define
 \begin{align}\label{times}
& T^\delta:={\Lambda}^{-1}\mm{ln}({\epsilon_0}/{\delta})>0,\quad\mbox{i.e.,}\;
 \delta e^{\Lambda T^\delta }=\epsilon_0,\\
&T^*:=\sup\left\{t\in(0,T^{\max})\left|~\sqrt{{\|\eta (\tau)\|_3^2+
\| u (\tau)\|_2^2+\beta(|d(\tau)|_3^2)}}\right.\right.\nonumber\\
&\qquad \qquad\quad \leqslant 2C_4\delta_0\mbox{ for any }\tau\in [0,t)\bigg\},\nonumber\\
&  T^{**}:=\sup\left\{t\in (0,T^{\max}) \left|~\left\|(\eta,u)(\tau)\right\|_{0}\leqslant 2 C_4\delta e^{\Lambda \tau}\mbox{ for any }\tau\in [0,t)\right\}.\right.\nonumber\\
&  T^{***}:=\sup\left\{t\in (0,T^{\max}) \left|~\int_0^t \|\eta\|^2_3 \mm{d}\tau\leqslant \Lambda / 2C_1\mbox{ for any }\tau\in [0,t)\right\}.\right.\nonumber
 \end{align}
 Noting that
\begin{align}
&\left.\sqrt{\|\eta(t)\|_3^2+\|u(t)\|_2^2+\beta( |d(t)|_3^2)}\right|_{t=0}\nonumber \\
&=\sqrt{\|\eta_0^\delta\|_3^2
+\|u_0^\delta\|_2^2+\beta( |d^\delta_0|_3^2)}\leqslant C_4\delta<2C_4\delta_0\leqslant \delta_2\mbox{ or }\delta_3,
\label{201809121553}
\end{align}
 thus $T^*>0$ by Proposition \ref{pro:0401nxd}.  Similarly, we also have  $T^{**}>0$.
Moreover, by Proposition \ref{pro:0401nxd} and Remark \ref{201811031411},
we can easily see that
 \begin{align}\label{0502n1}
&\sqrt{\|\eta (T^*)\|_3^2+
\| u (T^*)\|_2^2+\beta( |d(T^*)|_3^2)}=2C_4 \delta_0,\mbox{ if }T^*<\infty ,\\
\label{0502n111}  & \left\|(\eta,u) (T^{**}) \right\|_0
=2 C_4\delta e^{\Lambda T^{**}},\mbox{ if }T^{**}<T^{\max},\\
\label{0502n111asdfsafsa}  & \int_0^{T^{***}} \|\eta\|^2_3 \mm{d}\tau= \Lambda / 2C_1 ,\mbox{ if }T^{***}<T^{\max}.
\end{align}

We denote ${T}^{\min}:= \min\{T^\delta ,T^*,T^{**},T^{***}\}$.
 Noting that $(\eta,u)$ satisfies
$$
 \sup_{0\leqslant t< T^{\min}}\sqrt{\|\eta(t)\|_3^2 +\|u(t)\|_2^2+
\beta(|d(t)|_3^2)}\leqslant \delta_1,$$
thus,  for any $t\in (0,T^{\min})$, $(\eta,u)$ enjoys  \eqref{2016121521430850} and \eqref{2018008121027} with $(\eta^\delta_0,u^\delta_0)$ in place of $(\eta^0,u^0)$ by the third conclusion in Proposition \ref{pro:0401nxd}.
Noting that $\|(\eta,u)(t)\|_0\leqslant 2C_4 \delta e^{\Lambda t}$ on $(0,T^{\min})$, then we deduce from the estimate \eqref{2016121521430850} and \eqref{201809121553}  that, for all $t\in
(0, {T}^{\min})$,
\begin{align}
& \mathfrak{E}(t)+\beta( c\mathfrak{E}_L(t) +|d(t)|_3^2) +\int_0^t\mathfrak{O}(\tau)\mm{d}\tau\nonumber
\\
&\leqslant c_0\delta^2e^{2\Lambda t}+\Lambda\int_0^t\left(\mathfrak{E}(\tau)+\beta( c\mathfrak{E}_L(\tau)+ |d(\tau)|_3^2)
+\int_0^{\tau}\mathfrak{O}(s)\mm{d}s\right)\mm{d}\tau
 \label{0503} \end{align}
for some positive constant $c_0$.
   Applying Gronwall's inequality to the above estimate, we deduce that
$$
\mathfrak{E}(t)+\beta(  |d(t)|_3^2 )
 \lesssim \delta^2 \left(e^{2\Lambda t}+
 \Lambda\int_0^t e^{\Lambda (t+\tau)} \mm{d}\tau\right)\lesssim\delta^2 e^{2\Lambda t}.
 $$
 Putting the above estimate to  \eqref{0503}, we get
$$ \mathfrak{E}(t)+\beta( |d(t)|_3^2) +\int_0^t\mathfrak{O}(\tau)\mm{d}\tau \lesssim \delta^2 e^{2\Lambda t},$$
which, together with \eqref{2018008121027}, yields that, for some constant $C_6$,
\begin{align}
& \mathcal{E}(t)+\beta( |d(t)|_3^2 )+ \int_0^t(\|\nabla q_\tau\|_1^2+|  \llbracket  q_\tau \rrbracket |_{1/2}^2 )\mm{d}\tau
\leqslant (C_6\delta e^{\Lambda t})^2\leqslant C_6^2\epsilon_0^2,
\label{201702092114}\\
& \int_0^t \|\eta\|^2_3 \mm{d}\tau \leqslant \left(\frac{C_6}{\sqrt{2\Lambda}}\delta e^{\Lambda t}\right)^2 \leqslant C_6^2\epsilon_0^2/2\Lambda \mbox{ on }(0, {T}^{\min}).
\label{201811061630}
\end{align}

If $\bar{M}_3=0$, we define
 $$
\begin{aligned}
m_0:=& \min_{\tilde{\chi}^0=\tilde{\eta}^0, \tilde{u}^0} \{|\tilde{\chi}^0_3|_{L^1},\|\tilde{\chi}^0_{\mm{h}}\|_{L^1},\|\tilde{\chi}^0_3\|_{L^1}, \|\partial_3\tilde{\chi}^0_{\mm{h}}\|_{L^1},\|\partial_3\tilde{\chi}^0_3\|_{L^1},
 \|\mm{div}_{\mm{h}}\tilde{\chi}^0_{\mm{h}}\|_{L^1} \};
\end{aligned}$$
else $$
\begin{aligned}
m_0:=& \min_{\tilde{\chi}^0=\tilde{\eta}^0, \tilde{u}^0} \{|\tilde{\chi}^0_3|_{L^1},\|\tilde{\chi}^0_{\mm{h}}\|_{L^1},\|\tilde{\chi}^0_3\|_{L^1}, \|\partial_3\tilde{\chi}^0_{\mm{h}}\|_{L^1}, \|\partial_3\tilde{\chi}^0_3\|_{L^1},
 \|\mm{div}_{\mm{h}}\tilde{\chi}^0_{\mm{h}}\|_{L^1}, \\
&\qquad \qquad  \|\partial_{\bar{M}}\tilde{\chi}^0_{\mm{h}}\|_{L^1}
, \|\partial_{\bar{M}}\tilde{\chi}^0_3\|_{L^1}\} .
\end{aligned}$$
By \eqref{n05022052} and \eqref{n05022052xx}, $m_0>0$.
Now we define that
 \begin{equation}\label{defined}
\epsilon_0:=\min\left\{\frac{\theta}{C_6},\frac{C_4\delta_0}{C_6},
\frac{ \Lambda}{\sqrt{2 C_1} C_6 },\frac{C_4^2}{4C_5^2},
\frac{\delta_5}{C_6},
\frac{m_0^2}{4C_5^2}  \right\}>0,
 \end{equation}
 where $\theta$  and $C_5$ come  from Proposition \ref{lem:0401xx} with $\gamma=C_6$ and $\beta=\epsilon_0$.
Noting that $(\eta,u)$ satisfies \eqref{201702092114},  $\epsilon_0\leqslant \theta/{C_6}$ and $\beta=\epsilon_0$, then, by Proposition \ref{lem:0401xx} with $\gamma=C_6$, we immediately have \eqref{ereroe}--\eqref{201805291548} for any $t\in (0,T^{\min})$.
Consequently, we further have the  relation
\begin{equation}
\label{201702092227}
T^\delta =T^{\min}\neq T^*\mbox{ or }T^{**},
\end{equation}
which can be showed by contradiction as follows:

If $T^{\min}=T^*$, then   $T^*<\infty$. Noting that $\epsilon_0\leqslant C_4\delta_0/C_6 $, thus we deduce from \eqref{201702092114} that
\begin{equation*}\begin{aligned}
\sqrt{\|\eta (T^*)\|_3^2+\| u (T^*)\|_2^2+\beta( | d (T^*)|_3^2)}\leqslant C_4\delta_0< 2 C_4 \delta_0,
 \end{aligned} \end{equation*}
which contradicts \eqref{0502n1}. Hence, $T^{\min}\neq T^*$.

If $T^{\min}=T^{**}$, then   $T^{**}<T^*\leqslant T^{\max}$.
 Noting that $\sqrt{\epsilon_0}\leqslant C_4/2C_5$, we can deduce from  \eqref{0501},
 \eqref{ereroe}, \eqref{times} and the fact $\varepsilon_0\leqslant c_3^2/4C_5^2$  that \begin{equation*}\begin{aligned}
\| (\eta,u)  (T^{**})\|_{0}
&\leqslant \|(\eta^\mm{a},u^\mm{a})(T^{**})\|_{0}+\|(\eta^\mm{d}, u^\mm{d}) (T^{**})\|_{0}\leqslant  \delta e^{{\Lambda T^{**}}}(C_4+ C_5\sqrt{\delta e^{\Lambda T^{**}}})\\
&\leqslant  \delta e^{{\Lambda T^{**}}}(C_4+ C_5\sqrt{\epsilon_0})
\leqslant 3C_4 \delta e^{\Lambda T^{**}}/2<2C_4\delta e^{\Lambda T^{**}},
 \end{aligned} \end{equation*}
which contradicts \eqref{0502n111}. Hence, $T^{\min}\neq T^{**}$.

If $T^{\min}=T^{***}$, then   $T^{***}< T^{\max}$. Noting that $\epsilon_0\leqslant \Lambda^2/2 C_1 C_6^2 $, thus we deduce from \eqref{201811061630} that
\begin{equation*}\begin{aligned}
 \int_0^{T^{***}} \|\eta\|^2_3 \mm{d}\tau \leqslant   \Lambda / 4C_1< \Lambda / 2C_1,
 \end{aligned} \end{equation*}
which contradicts \eqref{0502n111asdfsafsa}. Hence, $T^{\min}\neq T^{***}$.
We immediately see that \eqref{201702092227} holds. This completes the verification of \eqref{201702092227}.

By the relation \eqref{201702092227} and the definition of $T^{\min}$, we see that $T^\delta<T^*\leqslant T^{\max}$. By \eqref{201702092114} and the fact $\epsilon_0\leqslant \delta_5/C_6$, we see that
$\sqrt{\|\eta(t)\|_3^2+\|u(t)\|_2^2} \leqslant \delta_5$. Then, by Proposition \ref{201809012320}, we see that $\eta \in H^{3,1}_{0,*}$. Hence $(\eta,u)\in C^0(\overline{I_T},H^{3,1}_{0,*}\times  H^2_0)$ for some $T\in (T^\delta,T^{\max})$,  and $(\eta,u,q)$ also enjoys the regularity mentioned in Remark \ref{201809121705}.

Noting that $\sqrt{\epsilon_0}\leqslant m_0/2C_5$ and \eqref{ereroe} holds for $t=T^\delta$, thus, using of \eqref{0501}, \eqref{ereroe}, \eqref{2018090119021}  and \eqref{times}, we can easily deduce the following instability relations:
$$
 \begin{aligned}
|\chi_3(T^\delta)|_{L^1}\geqslant & |\chi_3^{\mm{a}}(T^\delta)|_{L^1}^2-
|\chi_3^{\mm{d}}(T^\delta)|_{L^1}^2 \\
> &  \delta e^{\Lambda T^\delta }( |\tilde{\chi}^{0}_3|_{L^1}- C_5\sqrt{\delta e^{\Lambda T^\delta }})
\geqslant  (m_0 -C_5\sqrt{\epsilon_0})\epsilon_0 \geqslant m_0 \epsilon_0 /2,\\
\end{aligned}
$$
and
$$
 \begin{aligned}
\|\mathcal{A}_{3k}\partial_k\chi_3(T^\delta)
\|_{{L^1}}\geqslant &
\| \partial_3 \chi^{\mm{a}}_3(T^\delta) \|_{{L^1}}
-
\|\mathcal{A}_{3k}\partial_k\chi_3(T^\delta)
-\partial_3 \chi^{\mm{a}}_3(T^\delta)
\|_{{L^1}} \\
>&   \delta e^{\Lambda T^\delta }( \|\partial_3\tilde{\chi}^{0}_3\|_{L^1}- C_5\sqrt{\delta e^{\Lambda T^\delta }}) \geqslant m_0 \epsilon_0 /2,\\
\end{aligned}
$$
where $\chi=\eta$ or $u$. Similarly, we also  easily  verify that $(\eta,u)$ satisfies the rest instability relations in \eqref{201806012326} and \eqref{2018109121834}.
This completes the proof of  Theorem \ref{thm3} by taking $\epsilon:= m_0\epsilon_0 /8$.


\section{Proof of Theorem  \ref{201806012301xx}}\label{201806012303}

This section is devoted to the verification of Theorem  \ref{201806012301xx}.
Since $\eta\in H^{3,1}_{0,*}$, the definitions of $d$, $v$, $V$, $N$ and $\sigma$ in \eqref{usedopso} make sense.
Moreover, by Remark \ref{201809273138}, $v$, $V$, $N$, $\sigma$ and $d$   still are horizontally periodic functions.  Next we start to verify that $v$, $V$, $N$, $\sigma$ and $d$  satisfy the conclusions in Theorem \ref{201806012301xx} and \eqref{201808151856}--\eqref{201808151856sdfa} in Remark \ref{rem:201810011518}.
In what follows, we denote $\det \nabla_{y_{\mm{h} }} \zeta_{\mm{h}}(y_{\mm{h}},0,t)$ by $\Theta(y_{\mm{h}},t)$ for the simplicity.

\subsection{Homeomorphism}\label{201811191634}

Let $\zeta=\eta+y$, $\tilde{y}:=(y,t)$, $\tilde{y}_{\mm{h}} :=(y_{\mm{h}},t)$, $\tilde{x}=(x,t)$ and $\tilde{x}_{\mm{h}}=(x_{\mm{h}},t)$. Since $\eta(t)\in H^{3,1}_{0,*}$ and $\eta\in C(\overline{I_T},H^3)$, then we have
\begin{align}
\label{201806031606}
&\nabla_x \zeta^{-1} =(\nabla_y \zeta)^{-1}|_{y=\zeta^{-1} }=\mathcal{A}^{\mm{T}}|_{y=\zeta^{-1} }\mbox{ in }\Omega,\\
& \tilde{\zeta}: \overline{\Omega^{T}} \mapsto  \overline{\Omega^{T}} \mbox{ is bijective mapping}, \label{201809300836} \\
 & \tilde{\zeta}_\pm: X \mapsto \tilde{\zeta}_\pm(X)\mbox{ are bijective mappings}, \label{201809300834}\\
& \det \nabla_{\tilde{y}}\tilde{\zeta}(\tilde{y})=1\mbox{ in }\overline{\Omega_\pm^T}, \label{201810011535} \\
&  \bar{\zeta}(\tilde{y}_{\mm{h}}): Y \to  Y\mbox{ is a bijective mapping},\label{2018081518dsfgs5sdfa6}\\
&\det \nabla_{\tilde{y}_\mm{h}}\bar{\zeta}(\tilde{y}_{\mm{h}})=\Theta(y_{\mm{h}},t)\geqslant 1/2, \label{2018081518dsfgs5sdfa6xx}
 \end{align}
 where $X:=\overline{\Omega_\pm^{T}}$ or $\Omega_\pm^{T}$, and $Y:=\overline{\mathbb{R}^2_{T}} $ or ${\mathbb{R}^2_{T}}$.

We denote the inverse functions of $\tilde{\zeta}_\pm(\tilde{y})$, resp. $\bar{\zeta}(\tilde{y}_{\mm{h}})$ by $\tilde{\zeta}_{\pm}^{-1}(\tilde{x})$, resp. $\bar{\zeta}^{-1}(\tilde{x}_{\mm{h}})$.
Recalling the regularity
\begin{align}
&  \zeta\in C^0(\overline{I_T},H^3) \mbox{ and }\zeta_t=u\in C^0(\overline{I_T},H^2),\label{201810011540x}
\end{align}
thus, using embedding theorem $H^{k+2}(\Omega_\pm)\hookrightarrow C^{k}(\overline{\Omega}_\pm)$ for $k\geqslant 0$, we can get
\begin{align}
&\label{201809271654}
\tilde{\zeta}_\pm \in C^1(\overline{\Omega_\pm^T}),\\
& \label{201809271654xx}
\tilde{\zeta} \in C^1(\overline{\mathbb{R}^T_2}).
 \end{align}
By \eqref{n0101nn}$_5$, we further get \begin{equation}
\label{20180927165xx4}
\tilde{\zeta},\ u, \nabla_{y_{\mm{h}}} {\zeta} \in C^0(\overline{\Omega ^T}).
 \end{equation}

Following the argument of \eqref{201804151826}, we derive
 from \eqref{201809300836} and the continuity of $\tilde{\zeta}$ in \eqref{20180927165xx4}
 that
\begin{equation}
\label{201810012335}
\tilde{\zeta} (\tilde{y}) :\overline{\Omega ^T}\to  \overline{\Omega ^T}\mbox{ is a homeomorphism mapping}.
\end{equation}
  Similar to \eqref{2018031adsadfa21601xx} and \eqref{201810012335}, by \eqref{201809300834}, \eqref{201810011535} and the continuity of $\tilde{\zeta}_\pm$ in \eqref{201809271654}, we have
 \begin{align}
& \tilde{\zeta}_\pm(\tilde{y}):\overline{ \Omega_\pm^T} \to \tilde{\zeta}_\pm( \overline{\Omega_\pm^T})\mbox{ is a  homeomorphism  mapping},\label{201808151856sdfaxxxxx}\\
& \tilde{\zeta}_\pm(\tilde{y}): \Omega_\pm^T \to \tilde{\zeta}_\pm( {\Omega_\pm^T})\mbox{ is a }C^1\mbox{-diffeomorphic mapping}.\label{201808151856sdfaxx}
\end{align}
Moreover, $\nabla_{\tilde{x}}\tilde{\zeta}^{-1}_\pm
=(\nabla_{\tilde{y}}\tilde{\zeta}_\pm)^{-1}|_{\tilde{y}=\tilde{\zeta}_\pm^{-1}}$. In particular,
\begin{equation}
\label{2018060321637}
\partial_t \zeta^{-1}_\pm=-((\nabla_y \zeta_\pm)^{-1}u_\pm)|_{y=\zeta^{-1}_\pm}.
\end{equation}

Similar to \eqref{201806031606}, \eqref{201808151856sdfaxxxxx} and \eqref{201808151856sdfaxx}, we also derive from \eqref{2018081518dsfgs5sdfa6}--\eqref{2018081518dsfgs5sdfa6xx}  and \eqref{201809271654xx}  that
\begin{align}
&\bar{\zeta}(y_{\mm{h}},t): \overline{\mathbb{R}^2_{T}} \to   \overline{\mathbb{R}^2_{T}}\mbox{ is a homeomorphism mapping},\label{2018081518dsfgs56}\\
&\bar{\zeta}(y_{\mm{h}},t):  {\mathbb{R}^2_{T}} \to  {\mathbb{R}^2_{T}}\mbox{ is a }C^1\mbox{-diffeomorphic mapping},\label{2018081518dsfgsafdasfs56} \\
&\label{201809291022}
\nabla_{x_\mm{h}} (\zeta_{\mm{h}})^{-1}(x_{\mm{h}},t)=(\nabla_{y_\mm{h}} \zeta_{\mm{h}}(y_{\mm{h}},0,t))^{-1}|_{y_{\mm{h}}=(\zeta_{\mm{h}})^{-1}(x_{\mm{h}},t)}.
\end{align}
Summing up \eqref{201810012335}--\eqref{201808151856sdfaxx}, \eqref{2018081518dsfgs56} and \eqref{2018081518dsfgsafdasfs56}, we get \eqref{201808151856} and \eqref{201808151856sdfa}.

\subsection{Regularity of solutions of the mixed RT problem}\label{201811142125}

Let $a:=\|\eta\|_{C^0(\overline{I_T},H^3)}$. In what follows, the notation $C(a)$ denotes a generic positive  constant, which may vary from line to line, and depends on $a$, $\Omega$, and increases with respect to $a$.

By transformation of Lagrangian coordinate (i.e., $x=\zeta(y)$), regularity of $(\zeta,u)$ and \eqref{2011801941224538}, we can estimate that, for any given $t\geqslant 0$,
\begin{align}
\label{201811210953}
\|v(t)\|_{L^2(\Omega^{\mm{p}}(t))}^2= &
\int_{\zeta^{-1}(\mathcal{T}_{-l}^{\tau} )/\{y_3=0\}}|
u|^2\mm{d}y \leqslant 2\|u(t)\|_0^2.
\end{align}
Noting that $\partial_{x_i} v= (\mathcal{A}_{il}\partial_{y_l} u )|_{y=\zeta^{-1}(x)}$ for $1\leqslant i\leqslant 3$, we have
\begin{align}
\label{201811210954}
\|\partial_{x_i}v(t)\|_{L^2(\Omega^{\mm{p}}(t))}^2= &
\int_{\zeta^{-1}(\mathcal{T}_{-l}^{\tau} )/\{y_3=0\}}|
\mathcal{A}_{il}\partial_{y_l} u |^2\mm{d}y \leqslant \|u(t)\|_1^2 .
\end{align}

Similarly, we can further derive that, for any $1\leqslant i$, $j$, $k\leqslant 3$,
$$
\begin{aligned}
& \|\partial_{x_ j}\partial_{x_i}v(t)\|_{L^2(\Omega^{\mm{p}}(t))}< \|u(t)\|_2 ,
\\
& \|\partial_{x_ k}\partial_{x_ j}\partial_{x_i}v(t)\|_{L^2(\Omega^{\mm{p}}(t))}<\|u(t)\|_3
\end{aligned}$$
from the relations $$
\begin{aligned}
&\partial_{x_j}\partial_{x_i}v=(\mathcal{A}_{jm} \partial_{y_m}(\mathcal{A}_{il}\partial_{y_l} u ))|_{y=\zeta^{-1}}\\
&
 \partial_{x_k}\partial_{x_j}\partial_{x_i}v=(\mathcal{A}_{kn} \partial_n (\mathcal{A}_{jm} \partial_{y_m}(\mathcal{A}_{il}\partial_{y_l} u )))|_{y=\zeta^{-1}}.
\end{aligned}$$
Hence
\begin{align}
& \|v(t)\|_{H^2(\Omega^{\mm{p}}(t))}\leqslant C(a) \|u\|_2\mbox{ for any  }t\in \overline{I_T},\label{2018092137}\\
&\|v(t)\|_{H^3(\Omega^{\mm{p}}(t))}\leqslant C(a) \|u\|_3\mbox{ for a.e. }t\in I_T.\label{2018092137zxx} \end{align}

 Thanks to  the continuity of $(\nabla_y\zeta_\pm, u)$ in \eqref{201809271654} and \eqref{20180927165xx4}--\eqref{201808151856sdfaxxxxx}, we get   \eqref{2018090714874}.
Then \eqref{0103nxx}$_7$ and \eqref{0103nxx}$_8$ obviously hold due to \eqref{n0101nn}$_6$, \eqref{n0101nn}$_7$ and the continuity of $u$ in \eqref{2018090714874}.
Thus, putting \eqref{2018092137}, \eqref{0103nxx}$_7$, \eqref{0103nxx}$_8$ and the fact ``$\mm{div}v=0$ in $\Omega(t)$" together, we get
\begin{equation}
v\in H_\sigma^2(\Omega(t))\mbox{ for any  }t\in \overline{I_T}.\label{201811201620}
\end{equation}

 Similar to \eqref{2018092137}, we can also get
\begin{align}
&\|(V, N )(t)\|_{H^2 (\Omega^{\mm{p}}(t))}\leqslant C(a) \|\eta(t)\|_3,\label{201811141558} \\
&\|\nabla\sigma(t)\|_{H^{k} (\Omega^{\mm{p}}(t))}\leqslant C(a) \|\nabla q(t)\|_{k},\label{201811141558xx12} \\
& \|\sigma(t)\|_{H^{k+1}( \Omega^{\mm{p}}(t))}\leqslant C(a) \|q(t)\|_{k+1},
\label{201811141558xx}  \end{align}
where $k=0$ and $1$.
Thus \eqref{201808160842} holds by \eqref{201811201620}--\eqref{201811141558xx} with $k=0$.

By the definitions of $v$, $V$ and $N$, we can exploit \eqref{n0101nn}$_1$, \eqref{06132125}, \eqref{201806031606}, \eqref{2018060321637}  and the chain rule of differentiation to derive that
\begin{align}
& U_t+v\cdot\nabla U=\nabla vU,\label{201811142017} \\
& M_t+ v\cdot\nabla M =M\cdot \nabla v,\label{201811142017cc} \\
& \mm{div}M=0,\label{201811142017er}\\
& \mm{div}(M\otimes M)={M}\cdot \nabla M =\partial_{\bar{M}}^2\eta|_{y=\zeta^{-1}((x,t)},\label{201811142017df}
\end{align}
where $U:=V+I=\nabla\zeta(y,t)|_{y={\zeta}^{-1}(x,t)}$ and $M:=N+\bar{M}= \partial_{\bar{M}}\zeta(y)|_{y={\zeta}^{-1}(x,t)}$.
Hence $(v, V,N)$ satisfies \eqref{0103nxx}$_2$ and \eqref{0103nxx}$_3$ by \eqref{201811142017} and \eqref{201811142017cc}.

Recalling the definitions of $V$, $\sigma$ and $\mathcal{A}$, and using \eqref{201806031606}, \eqref{201811142017df} and  chain rule of differentiation, we can easily get from \eqref{n0101nn}$_1$--\eqref{n0101nn}$_3$ that $(v,V,N,q) $ satisfies \eqref{0103nxx}$_1$ and $\mm{div}v=0$, which, together with \eqref{201811142017er}, yields   \eqref{0103nxx}$_4$. Using \eqref{0103nxx}$_1$--\eqref{0103nxx}$_3$, the regularity of $(v,V,N,\sigma)$ and product estimate \eqref{fgestims}, we easily further get \eqref{201819292017}.

For any given $t \geqslant 0$, $ \Xi(x ,t):=((\zeta_{\mm{h}})^{-1}(x_{\mm{h}},t),x_3):\Omega_+\to \Omega_+$ is a bijective mapping. We denote the inverse function of $\Xi(x ,t)$ by $\Xi^{-1}(y ,t)$. Then $\Xi$ an $\Xi^{-1}$ enjoy the following properties:
$$
\begin{aligned}
& \Xi^{-1}(y ,t) = (\zeta_{\mm{h}}(y_{\mm{h}},0,t),y_3 ),\   \Xi^{-1}(\mathcal{T}_{-l}^\tau ,t)\subset  \mathcal{T}_{-l,\tau}^{1/2},\\
& \nabla_x \Xi(x ,t) =
\left(\frac{1}{\Theta(y_{\mm{h}},t) }
\left.\left( \begin{array}{ccc}
  \partial_2 \zeta_2 (y_{\mm{h}},t)    &    -\partial_2 \zeta_1 (y_{\mm{h}},t)  & 0 \\
  -\partial_1 \zeta_2 (y_{\mm{h}},t)  &    \partial_1 \zeta_1 (y_{\mm{h}},t)  & 0 \\
 0 & 0 &\Theta(y_{\mm{h}},t)
                                                    \end{array}\right)\right)
\right|_{y_{\mm{h}}=(\zeta_{\mm{h}})^{-1}(x_{\mm{h}},t)}
,\\
&1/2<  \Theta(y_{\mm{h}},t)  =\det  \nabla_y  \Xi^{-1}(y ,t)<C(a),
\end{aligned}
 $$
thus, similarly to \eqref{2018092137zxx}, we easily check that
$$\|\eta_3^+(y,t)|_{y=\Xi(x ,t)}\|_{H^3(\Omega_+)}\leqslant C(a) \|\eta(t)\|_3,$$
which, together with trace estimate,  yields
\begin{equation}
\label{201809292130}
|d(x_\mm{h},t)=(\eta_3^+(y,t)|_{y=\Xi(x ,t)})|_{x_3=0}|_{5/2}\leqslant C(a) \|\eta(t)\|_3.
\end{equation}
Moreover, by \eqref{201809291022}, we can compute out
\begin{align}
 \nabla_{x_{\mm{h}}}d= &
(\nabla_{x_\mm{h}}(\zeta_{\mm{h}})^{-1}(x_{\mm{h}},t))^{\mm{T}}
\nabla_{y_\mm{h}}\zeta_3
(y_{\mm{h}},0,t)|_{y_{\mm{h}}=(\zeta_{\mm{h}})^{-1}(x_{\mm{h}},t)} \nonumber \\
=& ((\nabla_{y_\mm{h}} \zeta_{\mm{h}}(y_{\mm{h}},0,t))^{-{\mm{T}}}
\nabla_{y_\mm{h}}\zeta_3(y_{\mm{h}},0,t))|_{y_{\mm{h}}=
(\zeta_{\mm{h}})^{-1}(x_{\mm{h}},t)}   \label{201810012320}\end{align}
and
\begin{align}
\partial_{x_{i}} \nabla_{x_{\mm{h}}} d= & (\partial_{y_{\mm{j}}} (\nabla_{y_\mm{h}} \zeta_{\mm{h}}(y_{\mm{h}},0,t))^{-{\mm{T}}}
\nabla_{y_\mm{h}}\zeta_3(y_{\mm{h}},0,t)) |_{y_{\mm{h}}
=(\zeta_{\mm{h}})^{-1}(x_{\mm{h}},t)}
\partial_{x_i}((\zeta_{\mm{h}})^{-1}(x_{\mm{h}},t))_j  .  \label{201810012320xx}
\end{align}

Similar to \eqref{201809292130}, we also have
$$|( 1+|\nabla_{x_{\mm{h}}} d(x_{\mm{h}},t)|^2 )^{-1/2}|_{3/2}\leqslant C(a),$$
which, together with \eqref{201809292130} and \eqref{06011948}, yields
\begin{equation}
\nu(x_{\mm{h}},t) \in H^{3/2}(\mathbb{T})\mbox{  for any }t\geqslant 0.
\label{201811201627}
\end{equation}
Moreover, making use of the continuity of $(\zeta,\nabla_{y_{\mm{h}}} {\zeta})$ in \eqref{20180927165xx4}, \eqref{2018081518dsfgs56} and \eqref{201810012320},  we see that
$\nabla_{\mm{h}}d$ is continuous on $\overline{\mathbb{R}^2_T}$, which yields  \eqref{2018090414874}.

Since, for any given $t\geqslant 0$,
\begin{equation}
\label{201810011624}
 \zeta(\cdot,t): \{y_3=0\}\to \{x_3=d(x_{\mm{h}},t)\}\mbox{ is a bijective mapping},
\end{equation}
then, for any given $x_{\mm{h}}$, there exists $y_{\mm{h}}$ such that
\begin{equation}
\label{xx201809282115}
x_{\mm{h}}:=\zeta_{\mm{h}}(y_{\mm{h}},0,t)=\eta_{\mm{h}}(y_{\mm{h}},0,t)\mbox{ and }
d(x_{\mm{h}},t)=\eta_3(y_{\mm{h}},0,t).
\end{equation} We can further compute out that
\begin{align}
u_3|_{y_3=0}=&\partial_t \eta_3(y_{\mm{h}},0,t)=\partial_t d (\zeta_{\mm{h}}(y_{\mm{h}},0,t),t)\nonumber \\
=& u_1|_{y_3=0} \partial_1d|_{x_{\mm{h}}=\zeta_{\mm{h}}(y_{\mm{h}},0,t)}+u_2 |_{y_3=0} \partial_2d|_{x_{\mm{h}}=\zeta_{\mm{h}}(y_{\mm{h}},0,t)} +d_t|_{x_{\mm{h}}=\zeta_{\mm{h}}(y_{\mm{h}},0,t)},
\label{201810011601}
\end{align}
and \begin{equation}
\label{201810011617}
v(x_{\mm{h}},d(x_{\mm{h}},t),t)=v(\zeta_{\mm{h}}(y_{\mm{h}},0,t),\zeta_3(y_{\mm{h}},0,t),t)
=u(y_{\mm{h}},0,t)|_{y_{\mm{h}}=(\zeta_{\mm{h}})^{-1}(x_{\mm{h}},t)}.
\end{equation}
Putting $y_{\mm{h}}= (\zeta_{\mm{h}})^{-1}(x_{\mm{h}},t)$ into \eqref{201810011601},
and then using \eqref{201810011617}, we get \eqref{0103nxx}$_5$ for any given $t\geqslant 0$.

Similar to \eqref{201811210953}  and \eqref{201811210954},
exploiting \eqref{201809291022} and \eqref{esmmdforinfty}, we have,
\begin{equation}
\label{201811211159}
\left|u_+((\zeta_{\mm{h}})^{-1}(x_{\mm{h}},t),0,t)\right|_1\lesssim C(a)| u_+(t)|_{1}.
\end{equation}
Noting that $\zeta_{\mm{h}}(y_{\mm{h}},0,t)$ satisfies \eqref{201803121601xx2108} with $k=3$ for given $t$, by \eqref{201809291022},  \eqref{esmmdforinfty} and the fact
$\Theta(y_{\mm{h}},t)\geqslant 1/2$, we have
$$\sup_{y_{\mm{h}}^1,y_{\mm{h}}^2\in\mathbb{R}^2}\frac{|y_{\mm{h}}^1-y_{\mm{h}}^2 |}{|\zeta_{\mm{h}}(y_{\mm{h}}^1,0,t)-\zeta_{\mm{h}}(y_{\mm{h}}^2,0,t)|}
=\sup_{x_{\mm{h}}^1,x_{\mm{h}}^2\in\mathbb{R}^2}\frac{|
(\zeta_{\mm{h}})^{-1}(x_{\mm{h}}^1,t)-(\zeta_{\mm{h}})^{-1}(x_{\mm{h}}^2,t) |}{|x_{\mm{h}}^1 -x_{\mm{h}}^2 |}\leqslant  C(a).$$
Thus, making use of \eqref{201809291022}, \eqref{esmmdforinfty},  \eqref{06041605xdsf} and the above estimate,
we further obtain that
\begin{align}
& \left|\partial_{x_i} u_+(( \zeta_{\mm{h}})^{-1}(x_{\mm{h}},t),0,t)\right|_{1/2}\nonumber \\
&\lesssim  C(a)\left|
((\nabla_{y_\mm{h}} \zeta_{\mm{h}}(y_{\mm{h}},0,t))^{-1})_{ji}
\partial_{y_j} u_+
 \right|_{1/2} \lesssim C(a)|u_+|_{3/2}, \label{201811211200}
\end{align}
where $(A)_{ij}$ denotes the $(i,j)$-th entry of matrix $A$.

Noting that, by \eqref{201810011624},
$$v_+(x_{\mm{h}}, d(x_{\mm{h}},t), t) = u_+((\zeta_{\mm{h}})^{-1}(x_{\mm{h}},t),0,t),$$
thus, by \eqref{201811211159}, \eqref{201811211200} and trace estimate, we get
\begin{align}
&|v_+(x_{\mm{h}}, d(x_{\mm{h}},t), t)|_{3/2}\nonumber \\
&= |u_+((\zeta_{\mm{h}})^{-1}(x_{\mm{h}},t),0,t) |_{3/2}
\leqslant C(a)\| u\|_2.
\label{2018090232130}
\end{align}
Thus we can derive from
  \eqref{0103nxx}$_4$, \eqref{201809292130}, \eqref{2018090232130}  and \eqref{06011948} that
 \begin{equation}
d_t\in H^{3/2}(\mathbb{T}).
\label{201811201629}
\end{equation}
We get \eqref{201809271720} by \eqref{201809292130}, \eqref{201811201627} and \eqref{201811201629}.

Using  relations \eqref{2018081810135} and \eqref{201810012320},
we easily compute out that
$$
 \begin{aligned} \mathcal{A}e_3|_{y_3=0}= & (-\partial_1\eta_3+\partial_1 \eta_2\partial_2\eta_3-\partial_1\eta_3\partial_2\eta_2,
-\partial_2\eta_3+\partial_1 \eta_3\partial_2\eta_1-\partial_1\eta_1\partial_2\eta_3,
\det\nabla_{y_\mm{h}}\zeta_{\mm{h}})^{\mm{T}}|_{y_3=0} \\
=&\Theta(y_{\mm{h}},t) (- \partial_{x_1}d, - \partial_{x_2}d, 1  )^{\mm{T}}|_{x_{\mm{h}}=\zeta_{\mm{h}}(y_{\mm{h}},0,t)},
\end{aligned} $$
which yields that
\begin{align}
\label{201809282114}
&\nu|_{x_{\mm{h}}=\zeta_{\mm{h}}(y_{\mm{h}},0,t)}=\left.\frac{(-\partial_{x_1} d,-\partial_{x_2} d,1)^{\mm{T}}}{\sqrt{1+
|\nabla_{x_\mm{h}}d|^2}} \right|_{x_{\mm{h}}=\zeta_{\mm{h}}(y_{\mm{h}},0,t)} = \left. \frac{\mathcal{A}e_3}{|\mathcal{A}e_3|}\right|_{y_3=0}.
\end{align}
Using \eqref{201809291022}, \eqref{201810012320} and \eqref{201810012320xx}, we derive that
\begin{align}
\label{20180928212154}
\mathcal{C}|_{x_{\mm{h}}=\zeta_{\mm{h}}(y_{\mm{h}},0,t)}
=H|_{y_3=0}.
\end{align}
In addition, we also check that, for any given $t\geqslant 0$,
  \begin{align}
& \left(\left(\mathcal{S}^{\mm{V},\mm{M}}(\sigma, v,V +I,N+\bar{M})|_{\Omega_\pm(t)}\right)|_{x_3=d(x_{\mm{h}},t)}\nu\right)|_{x_{\mm{h}}=
\zeta_{\mm{h}}(y_{\mm{h}},0,t)}\nonumber \\
&=((\mathcal{S}_{\mathcal{A}}(q,u,\eta )|_{\Omega_\pm}\mathcal{A}e_3-\bar{M}_3\partial_{\bar{M}}\eta_\pm)/|\mathcal{A}e_3|)|_{y_3=0}。
\label{201809282115}
\end{align}
Consequently, by \eqref{201810011624}, \eqref{xx201809282115} and \eqref{201809282114}--\eqref{201809282115},  we get  \eqref{0103nxx}$_6$ from \eqref{n0101nn}$_4$  for any given $t$.

\subsection{Higher regularity of $d$ for $\vartheta\neq 0$}

Next we further derive \eqref{201811141303}.
Multiplying  \eqref{0103nxx}$_6$ by $\nu$ yields that
\begin{equation}
\label{201811131906}
\llbracket \mathcal{S}^{\mm{V},\mm{M}}(\sigma,v,V +I,N+\bar{M})\nu\cdot \nu -g \rho d  \rrbracket =\vartheta \mathcal{C}    \mbox{ on }\mathbb{T}.
\end{equation}
Since $|d|_{5/2}\lesssim \|\eta\|_3$,
 thus, similar to \cite[Lemma 3.2]{TATNLTE1995},  following the argument of \cite[Theorem 4]{SVAUZNSE}, we derive from \eqref{201811131906} that, for sufficiently small $\delta$,
\begin{align}
 |d|_{7/2}\lesssim & |\llbracket   \mathcal{S}^{\mm{V},\mm{M}}(\sigma, v,V+I,N+\bar{M})\nu- g d \rho  \nu        \rrbracket|_{3/2}+|d|_{3/2}\nonumber \\
= &\left|( \llbracket \mathcal{S}_{\mathcal{A}}(q,u,\eta ) \mathcal{A}e_3-\bar{M}_3\partial_{\bar{M}}\eta \rrbracket /|\mathcal{A}e_3|) |_{y_{\mm{h}}=(\zeta_{\mm{h}})^{-1}(x_{\mm{h}},t)}\right|_{3/2}+|d|_{3/2}.\label{201811211302xx}
\end{align}

Following the argument of \eqref{2018090232130},
and further  exploiting \eqref{201811231301}, we have
\begin{equation}
| (\mathbb{D}_{\mathcal{A}}u_\pm)|_{y_3=0}|_{(\zeta_{\mm{h}})^{-1}(x_{\mm{h}},t)} |_{3/2}\leqslant
C(a)|\nabla u|_{3/2} \leqslant C(a)\|\nabla u\|_2.
\label{201811211349}
\end{equation}
Noting that, for sufficiently small $\delta$,
$$\mathcal{A}e_3/|\mathcal{A}e_3|_{3/2}+\left|1/|\mathcal{A}e_3|\right|_{3/2}\leqslant C(a),$$
then, using \eqref{06011948}, we derive from the above two estimates that
\begin{align}
& \left|  \llbracket\mathbb{D}_{\mathcal{A}}u  \mathcal{A}e_3/|\mathcal{A}e_3| \rrbracket  |_{y_{\mm{h}}=(\zeta_{\mm{h}})^{-1}(x_{\mm{h}},t)}\right|_{3/2} \leqslant C(a) \|\nabla u\|_2.
\nonumber
\end{align}
Similarly, we also have
\begin{align}
& \left|( \llbracket \mathcal{S}_{\mathcal{A}}(q,0,\eta ) \mathcal{A}e_3-\bar{M}_3\partial_{\bar{M}}\eta \rrbracket/|\mathcal{A}e_3|) |_{y_{\mm{h}}=(\zeta_{\mm{h}})^{-1}(x_{\mm{h}},t)}\right|_{3/2} \leqslant C(a)( |\llbracket q  \rrbracket|_{3/2}+\|\nabla \eta\|_2).\nonumber
\end{align}
Thus, putting \eqref{201809292130} and the above two estimates into \eqref{201811211302xx}  yields
\begin{align}\label{2018111401554}
 |d|_{7/2}\leqslant C(a) (\|(\eta,u)\|_3+ |\llbracket q \rrbracket|_{3/2} ).
\end{align}

Similar to \eqref{2018090232130}, we can estimate that
$$
|v_\pm(x_{\mm{h}},d(x_{\mm{h}},t))|_{5/2}\leqslant C(a)\|u\|_3\mbox{ for a.e. }t>0.
$$
Using \eqref{2018090232130}, \eqref{201811231301}, \eqref{201809271650} and the above estimate, we can derive from \eqref{0103nxx}$_5$ that
\begin{align}
|d_t|_{5/2}\lesssim & |v^3_+|_{5/2}+|\partial_{x_1} d v^1_++\partial_{x_2} d v^2_+|_{5/2}\nonumber  \\
\lesssim &  |v^3_+|_{5/2}+|d|_{7/2}| v_+|_{3/2}+|d|_{5/2}|v_+|_{5/2}\nonumber \\
\leqslant & C(a)( |d|_{7/2}\|u\|_{2}+\|u\|_{3}). \label{201811141556}
\end{align}
We further derive from \eqref{2018111401554} and \eqref{201811141556} that
\begin{align}
|d|_{3}^2\lesssim &  |d^0|_{3}^2+\int_0^t|d|_{7/2}|d_{\tau}|_{5/2}\mm{d}\tau \nonumber \\
\leqslant  & C(a)\left(|d^0|_{3}^2+\int_0^t
\left(\|u\|_{3}^2+ (1+\|u\|_{2}^2) (\|(\eta,u)\|_3^2+ |\llbracket q \rrbracket|_{3/2}^2  \right)\mm{d}\tau\right) .
\label{201811141623}
\end{align}
By the regularity of $(\eta,u,q)$ and \eqref{2018111401554}--\eqref{201811141623}, we immediately get \eqref{201811141303}.

In addition, exploiting \eqref{201811141623},
\eqref{esmmdforinfty} and \eqref{201811231301},  we derive from \eqref{20180928212154} that
\begin{align}
|H|_{1} \leqslant &C(a) |\mathcal{C}|_1  \leqslant C(a) |d|_3 . \label{201811141623xx}
\end{align}

\subsection{Instability relations and properties of initial data}

Finally we turn to deduce the last two conclusions in Theorem \ref{201806012301xx}.  Noting that $\zeta$ satisfies $\det \nabla \zeta=1$, \eqref{201803121601xx21082109},  \eqref{2011801941224538} and \eqref{201806031606}, thus, using transformations of integral variable, we derive the following instability relations  from the conclusions  ``$|d_3(T^\delta)|_{L^1}\geqslant 4\epsilon$" and ``$\|\mm{div}_{\mm{h}}\eta_{\mm{h}}(T^\delta)\|_{L^1}\geqslant 4\epsilon  $" in \eqref{201806012326}:
$$\begin{aligned}
|d(T^\delta)|_{L^1}= & \int_{(\zeta_{\mm{h}})^{-1}(\mathcal{T}, T^\delta)} |\eta_3(y_{\mm{h}},0,T^\delta)|\Theta(y_{\mm{h}},T^\delta)\mm{d}y_{\mm{h}}
\geqslant  \frac{1}{4} |\eta_3( T^\delta)|_0^2  \geqslant  \epsilon
\end{aligned}
$$
and
$$
\begin{aligned}
\|(V_{11}+V_{22})(T^\delta)\|_{L^1_\delta}=
\int_{\zeta^{-1}(\mathcal{T}_{-l}^{\tau}(T^\delta))/\{y_3=0\}}|\mm{div}_{\mm{h}}\eta_{\mm{h}}(T^\delta)|\mm{d}y
\geqslant  \frac{1}{2} \|\mm{div}_{\mm{h}}\eta_{\mm{h}}(T^\delta)\|_0^2 \geqslant 2{\epsilon}.
\end{aligned}
$$
Similarly, we can derive the rest instability relations in \eqref{201808160850}--\eqref{201805201416} from \eqref{201806012326}.

Noting that \eqref{0103nxx}$_4$ and \eqref{0103nxx}$_6$ also hold for $t=0$, thus we get \eqref{2018091019281x}--\eqref{xx2018091019261x}.
Noting that $\eta^0\in H^4$,  then we easily get \eqref{2018090414874xx} from  the definition of $(v^0,V^0,\sigma^0,d^0)$ by following the argument of  \eqref{201811141558}, \eqref{201811141558xx} and \eqref{201809292130}.
This completes the proof of Theorem \ref{201806012301xx}.

\section{Appendix}\label{20180816}

This Appendix is devoted to listing some mathematical analysis tools, which have been used in the previous sections. It should be noted that in this appendix we still adapt the simplified mathematical notations in Section \ref{introud}. For the simplicity, we still use the notation $a\lesssim b$ to mean that $a\leqslant cb$ for some constant $c>0$, where the positive constant $c$ may depend on the domains, and other given parameters locating lemmas.

\begin{lem}
\label{201806171834}
Embedding inequality (see \cite[Theorems 4.12 and 7.58]{ARAJJFF}): Let $D\subset \mathbb{R}^3$ be a bounded Lipschitz domain, then
\begin{align}
&\label{esmmdforinftfdsdy}\|f\|_{L^p(D)}\lesssim \| f\|_{H^1(D)} \mbox{ for }2\leqslant p\leqslant 6,\\
&\label{esmmdforinfty}\|f\|_{C^0(\bar{D})}= \|f\|_{L^\infty(D)}\lesssim\| f\|_{H^2(D)},\\
& \|\phi\|_{C(\mathbb{R}^2)}\lesssim |\phi|_{3/2}\mbox{ for any }\phi\in H^{3/2}(\mathbb{T})
\label{201811231301} ,\\
&\label{201809271650}
\|\phi \|_{W^{k,4}(\mathcal{T})}\lesssim |\phi|_{k+1}\mbox{ for any }\phi\in H^{k+1}(\mathbb{T}),\mbox{ where }k\geqslant 0\mbox{ is integer}.
\end{align}
\end{lem}
\begin{lem}
Interpolation inequality in $H^j$ (see \cite[5.2 Theorem]{ARAJJFF}): Let $D$ be a domain in $\mathbb{R}^N$ satisfying the cone condition, then,
 for any $0\leqslant j< i$, $\varepsilon>0$,
\begin{equation}
\|f\|_{H^j(D)}\lesssim\|f\|_{L^2(D)}^{1-\frac{j}{i}}\|f\|_{H^i(D)}^{\frac{j}{i}}\leqslant C(\varepsilon)\|f\|_{L^2(D)} +\varepsilon\|f\|_{H^i(D)},  \label{201807291850}
\end{equation}
where the constant $C(\varepsilon)$ depends on the domain  and $\varepsilon$, and Young's inequality has been used in the last
inequality in \eqref{201807291850}.
\end{lem}

\begin{lem}\label{xfsddfsf201805072234FRe}
Friedrichs's inequality (see \cite[Lemma 1.42]{NASII04}): Let $1\leqslant p<\infty$, and $D$ be a bounded Lipschitz domain. Let a set $\Gamma\subset \partial D$ be measurable with respect to the $(N-1)$-dimensional measure $\tilde{\mu}:=\mm{meas}_{N-1}$ defined on $\partial D$ and let $\mm{meas}_{N-1}(\Gamma)>0$. Then 
\begin{equation*}
\|w\|_{W^{1,p}(D)}\lesssim \|\nabla w\|_{L^p(D)}^2\end{equation*}
 for all $u\in W^{1,p}(D)$ satisfying that the trace of $u$ on $\Gamma$ is equal to $0$ a.e. with respect to the $(N-1)$-dimensional measure $\tilde{\mu}$.
\end{lem}

\begin{lem}\label{xfsddfsf201805072234Ponc}Poinc\'are's inequality (see \cite[Lemma 1.43]{GYTIAE2}):  Let $1\leqslant p<\infty$, and $D$ be a bounded Lipschitz domain. Then, for any $w\in W^{1,p}(D)$,
\begin{equation*}
\|w\|_{L^p(D)}^p\lesssim  \|\nabla u\|_{L^p(D)}^p+\left|\int_D u\mm{d}y\right|^p .\end{equation*}
 \end{lem}

\begin{lem}\label{xfsddfsf201805072234}
 Korn's inequality (see \cite[Lemma 10.7]{GYTIAE2}): Let $D$ be a bounded domain or $\Omega\!\!\!\!-$, then for any $w\in H_0^1(D)$,
\begin{equation*}\|w\|_{H^1(D)}^2\lesssim \|\mathbb{D} w \|_{L^2(D)}^2.\end{equation*}
\end{lem}

\begin{lem}\label{xfsddfsf2018050813379}
Stein extension theorem (see \cite[Section 5.24]{ARAJJFF}): If $D$ is a domain in $\mathbb{R}^N$ satisfying the strong local Lipschitz condition, then there exists a total extension operator for $D$, i.e.,  there exist a linear operator $E$ mapping $W^{m,p}(D)$ into $W^{m,p}(\mathbb{R}^N)$ for every $p>1$ and every integer $m \geqslant 0 $, and a constant $K=K(m,p)$ such that for every $u\in W^{m,p}(\Omega)$
\begin{enumerate}
  \item[(1)] $Eu(x)=u(x)$ a.e. in $\Omega$,
  \item[(2)] $\|Eu\|_{W^{m,p}(\mathbb{R}^N)}\leqslant K \|u\|_{W^{m,p}(\Omega)}$.
\end{enumerate}
\end{lem}

\begin{lem}
\label{xfsddfsf201805072254}
 Trace estimate: Let $k$, $m$ and $n$ be nonnegative integers, then
\begin{align}
\nonumber \|f|_{y_3=a}\|_{H^{1/2}(\mathcal{T}_{m,n})}\lesssim \|f\|_{1}\mbox{ for any }f\in H^{1} \mbox{ and }a\in (-l,\tau).
\end{align}
which can be derived from  \cite[Theorem 1.46]{NASII04} by Stein extension theorem and a cut-off technique.
\end{lem}
\begin{rem}
By  trace estimate, we further have
\begin{equation}
\label{201808051835}
|f|_{j+1/2}\lesssim \|f\|_{\underline{j},1}\mbox{ for any }f\in H^{j+1}.
\end{equation}
\end{rem}
\begin{lem}
Negative trace estimate:
 \begin{equation}\label{06201929201811041930}|u_3|_{-1/2}\lesssim \|u\|_0+\|\mm{div}u\|_0\;\;\;\mbox{ for any }u:=(u_1,u_2,u_3)\in H_0^1.
\end{equation}
\end{lem}
\begin{pf}
Estimate \eqref{06201929201811041930} can be derived by integration by parts and an inverse trace theorem  \cite[Lemma 1.47]{NASII04}. \hfill $\Box$
\end{pf}
\begin{lem}\label{xfsddfsf20180508}
Product estimates: Let $D\in \mathbb{R}^N$ be a bounded Lipschitz domain for $N=2$ or $3$. The functions $f$, $g$ are defined on $D$, $\phi$, $\varphi$
are defined on $\mathbb{T}$, and $\varpi$ is defined on $\Omega$.
\begin{align}
\label{fgestims}
&\begin{aligned}
&
 \|fg\|_{H^i(D)}\lesssim \left\{     \begin{array}{ll}
                      \|f\|_{H^1(\Omega)}\|g\|_{H^1(D)} & \hbox{ for }i=0;  \\
  \|f\|_{H^i(D)}\|g\|_{H^2(D)} & \hbox{ for }0\leqslant i\leqslant 2;  \\
                    \|f\|_{H^2(D)}\|g\|_{H^i(D)}+\|f\|_{H^i(D)}\|g\|_{H^2(D)}& \hbox{ for }3\leqslant i\leqslant 5.
                    \end{array}         \right.
 \end{aligned} \\
&\label{06011948} |\phi\varphi|_{r}\lesssim |\phi|_{s_1}
|\varphi|_{s_2}\mbox{ for }0\leqslant r\leqslant s_1\leqslant s_2\mbox{ and }s_1>1,\\
&\label{06041605} |\phi\varphi|_{1/2}\lesssim
 \|\phi\|_{W^{1,r}(\mathcal{T})}
|\varphi|_{1/2}   \mbox{ for any }r >2,
\\
 &\label{06041605xdsf} |\varpi\varphi|_{1/2}\lesssim \|\varpi\|_{2}
|\varphi|_{1/2},\\
&\label{0604221702}|\varpi\varphi|_{-1/2}\lesssim
\left\{
  \begin{array}{ll}
 \|\varpi \|_{2}| \varphi|_{-1/2}, \\
 |\varpi |_{r}| \varphi|_{-1/2}    \hbox{ for any }r>3/2.
  \end{array}
\right.
\end{align} \end{lem}
\begin{pf}
The product estimate \eqref{fgestims} can be derived by using
H\"older's inequality and the embedding inequalities \eqref{esmmdforinftfdsdy}--\eqref{esmmdforinfty}.
Estimate \eqref{06011948} can be founded in \cite[Lemma 10.1]{GYTIAE2},
estimate \eqref{06041605} can be obtained by following
the proof of \cite[Lemma 10.2]{GYTIAE2}, and estimate \eqref{06041605xdsf}  can be derived by using \eqref{201809271650} and \eqref{06041605}, which can be founded in  \cite[Theorem 7.58]{ARAJJFF}.
 Finally, using \eqref{06041605xdsf}, we can get \eqref{0604221702} by following
the proof of the third conclusion in  \cite[Lemma 10.1]{GYTIAE2}.
\hfill $\Box$
\end{pf}

\begin{lem}\label{xfsddfsf2018050813379safdadf} Difference quotients and weak derivatives: Let
$D$ be $\Omega\!\!\!\!-$, or $\mathbb{T}$, and $\|\cdot \|_{L^p(\mathbb{T})}:=\|\cdot \|_{L^p(\mathcal{T})}$.
\begin{enumerate}
 \item[(1)] Suppose  $1\leqslant p<\infty$ and $w\in W^{1,p}(D)$. Then
$\|D^{h}_{\mm{h}} w\|_{L^p(D) }\lesssim  \|\nabla_{\mm{h}}w\|_{L^p(D)}$.
 \item[(2)] Assume $1< p<\infty$,  $w\in L^p(D)$, and there exists a constant $c$ such that
$\|D^h_{\mm{h}} w\|_{L^p(D)}\leqslant c$. Then $\nabla_{\mm{h}}
w\in L^{p}(D)$ satisfies $\|\nabla_{\mm{h}}w\|_{L^p(D)}\leqslant c$ and
 $D^{-h_k}_{\mm{h}} w \rightharpoonup \nabla_{\mm{h}}w$ in $L^p(D)$ for some subsequence
$-h_k\to 0$.
\end{enumerate}
\end{lem}
\begin{pf}
Following the argument of  \cite[Theorem 3]{ELGP}, and use the periodicity of
$w$, we can easily get the desired conclusions. \hfill $\Box$
\end{pf}
\begin{lem}
\label{201807312051}
 Estimates in Sobelve--Sobodetskii spaces: If $f\in H^{1/2}(\mathbb{T})$, then
\begin{equation}
\label{201808010840}
\int_{\mathcal{T}}{\int_\mathcal{T}} \frac{|f(x)-f(y)|^2}{|x-y|^3}\mm{d}y\mm{d}x
\leqslant \int_{\mathcal{T}_{4,4}} |\mathfrak{D}_{\mf{h}}^{3/2} f|^2_0 \mm{d}\mf{h}\leqslant \int_{\mathcal{T}_{4,4}}\int_{\mathcal{T}_{6,6}}\frac{|f(x)-f(y)|^2}{|x-y|^3}\mm{d}y\mm{d}x.  \end{equation}
In addition,
\begin{align}
\label{201807312054}
&|f|_{1/2}^2\lesssim |f|_0^2+\int_{\mathcal{T}_{4,4}} |\mathfrak{D}_{\mf{h}}^{3/2} f|^2_0\mm{d}\mf{h}\lesssim |f|_{1/2}^2,\\
\label{201806291928}
 & \int_{\mathcal{T}_{4,4}}\|\mathfrak{D}_{\mf{h}}^{3/2} \varphi\|_0^2\mm{d}\mf{h}\lesssim \|\varphi\|_{1}^2\mbox{ for any }\varphi\in H^1,\\
& \int_{\mathcal{T}_{4,4}} |\mathfrak{D}_{\mf{h}}^{3/2} \phi|^2_{1/2}\mm{d}\mf{h}\lesssim |\phi|_{1}^2\mbox{ for any }\phi\in H^1(\mathbb{T}).  \label{201811132119}
\end{align}
\end{lem}
\begin{pf}
Noting that
$$\begin{aligned}
&\mathcal{T}\times \mathcal{T} \subset \{(y,x)\in \mathbb{R}^4~|~y\in \mathcal{T},\  x\in \mathcal{T}_{4,4}+y\} ,\\
&\{(y,x)\in\mathbb{R}^4~|~y\in \mathcal{T},\ x\in \mathcal{T}_{4,4}+y\}\subset \mathcal{T}_{4,4}\times  \mathcal{T}_{6,6} ,
\end{aligned}$$
thus, using the definition of operator $\mathfrak{D}_{\mf{h}}^{3/2} $, change of variable, Fubini's theorem and the periodicity of $f$, we can easily get \eqref{201808010840}. Exploiting \eqref{201808010840}, Fubini's theorem and trace estimate, we can further obtain \eqref{201807312054} and \eqref{201806291928}.

Now turn to prove \eqref{201811132119}. Since $\phi\in H^{1/2}(\mathbb{T})$, by Stein extension theorem, we have a extension function $\tilde{\phi}$ satisfying
$$ \tilde{\phi} \in H^{1}(\mathbb{R}^2),\ \tilde{\phi}=\phi\mbox{ a.e. in }\mathcal{T}_{6,6}\mbox{ and }\|\tilde{\phi}\|_{H^{1}(\mathbb{R}^2)}\leqslant c|\phi|_{1}$$
for some constant $c$.
Then we further have
 $$\int_{\mathcal{T}_{4,4}} |\mathfrak{D}_{\mf{h}}^{3/2} \phi|^2_{1/2}\mm{d}\mf{h}\leqslant
\int_{\mathbb{R}^2} |\mathfrak{D}_{\mf{h}}^{3/2} \tilde{\phi}|^2_{1/2}\mm{d}\mf{h}\lesssim \|\tilde{\phi}\|_{H^1(\mathbb{R}^2)}^2\lesssim |\phi|_1,$$
which yields \eqref{201811132119}.
\hfill $\Box$
\end{pf}

\begin{lem}\label{20180812}
Dual estimates: If $\varphi$, $\psi\in H^{1/2}$, and $\partial_{\mm{h}}\varphi \in L^1(\mathcal{T})$, then
\begin{equation}
\label{201808121426}
\left|\int_{\Sigma}\partial_{\mm{h}}\varphi \psi\mm{d}y_\mm{h}\right|\lesssim |\varphi|_{1/2} |\psi|_{1/2}.
\end{equation}Moreover, we have
\begin{equation}
\label{201808121247}
|\partial_{\mm{h}}\varphi|_{-1/2}\lesssim |\varphi|_{1/2}.
\end{equation}
\end{lem}
\begin{pf}

First, by using the Parseval's relation on the Torus (see \cite[Proposition 3.1.16]{grafakos2008classical}), H\"older's inequality in $l^2$, and Newton--Leibnitz formula, we can derive
 \begin{equation}
\label{201808121247ss}
\left|\int_{\Sigma}\partial_{\mm{h}}f\chi\mm{d}y_\mm{h}\right|\lesssim \|f\|_1 \|\chi\|_1
\end{equation}
for any $f$, $\chi\in H^1_0$ satisfying $\partial_{\mm{h}}f|_{\Sigma}\in L^1(\mathcal{T})$. Then we further use Stein extension theorem and  cut-off technique to verify that \eqref{201808121247ss} also holds for $f$, $g\in H^1$ satisfying $\partial_{\mm{h}}f|_{\Sigma}\in L^1(\mathcal{T})$.

Finally, since $\varphi$, $\psi\in H^{1/2}$, by an inverse trace theorem \cite[Lemma 1.47]{NASII04}, there exist a Lipschitz domain $D$, and two functions  $f$, $\chi\in H^1(D)$, such that  $f|_{\mathcal{T}}=\varphi$, $\chi|_\mathcal{T}=\psi$, $\mathcal{T}\subset \partial D$ and
$$\|f\|_{H^1(D)}\lesssim |\varphi|_{1/2}\mbox{ and }\|\chi\|_{H^1(D)}\lesssim |\psi|_{1/2}.$$
Consequently, by \eqref{201808121247ss}, we have \begin{equation*}
\left|\int_{\Sigma}\partial_{\mm{h}}\varphi\psi\mm{d}y_\mm{h}\right|\lesssim   |\varphi|_{1/2} |\psi|_{1/2},
\end{equation*}
which yields \eqref{201808121426} and \eqref{201808121247}.
 \hfill $\Box$
\end{pf}
\begin{lem}\label{201805241042}
 Existence theory of one-layer (steady) Stokes problem (see \cite[Lemma A.8]{WYJTIKCT}):
Denoting $D$ be a bounded domain of $C^{k+2}$-class, or $\Omega_\pm$, and $\partial\Omega_\pm:=\Sigma_0\cup \Sigma_\pm$. Let $\mu$ be constant, $k\geqslant 0$,  $f^{\mm{O},1}\in H^{k}$, $f^{\mm{O},2}\in H^{k+1}$ and $f^{\mm{O},3}\in H^{k+3/2}$ be given such that
\begin{equation*}
\int_{D} f^{\mm{O},2} \mm{d}x=\int_{\partial D}  f^{\mm{O},3}\vec{n}\mm{d}y_{\mm{h}},
\end{equation*}
where $\vec{n}$ denote the outer normal vector of $\partial D$,   then there exists a unique strong solution $ (u,q)\in   H^{k+2}(D) \times \underline{H}^{ k+1}(D)$, which solves
 \begin{equation}
 \left\{\begin{array}{ll}
\nabla q-\mu\Delta u =f^{\mm{O},1},\quad \mm{div}u=f^{\mm{O},2}&\mbox{ in }  D, \\[1mm]
 u=f^{\mm{O},3}  &\mbox{ on }\partial D.
\end{array}\right.
\label{201808311539}
\end{equation}
Moreover,  any solution  $ (u,q)\in   H^{k+2}(D) \times  {H}^{ k+1}(D)$ of \eqref{201808311539} enjoys the following estimate \begin{equation}
\label{xfsddfsf201705141252}
\|u\|_{H^{k+2}(D)}+\|  \nabla q\|_{H^k(D)}\lesssim \|f^{\mm{O},1}\|_{H^k(D)}+\|f^{\mm{O},2}\|_{H^{k+1}(D)}+\|f^{\mm{O},3}\|_{H^{k+3/2}(\partial D)}.
\end{equation}
\end{lem}
\begin{lem}
\label{xfsddfsf201805072212}
Existence  theory of a stratified (steady) Stokes problem (see \cite[Theorem 3.1]{WYJTIKCT})]: let $k\geqslant 0$, $f^{\mm{S},1}\in H^{k}$, $f^{\mm{S},2}\in H^{k+1}$, $f^{\mm{S},3}\in H^{k+1/2}$ and $f^{\mm{S},4}\in H^{k+3/2}$, then there exists a unique solution $(u,q)\in H^{k+2}\cap \underline{H}^{k+1}$ satisfying
\begin{align}
\left\{\begin{array}{ll}
\nabla q-\mu\Delta u =f^{\mm{S},1},\quad \mm{div}u=f^{\mm{S},2}&\mbox{ in }  \Omega, \\[1mm]
\llbracket u \rrbracket=0,\ \llbracket (q I-\mathbb{D}u )e_3 \rrbracket =f^{\mm{S},3} &\mbox{ on }\Sigma, \\
u=0 &\mbox{ on }\Sigma_-^+.
\end{array}\right.
\label{201808311541}
\end{align}
Moreover,  any solution  $(u,q)\in   H^{k+2} \times  {H}^{ k+1}$ of \eqref{201808311541} enjoys the following estimate
\begin{equation}
\label{2011805302036}
 \|(u ,q )\|_{\mm{S},k}   \lesssim \|f^{\mm{S},1}\|_{k}+\|f^{\mm{S},2}\|_{k+1}+|f^{\mm{S},3}|_{k+1/2} .
\end{equation}
\end{lem}
\begin{lem}
\label{xfsddfsf201805081359}
Inverse function theorem (see \cite[Theorem 9.24]{WalterRudin}): Let $N\geqslant 1$. Suppose $f$ is a $C^1$-mapping of an open set $E\subset \mathbb{R}^N$ into $\mathbb{R}^N$, $f'(a)$ is invertible for some $a\in E$, and $b=f(a)$. Then
\begin{enumerate}
  \item[(1)] there exist open sets $U$ and $V$ in $\mathbb{R}^N$ such that $a\in U$, $b\in V$, $f$ is one-to-one on $U$, and $f(U)=V$;
  \item[(2)] if $g$ is the inverse of $f$ (which exists, by (1)), defined in $V$ by
  $g(f(x))=x$ for $x\in U$, then $g\in C^1(V)$ and $\nabla g(y)=(\nabla f)^{-1}|_{x=f^{-1}(y)}$.
\end{enumerate}
\end{lem}

\begin{lem}
\label{xfsddfsf201803121937}
Global existence of inverse function:
Let $N\geqslant 2$, $\mathcal{D}\subseteq \mathbb{R}^N$ be an open set, $\zeta:\mathcal{D}\to \mathbb{R}^N$
belongs to $C^1(\mathcal{D})$, and $\nabla \zeta(x)$ be invertible for all $x\in \mathcal{D}$. Suppose that
$K$ is a compact subset of $\mathcal{D}$, the boundary $\partial K$ is connected, and $\zeta:\partial K\to \mathbb{R}^N$ is injective, Then $\zeta :K\to \mathbb{R}^N$ is injective.
\end{lem}
\begin{pf}
See Theorem 2 in \cite{Kestelman} or Corollary 2 in \cite{MGHCOLDUMJ}.
\end{pf}

\begin{lem}\label{201806091256}
Concerning threshold of coefficient of surface tension:
We have
\begin{equation}
\label{2018060102215680}
a:=\sup_{w\in H_{\sigma,3}^1 }
\frac{|w_3|_0^2}{
|\nabla_{\mm{h}}w_3 |_0^2}=  {\max\{L_1^2,L_2^2\} }.
\end{equation}
\end{lem}
\begin{pf}
Without loss of generality, we assume that $L_1^2={\max\{L_1^2,L_2^2\}}$. It should be noted that
\begin{equation}\label{201806101508}
|\nabla_{\mm{h}}w_3 |_{0}^2=0\mbox{ if and only if } w_3=0\mbox{ for any given }w\in H_{\sigma,\Sigma}^1.
\end{equation}
In fact, let $w\in H_\sigma^1$.
Since $\mm{div}w=0$, we have
$$-\int_{\mathcal{T}} w_3\mm{d}y_{\mm{h}}=\int_{\mathcal{T}\times (0,\tau)}\mm{div}w\mm{d}y=0.$$
Thus, using Pocare's inequality, we have
$$|w_3|_{0}\leqslant |\nabla_{\mm{h}}w_3|_{0},$$
which immediately implies the assertion \eqref{201806101508}.

(1) Now we prove  $a\geqslant L_1^2$. We choose a non-zero function $\psi\in H_0^2(-l,\tau)$. We denote
 $$\tilde{w}=(\psi'  (y_3)\cos( L^{-1}_1y_1 ),0,  L^{-1}_1\psi (y_3) \sin (L^{-1}_1 y_1 )),$$ then $\tilde{w} \in H^1_{\sigma,3}$ and
$$\begin{aligned}
& \frac{|\tilde{w}_3|^2_0} {|\nabla_{\mm{h}} \tilde{ w}_3|^2_0}
=  \frac{   \int_0^{2\pi L_1}
 \sin^2 ( L^{-1}_1 y_1 )\mm{d}y_1  }{L^{-2}_1  \int_0^{2\pi L_1}
  \cos^2(  L^{-1}_1 y_1) \mm{d}y_1  }=L_1^2,
\end{aligned} $$
which yields $a\geqslant L_1^2$.

(2)  Next we turn to the proof of $a\leqslant L_1^2$. Let $w\in H_{\sigma,3}^1$.  Then
 $|\nabla_{\mm{h}} w_3 |_{0}^2\neq 0$ by \eqref{201806101508}.
 Let $\hat{w}_3(\xi,y_3)$ be the horizontal Fourier transform of $w_3(y) $, i.e.,
$$ \hat{w}_3(\xi,y_3)=\int_{\mathcal{T}^0}w_3( y_{\mm{h}},y_3)e^{-\mm{i}y_{\mm{h}}\cdot\xi}\mm{d}x_{\mm{h}}, $$
where  $\xi=(\xi_1,\xi_2)$, then $\widehat{\partial_3  w_3} = \partial_{3} \widehat{w}_3$. We denote $\psi(\xi,y_3):= \psi_1(\xi,y_3) + \mm{i}\psi_2(\xi,y_3):=\hat{w}_3(\xi,y_3)$, where $\psi_1$ and $\psi_2$ are real functions.
Noting that $\psi(0)=0$, by  Parseval theorem  (see \cite[Proposition 3.1.16]{grafakos2008classical}), we have
 \begin{equation*}
\begin{aligned}
|\nabla_{\mm{h}}{w}_3 |^2_0=& \frac{1}{\pi^2 L_1L_2}\sum_{\xi\in (L^{-1}_1\mathbb{Z}\times L^{-1}_2\mathbb{Z})}
|\xi|^2|\psi(\xi,0)|^2
\geqslant
 L^{-1}_1 |w_3|^2_0,
\end{aligned}\end{equation*}
which immediately yields that $a\leqslant L^2_1$. The proof is complete.
\hfill
$\Box$
\end{pf}
\begin{rem}
Observing the derivation of ``$a\leqslant L^2_1$", it is easy to see that
$$
\sup_{w\in H_{0,3}^1 }
\frac{|w_3|_{0}^2}{
|\nabla_{\mm{h}}w_3|_{0}^2} \leqslant  {\max\{L_1^2,L_2^2\} }.$$
which, together with \eqref{2018060102215680} and the fact $ H_\sigma^1\subset H_0^1$, implies that
\begin{equation*}
 \sup_{w\in H_{0,3}^1 }
\frac{|w_3|_{0}^2}{
|\nabla_{\mm{h}}w_3|_{0 }^2}=  {\max\{L_1^2,L_2^2\} }.
\end{equation*}
Here we have defined that $H_{0,3}^1:=\{w\in H_{0}^1~|~w_3\in H^1(\mathbb{T}), \int_{\mathcal{T}}w_3(y_{\mm{h}},0)\mm{d}y_{\mm{h}}=0,\ w_3\neq 0\}$.\end{rem}

\begin{lem}\label{201810092222}
Dirichlet's approximation theorem: given $\alpha \in \mathbb{R}$ and any positive integer $N$, there exist integers $n$, $m$ with $1\leqslant n\leqslant N$ such that $|n\alpha- m|<1/N$.
\end{lem}
\begin{pf}
See \cite[Lemma 4.21]{LEEMJHNgrdtextsinn}.
\end{pf}
 \begin{lem}
\label{201806151836}
Let  $L_1$, $L_2$ be two positive constants, and $\bar{M}\in \mathbb{R}^3$ be a non-zero  vector.
If $\bar{M}_3=0$, then, for any given constant $a>0$, there always exists a non-zero function $w\subset H_{\sigma}^1$ such that
\begin{equation}
\label{201810052024}
\|\partial_{\bar{M}}w\|^2_{0}< a|w_3 |_{0}^2.
\end{equation}
\end{lem}
\begin{pf}
We only consider the case $\bar{M}_1\neq 0$, because the other case of ``$\bar{M}_1= 0$ and $\bar{M}_2\neq 0$" can be dealt with similarly.

Let  $\psi\in C_0^\infty(\mathbb{R})$ satisfying $\psi(0)\neq 0$ and
$\psi= 0 $ on $\mathbb{R}\setminus (a,b)$, where $(a,b)\subset (-l,\tau)$.
Denoting $ C_\sigma^\infty:=\{w\in C_0^\infty(\Omega_{-l}^\tau)~|~\mm{div}w=0\}$.

If $L_1 \bar{M}_2 / \bar{M}_1  L_2$ is a rational number, we can choose some integer $i$ such that
$  {\cos} (i \bar{M}_2y_1/\bar{M}_1L_2- iy_2/L_2)  $
 is a periodic function with  periodicity length $2\pi L_1$ in $y_1$-direction and with periodicity length $2\pi L_2$ in $y_2$-direction.
We further define
 \begin{align}
 w =& (0, -\psi'  {\cos} (i \bar{M}_2y_1/\bar{M}_1L_2- iy_2/L_2),  L^{-1}_2 i \psi  {\sin} (i \bar{M}_2y_1/\bar{M}_1L_2- iy_2/L_2)). \nonumber
\end{align}
Then $ w \in  H_\sigma^1$. Moreover, it is easy to check that
$|w_3|_{0}^2\neq 0$  and $\partial_{\bar{M}}w=0 $. Hence \eqref{201810052024} holds.

If $L_1 \bar{M}_2 / \bar{M}_1  L_2$ is not a rational number,  by
Dirichlet's approximation theorem (see Lemma \ref{201810092222}), there exists a sequence $\{(p_n,q_n)\}_{n=1}^\infty$
such that $p_n\geqslant 1$, $q_n$ are integers, and
$$|r_n|\to 0\mbox{ as }n\to \infty,$$
where we have defined that $r_n:=p_n\bar{M}_2 / L_2 - q_n \bar{M}_1 / L_1  $.
Now we define
 \begin{align}
 w =& (0, -\psi'  {\cos} (
q_n y_1 /L_1 - p_n  y_2/ L_2
 ),  L^{-1}_2 p_n \psi  {\sin} (
q_n y_1 /L_1 - p_n  y_2/ L_2 )). \nonumber
\end{align}
Then $ w \in  H_\sigma^1$. Moreover, it is easy to check that
$$
\begin{aligned}
& |w_3|_0^2= 8\pi^2 L_1 L^{-1}_2 p_n^2\psi^2 (0)  ,\\
& \|\partial_{\bar{M}}w\|_0^2=
 8\pi^2 L_1 L_2 r_n^2 ( \|\psi'   \|_{L^2(-l,\tau)}^2/2+L_2^{-2}p_n^2 \|\psi    \|_{L^2(-l,\tau)}^2).
\end{aligned}
$$ Hence \eqref{201810052024} holds for sufficiently large $n$.
\hfill$\Box$
\end{pf}

\vspace{4mm} \noindent\textbf{Acknowledgements.}
 The research of Fei Jiang was supported by NSFC (Grant No. 11671086), and the NSF of Fujian Province of China (Grant No. 2016J06001), and the research of Song Jiang by the Basic Research Program (2014CB745002)
and NSFC (Grant Nos. 11571046, 11631008).

\renewcommand\refname{References}
\renewenvironment{thebibliography}[1]{%
\section*{\refname}
\list{{\arabic{enumi}}}{\def\makelabel##1{\hss{##1}}\topsep=0mm
\parsep=0mm
\partopsep=0mm\itemsep=0mm
\labelsep=1ex\itemindent=0mm
\settowidth\labelwidth{\small[#1]}%
\leftmargin\labelwidth \advance\leftmargin\labelsep
\advance\leftmargin -\itemindent
\usecounter{enumi}}\small
\def\newblock{\ }
\sloppy\clubpenalty4000\widowpenalty4000
\sfcode`\.=1000\relax}{\endlist}
\bibliographystyle{model1b-num-names}

\begin{thebibliography}{55}
\expandafter\ifx\csname natexlab\endcsname\relax\def\natexlab#1{#1}\fi
\providecommand{\bibinfo}[2]{#2}
\ifx\xfnm\relax \def\xfnm[#1]{\unskip,\space#1}\fi
\bibitem[{Adams and John(2005)}]{ARAJJFF}
\bibinfo{author}{R.A. Adams}, \bibinfo{author}{J.F.F. John},
  \bibinfo{title}{{Sobolev Space}}, \bibinfo{publisher}{Academic Press: New
  York}, \bibinfo{year}{2005}.
\bibitem[{Boffetta et~al.(2010)Boffetta, Mazzino, Musacchio and
  Vozella}]{BGMAMSVLRTI}
\bibinfo{author}{G.~Boffetta}, \bibinfo{author}{A.~Mazzino},
  \bibinfo{author}{S.~Musacchio}, \bibinfo{author}{L.~Vozella},
  \bibinfo{title}{{ Rayleigh--Taylor instability in a viscoelastic binary
  fluid}}, \bibinfo{journal}{J. Fluid Mech.} \bibinfo{volume}{643}
  (\bibinfo{year}{2010}) \bibinfo{pages}{127--136}.
\bibitem[{Chandrasekhar(1961)}]{CSHHSCPO}
\bibinfo{author}{S.~Chandrasekhar}, \bibinfo{title}{{Hydrodynamic and
  Hydromagnetic Stability, The International Series of Monographs on Physics}},
  \bibinfo{publisher}{Oxford, Clarendon Press}, \bibinfo{year}{1961}.
\bibitem[{Duan et~al.(2015)Duan, Jiang and Yin}]{DRJFJSRS}
\bibinfo{author}{R.~Duan}, \bibinfo{author}{F.~Jiang}, \bibinfo{author}{J.P.
  Yin}, \bibinfo{title}{{ Rayleigh--Taylor instability for compressible
  rotating flows}}, \bibinfo{journal}{Acta Math. Sci. Engl. Ser.}
  \bibinfo{volume}{35} (\bibinfo{year}{2015}) \bibinfo{pages}{1359--1385}.
\bibitem[{Evans(1998)}]{ELGP}
\bibinfo{author}{L.C. Evans}, \bibinfo{title}{{Partial Differential
  Equations}}, \bibinfo{publisher}{American Mathematical Society},
  \bibinfo{address}{USA}, \bibinfo{year}{1998}.
\bibitem[{Friedlander et~al.(1997)Friedlander, Strauss and Vishik}]{FSSWVMNA}
\bibinfo{author}{S.~Friedlander}, \bibinfo{author}{W.~Strauss},
  \bibinfo{author}{M.~Vishik}, \bibinfo{title}{{Nonlinear instability in an
  ideal fluid}}, \bibinfo{journal}{Ann. Inst. H. Poincar\'e Anal. Non
  Lin\'eaire} \bibinfo{volume}{14} (\bibinfo{year}{1997})
  \bibinfo{pages}{187--209}.
\bibitem[{Grafakos(2008)}]{grafakos2008classical}
\bibinfo{author}{L.~Grafakos}, \bibinfo{title}{{ Classical fourier analysis
  (second edition)}}, \bibinfo{publisher}{Springer},
  \bibinfo{address}{Germany}, \bibinfo{year}{2008}.
\bibitem[{Guo et~al.(2007)Guo, Hallstrom and Spirn}]{GYHCSDDC}
\bibinfo{author}{Y.~Guo}, \bibinfo{author}{C.~Hallstrom},
  \bibinfo{author}{D.~Spirn}, \bibinfo{title}{{Dynamics near unstable,
  interfacial fluids}}, \bibinfo{journal}{Commun. Math. Phys.}
  \bibinfo{volume}{270} (\bibinfo{year}{2007}) \bibinfo{pages}{635--689}.
\bibitem[{Guo and Strauss(1995{\natexlab{a}})}]{GYSWIC}
\bibinfo{author}{Y.~Guo}, \bibinfo{author}{W.A. Strauss},
  \bibinfo{title}{{Instability of periodic BGK equilibria}},
  \bibinfo{journal}{Comm. Pure Appl. Math.} \bibinfo{volume}{48}
  (\bibinfo{year}{1995}{\natexlab{a}}) \bibinfo{pages}{861--894}.
\bibitem[{Guo and Strauss(1995{\natexlab{b}})}]{GYSWICNonlinea}
\bibinfo{author}{Y.~Guo}, \bibinfo{author}{W.A. Strauss},
  \bibinfo{title}{{Nonlinear instability of double-humped equilibria}},
  \bibinfo{journal}{Ann. Inst. H. Poinc\'are Anal. Non
  Lin$\mathrm{\acute{e}}$aire} \bibinfo{volume}{12}
  (\bibinfo{year}{1995}{\natexlab{b}}) \bibinfo{pages}{339--352}.
\bibitem[{Guo and Tice(2011{\natexlab{a}})}]{GYTI1}
\bibinfo{author}{Y.~Guo}, \bibinfo{author}{I.~Tice},
  \bibinfo{title}{{Compressible, inviscid Rayleigh--Taylor instability}},
  \bibinfo{journal}{Indiana Univ. Math. J.} \bibinfo{volume}{60}
  (\bibinfo{year}{2011}{\natexlab{a}}) \bibinfo{pages}{677--712}.
\bibitem[{Guo and Tice(2011{\natexlab{b}})}]{GYTI2}
\bibinfo{author}{Y.~Guo}, \bibinfo{author}{I.~Tice}, \bibinfo{title}{{Linear
  Rayleigh--Taylor instability for viscous, compressible fluids}},
  \bibinfo{journal}{SIAM J. Math. Anal.} \bibinfo{volume}{42}
  (\bibinfo{year}{2011}{\natexlab{b}}) \bibinfo{pages}{1688--1720}.
\bibitem[{Guo and Tice(2013{\natexlab{a}})}]{GYTIAE2}
\bibinfo{author}{Y.~Guo}, \bibinfo{author}{I.~Tice}, \bibinfo{title}{{Almost
  exponential decay of periodic viscous surface waves without surface
  tension}}, \bibinfo{journal}{Arch. Ration. Mech. Anal.} \bibinfo{volume}{207}
  (\bibinfo{year}{2013}{\natexlab{a}}) \bibinfo{pages}{459--531}.
\bibitem[{Guo and Tice(2013{\natexlab{b}})}]{GYTIDAP}
\bibinfo{author}{Y.~Guo}, \bibinfo{author}{I.~Tice}, \bibinfo{title}{{Decay of
  viscous surface waves without surface tension in horizontally infinite
  domains}}, \bibinfo{journal}{Anal. PDE} \bibinfo{volume}{6}
  (\bibinfo{year}{2013}{\natexlab{b}}) \bibinfo{pages}{1429--1533}.
\bibitem[{Guo and Tice(2013{\natexlab{c}})}]{GYTILW1}
\bibinfo{author}{Y.~Guo}, \bibinfo{author}{I.~Tice}, \bibinfo{title}{{Local
  well-posedness of the viscous surface wave problem without surface tension}},
  \bibinfo{journal}{Anal. PDE} \bibinfo{volume}{6}
  (\bibinfo{year}{2013}{\natexlab{c}}) \bibinfo{pages}{287--369}.
\bibitem[{Hester et~al.(1996)Hester, Stone, Scowen and et.al.}]{HJJSJMSPAWFPC}
\bibinfo{author}{J.J. Hester}, \bibinfo{author}{J.M. Stone},
  \bibinfo{author}{P.A. Scowen}, \bibinfo{author}{et.al.},
  \bibinfo{title}{{WFPC2 studies of the Crab Nebula. III. magnetic
  Rayleigh--Taylor instabilities and the origin of the filaments }},
  \bibinfo{journal}{Astrophys J.} \bibinfo{volume}{456} (\bibinfo{year}{1996})
  \bibinfo{pages}{225--233}.
\bibitem[{Hide(1955)}]{HRWP}
\bibinfo{author}{R.~Hide}, \bibinfo{title}{{Waves in a heavy, viscous,
  incompressible, electrically conducting fluid of variable density, in the
  presence of a magnetic field}}, \bibinfo{journal}{Proc. Roy. Soc. (London) A}
  \bibinfo{volume}{233} (\bibinfo{year}{1955}) \bibinfo{pages}{376--396}.
\bibitem[{Huang et~al.(2017)Huang, Jiang and Wang}]{HGJJJWWWOJM}
\bibinfo{author}{G.J. Huang}, \bibinfo{author}{J.~Jiang}, \bibinfo{author}{W.W.
  Wang}, \bibinfo{title}{{On the nonlinear Rayleigh--Taylor instability of
  nonhomogeneous incompressible viscoelastic fluids under $L^2$-norm}},
  \bibinfo{journal}{J. Math. Anal. Appl.} \bibinfo{volume}{455}
  (\bibinfo{year}{2017}) \bibinfo{pages}{873--904}.
\bibitem[{Hwang and Guo(2003)}]{HHJGY}
\bibinfo{author}{H.J. Hwang}, \bibinfo{author}{Y.~Guo}, \bibinfo{title}{{On the
  dynamical Rayleigh--Taylor instability}}, \bibinfo{journal}{Arch. Rational
  Mech. Anal.} \bibinfo{volume}{167} (\bibinfo{year}{2003})
  \bibinfo{pages}{235--253}.
\bibitem[{Isobe et~al.(2005)Isobe, Miyagoshi, Shibata and
  Yokoyama}]{IHMTSKYTFstru}
\bibinfo{author}{H.~Isobe}, \bibinfo{author}{T.~Miyagoshi},
  \bibinfo{author}{K.~Shibata}, \bibinfo{author}{T.~Yokoyama},
  \bibinfo{title}{{Filamentary structure on the Sun from the magnetic
  Rayleigh--Taylor instability}}, \bibinfo{journal}{Nature}
  \bibinfo{volume}{434} (\bibinfo{year}{2005}) \bibinfo{pages}{478--481}.
\bibitem[{Isobe et~al.(2006)Isobe, Miyagoshi, Shibata and Yokoyama}]{IHMTSKYTT}
\bibinfo{author}{H.~Isobe}, \bibinfo{author}{T.~Miyagoshi},
  \bibinfo{author}{K.~Shibata}, \bibinfo{author}{T.~Yokoyama},
  \bibinfo{title}{{Three-dimensional simulation of solar emerging flux using
  the earth simulator I. magnetic Rayleigh--Taylor instability at the top of
  the emerging flux as the origin of filamentary structure}},
  \bibinfo{journal}{Publ. Astron. Soc. Japan} \bibinfo{volume}{58}
  (\bibinfo{year}{2006}) \bibinfo{pages}{423--438}.
\bibitem[{Jang et~al.(2016)Jang, Tice and Wang}]{JJTIWYHCMP}
\bibinfo{author}{J.~Jang}, \bibinfo{author}{I.~Tice}, \bibinfo{author}{Y.J.
  Wang}, \bibinfo{title}{{The compressible viscous surface-internal wave
  problem: stability and vanishing surface tension limit}},
  \bibinfo{journal}{Commun. Math. Phys.} \bibinfo{volume}{343}
  (\bibinfo{year}{2016}) \bibinfo{pages}{1039--1113}.
\bibitem[{Jiang and Jiang(2014)}]{JFJSO2014}
\bibinfo{author}{F.~Jiang}, \bibinfo{author}{S.~Jiang}, \bibinfo{title}{{On
  instability and stability of three-dimensional gravity flows in a bounded
  domain}}, \bibinfo{journal}{Adv. Math.} \bibinfo{volume}{264}
  (\bibinfo{year}{2014}) \bibinfo{pages}{831--863}.
\bibitem[{Jiang and Jiang(2015)}]{JFJSJMFM}
\bibinfo{author}{F.~Jiang}, \bibinfo{author}{S.~Jiang}, \bibinfo{title}{{ On
  linear instability and stability of the Rayleigh--Taylor problem in
  magnetohydrodynamics}}, \bibinfo{journal}{J. Math. Fluid Mech.}
  \bibinfo{volume}{17} (\bibinfo{year}{2015}) \bibinfo{pages}{639--668}.
\bibitem[{Jiang and Jiang(2016)}]{JFJSSETEFP}
\bibinfo{author}{F.~Jiang}, \bibinfo{author}{S.~Jiang}, \bibinfo{title}{{
  Stabilizing effect of the equilibrium magnetic fields upon the Parker
  instability}}, \bibinfo{journal}{Under review}  (\bibinfo{year}{2016}).
\bibitem[{Jiang and Jiang(2018{\natexlab{a}})}]{JFJSJMFMOSERT}
\bibinfo{author}{F.~Jiang}, \bibinfo{author}{S.~Jiang}, \bibinfo{title}{{ On
  the stabilizing effect of the magnetic field in the magnetic Rayleigh--Taylor
  problem}}, \bibinfo{journal}{SIAM J. Math. Anal.} \bibinfo{volume}{50}
  (\bibinfo{year}{2018}{\natexlab{a}}) \bibinfo{pages}{491--540}.
\bibitem[{Jiang and Jiang(2018{\natexlab{b}})}]{JFJSOMITN}
\bibinfo{author}{F.~Jiang}, \bibinfo{author}{S.~Jiang}, \bibinfo{title}{{On
  magnetic inhibition theory in non-resistive magnetohydrodynamic fluids}},
  \bibinfo{journal}{arXiv:1803.00307v2}  (\bibinfo{year}{2018}{\natexlab{b}}).
\bibitem[{Jiang et~al.(2016{\natexlab{a}})Jiang, Jiang and Wang}]{JFJSWWWN}
\bibinfo{author}{F.~Jiang}, \bibinfo{author}{S.~Jiang}, \bibinfo{author}{W.W.
  Wang}, \bibinfo{title}{{Nonlinear Rayleigh--Taylor instability in
  nonhomogeneous incompressible viscous magnetohydrodynamic fluids}},
  \bibinfo{journal}{Discrete Contin. Dyn. Syst.-S} \bibinfo{volume}{9}
  (\bibinfo{year}{2016}{\natexlab{a}}) \bibinfo{pages}{1853--1898}.
\bibitem[{Jiang et~al.(2014)Jiang, Jiang and Wang}]{JFJSWWWOA}
\bibinfo{author}{F.~Jiang}, \bibinfo{author}{S.~Jiang}, \bibinfo{author}{Y.J.
  Wang}, \bibinfo{title}{{On the Rayleigh--Taylor instability for the
  incompressible viscous magnetohydrodynamic equations}},
  \bibinfo{journal}{Comm. Partial Differential Equations} \bibinfo{volume}{39}
  (\bibinfo{year}{2014}) \bibinfo{pages}{399--438}.
\bibitem[{Jiang et~al.(2017)Jiang, Jiang and Wu}]{JFJWGCOSdd}
\bibinfo{author}{F.~Jiang}, \bibinfo{author}{S.~Jiang}, \bibinfo{author}{G.C.
  Wu}, \bibinfo{title}{{ On stabilizing effect of elasticity in the
  Rayleigh--Taylor problem of stratified viscoelastic fluids}},
  \bibinfo{journal}{J. Funct. Anal.} \bibinfo{volume}{272}
  (\bibinfo{year}{2017}) \bibinfo{pages}{3763--3824}.
\bibitem[{Jiang et~al.(2016{\natexlab{b}})Jiang, Wu and Zhong}]{FJWGCZXOE}
\bibinfo{author}{F.~Jiang}, \bibinfo{author}{G.C. Wu},
  \bibinfo{author}{X.~Zhong}, \bibinfo{title}{{ On exponential stability of
  gravity driven viscoelastic flows}}, \bibinfo{journal}{J. Differential
  Equations} \bibinfo{volume}{260} (\bibinfo{year}{2016}{\natexlab{b}})
  \bibinfo{pages}{7498--7534}.
\bibitem[{Jun et~al.(1966)Jun, Norman and Stone}]{JBINMLSJMA}
\bibinfo{author}{B.I. Jun}, \bibinfo{author}{M.L. Norman},
  \bibinfo{author}{J.M. Stone}, \bibinfo{title}{{A numerical study of
  Rayleigh--Taylor instability in magnetic fluids}},
  \bibinfo{journal}{Astrophys J.} \bibinfo{volume}{453} (\bibinfo{year}{1966})
  \bibinfo{pages}{332--349}.
\bibitem[{Kestelman(1971)}]{Kestelman}
\bibinfo{author}{H.~Kestelman}, \bibinfo{title}{Mappings with non-vanishing
  jacobian}, \bibinfo{journal}{The American Mathematical Monthly}
  \bibinfo{volume}{78} (\bibinfo{year}{1971}) \bibinfo{pages}{662--663}.
\bibitem[{Kruskal and Schwarzschild(1954)}]{KMSMSP}
\bibinfo{author}{M.~Kruskal}, \bibinfo{author}{M.~Schwarzschild},
  \bibinfo{title}{Some instabilities of a completely ionized plasma},
  \bibinfo{journal}{Proc. Roy. Soc. (London) A} \bibinfo{volume}{233}
  (\bibinfo{year}{1954}) \bibinfo{pages}{348--360}.
\bibitem[{Le~Meur(2011)}]{LeMHVJ}
\bibinfo{author}{V.J.H. Le~Meur}, \bibinfo{title}{Well-posedness of surface
  wave equations above a viscoelastic fluid}, \bibinfo{journal}{J. Math. Fluid
  Mech.} \bibinfo{volume}{13} (\bibinfo{year}{2011}) \bibinfo{pages}{481--514}.
\bibitem[{Lee(2013)}]{LEEMJHNgrdtextsinn}
\bibinfo{author}{J.M. Lee}, \bibinfo{title}{Introduction to Smooth Manifolds},
  \bibinfo{publisher}{Springer}, \bibinfo{year}{2013}.
\bibitem[{Majda and Bertozzi(2002)}]{MAJBAL}
\bibinfo{author}{A.J. Majda}, \bibinfo{author}{A.L. Bertozzi},
  \bibinfo{title}{{Vorticity and incompressible flow}},
  \bibinfo{publisher}{Cambridge university press}, \bibinfo{year}{2002}.
\bibitem[{group~of Mathematics~Handbook(2016)}]{WRGMDSFA}
\bibinfo{author}{W.~group~of Mathematics~Handbook},
  \bibinfo{title}{{Mathematics Handbook}}, \bibinfo{publisher}{Higher Education
  Press}, \bibinfo{year}{2016}.
\bibitem[{Meistersand and Olech(1963)}]{MGHCOLDUMJ}
\bibinfo{author}{G.H. Meistersand}, \bibinfo{author}{C.~Olech},
  \bibinfo{title}{Locally one-to-one mappings and a classical theorem on
  schlicht functions}, \bibinfo{journal}{Duke Math. J.} \bibinfo{volume}{30}
  (\bibinfo{year}{1963}) \bibinfo{pages}{63--80}.
\bibitem[{Novotn{\`y} and Stra{\v{s}}kraba(2004)}]{NASII04}
\bibinfo{author}{A.~Novotn{\`y}}, \bibinfo{author}{I.~Stra{\v{s}}kraba},
  \bibinfo{title}{{Introduction to the Mathematical Theory of Compressible
  Flow}}, \bibinfo{publisher}{Oxford University Press, USA},
  \bibinfo{year}{2004}.
\bibitem[{Pr$\mathrm{\ddot{u}}$ess and Simonett(2010)}]{PJSGOI5x}
\bibinfo{author}{J.~Pr$\mathrm{\ddot{u}}$ess}, \bibinfo{author}{G.~Simonett},
  \bibinfo{title}{{On the Rayleigh--Taylor instability for the two-phase
  Navier--Stokes equations}}, \bibinfo{journal}{Indiana Univ. Math. J.}
  \bibinfo{volume}{59} (\bibinfo{year}{2010}) \bibinfo{pages}{1853--1871}.
\bibitem[{Rayleigh(1990)}]{RLIS}
\bibinfo{author}{L.~Rayleigh}, \bibinfo{title}{{Investigation of the character
  of the equilibrium of an in compressible heavy fluid of variable density}},
  \bibinfo{journal}{Scientific Paper, II}  (\bibinfo{year}{1990})
  \bibinfo{pages}{200--207}.
\bibitem[{Rudin(2004)}]{WalterRudin}
\bibinfo{author}{W.~Rudin}, \bibinfo{title}{{Principles of Mathematical
  Analysis (Third Edition)}}, \bibinfo{year}{2004}.
\bibitem[{Sharma and Sharma(1978)}]{SHKCSHAR}
\bibinfo{author}{R.C. Sharma}, \bibinfo{author}{K.C. Sharma}, \bibinfo{title}{{
  Rayleigh--Taylor instability of two viscoelastic superposed fluids}},
  \bibinfo{journal}{Acta Physica Academiae Scientiarum Hungaricae, Tomus}
  \bibinfo{volume}{45} (\bibinfo{year}{1978}) \bibinfo{pages}{213--220}.
\bibitem[{Solonnikov(1988)}]{SVAUZNSE}
\bibinfo{author}{V.A. Solonnikov}, \bibinfo{title}{Unsteady motion of a finite
  mass of fluid bounded by a free surface}, \bibinfo{journal}{Zap. Nauchn. Sem.
  LOMI 152, 137--157 (1986) (in Russian); English transl.: J. Sov. Math.}
  \bibinfo{volume}{40} (\bibinfo{year}{1988}) \bibinfo{pages}{672--686}.
\bibitem[{Tani and Tanaka(1995)}]{TATNLTE1995}
\bibinfo{author}{A.~Tani}, \bibinfo{author}{N.~Tanaka},
  \bibinfo{title}{{Large-time existence of surface waves in incompressible
  viscous fluids with or without surface tension}}, \bibinfo{journal}{Arch.
  Ration. Mech. Anal.} \bibinfo{volume}{130} (\bibinfo{year}{1995})
  \bibinfo{pages}{303--314}.
\bibitem[{Taylor(1950)}]{TGTP}
\bibinfo{author}{G.I. Taylor}, \bibinfo{title}{{The stability of liquid surface
  when accelerated in a direction perpendicular to their planes}},
  \bibinfo{journal}{Proc. Roy Soc. A} \bibinfo{volume}{201}
  (\bibinfo{year}{1950}) \bibinfo{pages}{192--196}.
\bibitem[{Wang(1994)}]{WJH}
\bibinfo{author}{J.~Wang}, \bibinfo{title}{Two-Dimensional Nonsteady Flows and
  Shock Waves (in Chinese)}, \bibinfo{publisher}{Science Press},
  \bibinfo{address}{Beijing, China}, \bibinfo{year}{1994}.
\bibitem[{Wang and Zhao(2018)}]{WWWYYZ2018}
\bibinfo{author}{W.W. Wang}, \bibinfo{author}{Y.Y. Zhao}, \bibinfo{title}{{ On
  the Rayleigh--Taylor instability in compressible viscoelastic fluids}},
  \bibinfo{journal}{J. Math. Anal. Appl.} \bibinfo{volume}{463}
  (\bibinfo{year}{2018}) \bibinfo{pages}{198--221}.
\bibitem[{Wang(2012)}]{WYC}
\bibinfo{author}{Y.~Wang}, \bibinfo{title}{{Critical magnetic number in the MHD
  Rayleigh--Taylor instability}}, \bibinfo{journal}{Journal of Math. Phys.}
  \bibinfo{volume}{53} (\bibinfo{year}{2012}) \bibinfo{pages}{073701}.
\bibitem[{Wang(2018)}]{WYTIVNMI}
\bibinfo{author}{Y.~Wang}, \bibinfo{title}{{ Sharp nonlinear stability
  criterion of viscous non-resistive MHD internal waves in 3D}},
  \bibinfo{journal}{Arch. Rational Mech. Anal.,
  https://doi.org/10.1007/s00205-018-1307-4}  (\bibinfo{year}{2018}).
\bibitem[{Wang and Tice(2012)}]{wang2011viscous}
\bibinfo{author}{Y.~Wang}, \bibinfo{author}{I.~Tice}, \bibinfo{title}{The
  viscous surface-internal wave problem: nonlinear rayleigh--taylor
  instability}, \bibinfo{journal}{Commun. P.D.E.} \bibinfo{volume}{37}
  (\bibinfo{year}{2012}) \bibinfo{pages}{1967--2028}.
\bibitem[{Wang et~al.(2014)Wang, Tice and Kim}]{WYJTIKCT}
\bibinfo{author}{Y.J. Wang}, \bibinfo{author}{I.~Tice},
  \bibinfo{author}{C.~Kim}, \bibinfo{title}{The viscous surface-internal wave
  problem: global well-posedness and decay}, \bibinfo{journal}{Arch. Rational
  Mech. Anal.} \bibinfo{volume}{212} (\bibinfo{year}{2014})
  \bibinfo{pages}{1--92}.
\bibitem[{Wilke(2017)}]{WMRTITNS2017}
\bibinfo{author}{M.~Wilke}, \bibinfo{title}{{Rayleigh-Taylor instability for
  the two-phase Navier--Stokes equations with surface tension in cylindrical
  domains}}, \bibinfo{journal}{arXiv:1703.05214v1 [math.AP] 15 Mar 2017}
  (\bibinfo{year}{2017}).
\bibitem[{Xu et~al.(2013)Xu, Zhang and Zhang}]{XLZPZZFGAR}
\bibinfo{author}{L.~Xu}, \bibinfo{author}{P.~Zhang}, \bibinfo{author}{Z.F.
  Zhang}, \bibinfo{title}{{Global solvability of a free boundary
  three-dimensional incompressible viscoelastic fluid system with surface
  tension}}, \bibinfo{journal}{Arch. Ration. Mech. Anal.} \bibinfo{volume}{208}
  (\bibinfo{year}{2013}) \bibinfo{pages}{753--803}.

\end{thebibliography}

\end{document}